
\ifx\shlhetal\undefinedcontrolsequence\let\shlhetal\relax\fi
\def\fmtname{AmS-TeX}

\def\fmtversion{2.1}
\catcode`\@=11
\ifx\amstexloaded@\relax\catcode`\@=\active
  \endinput\else\let\amstexloaded@\relax\fi
\newlinechar=`\^^J
\def\W@{\immediate\write\sixt@@n}
\def\CR@{\W@{^^J\fmtname - Version \fmtversion^^J%
COPYRIGHT 1985, 1990, 1991 - AMERICAN MATHEMATICAL SOCIETY^^J%
Use of this macro package is not restricted provided^^J%
each use is acknowledged upon publication.^^J}}
\CR@ \everyjob{\CR@}
\message{Loading definitions for}
\message{misc utility macros,}
\toksdef\toks@@=2
\long\def\rightappend@#1\to#2{\toks@{\\{#1}}\toks@@
 =\expandafter{#2}\xdef#2{\the\toks@@\the\toks@}\toks@{}\toks@@{}}
\def\alloclist@{}
\newif\ifalloc@
\def\showallocations{{\def\\{\immediate\write\m@ne}\alloclist@}\alloc@true}
\def\alloc@#1#2#3#4#5{\global\advance\count1#1by\@ne
 \ch@ck#1#4#2\allocationnumber=\count1#1
 \global#3#5=\allocationnumber
 \edef\next@{\string#5=\string#2\the\allocationnumber}%
 \expandafter\rightappend@\next@\to\alloclist@}
\newcount\count@@
\newcount\count@@@
\def\FN@{\futurelet\next}
\def\DN@{\def\next@}
\def\DNii@{\def\nextii@}
\def\RIfM@{\relax\ifmmode}
\def\RIfMIfI@{\relax\ifmmode\ifinner}
\def\setboxz@h{\setbox\z@\hbox}
\def\wdz@{\wd\z@}
\def\boxz@{\box\z@}
\def\setbox@ne{\setbox\@ne}
\def\wd@ne{\wd\@ne}
\def\iterate{\body\expandafter\iterate\else\fi}
\def\err@#1{\errmessage{AmS-TeX error: #1}}
\newhelp\defaulthelp@{Sorry, I already gave what help I could...^^J
Maybe you should try asking a human?^^J
An error might have occurred before I noticed any problems.^^J
``If all else fails, read the instructions.''}
\def\Err@{\errhelp\defaulthelp@\err@}
\def\eat@#1{}
\def\in@#1#2{\def\in@@##1#1##2##3\in@@{\ifx\in@##2\in@false\else\in@true\fi}%
 \in@@#2#1\in@\in@@}
\newif\ifin@
\def\space@.{\futurelet\space@\relax}
\space@. %
\newhelp\athelp@
{Only certain combinations beginning with @ make sense to me.^^J
Perhaps you wanted \string\@\space for a printed @?^^J
I've ignored the character or group after @.}
{\catcode`\~=\active 
 \lccode`\~=`\@ \lowercase{\gdef~{\FN@\at@}}}
\def\at@{\let\next@\at@@
 \ifcat\noexpand\next a\else\ifcat\noexpand\next0\else
 \ifcat\noexpand\next\relax\else
   \let\next\at@@@\fi\fi\fi
 \next@}
\def\at@@#1{\expandafter
 \ifx\csname\space @\string#1\endcsname\relax
  \expandafter\at@@@ \else
  \csname\space @\string#1\expandafter\endcsname\fi}
\def\at@@@#1{\errhelp\athelp@ \err@{\Invalid@@ @}}
\def\atdef@#1{\expandafter\def\csname\space @\string#1\endcsname}
\newhelp\defahelp@{If you typed \string\define\space cs instead of
\string\define\string\cs\space^^J
I've substituted an inaccessible control sequence so that your^^J
definition will be completed without mixing me up too badly.^^J
If you typed \string\define{\string\cs} the inaccessible control sequence^^J
was defined to be \string\cs, and the rest of your^^J
definition appears as input.}
\newhelp\defbhelp@{I've ignored your definition, because it might^^J
conflict with other uses that are important to me.}
\def\define{\FN@\define@}
\def\define@{\ifcat\noexpand\next\relax
 \expandafter\define@@\else\errhelp\defahelp@                               
 \err@{\string\define\space must be followed by a control
 sequence}\expandafter\def\expandafter\nextii@\fi}                          
\def\undefined@@@@@@@@@@{}
\def\preloaded@@@@@@@@@@{}
\def\next@@@@@@@@@@{}
\def\define@@#1{\ifx#1\relax\errhelp\defbhelp@                              
 \err@{\string#1\space is already defined}\DN@{\DNii@}\else
 \expandafter\ifx\csname\expandafter\eat@\string                            
 #1@@@@@@@@@@\endcsname\undefined@@@@@@@@@@\errhelp\defbhelp@
 \err@{\string#1\space can't be defined}\DN@{\DNii@}\else
 \expandafter\ifx\csname\expandafter\eat@\string#1\endcsname\relax          
 \global\let#1\undefined\DN@{\def#1}\else\errhelp\defbhelp@
 \err@{\string#1\space is already defined}\DN@{\DNii@}\fi
 \fi\fi\next@}

\def\predefine#1#2{\let#1#2}
\def\undefine#1{\let#1\undefined}
\message{page layout,}
\newdimen\captionwidth@
\captionwidth@\hsize
\advance\captionwidth@-1.5in
\def\pagewidth#1{\hsize#1\relax
 \captionwidth@\hsize\advance\captionwidth@-1.5in}
\def\pageheight#1{\vsize#1\relax}
\def\hcorrection#1{\advance\hoffset#1\relax}
\def\vcorrection#1{\advance\voffset#1\relax}
\message{accents/punctuation,}

\let\graveaccent\`
\let\acuteaccent\'
\let\tildeaccent\~
\let\hataccent\^
\let\underscore\_
\let\B\=
\let\D\.
\let\ic@\/
\def\/{\unskip\ic@}
\def\textfonti{\the\textfont\@ne}
\def\t#1#2{{\edef\next@{\the\font}\textfonti\accent"7F \next@#1#2}}
\def~{\unskip\nobreak\ \ignorespaces}
\def\.{.\spacefactor\@m}
\atdef@;{\leavevmode\null;}
\atdef@:{\leavevmode\null:}
\atdef@?{\leavevmode\null?}
\edef\@{\string @}
\def\relaxnext@{\let\next\relax}
\atdef@-{\relaxnext@\leavevmode
 \DN@{\ifx\next-\DN@-{\FN@\nextii@}\else
  \DN@{\leavevmode\hbox{-}}\fi\next@}%
 \DNii@{\ifx\next-\DN@-{\leavevmode\hbox{---}}\else
  \DN@{\leavevmode\hbox{--}}\fi\next@}%
 \FN@\next@}
\def\srdr@{\kern.16667em}
\def\drsr@{\kern.02778em}
\def\sldl@{\drsr@}
\def\dlsl@{\srdr@}
\atdef@"{\unskip\relaxnext@
 \DN@{\ifx\next\space@\DN@. {\FN@\nextii@}\else
  \DN@.{\FN@\nextii@}\fi\next@.}%
 \DNii@{\ifx\next`\DN@`{\FN@\nextiii@}\else
  \ifx\next\lq\DN@\lq{\FN@\nextiii@}\else
  \DN@####1{\FN@\nextiv@}\fi\fi\next@}%
 \def\nextiii@{\ifx\next`\DN@`{\sldl@``}\else\ifx\next\lq
  \DN@\lq{\sldl@``}\else\DN@{\dlsl@`}\fi\fi\next@}%
 \def\nextiv@{\ifx\next'\DN@'{\srdr@''}\else
  \ifx\next\rq\DN@\rq{\srdr@''}\else\DN@{\drsr@'}\fi\fi\next@}%
 \FN@\next@}

\def\textfontii{\the\textfont\tw@}
\def\lbrace@{\delimiter"4266308 }
\def\rbrace@{\delimiter"5267309 }
\def\{{\RIfM@\lbrace@\else{\textfontii f}\spacefactor\@m\fi}
\def\}{\RIfM@\rbrace@\else
 \let\@sf\empty\ifhmode\edef\@sf{\spacefactor\the\spacefactor}\fi
 {\textfontii g}\@sf\relax\fi}
\let\lbrace\{
\let\rbrace\}
\def\AmSTeX{{\textfontii A\kern-.1667em%
  \lower.5ex\hbox{M}\kern-.125emS}-\TeX}
\message{line and page breaks,}
\def\vmodeerr@#1{\Err@{\string#1\space not allowed between paragraphs}}
\def\mathmodeerr@#1{\Err@{\string#1\space not allowed in math mode}}
\def\linebreak{\RIfM@\mathmodeerr@\linebreak\else
 \ifhmode\unskip\unkern\break\else\vmodeerr@\linebreak\fi\fi}

\newskip\saveskip@
\def\allowlinebreak{\RIfM@\mathmodeerr@\allowlinebreak\else
 \ifhmode\saveskip@\lastskip\unskip
 \allowbreak\ifdim\saveskip@>\z@\hskip\saveskip@\fi
 \else\vmodeerr@\allowlinebreak\fi\fi}
\def\nolinebreak{\RIfM@\mathmodeerr@\nolinebreak\else
 \ifhmode\saveskip@\lastskip\unskip
 \nobreak\ifdim\saveskip@>\z@\hskip\saveskip@\fi
 \else\vmodeerr@\nolinebreak\fi\fi}
\def\newline{\relaxnext@
 \DN@{\RIfM@\expandafter\mathmodeerr@\expandafter\newline\else
  \ifhmode\ifx\next\par\else
  \expandafter\unskip\expandafter\null\expandafter\hfill\expandafter\break\fi
  \else
  \expandafter\vmodeerr@\expandafter\newline\fi\fi}%
 \FN@\next@}
\def\dmatherr@#1{\Err@{\string#1\space not allowed in display math mode}}
\def\nondmatherr@#1{\Err@{\string#1\space not allowed in non-display math
 mode}}
\def\onlydmatherr@#1{\Err@{\string#1\space allowed only in display math mode}}
\def\nonmatherr@#1{\Err@{\string#1\space allowed only in math mode}}
\def\mathbreak{\RIfMIfI@\break\else
 \dmatherr@\mathbreak\fi\else\nonmatherr@\mathbreak\fi}
\def\nomathbreak{\RIfMIfI@\nobreak\else
 \dmatherr@\nomathbreak\fi\else\nonmatherr@\nomathbreak\fi}
\def\allowmathbreak{\RIfMIfI@\allowbreak\else
 \dmatherr@\allowmathbreak\fi\else\nonmatherr@\allowmathbreak\fi}
\def\pagebreak{\RIfM@
 \ifinner\nondmatherr@\pagebreak\else\postdisplaypenalty-\@M\fi
 \else\ifvmode\removelastskip\break\else\vadjust{\break}\fi\fi}
\def\nopagebreak{\RIfM@
 \ifinner\nondmatherr@\nopagebreak\else\postdisplaypenalty\@M\fi
 \else\ifvmode\nobreak\else\vadjust{\nobreak}\fi\fi}
\def\nonvmodeerr@#1{\Err@{\string#1\space not allowed within a paragraph
 or in math}}
\def\vnonvmode@#1#2{\relaxnext@\DNii@{\ifx\next\par\DN@{#1}\else
 \DN@{#2}\fi\next@}%
 \ifvmode\DN@{#1}\else
 \DN@{\FN@\nextii@}\fi\next@}
\def\newpage{\vnonvmode@{\vfill\break}{\nonvmodeerr@\newpage}}
\def\smallpagebreak{\vnonvmode@\smallbreak{\nonvmodeerr@\smallpagebreak}}
\def\medpagebreak{\vnonvmode@\medbreak{\nonvmodeerr@\medpagebreak}}
\def\bigpagebreak{\vnonvmode@\bigbreak{\nonvmodeerr@\bigpagebreak}}
\def\NoBlackBoxes{\global\overfullrule\z@}
\def\BlackBoxes{\global\overfullrule5\p@}
\def\Invalid@#1{\def#1{\Err@{\Invalid@@\string#1}}}
\def\Invalid@@{Invalid use of }
\message{figures,}
\Invalid@\caption
\Invalid@\captionwidth
\newdimen\smallcaptionwidth@
\def\topspace{\mid@false\ins@}
\def\midspace{\mid@true\ins@}
\newif\ifmid@
\def\captionfont@{}
\def\ins@#1{\relaxnext@\allowbreak
 \smallcaptionwidth@\captionwidth@\gdef\thespace@{#1}%
 \DN@{\ifx\next\space@\DN@. {\FN@\nextii@}\else
  \DN@.{\FN@\nextii@}\fi\next@.}%
 \DNii@{\ifx\next\caption\DN@\caption{\FN@\nextiii@}%
  \else\let\next@\nextiv@\fi\next@}%
 \def\nextiv@{\vnonvmode@
  {\ifmid@\expandafter\midinsert\else\expandafter\topinsert\fi
   \vbox to\thespace@{}\endinsert}
  {\ifmid@\nonvmodeerr@\midspace\else\nonvmodeerr@\topspace\fi}}%
 \def\nextiii@{\ifx\next\captionwidth\expandafter\nextv@
  \else\expandafter\nextvi@\fi}%
 \def\nextv@\captionwidth##1##2{\smallcaptionwidth@##1\relax\nextvi@{##2}}%
 \def\nextvi@##1{\def\thecaption@{\captionfont@##1}%
  \DN@{\ifx\next\space@\DN@. {\FN@\nextvii@}\else
   \DN@.{\FN@\nextvii@}\fi\next@.}%
  \FN@\next@}%
 \def\nextvii@{\vnonvmode@
  {\ifmid@\expandafter\midinsert\else
  \expandafter\topinsert\fi\vbox to\thespace@{}\nobreak\smallskip
  \setboxz@h{\noindent\ignorespaces\thecaption@\unskip}%
  \ifdim\wdz@>\smallcaptionwidth@\centerline{\vbox{\hsize\smallcaptionwidth@
   \noindent\ignorespaces\thecaption@\unskip}}%
  \else\centerline{\boxz@}\fi\endinsert}
  {\ifmid@\nonvmodeerr@\midspace
  \else\nonvmodeerr@\topspace\fi}}%
 \FN@\next@}
\message{comments,}
\def\newcodes@{\catcode`\\12\catcode`\{12\catcode`\}12\catcode`\#12%
 \catcode`\%12\relax}
\def\oldcodes@{\catcode`\\0\catcode`\{1\catcode`\}2\catcode`\#6%
 \catcode`\%14\relax}
\def\comment{\newcodes@\endlinechar=10 \comment@}
{\lccode`\0=`\\
\lowercase{\gdef\comment@#1^^J{\comment@@#10endcomment\comment@@@}%
\gdef\comment@@#10endcomment{\FN@\comment@@@}%
\gdef\comment@@@#1\comment@@@{\ifx\next\comment@@@\let\next\comment@
 \else\def\next{\oldcodes@\endlinechar=`\^^M\relax}%
 \fi\next}}}
\def\pr@m@s{\ifx'\next\DN@##1{\prim@s}\else\let\next@\egroup\fi\next@}
\def\prime{{\null\prime@\null}}
\mathchardef\prime@="0230
\let\dsize\displaystyle

\let\ssize\scriptstyle

\message{math spacing,}
\def\,{\RIfM@\mskip\thinmuskip\relax\else\kern.16667em\fi}
\def\!{\RIfM@\mskip-\thinmuskip\relax\else\kern-.16667em\fi}
\let\thinspace\,
\let\negthinspace\!
\def\medspace{\RIfM@\mskip\medmuskip\relax\else\kern.222222em\fi}
\def\negmedspace{\RIfM@\mskip-\medmuskip\relax\else\kern-.222222em\fi}
\def\thickspace{\RIfM@\mskip\thickmuskip\relax\else\kern.27777em\fi}
\let\;\thickspace
\def\negthickspace{\RIfM@\mskip-\thickmuskip\relax\else
 \kern-.27777em\fi}
\atdef@,{\RIfM@\mskip.1\thinmuskip\else\leavevmode\null,\fi}
\atdef@!{\RIfM@\mskip-.1\thinmuskip\else\leavevmode\null!\fi}
\atdef@.{\RIfM@&&\else\leavevmode.\spacefactor3000 \fi}
\def\and{\DOTSB\;\mathchar"3026 \;}

\message{fractions,}
\def\frac#1#2{{#1\over#2}}

\newdimen\ex@
\ex@.2326ex
\Invalid@\thickness
\def\thickfrac{\relaxnext@
 \DN@{\ifx\next\thickness\let\next@\nextii@\else
 \DN@{\nextii@\thickness1}\fi\next@}%
 \DNii@\thickness##1##2##3{{##2\above##1\ex@##3}}%
 \FN@\next@}

\def\thickfracwithdelims#1#2{\relaxnext@\def\ldelim@{#1}\def\rdelim@{#2}%
 \DN@{\ifx\next\thickness\let\next@\nextii@\else
 \DN@{\nextii@\thickness1}\fi\next@}%
 \DNii@\thickness##1##2##3{{##2\abovewithdelims
 \ldelim@\rdelim@##1\ex@##3}}%
 \FN@\next@}

\def\:{\nobreak\hskip.1111em\mathpunct{}\nonscript\mkern-\thinmuskip{:}\hskip
 .3333emplus.0555em\relax}
\def\snug{\unskip\kern-\mathsurround}
\message{smash commands,}
\def\topsmash{\top@true\bot@false\smash@}
\def\botsmash{\top@false\bot@true\smash@}
\newif\iftop@
\newif\ifbot@
\def\smash{\top@true\bot@true\smash@}
\def\smash@{\RIfM@\expandafter\mathpalette\expandafter\mathsm@sh\else
 \expandafter\makesm@sh\fi}
\def\finsm@sh{\iftop@\ht\z@\z@\fi\ifbot@\dp\z@\z@\fi\leavevmode\boxz@}
\message{large operator symbols,}
\def\LimitsOnSums{\global\let\slimits@\displaylimits}
\def\NoLimitsOnSums{\global\let\slimits@\nolimits}
\LimitsOnSums
\mathchardef\coprod@="1360       \def\coprod{\DOTSB\coprod@\slimits@}
\mathchardef\bigvee@="1357       \def\bigvee{\DOTSB\bigvee@\slimits@}
\mathchardef\bigwedge@="1356     \def\bigwedge{\DOTSB\bigwedge@\slimits@}
\mathchardef\biguplus@="1355     \def\biguplus{\DOTSB\biguplus@\slimits@}
\mathchardef\bigcap@="1354       \def\bigcap{\DOTSB\bigcap@\slimits@}
\mathchardef\bigcup@="1353       \def\bigcup{\DOTSB\bigcup@\slimits@}
\mathchardef\prod@="1351         \def\prod{\DOTSB\prod@\slimits@}
\mathchardef\sum@="1350          \def\sum{\DOTSB\sum@\slimits@}
\mathchardef\bigotimes@="134E    \def\bigotimes{\DOTSB\bigotimes@\slimits@}
\mathchardef\bigoplus@="134C     \def\bigoplus{\DOTSB\bigoplus@\slimits@}
\mathchardef\bigodot@="134A      \def\bigodot{\DOTSB\bigodot@\slimits@}
\mathchardef\bigsqcup@="1346     \def\bigsqcup{\DOTSB\bigsqcup@\slimits@}
\message{integrals,}
\def\LimitsOnInts{\global\let\ilimits@\displaylimits}
\def\NoLimitsOnInts{\global\let\ilimits@\nolimits}
\NoLimitsOnInts
\def\int{\DOTSI\intop\ilimits@}
\def\oint{\DOTSI\ointop\ilimits@}
\def\intic@{\mathchoice{\hskip.5em}{\hskip.4em}{\hskip.4em}{\hskip.4em}}
\def\negintic@{\mathchoice
 {\hskip-.5em}{\hskip-.4em}{\hskip-.4em}{\hskip-.4em}}
\def\intkern@{\mathchoice{\!\!\!}{\!\!}{\!\!}{\!\!}}
\def\intdots@{\mathchoice{\plaincdots@}
 {{\cdotp}\mkern1.5mu{\cdotp}\mkern1.5mu{\cdotp}}
 {{\cdotp}\mkern1mu{\cdotp}\mkern1mu{\cdotp}}
 {{\cdotp}\mkern1mu{\cdotp}\mkern1mu{\cdotp}}}
\newcount\intno@
\def\iint{\DOTSI\intno@\tw@\FN@\ints@}
\def\iiint{\DOTSI\intno@\thr@@\FN@\ints@}
\def\iiiint{\DOTSI\intno@4 \FN@\ints@}
\def\idotsint{\DOTSI\intno@\z@\FN@\ints@}
\def\ints@{\findlimits@\ints@@}
\newif\iflimtoken@
\newif\iflimits@
\def\findlimits@{\limtoken@true\ifx\next\limits\limits@true
 \else\ifx\next\nolimits\limits@false\else
 \limtoken@false\ifx\ilimits@\nolimits\limits@false\else
 \ifinner\limits@false\else\limits@true\fi\fi\fi\fi}
\def\multint@{\int\ifnum\intno@=\z@\intdots@                                
 \else\intkern@\fi                                                          
 \ifnum\intno@>\tw@\int\intkern@\fi                                         
 \ifnum\intno@>\thr@@\int\intkern@\fi                                       
 \int}                                                                      
\def\multintlimits@{\intop\ifnum\intno@=\z@\intdots@\else\intkern@\fi
 \ifnum\intno@>\tw@\intop\intkern@\fi
 \ifnum\intno@>\thr@@\intop\intkern@\fi\intop}
\def\ints@@{\iflimtoken@                                                    
 \def\ints@@@{\iflimits@\negintic@\mathop{\intic@\multintlimits@}\limits    
  \else\multint@\nolimits\fi                                                
  \eat@}                                                                    
 \else                                                                      
 \def\ints@@@{\iflimits@\negintic@
  \mathop{\intic@\multintlimits@}\limits\else
  \multint@\nolimits\fi}\fi\ints@@@}
\def\LimitsOnNames{\global\let\nlimits@\displaylimits}
\def\NoLimitsOnNames{\global\let\nlimits@\nolimits@}
\LimitsOnNames
\def\nolimits@{\relaxnext@
 \DN@{\ifx\next\limits\DN@\limits{\nolimits}\else
  \let\next@\nolimits\fi\next@}%
 \FN@\next@}
\message{operator names,}
\def\newmcodes@{\mathcode`\'"27\mathcode`\*"2A\mathcode`\."613A%
 \mathcode`\-"2D\mathcode`\/"2F\mathcode`\:"603A }
\def\operatorname#1{\mathop{\newmcodes@\kern\z@\fam\z@#1}\nolimits@}
\def\operatornamewithlimits#1{\mathop{\newmcodes@\kern\z@\fam\z@#1}\nlimits@}
\def\qopname@#1{\mathop{\fam\z@#1}\nolimits@}
\def\qopnamewl@#1{\mathop{\fam\z@#1}\nlimits@}
\def\arccos{\qopname@{arccos}}
\def\arcsin{\qopname@{arcsin}}
\def\arctan{\qopname@{arctan}}
\def\arg{\qopname@{arg}}
\def\cos{\qopname@{cos}}
\def\cosh{\qopname@{cosh}}
\def\cot{\qopname@{cot}}
\def\coth{\qopname@{coth}}
\def\csc{\qopname@{csc}}
\def\deg{\qopname@{deg}}
\def\det{\qopnamewl@{det}}
\def\dim{\qopname@{dim}}
\def\exp{\qopname@{exp}}
\def\gcd{\qopnamewl@{gcd}}
\def\hom{\qopname@{hom}}
\def\inf{\qopnamewl@{inf}}
\def\injlim{\qopnamewl@{inj\,lim}}
\def\ker{\qopname@{ker}}
\def\lg{\qopname@{lg}}
\def\lim{\qopnamewl@{lim}}
\def\liminf{\qopnamewl@{lim\,inf}}
\def\limsup{\qopnamewl@{lim\,sup}}
\def\ln{\qopname@{ln}}
\def\log{\qopname@{log}}
\def\max{\qopnamewl@{max}}
\def\min{\qopnamewl@{min}}
\def\Pr{\qopnamewl@{Pr}}
\def\projlim{\qopnamewl@{proj\,lim}}
\def\sec{\qopname@{sec}}
\def\sin{\qopname@{sin}}
\def\sinh{\qopname@{sinh}}
\def\sup{\qopnamewl@{sup}}
\def\tan{\qopname@{tan}}
\def\tanh{\qopname@{tanh}}
\def\varinjlim{\mathop{\vtop{\ialign{##\crcr
 \hfil\rm lim\hfil\crcr\noalign{\nointerlineskip}\rightarrowfill\crcr
 \noalign{\nointerlineskip\kern-\ex@}\crcr}}}}
\def\varprojlim{\mathop{\vtop{\ialign{##\crcr
 \hfil\rm lim\hfil\crcr\noalign{\nointerlineskip}\leftarrowfill\crcr
 \noalign{\nointerlineskip\kern-\ex@}\crcr}}}}
\def\varliminf{\mathop{\underline{\vrule height\z@ depth.2exwidth\z@
 \hbox{\rm lim}}}}

\newdimen\buffer@
\buffer@\fontdimen13 \tenex
\newdimen\buffer
\buffer\buffer@

\def\ResetBuffer{\fontdimen13 \tenex\buffer@\global\buffer\buffer@}
\def\shave#1{\mathop{\hbox{$\m@th\fontdimen13 \tenex\z@                     
 \displaystyle{#1}$}}\fontdimen13 \tenex\buffer}

\message{multilevel sub/superscripts,}
\Invalid@\\
\def\Let@{\relax\iffalse{\fi\let\\=\cr\iffalse}\fi}
\Invalid@\vspace
\def\vspace@{\def\vspace##1{\crcr\noalign{\vskip##1\relax}}}
\def\multilimits@{\bgroup\vspace@\Let@
 \baselineskip\fontdimen10 \scriptfont\tw@
 \advance\baselineskip\fontdimen12 \scriptfont\tw@
 \lineskip\thr@@\fontdimen8 \scriptfont\thr@@
 \lineskiplimit\lineskip
 \vbox\bgroup\ialign\bgroup\hfil$\m@th\scriptstyle{##}$\hfil\crcr}
\def\Sb{_\multilimits@}
\def\endSb{\crcr\egroup\egroup\egroup}
\def\Sp{^\multilimits@}

\def\spreadlines#1{\RIfMIfI@\onlydmatherr@\spreadlines\else
 \openup#1\relax\fi\else\onlydmatherr@\spreadlines\fi}
\def\Mathstrut@{\copy\Mathstrutbox@}
\newbox\Mathstrutbox@
\setbox\Mathstrutbox@\null
\setboxz@h{$\m@th($}
\ht\Mathstrutbox@\ht\z@
\dp\Mathstrutbox@\dp\z@
\message{matrices,}
\newdimen\spreadmlines@
\def\spreadmatrixlines#1{\RIfMIfI@
 \onlydmatherr@\spreadmatrixlines\else
 \spreadmlines@#1\relax\fi\else\onlydmatherr@\spreadmatrixlines\fi}
\def\matrix{\null\,\vcenter\bgroup\Let@\vspace@
 \normalbaselines\openup\spreadmlines@\ialign
 \bgroup\hfil$\m@th##$\hfil&&\quad\hfil$\m@th##$\hfil\crcr
 \Mathstrut@\crcr\noalign{\kern-\baselineskip}}
\def\endmatrix{\crcr\Mathstrut@\crcr\noalign{\kern-\baselineskip}\egroup
 \egroup\,}
\def\format{\crcr\egroup\iffalse{\fi\ifnum`}=0 \fi\format@}
\newtoks\hashtoks@
\hashtoks@{#}
\def\format@#1\\{\def\preamble@{#1}%
 \def\l{$\m@th\the\hashtoks@$\hfil}%
 \def\c{\hfil$\m@th\the\hashtoks@$\hfil}%
 \def\r{\hfil$\m@th\the\hashtoks@$}%
 \edef\preamble@@{\preamble@}\ifnum`{=0 \fi\iffalse}\fi
 \ialign\bgroup\span\preamble@@\crcr}
\def\smallmatrix{\null\,\vcenter\bgroup\vspace@\Let@
 \baselineskip9\ex@\lineskip\ex@
 \ialign\bgroup\hfil$\m@th\scriptstyle{##}$\hfil&&\thickspace\hfil
 $\m@th\scriptstyle{##}$\hfil\crcr}
\def\endsmallmatrix{\crcr\egroup\egroup\,}

\newmuskip\dotsspace@
\dotsspace@1.5mu
\def\strip@#1 {#1}
\def\spacehdots#1\for#2{\multispan{#2}\xleaders
 \hbox{$\m@th\mkern\strip@#1 \dotsspace@.\mkern\strip@#1 \dotsspace@$}\hfill}
\def\hdotsfor#1{\spacehdots\@ne\for{#1}}
\def\multispan@#1{\omit\mscount#1\unskip\loop\ifnum\mscount>\@ne\sp@n\repeat}
\def\spaceinnerhdots#1\for#2\after#3{\multispan@{\strip@#2 }#3\xleaders
 \hbox{$\m@th\mkern\strip@#1 \dotsspace@.\mkern\strip@#1 \dotsspace@$}\hfill}
\def\innerhdotsfor#1\after#2{\spaceinnerhdots\@ne\for#1\after{#2}}
\def\cases{\bgroup\spreadmlines@\jot\left\{\,\matrix\format\l&\quad\l\\}
\def\endcases{\endmatrix\right.\egroup}
\message{multiline displays,}
\newif\ifinany@
\newif\ifinalign@
\newif\ifingather@
\def\strut@{\copy\strutbox@}
\newbox\strutbox@
\setbox\strutbox@\hbox{\vrule height8\p@ depth3\p@ width\z@}
\def\topaligned{\null\,\vtop\aligned@}
\def\botaligned{\null\,\vbox\aligned@}
\def\aligned{\null\,\vcenter\aligned@}
\def\aligned@{\bgroup\vspace@\Let@
 \ifinany@\else\openup\jot\fi\ialign
 \bgroup\hfil\strut@$\m@th\displaystyle{##}$&
 $\m@th\displaystyle{{}##}$\hfil\crcr}
\def\endaligned{\crcr\egroup\egroup}

\def\alignedat#1{\null\,\vcenter\bgroup\doat@{#1}\vspace@\Let@
 \ifinany@\else\openup\jot\fi\ialign\bgroup\span\preamble@@\crcr}
\newcount\atcount@
\def\doat@#1{\toks@{\hfil\strut@$\m@th
 \displaystyle{\the\hashtoks@}$&$\m@th\displaystyle
 {{}\the\hashtoks@}$\hfil}
 \atcount@#1\relax\advance\atcount@\m@ne                                    
 \loop\ifnum\atcount@>\z@\toks@=\expandafter{\the\toks@&\hfil$\m@th
 \displaystyle{\the\hashtoks@}$&$\m@th
 \displaystyle{{}\the\hashtoks@}$\hfil}\advance
  \atcount@\m@ne\repeat                                                     
 \xdef\preamble@{\the\toks@}\xdef\preamble@@{\preamble@}}

\def\gathered{\null\,\vcenter\bgroup\vspace@\Let@
 \ifinany@\else\openup\jot\fi\ialign
 \bgroup\hfil\strut@$\m@th\displaystyle{##}$\hfil\crcr}
\def\endgathered{\crcr\egroup\egroup}
\newif\iftagsleft@
\def\TagsOnLeft{\global\tagsleft@true}
\def\TagsOnRight{\global\tagsleft@false}
\TagsOnLeft
\newif\ifmathtags@
\def\TagsAsMath{\global\mathtags@true}
\def\TagsAsText{\global\mathtags@false}
\TagsAsText
\def\tagform@#1{\hbox{\rm(\ignorespaces#1\unskip)}}
\def\thetag{\leavevmode\tagform@}
\def\tag#1$${\iftagsleft@\leqno\else\eqno\fi                                
 \maketag@#1\maketag@                                                       
 $$}                                                                        
\def\maketag@{\FN@\maketag@@}
\def\maketag@@{\ifx\next"\expandafter\maketag@@@\else\expandafter\maketag@@@@
 \fi}
\def\maketag@@@"#1"#2\maketag@{\hbox{\rm#1}}                                
\def\maketag@@@@#1\maketag@{\ifmathtags@\tagform@{$\m@th#1$}\else
 \tagform@{#1}\fi}
\interdisplaylinepenalty\@M
\def\allowdisplaybreaks{\RIfMIfI@
 \onlydmatherr@\allowdisplaybreaks\else
 \interdisplaylinepenalty\z@\fi\else\onlydmatherr@\allowdisplaybreaks\fi}
\Invalid@\allowdisplaybreak
\Invalid@\displaybreak
\Invalid@\intertext
\def\allowdisplaybreak@{\def\allowdisplaybreak{\crcr\noalign{\allowbreak}}}
\def\displaybreak@{\def\displaybreak{\crcr\noalign{\break}}}
\def\intertext@{\def\intertext##1{\crcr\noalign{%
 \penalty\postdisplaypenalty \vskip\belowdisplayskip
 \vbox{\normalbaselines\noindent##1}%
 \penalty\predisplaypenalty \vskip\abovedisplayskip}}}
\newskip\centering@
\centering@\z@ plus\@m\p@
\def\align{\relax\ifingather@\DN@{\csname align (in
  \string\gather)\endcsname}\else
 \ifmmode\ifinner\DN@{\onlydmatherr@\align}\else
  \let\next@\align@\fi
 \else\DN@{\onlydmatherr@\align}\fi\fi\next@}
\newhelp\andhelp@
{An extra & here is so disastrous that you should probably exit^^J
and fix things up.}
\newif\iftag@
\newcount\and@
\def\align@{\inalign@true\inany@true
 \vspace@\allowdisplaybreak@\displaybreak@\intertext@
 \def\tag{\global\tag@true\ifnum\and@=\z@\DN@{&&}\else
          \DN@{&}\fi\next@}%
 \iftagsleft@\DN@{\csname align \endcsname}\else
  \DN@{\csname align \space\endcsname}\fi\next@}
\def\Tag@{\iftag@\else\errhelp\andhelp@\err@{Extra & on this line}\fi}
\newdimen\lwidth@
\newdimen\rwidth@
\newdimen\maxlwidth@
\newdimen\maxrwidth@
\newdimen\totwidth@
\def\measure@#1\endalign{\lwidth@\z@\rwidth@\z@\maxlwidth@\z@\maxrwidth@\z@
 \global\and@\z@                                                            
 \setbox@ne\vbox                                                            
  {\everycr{\noalign{\global\tag@false\global\and@\z@}}\Let@                
  \halign{\setboxz@h{$\m@th\displaystyle{\@lign##}$}
   \global\lwidth@\wdz@                                                     
   \ifdim\lwidth@>\maxlwidth@\global\maxlwidth@\lwidth@\fi                  
   \global\advance\and@\@ne                                                 
   &\setboxz@h{$\m@th\displaystyle{{}\@lign##}$}\global\rwidth@\wdz@        
   \ifdim\rwidth@>\maxrwidth@\global\maxrwidth@\rwidth@\fi                  
   \global\advance\and@\@ne                                                
   &\Tag@
   \eat@{##}\crcr#1\crcr}}
 \totwidth@\maxlwidth@\advance\totwidth@\maxrwidth@}                       
\def\displ@y@{\global\dt@ptrue\openup\jot
 \everycr{\noalign{\global\tag@false\global\and@\z@\ifdt@p\global\dt@pfalse
 \vskip-\lineskiplimit\vskip\normallineskiplimit\else
 \penalty\interdisplaylinepenalty\fi}}}
\def\black@#1{\noalign{\ifdim#1>\displaywidth
 \dimen@\prevdepth\nointerlineskip                                          
 \vskip-\ht\strutbox@\vskip-\dp\strutbox@                                   
 \vbox{\noindent\hbox to#1{\strut@\hfill}}
 \prevdepth\dimen@                                                          
 \fi}}
\expandafter\def\csname align \space\endcsname#1\endalign
 {\measure@#1\endalign\global\and@\z@                                       
 \ifingather@\everycr{\noalign{\global\and@\z@}}\else\displ@y@\fi           
 \Let@\tabskip\centering@                                                   
 \halign to\displaywidth
  {\hfil\strut@\setboxz@h{$\m@th\displaystyle{\@lign##}$}
  \global\lwidth@\wdz@\boxz@\global\advance\and@\@ne                        
  \tabskip\z@skip                                                           
  &\setboxz@h{$\m@th\displaystyle{{}\@lign##}$}
  \global\rwidth@\wdz@\boxz@\hfill\global\advance\and@\@ne                  
  \tabskip\centering@                                                       
  &\setboxz@h{\@lign\strut@\maketag@##\maketag@}
  \dimen@\displaywidth\advance\dimen@-\totwidth@
  \divide\dimen@\tw@\advance\dimen@\maxrwidth@\advance\dimen@-\rwidth@     
  \ifdim\dimen@<\tw@\wdz@\llap{\vtop{\normalbaselines\null\boxz@}}
  \else\llap{\boxz@}\fi                                                    
  \tabskip\z@skip                                                          
  \crcr#1\crcr                                                             
  \black@\totwidth@}}                                                      
\newdimen\lineht@
\expandafter\def\csname align \endcsname#1\endalign{\measure@#1\endalign
 \global\and@\z@
 \ifdim\totwidth@>\displaywidth\let\displaywidth@\totwidth@\else
  \let\displaywidth@\displaywidth\fi                                        
 \ifingather@\everycr{\noalign{\global\and@\z@}}\else\displ@y@\fi
 \Let@\tabskip\centering@\halign to\displaywidth
  {\hfil\strut@\setboxz@h{$\m@th\displaystyle{\@lign##}$}%
  \global\lwidth@\wdz@\global\lineht@\ht\z@                                 
  \boxz@\global\advance\and@\@ne
  \tabskip\z@skip&\setboxz@h{$\m@th\displaystyle{{}\@lign##}$}%
  \global\rwidth@\wdz@\ifdim\ht\z@>\lineht@\global\lineht@\ht\z@\fi         
  \boxz@\hfil\global\advance\and@\@ne
  \tabskip\centering@&\kern-\displaywidth@                                  
  \setboxz@h{\@lign\strut@\maketag@##\maketag@}%
  \dimen@\displaywidth\advance\dimen@-\totwidth@
  \divide\dimen@\tw@\advance\dimen@\maxlwidth@\advance\dimen@-\lwidth@
  \ifdim\dimen@<\tw@\wdz@
   \rlap{\vbox{\normalbaselines\boxz@\vbox to\lineht@{}}}\else
   \rlap{\boxz@}\fi
  \tabskip\displaywidth@\crcr#1\crcr\black@\totwidth@}}
\expandafter\def\csname align (in \string\gather)\endcsname
  #1\endalign{\vcenter{\align@#1\endalign}}
\Invalid@\endalign
\newif\ifxat@
\def\alignat{\RIfMIfI@\DN@{\onlydmatherr@\alignat}\else
 \DN@{\csname alignat \endcsname}\fi\else
 \DN@{\onlydmatherr@\alignat}\fi\next@}
\newif\ifmeasuring@
\newbox\savealignat@
\expandafter\def\csname alignat \endcsname#1#2\endalignat                   
 {\inany@true\xat@false
 \def\tag{\global\tag@true\count@#1\relax\multiply\count@\tw@
  \xdef\tag@{}\loop\ifnum\count@>\and@\xdef\tag@{&\tag@}\advance\count@\m@ne
  \repeat\tag@}%
 \vspace@\allowdisplaybreak@\displaybreak@\intertext@
 \displ@y@\measuring@true                                                   
 \setbox\savealignat@\hbox{$\m@th\displaystyle\Let@
  \attag@{#1}
  \vbox{\halign{\span\preamble@@\crcr#2\crcr}}$}%
 \measuring@false                                                           
 \Let@\attag@{#1}
 \tabskip\centering@\halign to\displaywidth
  {\span\preamble@@\crcr#2\crcr                                             
  \black@{\wd\savealignat@}}}                                               
\Invalid@\endalignat
\def\xalignat{\RIfMIfI@
 \DN@{\onlydmatherr@\xalignat}\else
 \DN@{\csname xalignat \endcsname}\fi\else
 \DN@{\onlydmatherr@\xalignat}\fi\next@}
\expandafter\def\csname xalignat \endcsname#1#2\endxalignat
 {\inany@true\xat@true
 \def\tag{\global\tag@true\def\tag@{}\count@#1\relax\multiply\count@\tw@
  \loop\ifnum\count@>\and@\xdef\tag@{&\tag@}\advance\count@\m@ne\repeat\tag@}%
 \vspace@\allowdisplaybreak@\displaybreak@\intertext@
 \displ@y@\measuring@true\setbox\savealignat@\hbox{$\m@th\displaystyle\Let@
 \attag@{#1}\vbox{\halign{\span\preamble@@\crcr#2\crcr}}$}%
 \measuring@false\Let@
 \attag@{#1}\tabskip\centering@\halign to\displaywidth
 {\span\preamble@@\crcr#2\crcr\black@{\wd\savealignat@}}}
\def\attag@#1{\let\Maketag@\maketag@\let\TAG@\Tag@                          
 \let\Tag@=0\let\maketag@=0
 \ifmeasuring@\def\llap@##1{\setboxz@h{##1}\hbox to\tw@\wdz@{}}%
  \def\rlap@##1{\setboxz@h{##1}\hbox to\tw@\wdz@{}}\else
  \let\llap@\llap\let\rlap@\rlap\fi                                         
 \toks@{\hfil\strut@$\m@th\displaystyle{\@lign\the\hashtoks@}$\tabskip\z@skip
  \global\advance\and@\@ne&$\m@th\displaystyle{{}\@lign\the\hashtoks@}$\hfil
  \ifxat@\tabskip\centering@\fi\global\advance\and@\@ne}
 \iftagsleft@
  \toks@@{\tabskip\centering@&\Tag@\kern-\displaywidth
   \rlap@{\@lign\maketag@\the\hashtoks@\maketag@}%
   \global\advance\and@\@ne\tabskip\displaywidth}\else
  \toks@@{\tabskip\centering@&\Tag@\llap@{\@lign\maketag@
   \the\hashtoks@\maketag@}\global\advance\and@\@ne\tabskip\z@skip}\fi      
 \atcount@#1\relax\advance\atcount@\m@ne
 \loop\ifnum\atcount@>\z@
 \toks@=\expandafter{\the\toks@&\hfil$\m@th\displaystyle{\@lign
  \the\hashtoks@}$\global\advance\and@\@ne
  \tabskip\z@skip&$\m@th\displaystyle{{}\@lign\the\hashtoks@}$\hfil\ifxat@
  \tabskip\centering@\fi\global\advance\and@\@ne}\advance\atcount@\m@ne
 \repeat                                                                    
 \xdef\preamble@{\the\toks@\the\toks@@}
 \xdef\preamble@@{\preamble@}
 \let\maketag@\Maketag@\let\Tag@\TAG@}                                      
\Invalid@\endxalignat
\def\xxalignat{\RIfMIfI@
 \DN@{\onlydmatherr@\xxalignat}\else\DN@{\csname xxalignat
  \endcsname}\fi\else
 \DN@{\onlydmatherr@\xxalignat}\fi\next@}
\expandafter\def\csname xxalignat \endcsname#1#2\endxxalignat{\inany@true
 \vspace@\allowdisplaybreak@\displaybreak@\intertext@
 \displ@y\setbox\savealignat@\hbox{$\m@th\displaystyle\Let@
 \xxattag@{#1}\vbox{\halign{\span\preamble@@\crcr#2\crcr}}$}%
 \Let@\xxattag@{#1}\tabskip\z@skip\halign to\displaywidth
 {\span\preamble@@\crcr#2\crcr\black@{\wd\savealignat@}}}
\def\xxattag@#1{\toks@{\tabskip\z@skip\hfil\strut@
 $\m@th\displaystyle{\the\hashtoks@}$&%
 $\m@th\displaystyle{{}\the\hashtoks@}$\hfil\tabskip\centering@&}%
 \atcount@#1\relax\advance\atcount@\m@ne\loop\ifnum\atcount@>\z@
 \toks@=\expandafter{\the\toks@&\hfil$\m@th\displaystyle{\the\hashtoks@}$%
  \tabskip\z@skip&$\m@th\displaystyle{{}\the\hashtoks@}$\hfil
  \tabskip\centering@}\advance\atcount@\m@ne\repeat
 \xdef\preamble@{\the\toks@\tabskip\z@skip}\xdef\preamble@@{\preamble@}}
\Invalid@\endxxalignat
\newdimen\gwidth@
\newdimen\gmaxwidth@
\def\gmeasure@#1\endgather{\gwidth@\z@\gmaxwidth@\z@\setbox@ne\vbox{\Let@
 \halign{\setboxz@h{$\m@th\displaystyle{##}$}\global\gwidth@\wdz@
 \ifdim\gwidth@>\gmaxwidth@\global\gmaxwidth@\gwidth@\fi
 &\eat@{##}\crcr#1\crcr}}}
\def\gather{\RIfMIfI@\DN@{\onlydmatherr@\gather}\else
 \ingather@true\inany@true\def\tag{&}%
 \vspace@\allowdisplaybreak@\displaybreak@\intertext@
 \displ@y\Let@
 \iftagsleft@\DN@{\csname gather \endcsname}\else
  \DN@{\csname gather \space\endcsname}\fi\fi
 \else\DN@{\onlydmatherr@\gather}\fi\next@}
\expandafter\def\csname gather \space\endcsname#1\endgather
 {\gmeasure@#1\endgather\tabskip\centering@
 \halign to\displaywidth{\hfil\strut@\setboxz@h{$\m@th\displaystyle{##}$}%
 \global\gwidth@\wdz@\boxz@\hfil&
 \setboxz@h{\strut@{\maketag@##\maketag@}}%
 \dimen@\displaywidth\advance\dimen@-\gwidth@
 \ifdim\dimen@>\tw@\wdz@\llap{\boxz@}\else
 \llap{\vtop{\normalbaselines\null\boxz@}}\fi
 \tabskip\z@skip\crcr#1\crcr\black@\gmaxwidth@}}
\newdimen\glineht@
\expandafter\def\csname gather \endcsname#1\endgather{\gmeasure@#1\endgather
 \ifdim\gmaxwidth@>\displaywidth\let\gdisplaywidth@\gmaxwidth@\else
 \let\gdisplaywidth@\displaywidth\fi\tabskip\centering@\halign to\displaywidth
 {\hfil\strut@\setboxz@h{$\m@th\displaystyle{##}$}%
 \global\gwidth@\wdz@\global\glineht@\ht\z@\boxz@\hfil&\kern-\gdisplaywidth@
 \setboxz@h{\strut@{\maketag@##\maketag@}}%
 \dimen@\displaywidth\advance\dimen@-\gwidth@
 \ifdim\dimen@>\tw@\wdz@\rlap{\boxz@}\else
 \rlap{\vbox{\normalbaselines\boxz@\vbox to\glineht@{}}}\fi
 \tabskip\gdisplaywidth@\crcr#1\crcr\black@\gmaxwidth@}}
\newif\ifctagsplit@
\def\CenteredTagsOnSplits{\global\ctagsplit@true}
\def\TopOrBottomTagsOnSplits{\global\ctagsplit@false}
\TopOrBottomTagsOnSplits
\def\split{\relax\ifinany@\let\next@\insplit@\else
 \ifmmode\ifinner\def\next@{\onlydmatherr@\split}\else
 \let\next@\outsplit@\fi\else
 \def\next@{\onlydmatherr@\split}\fi\fi\next@}
\def\insplit@{\global\setbox\z@\vbox\bgroup\vspace@\Let@\ialign\bgroup
 \hfil\strut@$\m@th\displaystyle{##}$&$\m@th\displaystyle{{}##}$\hfill\crcr}
\def\endsplit{\crcr\egroup\egroup\iftagsleft@\expandafter\lendsplit@\else
 \expandafter\rendsplit@\fi}
\def\rendsplit@{\global\setbox9 \vbox
 {\unvcopy\z@\global\setbox8 \lastbox\unskip}
 \setbox@ne\hbox{\unhcopy8 \unskip\global\setbox\tw@\lastbox
 \unskip\global\setbox\thr@@\lastbox}
 \global\setbox7 \hbox{\unhbox\tw@\unskip}
 \ifinalign@\ifctagsplit@                                                   
  \gdef\split@{\hbox to\wd\thr@@{}&
   \vcenter{\vbox{\moveleft\wd\thr@@\boxz@}}}
 \else\gdef\split@{&\vbox{\moveleft\wd\thr@@\box9}\crcr
  \box\thr@@&\box7}\fi                                                      
 \else                                                                      
  \ifctagsplit@\gdef\split@{\vcenter{\boxz@}}\else
  \gdef\split@{\box9\crcr\hbox{\box\thr@@\box7}}\fi
 \fi
 \split@}                                                                   
\def\lendsplit@{\global\setbox9\vtop{\unvcopy\z@}
 \setbox@ne\vbox{\unvcopy\z@\global\setbox8\lastbox}
 \setbox@ne\hbox{\unhcopy8\unskip\setbox\tw@\lastbox
  \unskip\global\setbox\thr@@\lastbox}
 \ifinalign@\ifctagsplit@                                                   
  \gdef\split@{\hbox to\wd\thr@@{}&
  \vcenter{\vbox{\moveleft\wd\thr@@\box9}}}
  \else                                                                     
  \gdef\split@{\hbox to\wd\thr@@{}&\vbox{\moveleft\wd\thr@@\box9}}\fi
 \else
  \ifctagsplit@\gdef\split@{\vcenter{\box9}}\else
  \gdef\split@{\box9}\fi
 \fi\split@}
\def\outsplit@#1$${\align\insplit@#1\endalign$$}
\newdimen\multlinegap@
\multlinegap@1em
\newdimen\multlinetaggap@
\multlinetaggap@1em
\def\MultlineGap#1{\global\multlinegap@#1\relax}
\def\multlinegap#1{\RIfMIfI@\onlydmatherr@\multlinegap\else
 \multlinegap@#1\relax\fi\else\onlydmatherr@\multlinegap\fi}
\def\nomultlinegap{\multlinegap{\z@}}
\def\multline{\RIfMIfI@
 \DN@{\onlydmatherr@\multline}\else
 \DN@{\multline@}\fi\else
 \DN@{\onlydmatherr@\multline}\fi\next@}
\newif\iftagin@
\def\tagin@#1{\tagin@false\in@\tag{#1}\ifin@\tagin@true\fi}
\def\multline@#1$${\inany@true\vspace@\allowdisplaybreak@\displaybreak@
 \tagin@{#1}\iftagsleft@\DN@{\multline@l#1$$}\else
 \DN@{\multline@r#1$$}\fi\next@}
\newdimen\mwidth@
\def\rmmeasure@#1\endmultline{%
 \def\shoveleft##1{##1}\def\shoveright##1{##1}
 \setbox@ne\vbox{\Let@\halign{\setboxz@h
  {$\m@th\@lign\displaystyle{}##$}\global\mwidth@\wdz@
  \crcr#1\crcr}}}
\newdimen\mlineht@
\newif\ifzerocr@
\newif\ifonecr@
\def\lmmeasure@#1\endmultline{\global\zerocr@true\global\onecr@false
 \everycr{\noalign{\ifonecr@\global\onecr@false\fi
  \ifzerocr@\global\zerocr@false\global\onecr@true\fi}}
  \def\shoveleft##1{##1}\def\shoveright##1{##1}%
 \setbox@ne\vbox{\Let@\halign{\setboxz@h
  {$\m@th\@lign\displaystyle{}##$}\ifonecr@\global\mwidth@\wdz@
  \global\mlineht@\ht\z@\fi\crcr#1\crcr}}}
\newbox\mtagbox@
\newdimen\ltwidth@
\newdimen\rtwidth@
\def\multline@l#1$${\iftagin@\DN@{\lmultline@@#1$$}\else
 \DN@{\setbox\mtagbox@\null\ltwidth@\z@\rtwidth@\z@
  \lmultline@@@#1$$}\fi\next@}
\def\lmultline@@#1\endmultline\tag#2$${%
 \setbox\mtagbox@\hbox{\maketag@#2\maketag@}
 \lmmeasure@#1\endmultline\dimen@\mwidth@\advance\dimen@\wd\mtagbox@
 \advance\dimen@\multlinetaggap@                                            
 \ifdim\dimen@>\displaywidth\ltwidth@\z@\else\ltwidth@\wd\mtagbox@\fi       
 \lmultline@@@#1\endmultline$$}
\def\lmultline@@@{\displ@y
 \def\shoveright##1{##1\hfilneg\hskip\multlinegap@}%
 \def\shoveleft##1{\setboxz@h{$\m@th\displaystyle{}##1$}%
  \setbox@ne\hbox{$\m@th\displaystyle##1$}%
  \hfilneg
  \iftagin@
   \ifdim\ltwidth@>\z@\hskip\ltwidth@\hskip\multlinetaggap@\fi
  \else\hskip\multlinegap@\fi\hskip.5\wd@ne\hskip-.5\wdz@##1}
  \halign\bgroup\Let@\hbox to\displaywidth
   {\strut@$\m@th\displaystyle\hfil{}##\hfil$}\crcr
   \hfilneg                                                                 
   \iftagin@                                                                
    \ifdim\ltwidth@>\z@                                                     
     \box\mtagbox@\hskip\multlinetaggap@                                    
    \else
     \rlap{\vbox{\normalbaselines\hbox{\strut@\box\mtagbox@}%
     \vbox to\mlineht@{}}}\fi                                               
   \else\hskip\multlinegap@\fi}                                             
\def\multline@r#1$${\iftagin@\DN@{\rmultline@@#1$$}\else
 \DN@{\setbox\mtagbox@\null\ltwidth@\z@\rtwidth@\z@
  \rmultline@@@#1$$}\fi\next@}
\def\rmultline@@#1\endmultline\tag#2$${\ltwidth@\z@
 \setbox\mtagbox@\hbox{\maketag@#2\maketag@}%
 \rmmeasure@#1\endmultline\dimen@\mwidth@\advance\dimen@\wd\mtagbox@
 \advance\dimen@\multlinetaggap@
 \ifdim\dimen@>\displaywidth\rtwidth@\z@\else\rtwidth@\wd\mtagbox@\fi
 \rmultline@@@#1\endmultline$$}
\def\rmultline@@@{\displ@y
 \def\shoveright##1{##1\hfilneg\iftagin@\ifdim\rtwidth@>\z@
  \hskip\rtwidth@\hskip\multlinetaggap@\fi\else\hskip\multlinegap@\fi}%
 \def\shoveleft##1{\setboxz@h{$\m@th\displaystyle{}##1$}%
  \setbox@ne\hbox{$\m@th\displaystyle##1$}%
  \hfilneg\hskip\multlinegap@\hskip.5\wd@ne\hskip-.5\wdz@##1}%
 \halign\bgroup\Let@\hbox to\displaywidth
  {\strut@$\m@th\displaystyle\hfil{}##\hfil$}\crcr
 \hfilneg\hskip\multlinegap@}
\def\endmultline{\iftagsleft@\expandafter\lendmultline@\else
 \expandafter\rendmultline@\fi}
\def\lendmultline@{\hfilneg\hskip\multlinegap@\crcr\egroup}
\def\rendmultline@{\iftagin@                                                
 \ifdim\rtwidth@>\z@                                                        
  \hskip\multlinetaggap@\box\mtagbox@                                       
 \else\llap{\vtop{\normalbaselines\null\hbox{\strut@\box\mtagbox@}}}\fi     
 \else\hskip\multlinegap@\fi                                                
 \hfilneg\crcr\egroup}
\def\bmod{\mskip-\medmuskip\mkern5mu\mathbin{\fam\z@ mod}\penalty900
 \mkern5mu\mskip-\medmuskip}
\def\pmod#1{\allowbreak\ifinner\mkern8mu\else\mkern18mu\fi
 ({\fam\z@ mod}\,\,#1)}
\def\pod#1{\allowbreak\ifinner\mkern8mu\else\mkern18mu\fi(#1)}
\def\mod#1{\allowbreak\ifinner\mkern12mu\else\mkern18mu\fi{\fam\z@ mod}\,\,#1}
\message{continued fractions,}
\newcount\cfraccount@
\def\cfrac{\bgroup\bgroup\advance\cfraccount@\@ne\strut
 \iffalse{\fi\def\\{\over\displaystyle}\iffalse}\fi}
\def\lcfrac{\bgroup\bgroup\advance\cfraccount@\@ne\strut
 \iffalse{\fi\def\\{\hfill\over\displaystyle}\iffalse}\fi}
\def\rcfrac{\bgroup\bgroup\advance\cfraccount@\@ne\strut\hfill
 \iffalse{\fi\def\\{\over\displaystyle}\iffalse}\fi}
\def\gloop@#1\repeat{\gdef\body{#1}\iterate}
\def\endcfrac{\gloop@\ifnum\cfraccount@>\z@\global\advance\cfraccount@\m@ne
 \egroup\hskip-\nulldelimiterspace\egroup\repeat}
\message{compound symbols,}
\def\binrel@#1{\setboxz@h{\thinmuskip0mu
  \medmuskip\m@ne mu\thickmuskip\@ne mu$#1\m@th$}%
 \setbox@ne\hbox{\thinmuskip0mu\medmuskip\m@ne mu\thickmuskip
  \@ne mu${}#1{}\m@th$}%
 \setbox\tw@\hbox{\hskip\wd@ne\hskip-\wdz@}}
\def\overset#1\to#2{\binrel@{#2}\ifdim\wd\tw@<\z@
 \mathbin{\mathop{\kern\z@#2}\limits^{#1}}\else\ifdim\wd\tw@>\z@
 \mathrel{\mathop{\kern\z@#2}\limits^{#1}}\else
 {\mathop{\kern\z@#2}\limits^{#1}}{}\fi\fi}
\def\underset#1\to#2{\binrel@{#2}\ifdim\wd\tw@<\z@
 \mathbin{\mathop{\kern\z@#2}\limits_{#1}}\else\ifdim\wd\tw@>\z@
 \mathrel{\mathop{\kern\z@#2}\limits_{#1}}\else
 {\mathop{\kern\z@#2}\limits_{#1}}{}\fi\fi}
\def\oversetbrace#1\to#2{\overbrace{#2}^{#1}}
\def\undersetbrace#1\to#2{\underbrace{#2}_{#1}}
\def\sideset#1\and#2\to#3{%
 \setbox@ne\hbox{$\dsize{\vphantom{#3}}#1{#3}\m@th$}%
 \setbox\tw@\hbox{$\dsize{#3}#2\m@th$}%
 \hskip\wd@ne\hskip-\wd\tw@\mathop{\hskip\wd\tw@\hskip-\wd@ne
  {\vphantom{#3}}#1{#3}#2}}
\def\rightarrowfill@#1{\setboxz@h{$#1-\m@th$}\ht\z@\z@
  $#1\m@th\copy\z@\mkern-6mu\cleaders
  \hbox{$#1\mkern-2mu\box\z@\mkern-2mu$}\hfill
  \mkern-6mu\mathord\rightarrow$}
\def\leftarrowfill@#1{\setboxz@h{$#1-\m@th$}\ht\z@\z@
  $#1\m@th\mathord\leftarrow\mkern-6mu\cleaders
  \hbox{$#1\mkern-2mu\copy\z@\mkern-2mu$}\hfill
  \mkern-6mu\box\z@$}
\def\leftrightarrowfill@#1{\setboxz@h{$#1-\m@th$}\ht\z@\z@
  $#1\m@th\mathord\leftarrow\mkern-6mu\cleaders
  \hbox{$#1\mkern-2mu\box\z@\mkern-2mu$}\hfill
  \mkern-6mu\mathord\rightarrow$}
\def\overrightarrow{\mathpalette\overrightarrow@}
\def\overrightarrow@#1#2{\vbox{\ialign{##\crcr\rightarrowfill@#1\crcr
 \noalign{\kern-\ex@\nointerlineskip}$\m@th\hfil#1#2\hfil$\crcr}}}

\def\overleftarrow{\mathpalette\overleftarrow@}
\def\overleftarrow@#1#2{\vbox{\ialign{##\crcr\leftarrowfill@#1\crcr
 \noalign{\kern-\ex@\nointerlineskip}$\m@th\hfil#1#2\hfil$\crcr}}}
\def\overleftrightarrow{\mathpalette\overleftrightarrow@}
\def\overleftrightarrow@#1#2{\vbox{\ialign{##\crcr\leftrightarrowfill@#1\crcr
 \noalign{\kern-\ex@\nointerlineskip}$\m@th\hfil#1#2\hfil$\crcr}}}
\def\underrightarrow{\mathpalette\underrightarrow@}
\def\underrightarrow@#1#2{\vtop{\ialign{##\crcr$\m@th\hfil#1#2\hfil$\crcr
 \noalign{\nointerlineskip}\rightarrowfill@#1\crcr}}}

\def\underleftarrow{\mathpalette\underleftarrow@}
\def\underleftarrow@#1#2{\vtop{\ialign{##\crcr$\m@th\hfil#1#2\hfil$\crcr
 \noalign{\nointerlineskip}\leftarrowfill@#1\crcr}}}
\def\underleftrightarrow{\mathpalette\underleftrightarrow@}
\def\underleftrightarrow@#1#2{\vtop{\ialign{##\crcr$\m@th\hfil#1#2\hfil$\crcr
 \noalign{\nointerlineskip}\leftrightarrowfill@#1\crcr}}}
\message{various kinds of dots,}
\let\DOTSI\relax
\let\DOTSB\relax

\newif\ifmath@
{\uccode`7=`\\ \uccode`8=`m \uccode`9=`a \uccode`0=`t \uccode`!=`h
 \uppercase{\gdef\math@#1#2#3#4#5#6\math@{\global\math@false\ifx 7#1\ifx 8#2%
 \ifx 9#3\ifx 0#4\ifx !#5\xdef\meaning@{#6}\global\math@true\fi\fi\fi\fi\fi}}}
\newif\ifmathch@
{\uccode`7=`c \uccode`8=`h \uccode`9=`\"
 \uppercase{\gdef\mathch@#1#2#3#4#5#6\mathch@{\global\mathch@false
  \ifx 7#1\ifx 8#2\ifx 9#5\global\mathch@true\xdef\meaning@{9#6}\fi\fi\fi}}}
\newcount\classnum@
\def\getmathch@#1.#2\getmathch@{\classnum@#1 \divide\classnum@4096
 \ifcase\number\classnum@\or\or\gdef\thedots@{\dotsb@}\or
 \gdef\thedots@{\dotsb@}\fi}
\newif\ifmathbin@
{\uccode`4=`b \uccode`5=`i \uccode`6=`n
 \uppercase{\gdef\mathbin@#1#2#3{\relaxnext@
  \DNii@##1\mathbin@{\ifx\space@\next\global\mathbin@true\fi}%
 \global\mathbin@false\DN@##1\mathbin@{}%
 \ifx 4#1\ifx 5#2\ifx 6#3\DN@{\FN@\nextii@}\fi\fi\fi\next@}}}
\newif\ifmathrel@
{\uccode`4=`r \uccode`5=`e \uccode`6=`l
 \uppercase{\gdef\mathrel@#1#2#3{\relaxnext@
  \DNii@##1\mathrel@{\ifx\space@\next\global\mathrel@true\fi}%
 \global\mathrel@false\DN@##1\mathrel@{}%
 \ifx 4#1\ifx 5#2\ifx 6#3\DN@{\FN@\nextii@}\fi\fi\fi\next@}}}
\newif\ifmacro@
{\uccode`5=`m \uccode`6=`a \uccode`7=`c
 \uppercase{\gdef\macro@#1#2#3#4\macro@{\global\macro@false
  \ifx 5#1\ifx 6#2\ifx 7#3\global\macro@true
  \xdef\meaning@{\macro@@#4\macro@@}\fi\fi\fi}}}
\def\macro@@#1->#2\macro@@{#2}
\newif\ifDOTS@
\newcount\DOTSCASE@
{\uccode`6=`\\ \uccode`7=`D \uccode`8=`O \uccode`9=`T \uccode`0=`S
 \uppercase{\gdef\DOTS@#1#2#3#4#5{\global\DOTS@false\DN@##1\DOTS@{}%
  \ifx 6#1\ifx 7#2\ifx 8#3\ifx 9#4\ifx 0#5\let\next@\DOTS@@\fi\fi\fi\fi\fi
  \next@}}}
{\uccode`3=`B \uccode`4=`I \uccode`5=`X
 \uppercase{\gdef\DOTS@@#1{\relaxnext@
  \DNii@##1\DOTS@{\ifx\space@\next\global\DOTS@true\fi}%
  \DN@{\FN@\nextii@}%
  \ifx 3#1\global\DOTSCASE@\z@\else
  \ifx 4#1\global\DOTSCASE@\@ne\else
  \ifx 5#1\global\DOTSCASE@\tw@\else\DN@##1\DOTS@{}%
  \fi\fi\fi\next@}}}
\newif\ifnot@
{\uccode`5=`\\ \uccode`6=`n \uccode`7=`o \uccode`8=`t
 \uppercase{\gdef\not@#1#2#3#4{\relaxnext@
  \DNii@##1\not@{\ifx\space@\next\global\not@true\fi}%
 \global\not@false\DN@##1\not@{}%
 \ifx 5#1\ifx 6#2\ifx 7#3\ifx 8#4\DN@{\FN@\nextii@}\fi\fi\fi
 \fi\next@}}}
\newif\ifkeybin@
\def\keybin@{\keybin@true
 \ifx\next+\else\ifx\next=\else\ifx\next<\else\ifx\next>\else\ifx\next-\else
 \ifx\next*\else\ifx\next:\else\keybin@false\fi\fi\fi\fi\fi\fi\fi}
\def\dots{\RIfM@\expandafter\mdots@\else\expandafter\tdots@\fi}
\def\tdots@{\unskip\relaxnext@
 \DN@{$\m@th\mathinner{\ldotp\ldotp\ldotp}\,
   \ifx\next,\,$\else\ifx\next.\,$\else\ifx\next;\,$\else\ifx\next:\,$\else
   \ifx\next?\,$\else\ifx\next!\,$\else$ \fi\fi\fi\fi\fi\fi}%
 \ \FN@\next@}
\def\mdots@{\FN@\mdots@@}
\def\mdots@@{\gdef\thedots@{\dotso@}
 \ifx\next\boldkey\gdef\thedots@\boldkey{\boldkeydots@}\else                
 \ifx\next\boldsymbol\gdef\thedots@\boldsymbol{\boldsymboldots@}\else       
 \ifx,\next\gdef\thedots@{\dotsc}
 \else\ifx\not\next\gdef\thedots@{\dotsb@}
 \else\keybin@
 \ifkeybin@\gdef\thedots@{\dotsb@}
 \else\xdef\meaning@{\meaning\next..........}\xdef\meaning@@{\meaning@}
  \expandafter\math@\meaning@\math@
  \ifmath@
   \expandafter\mathch@\meaning@\mathch@
   \ifmathch@\expandafter\getmathch@\meaning@\getmathch@\fi                 
  \else\expandafter\macro@\meaning@@\macro@                                 
  \ifmacro@                                                                
   \expandafter\not@\meaning@\not@\ifnot@\gdef\thedots@{\dotsb@}
  \else\expandafter\DOTS@\meaning@\DOTS@
  \ifDOTS@
   \ifcase\number\DOTSCASE@\gdef\thedots@{\dotsb@}%
    \or\gdef\thedots@{\dotsi}\else\fi                                      
  \else\expandafter\math@\meaning@\math@                                   
  \ifmath@\expandafter\mathbin@\meaning@\mathbin@
  \ifmathbin@\gdef\thedots@{\dotsb@}
  \else\expandafter\mathrel@\meaning@\mathrel@
  \ifmathrel@\gdef\thedots@{\dotsb@}
  \fi\fi\fi\fi\fi\fi\fi\fi\fi\fi\fi\fi
 \thedots@}
\def\plainldots@{\mathinner{\ldotp\ldotp\ldotp}}
\def\plaincdots@{\mathinner{\cdotp\cdotp\cdotp}}
\def\dotsi{\!\plaincdots@}
\let\dotsb@\plaincdots@
\newif\ifextra@
\newif\ifrightdelim@
\def\rightdelim@{\global\rightdelim@true                                    
 \ifx\next)\else                                                            
 \ifx\next]\else
 \ifx\next\rbrack\else
 \ifx\next\}\else
 \ifx\next\rbrace\else
 \ifx\next\rangle\else
 \ifx\next\rceil\else
 \ifx\next\rfloor\else
 \ifx\next\rgroup\else
 \ifx\next\rmoustache\else
 \ifx\next\right\else
 \ifx\next\bigr\else
 \ifx\next\biggr\else
 \ifx\next\Bigr\else                                                        
 \ifx\next\Biggr\else\global\rightdelim@false
 \fi\fi\fi\fi\fi\fi\fi\fi\fi\fi\fi\fi\fi\fi\fi}
\def\extra@{%
 \global\extra@false\rightdelim@\ifrightdelim@\global\extra@true            
 \else\ifx\next$\global\extra@true                                          
 \else\xdef\meaning@{\meaning\next..........}
 \expandafter\macro@\meaning@\macro@\ifmacro@                               
 \expandafter\DOTS@\meaning@\DOTS@
 \ifDOTS@
 \ifnum\DOTSCASE@=\tw@\global\extra@true                                    
 \fi\fi\fi\fi\fi}
\newif\ifbold@
\def\dotso@{\relaxnext@
 \ifbold@
  \let\next\delayed@
  \DNii@{\extra@\plainldots@\ifextra@\,\fi}%
 \else
  \DNii@{\DN@{\extra@\plainldots@\ifextra@\,\fi}\FN@\next@}%
 \fi
 \nextii@}
\def\extrap@#1{%
 \ifx\next,\DN@{#1\,}\else
 \ifx\next;\DN@{#1\,}\else
 \ifx\next.\DN@{#1\,}\else\extra@
 \ifextra@\DN@{#1\,}\else
 \let\next@#1\fi\fi\fi\fi\next@}
\def\ldots{\DN@{\extrap@\plainldots@}%
 \FN@\next@}
\def\cdots{\DN@{\extrap@\plaincdots@}%
 \FN@\next@}

\def\dotsc{\relaxnext@
 \DN@{\ifx\next;\plainldots@\,\else
  \ifx\next.\plainldots@\,\else\extra@\plainldots@
  \ifextra@\,\fi\fi\fi}%
 \FN@\next@}
\def\cdot{\mathchar"2201 }
\def\longrightarrow{\DOTSB\relbar\joinrel\rightarrow}

\def\mapsto{\DOTSB\mapstochar\rightarrow}

\def\hookrightarrow{\DOTSB\lhook\joinrel\rightarrow}

\message{special superscripts,}
\def\dddot#1{{\mathop{#1}\limits^{\vbox to-1.4\ex@{\kern-\tw@\ex@
 \hbox{\rm...}\vss}}}}
\def\ddddot#1{{\mathop{#1}\limits^{\vbox to-1.4\ex@{\kern-\tw@\ex@
 \hbox{\rm....}\vss}}}}
\def\sphat{^{\mathchoice{}{}%
 {\,\,\botsmash{\hbox{\lower4\ex@\hbox{$\m@th\widehat{\null}$}}}}%
 {\,\botsmash{\hbox{\lower3\ex@\hbox{$\m@th\hat{\null}$}}}}}}

\def\spacute{^{\!\botsmash{\hbox{\lower\@ne ex\hbox{\'{}}}}}}
\def\spgrave{^{\mathchoice{}{}{}{\!}%
 \botsmash{\hbox{\lower\@ne ex\hbox{\`{}}}}}}
\def\spdot{^{\hbox{\raise\ex@\hbox{\rm.}}}}
\def\spddot{^{\hbox{\raise\ex@\hbox{\rm..}}}}
\def\spdddot{^{\hbox{\raise\ex@\hbox{\rm...}}}}
\def\spddddot{^{\hbox{\raise\ex@\hbox{\rm....}}}}
\def\spbreve{^{\!\botsmash{\hbox{\lower4\ex@\hbox{\u{}}}}}}

\message{\string\text,}
\def\textonlyfont@#1#2{\def#1{\RIfM@
 \Err@{Use \string#1\space only in text}\else#2\fi}}
\textonlyfont@\rm\tenrm
\textonlyfont@\it\tenit
\textonlyfont@\sl\tensl
\textonlyfont@\bf\tenbf
\def\oldnos#1{\RIfM@{\mathcode`\,="013B \fam\@ne#1}\else
 \leavevmode\hbox{$\m@th\mathcode`\,="013B \fam\@ne#1$}\fi}
\def\text{\RIfM@\expandafter\text@\else\expandafter\text@@\fi}
\def\text@@#1{\leavevmode\hbox{#1}}
\def\mathhexbox@#1#2#3{\text{$\m@th\mathchar"#1#2#3$}}
\def\dag{{\mathhexbox@279}}
\def\ddag{{\mathhexbox@27A}}
\def\S{{\mathhexbox@278}}
\def\P{{\mathhexbox@27B}}
\newif\iffirstchoice@
\firstchoice@true
\def\text@#1{\mathchoice
 {\hbox{\everymath{\displaystyle}\def\textfonti{\the\textfont\@ne}%
  \def\textfontii{\the\textfont\tw@}\textdef@@ T#1}}
 {\hbox{\firstchoice@false
  \everymath{\textstyle}\def\textfonti{\the\textfont\@ne}%
  \def\textfontii{\the\textfont\tw@}\textdef@@ T#1}}
 {\hbox{\firstchoice@false
  \everymath{\scriptstyle}\def\textfonti{\the\scriptfont\@ne}%
  \def\textfontii{\the\scriptfont\tw@}\textdef@@ S\rm#1}}
 {\hbox{\firstchoice@false
  \everymath{\scriptscriptstyle}\def\textfonti
  {\the\scriptscriptfont\@ne}%
  \def\textfontii{\the\scriptscriptfont\tw@}\textdef@@ s\rm#1}}}
\def\textdef@@#1{\textdef@#1\rm\textdef@#1\bf\textdef@#1\sl\textdef@#1\it}
\def\rmfam{0}
\def\textdef@#1#2{%
 \DN@{\csname\expandafter\eat@\string#2fam\endcsname}%
 \if S#1\edef#2{\the\scriptfont\next@\relax}%
 \else\if s#1\edef#2{\the\scriptscriptfont\next@\relax}%
 \else\edef#2{\the\textfont\next@\relax}\fi\fi}
\scriptfont\itfam\tenit \scriptscriptfont\itfam\tenit
\scriptfont\slfam\tensl \scriptscriptfont\slfam\tensl
\newif\iftopfolded@
\newif\ifbotfolded@
\def\topfoldedtext{\topfolded@true\botfolded@false\foldedtext@}
\def\botfoldedtext{\botfolded@true\topfolded@false\foldedtext@}
\def\foldedtext{\topfolded@false\botfolded@false\foldedtext@}
\Invalid@\foldedwidth
\def\foldedtext@{\relaxnext@
 \DN@{\ifx\next\foldedwidth\let\next@\nextii@\else
  \DN@{\nextii@\foldedwidth{.3\hsize}}\fi\next@}%
 \DNii@\foldedwidth##1##2{\setbox\z@\vbox
  {\normalbaselines\hsize##1\relax
  \tolerance1600 \noindent\ignorespaces##2}\ifbotfolded@\boxz@\else
  \iftopfolded@\vtop{\unvbox\z@}\else\vcenter{\boxz@}\fi\fi}%
 \FN@\next@}
\message{math font commands,}
\def\bold{\RIfM@\expandafter\bold@\else
 \expandafter\nonmatherr@\expandafter\bold\fi}
\def\bold@#1{{\bold@@{#1}}}
\def\bold@@#1{\fam\bffam\relax#1}
\def\slanted{\RIfM@\expandafter\slanted@\else
 \expandafter\nonmatherr@\expandafter\slanted\fi}
\def\slanted@#1{{\slanted@@{#1}}}
\def\slanted@@#1{\fam\slfam\relax#1}
\def\roman{\RIfM@\expandafter\roman@\else
 \expandafter\nonmatherr@\expandafter\roman\fi}
\def\roman@#1{{\roman@@{#1}}}
\def\roman@@#1{\fam\rmfam\relax#1}
\def\italic{\RIfM@\expandafter\italic@\else
 \expandafter\nonmatherr@\expandafter\italic\fi}
\def\italic@#1{{\italic@@{#1}}}
\def\italic@@#1{\fam\itfam\relax#1}
\def\Cal{\RIfM@\expandafter\Cal@\else
 \expandafter\nonmatherr@\expandafter\Cal\fi}
\def\Cal@#1{{\Cal@@{#1}}}
\def\Cal@@#1{\noaccents@\fam\tw@#1}
\mathchardef\Gamma="0000
\mathchardef\Delta="0001
\mathchardef\Theta="0002
\mathchardef\Lambda="0003
\mathchardef\Xi="0004
\mathchardef\Pi="0005
\mathchardef\Sigma="0006
\mathchardef\Upsilon="0007
\mathchardef\Phi="0008
\mathchardef\Psi="0009
\mathchardef\Omega="000A
\mathchardef\varGamma="0100
\mathchardef\varDelta="0101
\mathchardef\varTheta="0102
\mathchardef\varLambda="0103
\mathchardef\varXi="0104
\mathchardef\varPi="0105
\mathchardef\varSigma="0106
\mathchardef\varUpsilon="0107
\mathchardef\varPhi="0108
\mathchardef\varPsi="0109
\mathchardef\varOmega="010A
\let\alloc@@\alloc@
\def\hexnumber@#1{\ifcase#1 0\or 1\or 2\or 3\or 4\or 5\or 6\or 7\or 8\or
 9\or A\or B\or C\or D\or E\or F\fi}
\def\loadmsam{%
 \font@\tenmsa=msam10
 \font@\sevenmsa=msam7
 \font@\fivemsa=msam5
 \alloc@@8\fam\chardef\sixt@@n\msafam
 \textfont\msafam=\tenmsa
 \scriptfont\msafam=\sevenmsa
 \scriptscriptfont\msafam=\fivemsa
 \edef\next{\hexnumber@\msafam}%
 \mathchardef\dabar@"0\next39
 \edef\dashrightarrow{\mathrel{\dabar@\dabar@\mathchar"0\next4B}}%
 \edef\dashleftarrow{\mathrel{\mathchar"0\next4C\dabar@\dabar@}}%
 \let\dasharrow\dashrightarrow
 \edef\ulcorner{\delimiter"4\next70\next70 }%
 \edef\urcorner{\delimiter"5\next71\next71 }%
 \edef\llcorner{\delimiter"4\next78\next78 }%
 \edef\lrcorner{\delimiter"5\next79\next79 }%
 \edef\yen{{\noexpand\mathhexbox@\next55}}%
 \edef\checkmark{{\noexpand\mathhexbox@\next58}}%
 \edef\circledR{{\noexpand\mathhexbox@\next72}}%
 \edef\maltese{{\noexpand\mathhexbox@\next7A}}%
 \global\let\loadmsam\empty}%
\def\loadmsbm{%
 \font@\tenmsb=msbm10 \font@\sevenmsb=msbm7 \font@\fivemsb=msbm5
 \alloc@@8\fam\chardef\sixt@@n\msbfam
 \textfont\msbfam=\tenmsb
 \scriptfont\msbfam=\sevenmsb \scriptscriptfont\msbfam=\fivemsb
 \global\let\loadmsbm\empty
 }
\def\widehat#1{\ifx\undefined\msbfam \DN@{362}%
  \else \setboxz@h{$\m@th#1$}%
    \edef\next@{\ifdim\wdz@>\tw@ em%
        \hexnumber@\msbfam 5B%
      \else 362\fi}\fi
  \mathaccent"0\next@{#1}}
\def\widetilde#1{\ifx\undefined\msbfam \DN@{365}%
  \else \setboxz@h{$\m@th#1$}%
    \edef\next@{\ifdim\wdz@>\tw@ em%
        \hexnumber@\msbfam 5D%
      \else 365\fi}\fi
  \mathaccent"0\next@{#1}}
\message{\string\newsymbol,}
\def\newsymbol#1#2#3#4#5{\define#1{}%
  \count@#2\relax \advance\count@\m@ne 
 \ifcase\count@
   \ifx\undefined\msafam\loadmsam\fi \let\next@\msafam
 \or \ifx\undefined\msbfam\loadmsbm\fi \let\next@\msbfam
 \else  \Err@{\Invalid@@\string\newsymbol}\let\next@\tw@\fi
 \mathchardef#1="#3\hexnumber@\next@#4#5\space}

\def\Bbb{\RIfM@\expandafter\Bbb@\else
 \expandafter\nonmatherr@\expandafter\Bbb\fi}
\def\Bbb@#1{{\Bbb@@{#1}}}
\def\Bbb@@#1{\noaccents@\fam\msbfam\relax#1}
\message{bold Greek and bold symbols,}
\def\loadbold{%
 \font@\tencmmib=cmmib10 \font@\sevencmmib=cmmib7 \font@\fivecmmib=cmmib5
 \skewchar\tencmmib'177 \skewchar\sevencmmib'177 \skewchar\fivecmmib'177
 \alloc@@8\fam\chardef\sixt@@n\cmmibfam
 \textfont\cmmibfam\tencmmib
 \scriptfont\cmmibfam\sevencmmib \scriptscriptfont\cmmibfam\fivecmmib
 \font@\tencmbsy=cmbsy10 \font@\sevencmbsy=cmbsy7 \font@\fivecmbsy=cmbsy5
 \skewchar\tencmbsy'60 \skewchar\sevencmbsy'60 \skewchar\fivecmbsy'60
 \alloc@@8\fam\chardef\sixt@@n\cmbsyfam
 \textfont\cmbsyfam\tencmbsy
 \scriptfont\cmbsyfam\sevencmbsy \scriptscriptfont\cmbsyfam\fivecmbsy
 \let\loadbold\empty
}
\def\boldnotloaded#1{\Err@{\ifcase#1\or First\else Second\fi
       bold symbol font not loaded}}
\def\mathchari@#1#2#3{\ifx\undefined\cmmibfam
    \boldnotloaded@\@ne
  \else\mathchar"#1\hexnumber@\cmmibfam#2#3\space \fi}
\def\mathcharii@#1#2#3{\ifx\undefined\cmbsyfam
    \boldnotloaded\tw@
  \else \mathchar"#1\hexnumber@\cmbsyfam#2#3\space\fi}
\edef\bffam@{\hexnumber@\bffam}
\def\boldkey#1{\ifcat\noexpand#1A%
  \ifx\undefined\cmmibfam \boldnotloaded\@ne
  \else {\fam\cmmibfam#1}\fi
 \else
 \ifx#1!\mathchar"5\bffam@21 \else
 \ifx#1(\mathchar"4\bffam@28 \else\ifx#1)\mathchar"5\bffam@29 \else
 \ifx#1+\mathchar"2\bffam@2B \else\ifx#1:\mathchar"3\bffam@3A \else
 \ifx#1;\mathchar"6\bffam@3B \else\ifx#1=\mathchar"3\bffam@3D \else
 \ifx#1?\mathchar"5\bffam@3F \else\ifx#1[\mathchar"4\bffam@5B \else
 \ifx#1]\mathchar"5\bffam@5D \else
 \ifx#1,\mathchari@63B \else
 \ifx#1-\mathcharii@200 \else
 \ifx#1.\mathchari@03A \else
 \ifx#1/\mathchari@03D \else
 \ifx#1<\mathchari@33C \else
 \ifx#1>\mathchari@33E \else
 \ifx#1*\mathcharii@203 \else
 \ifx#1|\mathcharii@06A \else
 \ifx#10\bold0\else\ifx#11\bold1\else\ifx#12\bold2\else\ifx#13\bold3\else
 \ifx#14\bold4\else\ifx#15\bold5\else\ifx#16\bold6\else\ifx#17\bold7\else
 \ifx#18\bold8\else\ifx#19\bold9\else
  \Err@{\string\boldkey\space can't be used with #1}%
 \fi\fi\fi\fi\fi\fi\fi\fi\fi\fi\fi\fi\fi\fi\fi
 \fi\fi\fi\fi\fi\fi\fi\fi\fi\fi\fi\fi\fi\fi}
\def\boldsymbol#1{%
 \DN@{\Err@{You can't use \string\boldsymbol\space with \string#1}#1}%
 \ifcat\noexpand#1A%
   \let\next@\relax
   \ifx\undefined\cmmibfam \boldnotloaded\@ne
   \else {\fam\cmmibfam#1}\fi
 \else
  \xdef\meaning@{\meaning#1.........}%
  \expandafter\math@\meaning@\math@
  \ifmath@
   \expandafter\mathch@\meaning@\mathch@
   \ifmathch@
    \expandafter\boldsymbol@@\meaning@\boldsymbol@@
   \fi
  \else
   \expandafter\macro@\meaning@\macro@
   \expandafter\delim@\meaning@\delim@
   \ifdelim@
    \expandafter\delim@@\meaning@\delim@@
   \else
    \boldsymbol@{#1}%
   \fi
  \fi
 \fi
 \next@}
\def\mathhexboxii@#1#2{\ifx\undefined\cmbsyfam
    \boldnotloaded\tw@
  \else \mathhexbox@{\hexnumber@\cmbsyfam}{#1}{#2}\fi}
\def\boldsymbol@#1{\let\next@\relax\let\next#1%
 \ifx\next\cdot\mathcharii@201 \else
 \ifx\next\prime{{\null\mathcharii@030 \null}}\else
 \ifx\next\lbrack\mathchar"4\bffam@5B \else
 \ifx\next\rbrack\mathchar"5\bffam@5D \else
 \ifx\next\{\mathcharii@466 \else
 \ifx\next\lbrace\mathcharii@466 \else
 \ifx\next\}\mathcharii@567 \else
 \ifx\next\rbrace\mathcharii@567 \else
 \ifx\next\surd{{\mathcharii@170}}\else
 \ifx\next\S{{\mathhexboxii@78}}\else
 \ifx\next\P{{\mathhexboxii@7B}}\else
 \ifx\next\dag{{\mathhexboxii@79}}\else
 \ifx\next\ddag{{\mathhexboxii@7A}}\else
 \DN@{\Err@{You can't use \string\boldsymbol\space with \string#1}#1}%
 \fi\fi\fi\fi\fi\fi\fi\fi\fi\fi\fi\fi\fi}
\def\boldsymbol@@#1.#2\boldsymbol@@{\classnum@#1 \count@@@\classnum@        
 \divide\classnum@4096 \count@\classnum@                                    
 \multiply\count@4096 \advance\count@@@-\count@ \count@@\count@@@           
 \divide\count@@@\@cclvi \count@\count@@                                    
 \multiply\count@@@\@cclvi \advance\count@@-\count@@@                       
 \divide\count@@@\@cclvi                                                    
 \multiply\classnum@4096 \advance\classnum@\count@@                         
 \ifnum\count@@@=\z@                                                        
  \count@"\bffam@ \multiply\count@\@cclvi
  \advance\classnum@\count@
  \DN@{\mathchar\number\classnum@}%
 \else
  \ifnum\count@@@=\@ne                                                      
   \ifx\undefined\cmmibfam \DN@{\boldnotloaded\@ne}%
   \else \count@\cmmibfam \multiply\count@\@cclvi
     \advance\classnum@\count@
     \DN@{\mathchar\number\classnum@}\fi
  \else
   \ifnum\count@@@=\tw@                                                    
     \ifx\undefined\cmbsyfam
       \DN@{\boldnotloaded\tw@}%
     \else
       \count@\cmbsyfam \multiply\count@\@cclvi
       \advance\classnum@\count@
       \DN@{\mathchar\number\classnum@}%
     \fi
  \fi
 \fi
\fi}
\newif\ifdelim@
\newcount\delimcount@
{\uccode`6=`\\ \uccode`7=`d \uccode`8=`e \uccode`9=`l
 \uppercase{\gdef\delim@#1#2#3#4#5\delim@
  {\delim@false\ifx 6#1\ifx 7#2\ifx 8#3\ifx 9#4\delim@true
   \xdef\meaning@{#5}\fi\fi\fi\fi}}}
\def\delim@@#1"#2#3#4#5#6\delim@@{\if#32%
\let\next@\relax
 \ifx\undefined\cmbsyfam \boldnotloaded\@ne
 \else \mathcharii@#2#4#5\space \fi\fi}
\def\vert{\delimiter"026A30C }
\def\Vert{\delimiter"026B30D }
\let\|\Vert
\def\backslash{\delimiter"026E30F }
\def\boldkeydots@#1{\bold@true\let\next=#1\let\delayed@=#1\mdots@@
 \boldkey#1\bold@false}  
\def\boldsymboldots@#1{\bold@true\let\next#1\let\delayed@#1\mdots@@
 \boldsymbol#1\bold@false}
\message{Euler fonts,}

\def\frak{\mathfont@\frak}

\def\loadmathfont#1{%
   \expandafter\font@\csname ten#1\endcsname=#110
   \expandafter\font@\csname seven#1\endcsname=#17
   \expandafter\font@\csname five#1\endcsname=#15
   \edef\next{\noexpand\alloc@@8\fam\chardef\sixt@@n
     \expandafter\noexpand\csname#1fam\endcsname}%
   \next
   \textfont\csname#1fam\endcsname \csname ten#1\endcsname
   \scriptfont\csname#1fam\endcsname \csname seven#1\endcsname
   \scriptscriptfont\csname#1fam\endcsname \csname five#1\endcsname
   \expandafter\def\csname #1\expandafter\endcsname\expandafter{%
      \expandafter\mathfont@\csname#1\endcsname}%
 \expandafter\gdef\csname load#1\endcsname{}%
}
\def\mathfont@#1{\RIfM@\expandafter\mathfont@@\expandafter#1\else
  \expandafter\nonmatherr@\expandafter#1\fi}
\def\mathfont@@#1#2{{\mathfont@@@#1{#2}}}
\def\mathfont@@@#1#2{\noaccents@
   \fam\csname\expandafter\eat@\string#1fam\endcsname
   \relax#2}
\message{math accents,}
\def\accentclass@{7}
\def\noaccents@{\def\accentclass@{0}}
\def\makeacc@#1#2{\def#1{\mathaccent"\accentclass@#2 }}
\makeacc@\hat{05E}
\makeacc@\check{014}
\makeacc@\tilde{07E}
\makeacc@\acute{013}
\makeacc@\grave{012}
\makeacc@\dot{05F}
\makeacc@\ddot{07F}
\makeacc@\breve{015}
\makeacc@\bar{016}

\newcount\skewcharcount@
\newcount\familycount@
\def\theskewchar@{\familycount@\@ne
 \global\skewcharcount@\the\skewchar\textfont\@ne                           
 \ifnum\fam>\m@ne\ifnum\fam<16
  \global\familycount@\the\fam\relax
  \global\skewcharcount@\the\skewchar\textfont\the\fam\relax\fi\fi          
 \ifnum\skewcharcount@>\m@ne
  \ifnum\skewcharcount@<128
  \multiply\familycount@256
  \global\advance\skewcharcount@\familycount@
  \global\advance\skewcharcount@28672
  \mathchar\skewcharcount@\else
  \global\skewcharcount@\m@ne\fi\else
 \global\skewcharcount@\m@ne\fi}                                            
\newcount\pointcount@
\def\getpoints@#1.#2\getpoints@{\pointcount@#1 }
\newdimen\accentdimen@
\newcount\accentmu@
\def\dimentomu@{\multiply\accentdimen@ 100
 \expandafter\getpoints@\the\accentdimen@\getpoints@
 \multiply\pointcount@18
 \divide\pointcount@\@m
 \global\accentmu@\pointcount@}
\def\Makeacc@#1#2{\def#1{\RIfM@\DN@{\mathaccent@
 {"\accentclass@#2 }}\else\DN@{\nonmatherr@{#1}}\fi\next@}}
\def\unbracefonts@{\let\Cal@\Cal@@\let\roman@\roman@@\let\bold@\bold@@
 \let\slanted@\slanted@@}
\def\mathaccent@#1#2{\ifnum\fam=\m@ne\xdef\thefam@{1}\else
 \xdef\thefam@{\the\fam}\fi                                                 
 \accentdimen@\z@                                                           
 \setboxz@h{\unbracefonts@$\m@th\fam\thefam@\relax#2$}
 \ifdim\accentdimen@=\z@\DN@{\mathaccent#1{#2}}
  \setbox@ne\hbox{\unbracefonts@$\m@th\fam\thefam@\relax#2\theskewchar@$}
  \setbox\tw@\hbox{$\m@th\ifnum\skewcharcount@=\m@ne\else
   \mathchar\skewcharcount@\fi$}
  \global\accentdimen@\wd@ne\global\advance\accentdimen@-\wdz@
  \global\advance\accentdimen@-\wd\tw@                                     
  \global\multiply\accentdimen@\tw@
  \dimentomu@\global\advance\accentmu@\@ne                                 
 \else\DN@{{\mathaccent#1{#2\mkern\accentmu@ mu}%
    \mkern-\accentmu@ mu}{}}\fi                                             
 \next@}\Makeacc@\Hat{05E}
\Makeacc@\Check{014}
\Makeacc@\Tilde{07E}
\Makeacc@\Acute{013}
\Makeacc@\Grave{012}
\Makeacc@\Dot{05F}
\Makeacc@\Ddot{07F}
\Makeacc@\Breve{015}
\Makeacc@\Bar{016}
\def\Vec{\RIfM@\DN@{\mathaccent@{"017E }}\else
 \DN@{\nonmatherr@\Vec}\fi\next@}
\def\accentedsymbol#1#2{\csname newbox\expandafter\endcsname
  \csname\expandafter\eat@\string#1@box\endcsname
 \expandafter\setbox\csname\expandafter\eat@
  \string#1@box\endcsname\hbox{$\m@th#2$}\define
  #1{\copy\csname\expandafter\eat@\string#1@box\endcsname{}}}
\message{roots,}
\def\sqrt#1{\radical"270370 {#1}}
\let\underline@\underline
\let\overline@\overline
\def\underline#1{\underline@{#1}}
\def\overline#1{\overline@{#1}}
\Invalid@\leftroot
\Invalid@\uproot
\newcount\uproot@
\newcount\leftroot@
\def\root{\relaxnext@
  \DN@{\ifx\next\uproot\let\next@\nextii@\else
   \ifx\next\leftroot\let\next@\nextiii@\else
   \let\next@\plainroot@\fi\fi\next@}%
  \DNii@\uproot##1{\uproot@##1\relax\FN@\nextiv@}%
  \def\nextiv@{\ifx\next\space@\DN@. {\FN@\nextv@}\else
   \DN@.{\FN@\nextv@}\fi\next@.}%
  \def\nextv@{\ifx\next\leftroot\let\next@\nextvi@\else
   \let\next@\plainroot@\fi\next@}%
  \def\nextvi@\leftroot##1{\leftroot@##1\relax\plainroot@}%
   \def\nextiii@\leftroot##1{\leftroot@##1\relax\FN@\nextvii@}%
  \def\nextvii@{\ifx\next\space@
   \DN@. {\FN@\nextviii@}\else
   \DN@.{\FN@\nextviii@}\fi\next@.}%
  \def\nextviii@{\ifx\next\uproot\let\next@\nextix@\else
   \let\next@\plainroot@\fi\next@}%
  \def\nextix@\uproot##1{\uproot@##1\relax\plainroot@}%
  \bgroup\uproot@\z@\leftroot@\z@\FN@\next@}
\def\plainroot@#1\of#2{\setbox\rootbox\hbox{$\m@th\scriptscriptstyle{#1}$}%
 \mathchoice{\r@@t\displaystyle{#2}}{\r@@t\textstyle{#2}}
 {\r@@t\scriptstyle{#2}}{\r@@t\scriptscriptstyle{#2}}\egroup}
\def\r@@t#1#2{\setboxz@h{$\m@th#1\sqrt{#2}$}%
 \dimen@\ht\z@\advance\dimen@-\dp\z@
 \setbox@ne\hbox{$\m@th#1\mskip\uproot@ mu$}\advance\dimen@ 1.667\wd@ne
 \mkern-\leftroot@ mu\mkern5mu\raise.6\dimen@\copy\rootbox
 \mkern-10mu\mkern\leftroot@ mu\boxz@}
\def\boxed#1{\setboxz@h{$\m@th\displaystyle{#1}$}\dimen@.4\ex@
 \advance\dimen@3\ex@\advance\dimen@\dp\z@
 \hbox{\lower\dimen@\hbox{%
 \vbox{\hrule height.4\ex@
 \hbox{\vrule width.4\ex@\hskip3\ex@\vbox{\vskip3\ex@\boxz@\vskip3\ex@}%
 \hskip3\ex@\vrule width.4\ex@}\hrule height.4\ex@}%
 }}}
\message{commutative diagrams,}
\let\ampersand@\relax
\newdimen\minaw@
\minaw@11.11128\ex@
\newdimen\minCDaw@
\minCDaw@2.5pc
\def\minCDarrowwidth#1{\RIfMIfI@\onlydmatherr@\minCDarrowwidth
 \else\minCDaw@#1\relax\fi\else\onlydmatherr@\minCDarrowwidth\fi}
\newif\ifCD@
\def\CD{\bgroup\vspace@\relax\let\ampersand@&\iffalse}\fi
 \CD@true\vcenter\bgroup\Let@\tabskip\z@skip\baselineskip20\ex@
 \lineskip3\ex@\lineskiplimit3\ex@\halign\bgroup
 &\hfill$\m@th##$\hfill\crcr}
\def\endCD{\crcr\egroup\egroup\egroup}
\newdimen\bigaw@
\atdef@>#1>#2>{\ampersand@                                                  
 \setboxz@h{$\m@th\ssize\;{#1}\;\;$}
 \setbox@ne\hbox{$\m@th\ssize\;{#2}\;\;$}
 \setbox\tw@\hbox{$\m@th#2$}
 \ifCD@\global\bigaw@\minCDaw@\else\global\bigaw@\minaw@\fi                 
 \ifdim\wdz@>\bigaw@\global\bigaw@\wdz@\fi
 \ifdim\wd@ne>\bigaw@\global\bigaw@\wd@ne\fi                                
 \ifCD@\enskip\fi                                                           
 \ifdim\wd\tw@>\z@
  \mathrel{\mathop{\hbox to\bigaw@{\rightarrowfill@\displaystyle}}%
    \limits^{#1}_{#2}}
 \else\mathrel{\mathop{\hbox to\bigaw@{\rightarrowfill@\displaystyle}}%
    \limits^{#1}}\fi                                                        
 \ifCD@\enskip\fi                                                          
 \ampersand@}                                                              
\atdef@<#1<#2<{\ampersand@\setboxz@h{$\m@th\ssize\;\;{#1}\;$}%
 \setbox@ne\hbox{$\m@th\ssize\;\;{#2}\;$}\setbox\tw@\hbox{$\m@th#2$}%
 \ifCD@\global\bigaw@\minCDaw@\else\global\bigaw@\minaw@\fi
 \ifdim\wdz@>\bigaw@\global\bigaw@\wdz@\fi
 \ifdim\wd@ne>\bigaw@\global\bigaw@\wd@ne\fi
 \ifCD@\enskip\fi
 \ifdim\wd\tw@>\z@
  \mathrel{\mathop{\hbox to\bigaw@{\leftarrowfill@\displaystyle}}%
       \limits^{#1}_{#2}}\else
  \mathrel{\mathop{\hbox to\bigaw@{\leftarrowfill@\displaystyle}}%
       \limits^{#1}}\fi
 \ifCD@\enskip\fi\ampersand@}
\begingroup
 \catcode`\~=\active \lccode`\~=`\@
 \lowercase{%
  \global\atdef@)#1)#2){~>#1>#2>}
  \global\atdef@(#1(#2({~<#1<#2<}}
\endgroup
\atdef@ A#1A#2A{\llap{$\m@th\vcenter{\hbox
 {$\ssize#1$}}$}\Big\uparrow\rlap{$\m@th\vcenter{\hbox{$\ssize#2$}}$}&&}
\atdef@ V#1V#2V{\llap{$\m@th\vcenter{\hbox
 {$\ssize#1$}}$}\Big\downarrow\rlap{$\m@th\vcenter{\hbox{$\ssize#2$}}$}&&}
\atdef@={&\enskip\mathrel
 {\vbox{\hrule width\minCDaw@\vskip3\ex@\hrule width
 \minCDaw@}}\enskip&}
\atdef@|{\Big\Vert&&}
\atdef@\vert{\Big\Vert&&}
\def\pretend#1\haswidth#2{\setboxz@h{$\m@th\scriptstyle{#2}$}\hbox
 to\wdz@{\hfill$\m@th\scriptstyle{#1}$\hfill}}
\message{poor man's bold,}
\def\pmb{\RIfM@\expandafter\mathpalette\expandafter\pmb@\else
 \expandafter\pmb@@\fi}
\def\pmb@@#1{\leavevmode\setboxz@h{#1}%
   \dimen@-\wdz@
   \kern-.5\ex@\copy\z@
   \kern\dimen@\kern.25\ex@\raise.4\ex@\copy\z@
   \kern\dimen@\kern.25\ex@\box\z@
}
\def\binrel@@#1{\ifdim\wd2<\z@\mathbin{#1}\else\ifdim\wd\tw@>\z@
 \mathrel{#1}\else{#1}\fi\fi}
\newdimen\pmbraise@
\def\pmb@#1#2{\setbox\thr@@\hbox{$\m@th#1{#2}$}%
 \setbox4\hbox{$\m@th#1\mkern.5mu$}\pmbraise@\wd4\relax
 \binrel@{#2}%
 \dimen@-\wd\thr@@
   \binrel@@{%
   \mkern-.8mu\copy\thr@@
   \kern\dimen@\mkern.4mu\raise\pmbraise@\copy\thr@@
   \kern\dimen@\mkern.4mu\box\thr@@
}}
\def\documentstyle#1{\W@{}\input #1.sty\relax}
\message{syntax check,}
\font\dummyft@=dummy
\fontdimen1 \dummyft@=\z@
\fontdimen2 \dummyft@=\z@
\fontdimen3 \dummyft@=\z@
\fontdimen4 \dummyft@=\z@
\fontdimen5 \dummyft@=\z@
\fontdimen6 \dummyft@=\z@
\fontdimen7 \dummyft@=\z@
\fontdimen8 \dummyft@=\z@
\fontdimen9 \dummyft@=\z@
\fontdimen10 \dummyft@=\z@
\fontdimen11 \dummyft@=\z@
\fontdimen12 \dummyft@=\z@
\fontdimen13 \dummyft@=\z@
\fontdimen14 \dummyft@=\z@
\fontdimen15 \dummyft@=\z@
\fontdimen16 \dummyft@=\z@
\fontdimen17 \dummyft@=\z@
\fontdimen18 \dummyft@=\z@
\fontdimen19 \dummyft@=\z@
\fontdimen20 \dummyft@=\z@
\fontdimen21 \dummyft@=\z@
\fontdimen22 \dummyft@=\z@
\def\fontlist@{\\{\tenrm}\\{\sevenrm}\\{\fiverm}\\{\teni}\\{\seveni}%
 \\{\fivei}\\{\tensy}\\{\sevensy}\\{\fivesy}\\{\tenex}\\{\tenbf}\\{\sevenbf}%
 \\{\fivebf}\\{\tensl}\\{\tenit}}
\def\font@#1=#2 {\rightappend@#1\to\fontlist@\font#1=#2 }
\def\dodummy@{{\def\\##1{\global\let##1\dummyft@}\fontlist@}}
\def\nopages@{\output{\setbox\z@\box\@cclv \deadcycles\z@}%
 \alloc@5\toks\toksdef\@cclvi\output}
\let\galleys\nopages@
\newif\ifsyntax@
\newcount\countxviii@
\def\syntax{\syntax@true\dodummy@\countxviii@\count18
 \loop\ifnum\countxviii@>\m@ne\textfont\countxviii@=\dummyft@
 \scriptfont\countxviii@=\dummyft@\scriptscriptfont\countxviii@=\dummyft@
 \advance\countxviii@\m@ne\repeat                                           
 \dummyft@\tracinglostchars\z@\nopages@\frenchspacing\hbadness\@M}
\def\first@#1#2\end{#1}
\def\printoptions{\W@{Do you want S(yntax check),
  G(alleys) or P(ages)?}%
 \message{Type S, G or P, followed by <return>: }%
 \begingroup 
 \endlinechar\m@ne 
 \read\m@ne to\ans@
 \edef\ans@{\uppercase{\def\noexpand\ans@{%
   \expandafter\first@\ans@ P\end}}}%
 \expandafter\endgroup\ans@
 \if\ans@ P
 \else \if\ans@ S\syntax
 \else \if\ans@ G\galleys
 \else\message{? Unknown option: \ans@; using the `pages' option.}%
 \fi\fi\fi}
\def\alloc@#1#2#3#4#5{\global\advance\count1#1by\@ne
 \ch@ck#1#4#2\allocationnumber=\count1#1
 \global#3#5=\allocationnumber
 \ifalloc@\wlog{\string#5=\string#2\the\allocationnumber}\fi}
\def\document{\def\alloclist@{}\def\fontlist@{}}
\let\enddocument\bye

\let\proclaim\undefined
\let\footnote\undefined
\let\=\undefined
\let\>\undefined

\catcode`\@=\active
\message{... finished}

\expandafter\ifx\csname mathdefs.tex\endcsname\relax
  \expandafter\gdef\csname mathdefs.tex\endcsname{}
\else \message{Hey!  Apparently you were trying to
  \string\input{mathdefs.tex} twice.   This does not make sense.} 
\errmessage{Please edit your file (probably \jobname.tex) and remove
any duplicate ``\string\input'' lines}\endinput\fi




\catcode`\X=12\catcode`\@=11

\def\n@wcount{\alloc@0\count\countdef\insc@unt}
\def\n@wwrite{\alloc@7\write\chardef\sixt@@n}
\def\n@wread{\alloc@6\read\chardef\sixt@@n}
\def\r@s@t{\relax}\def\v@idline{\par}\def\@mputate#1/{#1}
\def\l@c@l#1X{\firstpart.#1}\def\gl@b@l#1X{#1}\def\t@d@l#1X{{}}

\def\crossrefs#1{\ifx\all#1\let\tr@ce=\all\else\def\tr@ce{#1,}\fi
   \n@wwrite\cit@tionsout\openout\cit@tionsout=\jobname.cit 
   \write\cit@tionsout{\tr@ce}\expandafter\setfl@gs\tr@ce,}
\def\setfl@gs#1,{\def\@{#1}\ifx\@\empty\let\next=\relax
   \else\let\next=\setfl@gs\expandafter\xdef
   \csname#1tr@cetrue\endcsname{}\fi\next}
\def\m@ketag#1#2{\expandafter\n@wcount\csname#2tagno\endcsname
     \csname#2tagno\endcsname=0\let\tail=\all\xdef\all{\tail#2,}
   \ifx#1\l@c@l\let\tail=\r@s@t\xdef\r@s@t{\csname#2tagno\endcsname=0\tail}\fi
   \expandafter\gdef\csname#2cite\endcsname##1{\expandafter
     \ifx\csname#2tag##1\endcsname\relax?\else\csname#2tag##1\endcsname\fi
     \expandafter\ifx\csname#2tr@cetrue\endcsname\relax\else
     \write\cit@tionsout{#2tag ##1 cited on page \folio.}\fi}
   \expandafter\gdef\csname#2page\endcsname##1{\expandafter
     \ifx\csname#2page##1\endcsname\relax?\else\csname#2page##1\endcsname\fi
     \expandafter\ifx\csname#2tr@cetrue\endcsname\relax\else
     \write\cit@tionsout{#2tag ##1 cited on page \folio.}\fi}
   \expandafter\gdef\csname#2tag\endcsname##1{\expandafter
      \ifx\csname#2check##1\endcsname\relax
      \expandafter\xdef\csname#2check##1\endcsname{}%
      \else\immediate\write16{Warning: #2tag ##1 used more than once.}\fi
      \multit@g{#1}{#2}##1/X%
      \write\t@gsout{#2tag ##1 assigned number \csname#2tag##1\endcsname\space
      on page \number\count0.}%
   \csname#2tag##1\endcsname}}

\def\multit@g#1#2#3/#4X{\def\t@mp{#4}\ifx\t@mp\empty%
      \global\advance\csname#2tagno\endcsname by 1 
      \expandafter\xdef\csname#2tag#3\endcsname
      {#1\number\csname#2tagno\endcsnameX}%
   \else\expandafter\ifx\csname#2last#3\endcsname\relax
      \expandafter\n@wcount\csname#2last#3\endcsname
      \global\advance\csname#2tagno\endcsname by 1 
      \expandafter\xdef\csname#2tag#3\endcsname
      {#1\number\csname#2tagno\endcsnameX}
      \write\t@gsout{#2tag #3 assigned number \csname#2tag#3\endcsname\space
      on page \number\count0.}\fi
   \global\advance\csname#2last#3\endcsname by 1
   \def\t@mp{\expandafter\xdef\csname#2tag#3/}%
   \expandafter\t@mp\@mputate#4\endcsname
   {\csname#2tag#3\endcsname\lastpart{\csname#2last#3\endcsname}}\fi}
\def\t@gs#1{\def\all{}\m@ketag#1e\m@ketag#1s\m@ketag\t@d@l p
\let\realscite\scite
\let\realstag\stag
   \m@ketag\gl@b@l r \n@wread\t@gsin
   \openin\t@gsin=\jobname.tgs \re@der \closein\t@gsin
   \n@wwrite\t@gsout\openout\t@gsout=\jobname.tgs }
\outer\def\localtags{\t@gs\l@c@l}
\outer\def\globaltags{\t@gs\gl@b@l}
\outer\def\newlocaltag#1{\m@ketag\l@c@l{#1}}
\outer\def\newglobaltag#1{\m@ketag\gl@b@l{#1}}

\newif\ifpr@ 
\def\m@kecs #1tag #2 assigned number #3 on page #4.%
   {\expandafter\gdef\csname#1tag#2\endcsname{#3}
   \expandafter\gdef\csname#1page#2\endcsname{#4}
   \ifpr@\expandafter\xdef\csname#1check#2\endcsname{}\fi}
\def\re@der{\ifeof\t@gsin\let\next=\relax\else
   \read\t@gsin to\t@gline\ifx\t@gline\v@idline\else
   \expandafter\m@kecs \t@gline\fi\let \next=\re@der\fi\next}
\def\pretags#1{\pr@true\pret@gs#1,,}
\def\pret@gs#1,{\def\@{#1}\ifx\@\empty\let\n@xtfile=\relax
   \else\let\n@xtfile=\pret@gs \openin\t@gsin=#1.tgs \message{#1} \re@der 
   \closein\t@gsin\fi \n@xtfile}

\newcount\sectno\sectno=0\newcount\subsectno\subsectno=0
\newif\ifultr@local \def\ultralocal{\ultr@localtrue}
\def\firstpart{\number\sectno}
\def\lastpart#1{\ifcase#1 \or a\or b\or c\or d\or e\or f\or g\or h\or 
   i\or k\or l\or m\or n\or o\or p\or q\or r\or s\or t\or u\or v\or w\or 
   x\or y\or z \fi}

\def\resetall{\global\advance\sectno by 1\subsectno=0
   \gdef\firstpart{\number\sectno}\r@s@t}
\def\resetsub{\global\advance\subsectno by 1
   \gdef\firstpart{\number\sectno.\number\subsectno}\r@s@t}
\def\newsection#1\par{\resetall\vskip0pt plus.3\vsize\penalty-250
   \vskip0pt plus-.3\vsize\bigskip\bigskip
   \message{#1}\leftline{\bf#1}\nobreak\bigskip}
\def\subsection#1\par{\ifultr@local\resetsub\fi
   \vskip0pt plus.2\vsize\penalty-250\vskip0pt plus-.2\vsize
   \bigskip\smallskip\message{#1}\leftline{\bf#1}\nobreak\medskip}


\newdimen\marginshift

\newdimen\margindelta
\newdimen\marginmax
\newdimen\marginmin

\def\margininit{       
\marginmax=3 true cm                  
				      
\margindelta=0.1 true cm              
\marginmin=0.1true cm                 
\marginshift=\marginmin
}    

\def\t@gsjj#1,{\def\@{#1}\ifx\@\empty\let\next=\relax\else\let\next=\t@gsjj
   \def\@@{p}\ifx\@\@@\else
   \expandafter\gdef\csname#1cite\endcsname##1{\citejj{##1}}
   \expandafter\gdef\csname#1page\endcsname##1{?}
   \expandafter\gdef\csname#1tag\endcsname##1{\tagjj{##1}}\fi\fi\next}
\newif\ifshowstuffinmargin
\showstuffinmarginfalse
\def\jjtags{\ifx\shlhetal\relax 
  \else
\ifx\shlhetal\undefinedcontrolseq
\else
\showstuffinmargintrue
\ifx\all\relax\else\expandafter\t@gsjj\all,\fi\fi \fi
}

\def\tagjj#1{\realstag{#1}\mginpar{\zeigen{#1}}}
\def\citejj#1{\rechnen{#1}\mginpar{\zeigen{#1}}}     

\def\rechnen#1{\expandafter\ifx\csname stag#1\endcsname\relax ??\else
                           \csname stag#1\endcsname\fi}

\newdimen\theight

\def\marginfont{\sevenrm}

\def\trymarginbox#1{\setbox0=\hbox{\marginfont\hskip\marginshift #1}%
		\global\marginshift\wd0 
		\global\advance\marginshift\margindelta}

\def \mginpar#1{%
\ifvmode\setbox0\hbox to \hsize{\hfill\rlap{\marginfont\quad#1}}%
\ht0 0cm
\dp0 0cm
\box0\vskip-\baselineskip
\else 
             \vadjust{\trymarginbox{#1}%
		\ifdim\marginshift>\marginmax \global\marginshift\marginmin
			\trymarginbox{#1}%
                \fi
             \theight=\ht0
             \advance\theight by \dp0    \advance\theight by \lineskip
             \kern -\theight \vbox to \theight{\rightline{\rlap{\box0}}%
\vss}}\fi}


\def\t@gsoff#1,{\def\@{#1}\ifx\@\empty\let\next=\relax\else\let\next=\t@gsoff
   \def\@@{p}\ifx\@\@@\else
   \expandafter\gdef\csname#1cite\endcsname##1{\zeigen{##1}}
   \expandafter\gdef\csname#1page\endcsname##1{?}
   \expandafter\gdef\csname#1tag\endcsname##1{\zeigen{##1}}\fi\fi\next}
\def\verbatimtags{\showstuffinmarginfalse
\ifx\all\relax\else\expandafter\t@gsoff\all,\fi}
\def\zeigen#1{\hbox{$\langle$}#1\hbox{$\rangle$}}

\def\margincite#1{\ifshowstuffinmargin\mginpar{\zeigen{#1}}\fi}

\def\(#1){\edef\dot@g{\ifmmode\ifinner(\hbox{\noexpand\etag{#1}})
   \else\noexpand\eqno(\hbox{\noexpand\etag{#1}})\fi
   \else(\noexpand\ecite{#1})\fi}\dot@g}

\newif\ifbr@ck
\def\eat#1{}
\def\[#1]{\br@cktrue[\br@cket#1'X]}
\def\br@cket#1'#2X{\def\temp{#2}\ifx\temp\empty\let\next\eat
   \else\let\next\br@cket\fi
   \ifbr@ck\br@ckfalse\br@ck@t#1,X\else\br@cktrue#1\fi\next#2X}
\def\br@ck@t#1,#2X{\def\temp{#2}\ifx\temp\empty\let\neext\eat
   \else\let\neext\br@ck@t\def\temp{,}\fi
   \def\teemp{#1}\ifx\teemp\empty\else\rcite{#1}\fi\temp\neext#2X}
\def\resetbr@cket{\gdef\[##1]{[\rtag{##1}]}}
\def\references{\resetbr@cket\newsection References\par}

\newtoks\symb@ls\newtoks\s@mb@ls\newtoks\p@gelist\n@wcount\ftn@mber
    \ftn@mber=1\newif\ifftn@mbers\ftn@mbersfalse\newif\ifbyp@ge\byp@gefalse
\def\defm@rk{\ifftn@mbers\n@mberm@rk\else\symb@lm@rk\fi}
\def\n@mberm@rk{\xdef\m@rk{{\the\ftn@mber}}%
    \global\advance\ftn@mber by 1 }
\def\rot@te#1{\let\temp=#1\global#1=\expandafter\r@t@te\the\temp,X}
\def\r@t@te#1,#2X{{#2#1}\xdef\m@rk{{#1}}}
\def\b@@st#1{{$^{#1}$}}\def\str@p#1{#1}
\def\symb@lm@rk{\ifbyp@ge\rot@te\p@gelist\ifnum\expandafter\str@p\m@rk=1 
    \s@mb@ls=\symb@ls\fi\write\f@nsout{\number\count0}\fi \rot@te\s@mb@ls}
\def\byp@ge{\byp@getrue\n@wwrite\f@nsin\openin\f@nsin=\jobname.fns 
    \n@wcount\currentp@ge\currentp@ge=0\p@gelist={0}
    \re@dfns\closein\f@nsin\rot@te\p@gelist
    \n@wread\f@nsout\openout\f@nsout=\jobname.fns }
\def\m@kelist#1X#2{{#1,#2}}
\def\re@dfns{\ifeof\f@nsin\let\next=\relax\else\read\f@nsin to \f@nline
    \ifx\f@nline\v@idline\else\let\t@mplist=\p@gelist
    \ifnum\currentp@ge=\f@nline
    \global\p@gelist=\expandafter\m@kelist\the\t@mplistX0
    \else\currentp@ge=\f@nline
    \global\p@gelist=\expandafter\m@kelist\the\t@mplistX1\fi\fi
    \let\next=\re@dfns\fi\next}
\def\symbols#1{\symb@ls={#1}\s@mb@ls=\symb@ls} 
\def\bigsymbol{\textstyle}
\symbols{\bigsymbol\ast,\dagger,\ddagger,\sharp,\flat,\natural,\star}
\def\ftnumbers{\ftn@mberstrue} \def\ftsymbols{\ftn@mbersfalse}
\def\paginal{\byp@ge} \def\resetftnumbers{\ftn@mber=1}
\def\ftnote#1{\defm@rk\expandafter\expandafter\expandafter\footnote
    \expandafter\b@@st\m@rk{#1}}

\long\def\jump#1\endjump{}
\def\ssum{\mathop{\lower .1em\hbox{$\textstyle\Sigma$}}\nolimits}

\def\qed{\nobreak\kern 1em \vrule height .5em width .5em depth 0em}
\def\newneq{\hbox{\rlap{\hbox to 1\wd9{\hss$=$\hss}}\raise .1em 
   \hbox to 1\wd9{\hss$\scriptscriptstyle/$\hss}}}
\def\subsetne{\setbox9 = \hbox{$\subset$}\mathrel{\hbox{\rlap
   {\lower .4em \newneq}\raise .13em \hbox{$\subset$}}}}
\def\supsetne{\setbox9 = \hbox{$\subset$}\mathrel{\hbox{\rlap
   {\lower .4em \newneq}\raise .13em \hbox{$\supset$}}}}

\def\vbar{\mathchoice{\vrule height6.3ptdepth-.5ptwidth.8pt\kern-.8pt}
   {\vrule height6.3ptdepth-.5ptwidth.8pt\kern-.8pt}
   {\vrule height4.1ptdepth-.35ptwidth.6pt\kern-.6pt}
   {\vrule height3.1ptdepth-.25ptwidth.5pt\kern-.5pt}}
\def\f@dge{\mathchoice{}{}{\mkern.5mu}{\mkern.8mu}}
\def\b@c#1#2{{\rm \mkern#2mu\vbar\mkern-#2mu#1}}
\def\b@b#1{{\rm I\mkern-3.5mu #1}}
\def\b@a#1#2{{\rm #1\mkern-#2mu\f@dge #1}}
\def\bb#1{{\count4=`#1 \advance\count4by-64 \ifcase\count4\or\b@a A{11.5}\or
   \b@b B\or\b@c C{5}\or\b@b D\or\b@b E\or\b@b F \or\b@c G{5}\or\b@b H\or
   \b@b I\or\b@c J{3}\or\b@b K\or\b@b L \or\b@b M\or\b@b N\or\b@c O{5} \or
   \b@b P\or\b@c Q{5}\or\b@b R\or\b@a S{8}\or\b@a T{10.5}\or\b@c U{5}\or
   \b@a V{12}\or\b@a W{16.5}\or\b@a X{11}\or\b@a Y{11.7}\or\b@a Z{7.5}\fi}}

\catcode`\X=11 \catcode`\@=12




\let\thischap\jobname

\def\partof#1{\csname returnthe#1part\endcsname}
\def\chapof#1{\csname returnthe#1chap\endcsname}

\def\setchapter#1,#2,#3.{%
  \expandafter\def\csname returnthe#1part\endcsname{#2}%
  \expandafter\def\csname returnthe#1chap\endcsname{#3}%
}

\setchapter 300a,A,I.
\setchapter 300b,A,II.
\setchapter 300c,A,III.
\setchapter 300d,A,IV.
\setchapter 300e,A,V.
\setchapter 300f,A,VI.
\setchapter 300g,A,VII.
\setchapter   88,B,I.
\setchapter  600,B,II.
\setchapter  705,B,III.

\def\cprefix#1{
\edef\theotherpart{\partof{#1}}\edef\theotherchap{\chapof{#1}}%
\ifx\theotherpart\thispart
   \ifx\theotherchap\thischap 
    \else 
     \theotherchap%
    \fi
   \else 
     \theotherpart.\theotherchap\fi}

\def\sectioncite[#1]#2{%
     \cprefix{#2}#1}

\edef\thispart{\partof{\thischap}}
\edef\thischap{\chapof{\thischap}}


\def\spuriousreset{}


\expandafter\ifx\csname citeadd.tex\endcsname\relax
\expandafter\gdef\csname citeadd.tex\endcsname{}
\else \message{Hey!  Apparently you were trying to
\string\input{citeadd.tex} twice.   This does not make sense.} 
\errmessage{Please edit your file (probably \jobname.tex) and remove
any duplicate ``\string\input'' lines}\endinput\fi

\sectno=-1   
\newbox\noforkbox \newdimen\forklinewidth
\forklinewidth=0.3pt   
\setbox0\hbox{$\textstyle\bigcup$}
\setbox1\hbox to \wd0{\hfil\vrule width \forklinewidth depth \dp0
                        height \ht0 \hfil}
\wd1=0 cm
\setbox\noforkbox\hbox{\box1\box0\relax}
  \def\unionstick{\mathop{\copy\noforkbox}\limits}
\def\nonfork#1#2_#3{#1\unionstick_{\textstyle #3}#2}
\def\nonforkin#1#2_#3^#4{#1\unionstick_{\textstyle #3}^{\textstyle #4}#2}     
%
\setbox0\hbox{$\textstyle\bigcup$}
\setbox1\hbox to \wd0{\hfil{\sl /\/}\hfil}
\setbox2\hbox to \wd0{\hfil\vrule height \ht0 depth \dp0 width
                                \forklinewidth\hfil}
\wd1=0cm
\wd2=0cm
\newbox\doesforkbox
\setbox\doesforkbox\hbox{\box1\box0\relax}
\def\nunionstick{\mathop{\copy\doesforkbox}\limits}

\def\fork#1#2_#3{#1\nunionstick_{\textstyle #3}#2}
\def\forkin#1#2_#3^#4{#1\nunionstick_{\textstyle #3}^{\textstyle #4}#2}     
\localtags
\jjtags
\NoBlackBoxes
\define\mr{\medskip\roster}
\define\sn{\smallskip\noindent}
\define\mn{\medskip\noindent}
\define\lfoot{\lfloor}
\define\rfoot{\rfloor}
\define\vep{\varepsilon}
\define\bn{\bigskip\noindent}
\define\ub{\underbar}
\define\wilog{\text{without loss of generality}}
\define\ermn{\endroster\medskip\noindent}

\define\dbcu{\dsize\bigcup}
\define \nl{\newline}
\define \hra{\hookrightarrow}
\magnification=\magstep 1
\documentstyle{amsppt}

{    
\catcode`@11

\ifx\alicetwothousandloaded@\relax
  \endinput\else\global\let\alicetwothousandloaded@\relax\fi

\gdef\subjclass{\let\savedef@\subjclass
 \def\subjclass##1\endsubjclass{\let\subjclass\savedef@
   \toks@{\def\usualspace{{\rm\enspace}}\eightpoint}%
   \toks@@{##1\unskip.}%
   \edef\thesubjclass@{\the\toks@
     \frills@{{\noexpand\rm2000 {\noexpand\it Mathematics Subject
       Classification}.\noexpand\enspace}}%
     \the\toks@@}}%
  \nofrillscheck\subjclass}
} 


\expandafter\ifx\csname alice2jlem.tex\endcsname\relax
  \expandafter\xdef\csname alice2jlem.tex\endcsname{\the\catcode`@}
\else \message{Hey!  Apparently you were trying to
\string\input{alice2jlem.tex}  twice.   This does not make sense.}
\errmessage{Please edit your file (probably \jobname.tex) and remove
any duplicate ``\string\input'' lines}\endinput\fi

\expandafter\ifx\csname bib4plain.tex\endcsname\relax
  \expandafter\gdef\csname bib4plain.tex\endcsname{}
\else \message{Hey!  Apparently you were trying to \string\input
  bib4plain.tex twice.   This does not make sense.}
\errmessage{Please edit your file (probably \jobname.tex) and remove
any duplicate ``\string\input'' lines}\endinput\fi

\def\renewcommand{\newcommand}	       
\edef\cite{\the\catcode`@}%
\catcode`@ = 11
\let\@oldatcatcode = \cite
\chardef\@letter = 11
\chardef\@other = 12
%
%
%
%
\def\@innerdef#1#2{\edef#1{\expandafter\noexpand\csname #2\endcsname}}%
%
%
\@innerdef\@innernewcount{newcount}%
\@innerdef\@innernewdimen{newdimen}%
\@innerdef\@innernewif{newif}%
\@innerdef\@innernewwrite{newwrite}%
%
%
%
\def\@gobble#1{}%
%
%
%
\ifx\inputlineno\@undefined
   \let\@linenumber = \empty 
\else
   \def\@linenumber{\the\inputlineno:\space}%
\fi
%
%
%
\def\@futurenonspacelet#1{\def\cs{#1}%
   \afterassignment\@stepone\let\@nexttoken=
}%
\begingroup 
\def\\{\global\let\@stoken= }%
\\ 
\endgroup
\def\@stepone{\expandafter\futurelet\cs\@steptwo}%
\def\@steptwo{\expandafter\ifx\cs\@stoken\let\@@next=\@stepthree
   \else\let\@@next=\@nexttoken\fi \@@next}%
\def\@stepthree{\afterassignment\@stepone\let\@@next= }%
%
%
%
\def\@getoptionalarg#1{%
   \let\@optionaltemp = #1%
   \let\@optionalnext = \relax
   \@futurenonspacelet\@optionalnext\@bracketcheck
}%
%
%
\def\@bracketcheck{%
   \ifx [\@optionalnext
      \expandafter\@@getoptionalarg
   \else
      \let\@optionalarg = \empty
      \expandafter\@optionaltemp
   \fi
}%
\def\@@getoptionalarg[#1]{%
   \def\@optionalarg{#1}%
   \@optionaltemp
}%
%
%
%
\def\@nnil{\@nil}%
\def\@fornoop#1\@@#2#3{}%
\def\@for#1:=#2\do#3{%
   \edef\@fortmp{#2}%
   \ifx\@fortmp\empty \else
      \expandafter\@forloop#2,\@nil,\@nil\@@#1{#3}%
   \fi
}%
\def\@forloop#1,#2,#3\@@#4#5{\def#4{#1}\ifx #4\@nnil \else
       #5\def#4{#2}\ifx #4\@nnil \else#5\@iforloop #3\@@#4{#5}\fi\fi
}%
\def\@iforloop#1,#2\@@#3#4{\def#3{#1}\ifx #3\@nnil
       \let\@nextwhile=\@fornoop \else
      #4\relax\let\@nextwhile=\@iforloop\fi\@nextwhile#2\@@#3{#4}%
}%
%
%
%
\@innernewif\if@fileexists
\def\@testfileexistence{\@getoptionalarg\@finishtestfileexistence}%
\def\@finishtestfileexistence#1{%
   \begingroup
      \def\extension{#1}%
      \immediate\openin0 =
         \ifx\@optionalarg\empty\jobname\else\@optionalarg\fi
         \ifx\extension\empty \else .#1\fi
         \space
      \ifeof 0
         \global\@fileexistsfalse
      \else
         \global\@fileexiststrue
      \fi
      \immediate\closein0
   \endgroup
}%
%
%
%
%
\def\bibliographystyle#1{%
   \@readauxfile
   \@writeaux{\string\bibstyle{#1}}%
}%
\let\bibstyle = \@gobble
%
%
\let\bblfilebasename = \jobname
\def\bibliography#1{%
   \@readauxfile
   \@writeaux{\string\bibdata{#1}}%
   \@testfileexistence[\bblfilebasename]{bbl}%
   \if@fileexists
      \nobreak
      \@readbblfile
   \fi
}%
\let\bibdata = \@gobble
%
%
\def\nocite#1{%
   \@readauxfile
   \@writeaux{\string\citation{#1}}%
}%
\@innernewif\if@notfirstcitation
%
%
\def\cite{\@getoptionalarg\@cite}%
%
%
\def\@cite#1{%
   \let\@citenotetext = \@optionalarg
   \printcitestart
   \nocite{#1}%
   \@notfirstcitationfalse
   \@for \@citation :=#1\do
   {%
      \expandafter\@onecitation\@citation\@@
   }%
   \ifx\empty\@citenotetext\else
      \printcitenote{\@citenotetext}%
   \fi
   \printcitefinish
}%
\newif\ifweareinprivate
\weareinprivatetrue
\ifx\shlhetal\undefinedcontrolseq\weareinprivatefalse\fi
\ifx\shlhetal\relax\weareinprivatefalse\fi
\def\@onecitation#1\@@{%
   \if@notfirstcitation
      \printbetweencitations
   \fi
   \expandafter \ifx \csname\@citelabel{#1}\endcsname \relax
      \if@citewarning
         \message{\@linenumber Undefined citation `#1'.}%
      \fi
     \ifweareinprivate
      \expandafter\gdef\csname\@citelabel{#1}\endcsname{%
\strut 
\vadjust{\vskip-\dp\strutbox
\vbox to 0pt{\vss\parindent0cm \leftskip=\hsize 
\advance\leftskip3mm
\advance\hsize 4cm\strut\openup-4pt 
\rightskip 0cm plus 1cm minus 0.5cm ?  #1 ?\strut}}
         {\tt
            \escapechar = -1
            \nobreak\hskip0pt\pfeilsw
            \expandafter\string\csname#1\endcsname
             \pfeilso
            \nobreak\hskip0pt
         }%
      }%
     \else  
      \expandafter\gdef\csname\@citelabel{#1}\endcsname{%
            {\tt\expandafter\string\csname#1\endcsname}
      }%
     \fi  
   \fi
   \csname\@citelabel{#1}\endcsname
   \@notfirstcitationtrue
}%
%
%
\def\@citelabel#1{b@#1}%
%
%
\def\@citedef#1#2{\expandafter\gdef\csname\@citelabel{#1}\endcsname{#2}}%
%
%
%
\def\@readbblfile{%
   \ifx\@itemnum\@undefined
      \@innernewcount\@itemnum
   \fi
   \begingroup
      \def\begin##1##2{%
         \setbox0 = \hbox{\biblabelcontents{##2}}%
         \biblabelwidth = \wd0
      }%
      \def\end##1{}
      %
      %
      \@itemnum = 0
      \def\bibitem{\@getoptionalarg\@bibitem}%
      \def\@bibitem{%
         \ifx\@optionalarg\empty
            \expandafter\@numberedbibitem
         \else
            \expandafter\@alphabibitem
         \fi
      }%
      \def\@alphabibitem##1{%
         \expandafter \xdef\csname\@citelabel{##1}\endcsname {\@optionalarg}%
         \ifx\biblabelprecontents\@undefined
            \let\biblabelprecontents = \relax
         \fi
         \ifx\biblabelpostcontents\@undefined
            \let\biblabelpostcontents = \hss
         \fi
         \@finishbibitem{##1}%
      }%
      \def\@numberedbibitem##1{%
         \advance\@itemnum by 1
         \expandafter \xdef\csname\@citelabel{##1}\endcsname{\number\@itemnum}%
         \ifx\biblabelprecontents\@undefined
            \let\biblabelprecontents = \hss
         \fi
         \ifx\biblabelpostcontents\@undefined
            \let\biblabelpostcontents = \relax
         \fi
         \@finishbibitem{##1}%
      }%
      \def\@finishbibitem##1{%
         \biblabelprint{\csname\@citelabel{##1}\endcsname}%
         \@writeaux{\string\@citedef{##1}{\csname\@citelabel{##1}\endcsname}}%
         \ignorespaces
      }%
      %
      %
      \let\em = \bblem
      \let\newblock = \bblnewblock
      \let\sc = \bblsc
      \frenchspacing
      \clubpenalty = 4000 \widowpenalty = 4000
      \tolerance = 10000 \hfuzz = .5pt
      \everypar = {\hangindent = \biblabelwidth
                      \advance\hangindent by \biblabelextraspace}%
      \bblrm
      \parskip = 1.5ex plus .5ex minus .5ex
      \biblabelextraspace = .5em
      \bblhook
      \input \bblfilebasename.bbl
   \endgroup
}%
%
%
\@innernewdimen\biblabelwidth
\@innernewdimen\biblabelextraspace
%
%
%
\def\biblabelprint#1{%
   \noindent
   \hbox to \biblabelwidth{%
      \biblabelprecontents
      \biblabelcontents{#1}%
      \biblabelpostcontents
   }%
   \kern\biblabelextraspace
}%
%
%
%
\def\biblabelcontents#1{{\bblrm [#1]}}%
%
%
\def\bblrm{\rm}%
%
%
\def\bblem{\it}%
%
%
\def\bblsc{\ifx\@scfont\@undefined
              \font\@scfont = cmcsc10
           \fi
           \@scfont
}%
%
%
\def\bblnewblock{\hskip .11em plus .33em minus .07em }%
%
%
\let\bblhook = \empty
%
%
%
\def\printcitestart{[}
\def\printcitefinish{]}
\def\printbetweencitations{, }
\def\printcitenote#1{, #1}
%
%
%
\let\citation = \@gobble
%
%
%
\@innernewcount\@numparams
%
%
\def\newcommand#1{%
   \def\@commandname{#1}%
   \@getoptionalarg\@continuenewcommand
}%
%
%
\def\@continuenewcommand{%
   \@numparams = \ifx\@optionalarg\empty 0\else\@optionalarg \fi \relax
   \@newcommand
}%
%
%
\def\@newcommand#1{%
   \def\@startdef{\expandafter\edef\@commandname}%
   \ifnum\@numparams=0
      \let\@paramdef = \empty
   \else
      \ifnum\@numparams>9
         \errmessage{\the\@numparams\space is too many parameters}%
      \else
         \ifnum\@numparams<0
            \errmessage{\the\@numparams\space is too few parameters}%
         \else
            \edef\@paramdef{%
               \ifcase\@numparams
                  \empty  No arguments.
               \or ####1%
               \or ####1####2%
               \or ####1####2####3%
               \or ####1####2####3####4%
               \or ####1####2####3####4####5%
               \or ####1####2####3####4####5####6%
               \or ####1####2####3####4####5####6####7%
               \or ####1####2####3####4####5####6####7####8%
               \or ####1####2####3####4####5####6####7####8####9%
               \fi
            }%
         \fi
      \fi
   \fi
   \expandafter\@startdef\@paramdef{#1}%
}%
%
%
%
%
\def\@readauxfile{%
   \if@auxfiledone \else 
      \global\@auxfiledonetrue
      \@testfileexistence{aux}%
      \if@fileexists
         \begingroup
            \endlinechar = -1
            \catcode`@ = 11
            \input \jobname.aux
         \endgroup
      \else
         \message{\@undefinedmessage}%
         \global\@citewarningfalse
      \fi
      \immediate\openout\@auxfile = \jobname.aux
   \fi
}%
%
%
\newif\if@auxfiledone
\ifx\noauxfile\@undefined \else \@auxfiledonetrue\fi
%
%
%
%
\@innernewwrite\@auxfile
\def\@writeaux#1{\ifx\noauxfile\@undefined \write\@auxfile{#1}\fi}%
%
%
%
\ifx\@undefinedmessage\@undefined
   \def\@undefinedmessage{No .aux file; I won't give you warnings about
                          undefined citations.}%
\fi
%
%
\@innernewif\if@citewarning
\ifx\noauxfile\@undefined \@citewarningtrue\fi
%
%
%
\catcode`@ = \@oldatcatcode

\def\pfeilso{\leavevmode
            \vrule width 1pt height9pt depth 0pt\relax
           \vrule width 1pt height8.7pt depth 0pt\relax
           \vrule width 1pt height8.3pt depth 0pt\relax
           \vrule width 1pt height8.0pt depth 0pt\relax
           \vrule width 1pt height7.7pt depth 0pt\relax
            \vrule width 1pt height7.3pt depth 0pt\relax
            \vrule width 1pt height7.0pt depth 0pt\relax
            \vrule width 1pt height6.7pt depth 0pt\relax
            \vrule width 1pt height6.3pt depth 0pt\relax
            \vrule width 1pt height6.0pt depth 0pt\relax
            \vrule width 1pt height5.7pt depth 0pt\relax
            \vrule width 1pt height5.3pt depth 0pt\relax
            \vrule width 1pt height5.0pt depth 0pt\relax
            \vrule width 1pt height4.7pt depth 0pt\relax
            \vrule width 1pt height4.3pt depth 0pt\relax
            \vrule width 1pt height4.0pt depth 0pt\relax
            \vrule width 1pt height3.7pt depth 0pt\relax
            \vrule width 1pt height3.3pt depth 0pt\relax
            \vrule width 1pt height3.0pt depth 0pt\relax
            \vrule width 1pt height2.7pt depth 0pt\relax
            \vrule width 1pt height2.3pt depth 0pt\relax
            \vrule width 1pt height2.0pt depth 0pt\relax
            \vrule width 1pt height1.7pt depth 0pt\relax
            \vrule width 1pt height1.3pt depth 0pt\relax
            \vrule width 1pt height1.0pt depth 0pt\relax
            \vrule width 1pt height0.7pt depth 0pt\relax
            \vrule width 1pt height0.3pt depth 0pt\relax}

\def\pfeilsw{ \leavevmode 
            \vrule width 1pt height0.3pt depth 0pt\relax
            \vrule width 1pt height0.7pt depth 0pt\relax
            \vrule width 1pt height1.0pt depth 0pt\relax
            \vrule width 1pt height1.3pt depth 0pt\relax
            \vrule width 1pt height1.7pt depth 0pt\relax
            \vrule width 1pt height2.0pt depth 0pt\relax
            \vrule width 1pt height2.3pt depth 0pt\relax
            \vrule width 1pt height2.7pt depth 0pt\relax
            \vrule width 1pt height3.0pt depth 0pt\relax
            \vrule width 1pt height3.3pt depth 0pt\relax
            \vrule width 1pt height3.7pt depth 0pt\relax
            \vrule width 1pt height4.0pt depth 0pt\relax
            \vrule width 1pt height4.3pt depth 0pt\relax
            \vrule width 1pt height4.7pt depth 0pt\relax
            \vrule width 1pt height5.0pt depth 0pt\relax
            \vrule width 1pt height5.3pt depth 0pt\relax
            \vrule width 1pt height5.7pt depth 0pt\relax
            \vrule width 1pt height6.0pt depth 0pt\relax
            \vrule width 1pt height6.3pt depth 0pt\relax
            \vrule width 1pt height6.7pt depth 0pt\relax
            \vrule width 1pt height7.0pt depth 0pt\relax
            \vrule width 1pt height7.3pt depth 0pt\relax
            \vrule width 1pt height7.7pt depth 0pt\relax
            \vrule width 1pt height8.0pt depth 0pt\relax
            \vrule width 1pt height8.3pt depth 0pt\relax
            \vrule width 1pt height8.7pt depth 0pt\relax
            \vrule width 1pt height9pt depth 0pt\relax
      }


\def\widestnumber#1#2{}

\def\citewarning#1{\ifx\shlhetal\relax 
    \else
    \par{#1}\par
    \fi
}

\def\rm{\fam0 \tenrm}

\def\fakesubhead#1\endsubhead{\bigskip\noindent{\bf#1}\par}



%
%
%

%

\font\textrsfs=rsfs10
\font\scriptrsfs=rsfs7
\font\scriptscriptrsfs=rsfs5

\newfam\rsfsfam
\textfont\rsfsfam=\textrsfs
\scriptfont\rsfsfam=\scriptrsfs
\scriptscriptfont\rsfsfam=\scriptscriptrsfs

\edef\oldcatcodeofat{\the\catcode`\@}
\catcode`\@11

\def\Cal@@#1{\noaccents@ \fam \rsfsfam #1}

\catcode`\@\oldcatcodeofat


\expandafter\ifx \csname margininit\endcsname \relax\else\margininit\fi

\long\def\red#1\endred{}
\long\def\green#1\endgreen{}
\long\def\blue#1\endblue{}

\def\endred{ \unmatched endred! }
\def\endgreen{ \unmatched endgreen! }
\def\endblue{ \unmatched endblue! }

\ifx\latexcolors\undefinedcs\def\latexcolors{}\fi

\def\emptycs{}
\def\evaluatelatexcolors{%
        \ifx\latexcolors\emptycs\else
        \expandafter\xxevaluate\latexcolors\xxfertig\evaluatelatexcolors\fi}
\def\xxevaluate#1,#2\xxfertig{\setupthiscolor{#1}%
        \def\latexcolors{#2}}

\font\smallfont=cmsl7
\def\rutgerscolor{\ifmmode\else\endgraf\fi\smallfont
\advance\leftskip0.5cm\relax}
\def\setupthiscolor#1{\edef\tmptmpcs{\noexpand\bgroup\noexpand\rutgerscolor
\noexpand\def\noexpand\currentcolor{#1}%
\noexpand}%
\expandafter\let\csname#1\endcsname\tmptmpcs
\def\tmptmpcs{\checkColorUnmatched{#1}\popthecolor}
\expandafter\let\csname end#1\endcsname\tmptmpcs}

\def\checkColorUnmatched#1{\def\expectcolor{#1}%
    \ifx\expectcolor\currentcolor   
    \else \edef\failhere{\noexpand\tryingToClose '\currentcolor' with end\expectcolor}\failhere\fi}

\def\currentcolor{???}

\def\popthecolor{\ifmmode\else\endgraf\fi\egroup}

\expandafter\def\csname#1\endcsname{}

\evaluatelatexcolors

 \let\outerhead\head
 \def\head{\innerhead}
 \let\innerhead\outerhead

 \let\outersubhead\subhead
 \def\subhead{\innersubhead}
 \let\innersubhead\outersubhead

 \let\outersubsubhead\subsubhead
 \def\subsubhead{\innersubsubhead}
 \let\innersubsubhead\outersubsubhead

 \def\proclaim{\innerproclaim}
 \let\innerproclaim\outerproclaim

 %
 %
 %
 %

\def\demo#1{\medskip\noindent{\it #1.\/}}
\def\enddemo{\smallskip}

\def\remark#1{\medskip\noindent{\it #1.\/}}
\def\endremark{\smallskip}

\pageheight{8.5truein}
\topmatter
\title{Zero one laws for graphs with edge probabilities decaying with
distance, Part II} \endtitle
\rightheadtext{Zero one laws, etc., Part II}
\author {Saharon Shelah \thanks {\null\newline 
We thank John Baldwin and Shmuel Lifsches and Cigden Gencer and Alon Siton 
for helping in various ways and stages to make the paper more 
user friendly. \null\newline
The research partially supported by the United States -- Israel Binational
Science Foundation \null\newline
Publication 517 \null\newline
I would like to thank Alice Leonhardt for the beautiful typing. \null\newline
} \endthanks} \endauthor 

\affil{The Hebrew University of Jerusalem \\
 Einstein Institute of Mathematics \\
 Edmond J. Safra Campus, Givat Ram \\
 Jerusalem 91904, Israel
 \medskip
 Department of Mathematics \\
 Hill Center-Busch Campus \\
 Rutgers, The State University of New Jersey \\ 
 110 Frelinghuysen Road \\
 Piscataway, NJ 08854-8019  USA} \endaffil
\endtopmatter
\document

\newpage

\head {Introduction} \endhead  \resetall 
 \spuriousreset
\bigskip

This continues \cite{Sh:467} which is Part I and will be denoted here
by [I], background and a description of the results are
given in [I,\S0]; as this is the second part, our sections are named
\S4 - \S7 and not \S1 - \S4.

Recall that we fix an irrational $\alpha \in (0,1)_{\Bbb R}$ and the
random graph ${\Cal M}_n = {\Cal M}^0_n$ is drawn as follows
\mr
\item "{$(a)$}"  its set of elements is $[n] = \{1,\dotsc,n\}$
\sn
\item "{$(b)$}"  for $i <j$ in $[n]$ the probability of $\{i,j\}$
being an edge is $p_{|i-j|}$ where $p_\ell = 1/\ell^\alpha$ if 
$\ell >1$ as is $1/2^\alpha$ if $\ell=1$ or just
\footnote{originally we assume the former but actually also the latter works}
$p_\ell = 1/\ell^\alpha$ for $\ell > 1$
\sn
\item "{$(c)$}"  the drawing for the edges are independent
\sn
\item "{$(d)$}"  ${\Cal K}_n$ is the set of possible values of ${\Cal
M}_n,{\Cal K}$ is the class of graphs.
\ermn
Out main interest is to prove the 0-1 laws (for first order logic) for
this 0-1 context, but also to analyze the limit theory.
\bn
We can now explain our intentions.
\sn
\ub{Zero Step}:  We define relations $<^*_x$ on the class of graphs with
no apparent relation to the probability side.
\bn
\ub{First Step}:  We can prove that these $<^*_x$ have the formal
properties of $<_x$, like $<^*_i$ is a partial order etc., this is
done in \S4, e.g., in \scite{1.17} above. \nl
Remember from [I,\S1] that $A<_a B \Leftrightarrow$ for random enough
${\Cal M}_n$ and $f:A \hookrightarrow {\Cal M}_n$, the maximal number 
of pairwise disjoint $g\supseteq f$ satisfying $g:B \rightarrow 
{\Cal M}_n$ is $< n^\varepsilon$ (for every fixed $\varepsilon$).
\sn
\ub{Second Step}:   We shall start dealing with the two 
version of $<_a$: the $<_a$ from [I,\S1] and $<^*_a$ defined in 
\scite{1.11}(5) below.  We intend to prove:
\mr
\item "{$(*)$}"   $A<^*_a B\quad\Rightarrow\quad A<_a B$.
\ermn
For this it suffices to show that for every $f:A \hookrightarrow {\Cal M}_n$ 
and positive real $\varepsilon$, the expected value of the following is 
$\leq 1/n^{\varepsilon}$: the number of extensions
$g:B\rightarrow {\Cal M}_n$ of $f$ satisfying ``the sets 
Rang$(g\restriction (B \setminus A))$, $f(A)$ are with 
distance $\geq n^\varepsilon$". Then, the expected value of the 
number of $k$-tuples of such (pairwise) disjoint $g$ is 
$\leq \frac{1}{n^{k\varepsilon}}$. So if $k \varepsilon>|A|$, the 
expected value of the number of functions $f$ with $k$ pairwise 
disjoint such extensions $g$ is $< \frac{1}{n^{k\varepsilon-|A|}}$.  Hence for 
random enough ${\Cal M}_n$, for every $f: A \hookrightarrow {\Cal
M}_n$ there are no such $k$-tuples of pairwise disjoint $g$'s.
This will help to prove that $<^*_i = <_i$.  We do this and more
probability arguments in \S5.
\bn
\ub{Third step}:  We deduce from \S5 that $<^*_x = <_x$ for all
relevant $x$ and prove that the context is weakly nice.  We then work
somwhat more to prove the existential part of nice (the simple goodness (see 
Definition [I,2.12,(1)]) of appropriate candidate).
I.e. we first prove ``weakly niceness'' by proving that $A<^*_i B$
implies $(A, B)$ satisfy the demand for $<_i$ of [I,\S1], and in a strong
way the parallel thing for $\leq_s$. Those involve probability
estimation, i.e., quoting \S5.  
But we need more: sufficient conditions for appropriate tuples to
be simply good and this is the first part of \S6.
\bn
\ub{Fourth step}:  This is the universal part from niceness. 
This does not involve any probability,
just weight computations (and previous stages), in other words, purely
model theoretic investigation of the ``limit" theory.
By the ``universal part of
nice'' we mean (A) of [I,2.13,(1)] which includes:
\mr
\item "{${{}}$}"  if $\bar a \in {}^k({\Cal M}_n),b \in {\Cal M}_n$ 
\ub{then} there are $m_1<m$, $B\subseteq c \ell^{k,m_1}(\bar a)$ 
such that $\bar{a}\subseteq B$ and
\endroster
$$
c \ell^k(B) \nonforkin{}{}_{B}^{{\Cal M}_n}
(c \ell^k(\bar a b,{\Cal M}_n) \backslash c \ell^k(B,{\Cal M})) \cup B.
$$
\mn
This is done in the latter part of \S6.
\newpage

\head {\S4 Applications} \endhead  \resetall \sectno=4
 \spuriousreset
\bigskip

We intend to apply the general theorems 
(Lemmas [I,2.17,2.19]), to our problem.  
That is, we try to answer: does the main 
context ${\Cal M}^0_n$ with $p_i=1/i^{\alpha}$ for
$i>1$ satisfies the 0-1 law?  So here our irrational number $\alpha \in
(0,1)_{\Bbb R}$ is fixed. We work in Main Context (see \scite{1.1} 
below, the other one, ${\Cal M}^1_n$, would work out as well, see \S7).
\bigskip

\demo{\stag{1.1} Context}  A particular case of [I,1.1]:
$p_i=1/{i^\alpha}$ for $i>1$, $p_1=p_2$ (where $\alpha \in (0,1)_{\Bbb R}$
is a fix irrational) and the $n$-th random structure is 
${\Cal M}_n = {\Cal M}^0_n = ([n],R)$ 
(i.e. only the graph with the probability of $\{i,j\}$ being $p_{|i-j|}$).
\enddemo
\bn
\ub{\stag{1.2} Fact}  1) For any graph $A$

$$
1 = \lim_n \text{ Prob}(A \text{ is embeddable into } {\Cal M}_n).
$$
\mn
2) Moreover \footnote{Actually also ``$\ge c n$'' works for $c \in 
\Bbb R^{>0}$ depending on $A$ only.} for every $\varepsilon>0$

$$
1= \lim_n \text{ Prob}(A \text{ has } \geq 
n^{1-\varepsilon} \text{ disjoint copies in } {\Cal M}_n).
$$
\mn
This is easy, still, before proving it, note that since by our definition of
the closure $A\subseteq c \ell^{m,k}(\emptyset,{\Cal M}_n)$ implies
that $A$ has $< n^\varepsilon$ embeddings into ${\Cal M}_n$ we get:
\demo{\stag{1.3} Conclusion}   $\langle c \ell^{m,k}_{{\Cal M}_n}
(\emptyset): n<\omega\rangle$ satisfies the 0--1 law (being a 
sequence of empty models).
\enddemo
\bn
Hence (see [I,Def.1.4,Conclusion 2.19])
\demo{\stag{1.4} Conclusion}  ${\Cal K}_\infty = {\Cal K}$ and 
for our main theorem it suffices to prove \ub{simple almost niceness}
of ${\frak K}$ (see Def.[I,2.13]).

(Now \scite{1.3} explicate one part of 
what in fact we always meant by ``random enough'' in previous discussions.)
\enddemo
\bigskip

\demo{Proof of \scite{1.2}}  Let the nodes of $A$ be 
$\{a_0,\ldots,a_{k-1}\}$. Let the event ${\Cal E}_r^n$ be:

$$
a_\ell \mapsto 2rk+ 2 \ell \text{ is an embedding of } A \text{ into }
{\Cal M}_n.
$$
\mn
The point of this is that for various values of $r$ these tries are going to
speak on pairwise disjoint sets of nodes, so we get independent events.
\enddemo
\bn
Now \nl
\ub{\stag{1.5} Subfact}  Prob$({\Cal E}^n_r) = q>0$ 
(i.e. $>0$ but it does \ub{not} depend on $n,r$).
\sn
(Note: this is not true in a close context where the probability of
$\{i,j\}$ being an edge when $i \ne j$ is
$1/n^{\alpha}+1/2^{|i-j|}$, as in that case the probability depends 
on $n$. But still, we can have $\geq q>0$ which suffices); where:

$$
q = \dsize \prod_{\ell<m<k,\{\ell,m\} \text{ edge}} \,
1/(2(m-\ell))^{\alpha} \times
\dsize \prod_{\ell<m<k,\{\ell,m\}\text{not an edge}}\,
\left(1-\frac{1}{(2(m-\ell))^{\alpha}}\right).
$$
\mn
(What we need is that all the relevant edges have probability $>0$,
$<1$. Note: if we have retained $p=1/i^{\alpha}$ this is false 
for the pairs $(i, i+1)$, so we have changed $p_1$. Anyway, in our case we
multiplied by $2$ to avoid this (in the definition of the event)).   
For the second case, (the probability of edge being $1/n^\alpha + 1/2^{(i-j)}$)

$$
q \geq \dsize \prod_{\ell<m<k,\{\ell,m\} \text{ edge}}
\frac{1}{2^{|m-\ell|}}\times 
\dsize \prod_{\ell<m<k,\{\ell,m\}\text{ not an edge}}
(1-\frac{1}{(3/2)^{|m-\ell|}}).
$$
\mn
So these Prob$({\Cal E}^n_r)$ have a positive lower 
bound which does not depend on $r$.

Also the events 
${\Cal E}^n_0,\ldots,{\Cal E}^n_{[\frac{n}{2k}]-1}$ are independent.
So the probability that they all fail is

$$
\dsize \prod_{i<\lfloor\frac{n}{2k}\rfloor}(1- \text{Prob}({\Cal E}^n_i)) \le
\dsize \prod_i(1-q)\leq(1-q)^{\frac{n}{2k}}
$$
\mn
which goes to $0$ quite fast. The ``moreover'' is left to the 
reader. \hfill$\square_{\scite{1.4}}$\margincite{1.4}
\bigskip

\definition{\stag{1.6} Definition}  1) Let

$$
\align
{\Cal T} = \{(A,B,\lambda):&A \subseteq B 
\text{ graphs (generally: models from } {\Cal K}) \text{ and }\\
 &\lambda \text{ an equivalence relation on }B\setminus A\}
\endalign
$$
\mn
We may write $(A,B,\lambda)$ instead of $(A,B,\lambda\restriction(B\setminus
A))$. \nl
2) We say that $X\subseteq B$ is $\lambda$-closed if:

$$
x\in X \text{ and } x \in B \cap \text{ Dom}(\lambda) 
\text{ implies } x/\lambda \subseteq X.
$$
\mn
3) $A \leq^* B$ if\ \footnote{Note: this is in our present
specific context, so this definition does not apply to \S1, \S2,
\S3,\S7; in fact, in \S7 we give a different definition for a different
context.} 
$A \leq  B \in {\Cal K}_\infty$ (clearly $\leq^*$ is a partial order).
\enddefinition
\bn
\ub{Story}:

We would like to ask for any given copy of $A$ in ${\Cal M}_n$, is there a copy
of $B$ above it, and how many, we hope for a dichotomy: i.e. usually none,
always few \ub{or} always many. The point of $\lambda$ is to take distance
into account, because  for our present distribution 
being near is important, $b_1 \lambda b_2$ will
indicate that $b_1$ and $b_2$ are near. Note that being near is not
transitive, but ``luck'' helps us, we will succeed to ``pretend'' it is. We
will look at many candidates for a copy of 
$B\setminus A$ and compute the expected
value.  We would like to show that saying ``variance small'' says that the true
value is near the expected value.
\bigskip

\definition{\stag{1.7} Definition}  1) For $(A,B,\lambda) \in {\Cal T}$
let

$$
\bold v(A,B,\lambda) = \bold v_{\lambda}(A,B) = |(B\setminus A)/\lambda|
$$
\mn
be the number of $\lambda$--equivalence classes in $B \setminus A$
($\bold v$ stands for vertices).
\enddefinition
\bn
(This measures degrees of freedom in choosing candidates for $B$ over
a given copy of $A$.)

2) Let

$$
\bold e(A,B,\lambda) = \bold e_\lambda(A,B) = |e_\lambda(A,B)| 
\text{ where}
$$

$$
e_\lambda(A,B) = \{e:e\text{ an edge of }B,e \nsubseteq A,
\text{ and } e \nsubseteq x/\lambda \text{ for } x \in B \setminus A\}.
$$
\mn
[This measures the number of ``expensive", ``long'' edges ($\bold e$ 
stands for edges).]
\bn
\ub{Story}:

$\bold v$ larger means that there are more candidates for $B$,

$\bold e$ larger means that the probability per candidate is smaller.
\bigskip

\definition{\stag{1.8} Definition}  1) For 
$(A,B,\lambda)\in {\Cal T}$ and our given irrational 
$\alpha \in (0,1)_{\Bbb R}$ we define ($\bold w$ stands for weight)

$$
\bold w(A,B,\lambda) = \bold w_\lambda(A,B) = 
\bold v_\lambda(A,B) - \alpha \bold e_\lambda(A,B).
$$
\mn
2) Let

$$
\align
\Xi(A,B) =: \{\lambda:&(A,B,\lambda) \in {\Cal T},
\text{ and if } C \subseteq B\setminus A \text{ is a nonempty} \\
  &\lambda \text{-closed set then } \bold w_\lambda(A, C\cup A)>0\}.
\endalign
$$
\mn
3) If $A \le^* B$ then we let $\xi(A,B) = \text{ Max}\{\bold
w_\lambda(A,B):\lambda \in \Xi(A,B)\}$.
\enddefinition
\bigskip

\demo{\stag{1.9} Observation}  1) $(A,B,\lambda) \in 
{\Cal T} \and  A \neq B \Rightarrow \bold w_\lambda(A,B)\neq 0$. \nl
2) If $A\leq^* B\leq^* C$ (hence $A \le^* C$) and 
$(A,C,\lambda)\in {\Cal T}$ and $B$ is
$\lambda$-closed \ub{then}
\mr
\item "{$(a)$}"   $(A,B,\lambda\restriction (B\setminus A))\in {\Cal T}$, 
\sn
\item "{$(b)$}"   $(B,C,\lambda\restriction (C\setminus B))\in {\Cal
T}$ 
\sn
\item "{$(c)$}"   $\bold w_\lambda(A,C) = 
\bold w_{\lambda\restriction (B\setminus A)}(A,B)+ 
\bold w_{\lambda\restriction (C\setminus B)}(B,C)$
\sn
\item "{$(d)$}"   similarly for $\bold v$ and $\bold e$.
\ermn
3) Note that \scite{1.9}(2) legitimizes our writing $\lambda$ instead of
$\lambda\restriction (C\setminus A)$ or $\lambda\restriction
(B\setminus (C\cup A))$ when $(A,B,\lambda) \in {\Cal T}$ and 
$C$ is a $\lambda$-closed subset of $B$. Thus we
may write, e.g., $\bold w_\lambda(A\cup C,B)$ for 
$\bold w(A\cup C,B, \lambda\restriction (B\setminus A\setminus
C))$. \nl
4)  If $(A,B,\lambda) \in {\Cal T}$ and $D\subseteq B \setminus A$ and
$D^+=\bigcup \{x/\lambda: x\in D \}$ then $\bold w_{\lambda\restriction D^+}
(A,A \cup D^+) \leq \bold w_{\lambda\restriction D}(A,A \cup D)$ and
$D^+$ is $\lambda$-closed.
\enddemo
\bigskip

\demo{Proof}   1) As $\alpha$ is irrational and $\bold v_\lambda(A,B)$
is not zero. \nl
2) Clauses  (a), (b) are totally immediate, and for a
proof of clauses (c), (d) see the proof of \scite{1.16} below. \nl
3)  Left to the reader. \nl
4)  Clearly by the choice of $D^+$ we have
$\bold v_{\lambda\restriction D^+}(A,A\cup D^+) =
|D^+/\lambda| = |D/(\lambda \restriction D)| = 
\bold v_{\lambda\restriction D}(A,A\cup D)$ and $\bold e_{\lambda
\restriction D^+}(A,A \cup D^+) \ge \bold e_{\lambda \restriction
D}(A,A \cup D^+)$ hence $\bold w_{\lambda \restriction D^+}(A,A \cup
D^+) \le \bold w_{\lambda \restriction D}(A,A \cup D^+)$. 
\hfill$\square_{\scite{1.9}}$\margincite{1.9}
\enddemo
\bn
\ub{\stag{1.10} Discussion}:  Note: $\bold w_\lambda(A,B)$ measures in a 
sense the expected value of the number of copies of $B$ over a 
given copy of $A$ with $\lambda$ saying when one node is ``near to''
another.  Of
course, when $\lambda$ is the identity this degenerates to 
the definition in \cite{ShSp:304}. \nl
We would like to characterize $\leq_i$ and $\leq_s$ 
(from Definition [I,1.4,(3)] and Definition [I,1.4,(4)]), 
using $\bold w$ and to prove that
they are O.K. (meaning that they form a nice context). Looking at
the expected behaviour, we attempt to give an ``effective'' definition
(depending on $\alpha$ only). \nl
All of this, of course, just says what the intention of these
relations and functions is
(i.e. $<^*_i$, $<^*_s$, $<^*_{pr}$ and $\bold v$, $\bold e$, $\bold
w$ below); we still will not prove anything on the connections to
$\le_i,\le_s,\le_{\text{pr}}$. We may view it differently: 
We are, for our fix $\alpha$, defining $\bold w_{\lambda}(A,B)$ 
and investigating the 
$\leq^*_i,\leq^*_s,\leq^*_{pr}$ defined below per se ignoring the
probability side.
\bigskip

\definition{\stag{1.11} Definition}  1) $A\leq B$ means
$A$ is a submodel of $B$, and remember that by Definition \scite{1.6}(3),
$A\leq^* B$ means
\footnote{Note: this is in our present specific context, so 
this definition does not apply to \S1, \S2,
\S3, \S7; in fact in \S7 we give a different definition for a different
context.}  $A\leq B \in {\Cal K}_\infty$. \nl
2)  $A <^*_c B$ \ub{if} $A<^* B$ and for every $\lambda$, we have

$$
(A,B,\lambda) \in {\Cal T} \Rightarrow \bold w_\lambda(A,B)<0, 
$$
\mn
3)  $A\leq^*_i B$ \ub{if} $A\leq^* B$ and for every $A'$ we have

$$
A\leq^* A'<^* B\Rightarrow A'<^*_c B.
$$
\mn
Of course, $A \le^*_i B$ means $A \le^*_i B \and A \ne B$. \nl
4) $A\leq^*_s B$ \ub{if} $A\leq^* B$ and for no $A'$ do we have

$$
A<^*_i A'\leq^* B,
$$
\mn
Of course, $A <^*_s B$ means $A \le^*_s B \and A \ne B$. \nl
5) $A <^*_a B$ \ub{if} $A \leq^* B$, $\neg (A<^*_s B)$
(i.e. $A\leq^*B$ and there is $A' \subseteq B \setminus A$ such that
$A <^*_i A \cup A' \le^* B$), \nl
6)  $A<^*_{pr} B$ \ub{if} $A\leq^* B$ and $A<^*_s B$ but for no $C$
do we have $A<^*_s C<^*_s B$.
\enddefinition
\bigskip

\remark{\stag{1.12} Remark}  We \ub{intend} to prove that 
usually $\leq^*_x = \leq_x$ but it will take time. 
\endremark
\bigskip

\proclaim{\stag{1.13} Lemma}  Suppose $A'<^* B,(A',B,\lambda) \in 
{\Cal T}$ and $\bold w_\lambda(A',B)>0$.
\ub{Then} there is $A''$ satisfying $A'\leq^* A''<^* B$ such that
$A''$ is $\lambda$-closed and
\mr
\item "{$(*)_1$}"  $[A'',B,\lambda]$ we have 
$\bold w_\lambda(A'',B)>0$ and \ub{if} 
$C\subseteq B\setminus A''$, $C\neq\{\emptyset, B\setminus A''\}$ and $C$ is
$\lambda$--closed \ub{then} $\bold w_\lambda(A'',A''\cup C)>0$ and 
$\bold w_\lambda(A''\cup C,B)<0$.
\endroster
\endproclaim
\bigskip

\demo{Proof}   Let $C'$ be a maximal $\lambda$-closed 
subset of $B\setminus A'$ such that $\bold w_\lambda(A'\cup C',B)>0$. 
Such a $C'$ exists since $C'=\emptyset$ is as required and 
$B$ is finite. Let $A''=A'\cup C'$.  Since $C'$ is
$\lambda$-closed, $B \backslash A''$ is $\lambda$-closed and
$(A'',B,\lambda\restriction (B\setminus A''))
\in {\Cal T}$ and clearly $\bold w_\lambda(A'',B)>0$. Now
suppose $D\subseteq B\setminus A''$ is $\lambda$-closed, 
$D\notin\{\emptyset,B\setminus A''\}$. By the maximality of 
$C'$, $\bold w_\lambda(A''\cup D,B) < 0$.  Now (by \scite{1.9}(2)(c))

$$
\bold w_\lambda(A'',B) = \bold w_\lambda(A'',A''\cup D)
+ \bold w_\lambda(A''\cup D,B).
$$
\mn
and the left term is positive by the choice of $C'$ and $A''$, but
the right term is negative by the previous sentence so together  
we conclude $\bold w_\lambda(A'',A''\cup D)>0$. \nl
${{}}$  \hfill$\square_{\scite{1.13}}$\margincite{1.13}
\enddemo
\bigskip

\proclaim{\stag{1.14} Claim}  Assume $A <^* B$. 
The following statements are equivalent:
\mr
\widestnumber\item{$(iii)$}
\item "{$(i)$}"     $A<^*_i B$,
\sn
\item "{$(ii)$}"    for no $A'$ and $\lambda$ do we have:
{\roster
\itemitem{ $(*)_2$ }  $=(*)_{2}[A,A',B,\lambda]$ we have
$A \leq^* A'<^* B,
(A',B,\lambda) \in {\Cal T}$ and $\bold w_\lambda(A',B)>0$,
\endroster}
\item "{$(iii)$}"   for no $A'$, $\lambda$ do we have:
{\roster
\itemitem{ $(*)_3$ }  $=(*)_{3}[A,A',B,\lambda]$ we have $A\leq^* A'<^* B,
(A',B,\lambda) \in {\Cal T},\bold w_\lambda(A',B)>0$ and 
$(*)_1[A',B,\lambda]$ of \scite{1.13}.
\endroster}
\endroster
\endproclaim
\bigskip

\demo{Proof}   For the equivalence of 
the first and the second clauses read Definition \scite{1.11}(2), 
\scite{1.11}(3) (remembering \scite{1.9}(1)). Trivially $(*)_3
\Rightarrow (*)_2$ and hence the second clause implies the third one. Now
we will see that $(iii)\Rightarrow (ii)$. So suppose $\neg(ii)$, so let this
be exemplified by $A'$,$\lambda$ i.e. they satisfy $(\ast)_2$ so by
\scite{1.13} there is $A''$ such that $A'\leq^* A^{''}<^* B$ and
$(\ast)_1[A'',B,\lambda]$ of \scite{1.13} holds. So $A'',\lambda$
exemplified that $\neg(iii)$ holds.  \nl
${{}}$   \hfill$\square_{\scite{1.14}}$\margincite{1.14}
\enddemo
\bigskip

\demo{\stag{1.15} Observation}  1) If $(*)_3 [A,A',B,\lambda]$
from \scite{1.14}(iii) holds, \ub{then} we have: if 
$C\subseteq B\setminus A'$ is
$\lambda$-closed non-empty then $\bold w(A',A'\cup C,\lambda
\restriction C)>0$. \nl
[Why? If $C\neq B\setminus A'$ this is stated explicitly, otherwise
this means $\bold w(A',B,\lambda)>0$ which holds.] 
\nl
2) In $(*)_3$ of \scite{1.14}(iii), i.e., \scite{1.13}$(*)_1[A',B,\lambda]$,
we can allow any $\lambda$--closed $C\subseteq B\setminus A'$ if we make the
inequalities non-strict.  [Why?  If $C= \emptyset$ then 
$\bold w_{\lambda}(A',A'\cup C) = \bold w_{\lambda}(A',A') = 0$,
$\bold w_{\lambda}(A'\cup C,B) = \bold w_{\lambda}(A',B) > 0$. 
If $C = B\setminus A'$ then 
$\bold w_{\lambda}(A',A'\cup C) = \bold w_\lambda(A',B)>0$ and 
$\bold w_{\lambda}(A'\cup C,B) = \bold w_{\lambda}(B,B)=0$. 
Lastly if $C\notin \{\emptyset, B \setminus A'\}$ we use
\scite{1.13}$(*)_1[A',B,\lambda]$ itself.] \nl
3) If $(A,B,\lambda) \in {\Cal T},A' \leq^* A,B' \leq^* B,A' \leq^*B'$ and
$B \setminus A=B' \setminus A'$ \ub{then}
$(A',B',\lambda) \in {\Cal T},\bold w(A',B',\lambda) \geq 
\bold w(A,B,\lambda)$ also $\bold e(A',B',\lambda) \leq \bold e
(A,B,\lambda),\bold v(A',B',\lambda)= \bold v(A,B,\lambda)$. \nl
4) In (3) if in addition $\nonforkin {A}{B'}_{A'}^{M}$, i.e., no edge
$\{x,y\}$ with $x \in A \backslash A',y \in B' \backslash A'$ \ub{then} the
equalities hold.
\enddemo
\bigskip

\proclaim{\stag{1.16} Claim}  $A \leq_s^* B$ \ub{if and only if} 
either $A=B$ or for some $\lambda$ we have:

$(A,B,\lambda) \in {\Cal T}$ and $\bold w_\lambda(A,B)>0$, moreover

for every nonempty $\lambda$-closed $C\subseteq B\setminus A$, we have
$\bold w(A,A\cup C,\lambda\restriction C)>0$, that is $\Xi(A,B) \ne \emptyset$.
\endproclaim
\bigskip

\demo{Proof}  The \ub{only if} direction:

So we have $A\leq_s^* B$. If $A=B$ we are done: the left side holds 
as its first possibility is $A = B$.  So assume $A<_s^* B$. 
Let $C$ be minimal such that $A\leq^* C\leq^* B$ and for some $\lambda_0$
the triple $(C,B,\lambda_0) \in {\Cal T}$ satisfies: for every non empty
$\lambda_0$-closed $C'\subseteq B\setminus C$ we have $\bold w(C,C\cup C',
\lambda_0\restriction C')> 0$ (exists as $C=B$ is O.K. as there is no such
$C'$).  By \scite{1.9}(4): for every non empty $C' \subseteq B
\backslash C$ we have $\bold w(C,C \cup C',\lambda_0 \restriction C') > 0$
hence $\neg(C <^*_i C \cup C')$ by $(i) \Leftrightarrow (ii)$ of \scite{1.14}.
If $C=A$ we have finished by the definition of $\le^*_s$.  
Otherwise, the hypothesis $A \leq^*_s B$ implies that $\neg
(A<^*_i C)$, hence by \scite{1.14} the third clause $(iii)$ fails
which means that (recalling \scite{1.11}(1)), for some 
$C'$, $\lambda_1$ we have $A\leq^* C'<^* C$, $(C',C,\lambda_1) \in
{\Cal T},\bold w_{\lambda_1}(C',C) > 0$ 
and for every $\lambda_1$-closed $D\subseteq C\setminus C'$ 
satisfying $D\neq\{\emptyset,C\setminus C'\}$ we have

$$
\bold w(C',C'\cup D,\lambda_1\restriction D)>0,
$$
and

$$
\bold w\bigl(C'\cup D,C,\lambda_1 \restriction (C\setminus 
C'\setminus D)\bigr)<0.
$$
\mn
Define an equivalence relation $\lambda$ on $B\setminus C'$: an equivalence
class of $\lambda$ is an equivalence class of $\lambda_0$ or an equivalence
class of $\lambda_1$.

We shall show that $(C',B,\lambda)$ satisfies the requirement above on
$C$, thus contradicting the
minimality of $C$.  Clearly $A\leq^* C'\leq^* B$. So let
$D\subseteq B\setminus C'$ be $\lambda$-closed and
we define $D_0=D\cap (B\setminus C)$, $D_1=D\cap 
(C\setminus C')$. Clearly $D_0$ is $\lambda_0$-closed so 
$\bold w(C,C\cup D_0,\lambda \restriction D_0) \geq 0$ 
(see \scite{1.15}(2)), and $D_1$ is $\lambda_1$-closed so
$\bold w(C',C'\cup D_1,\lambda \restriction D_1)\geq 0$ 
(this follows from: for every $\lambda_1$-closed 
$D\subseteq C\setminus C'$ satisfying $D\neq \{\emptyset,C \setminus
C'\}$ we have $\bold w_\lambda(C',C'\cup D,\lambda_1 \restriction
D)>0$ and  by \scite{1.15}(2)). Now (in the last line we change $C'$
to $C$ twice), by \scite{1.15}(3) we will get

$$
\align
\bold v(C',C'\cup D,\lambda)=
&|D/\lambda|=|D_1/\lambda_1|+|D_0/\lambda_0|\\ 
   &= \bold v(C',C'\cup D_1,\lambda\restriction D_1) +
\bold v(C' \cup D_1,C' \cup D_1 \cup D_0,\lambda \restriction D_0)\\
   &= \bold v(C',C'\cup D_1,\lambda\restriction D_1)
+ \bold v(C,C\cup D_0,\lambda \restriction D_0),
\endalign
$$
\mn
and (using \scite{1.11}(3)):

$$
\align
\bold e(C',C'\cup D,\lambda) &= \bold e(C',C'\cup D_1,\lambda
\restriction D_1) \\
  &+ \bold e(C'\cup D_1,C'\cup D_1\cup D_0,\lambda\restriction D_0)\\
  &\le \bold e(C',C'\cup D_1,\lambda\restriction D_1) +
\bold e(C,C\cup D_0,\lambda \restriction D_0),
\endalign
$$
\mn
and hence

$$
\align
\bold w(C',C'\cup D,\lambda) &= \bold v(C',C'\cup D,\lambda) -
\alpha \bold e(C',C'\cup D,\lambda) \\
  & = \bold v(C',C'\cup D_1,\lambda\restriction D_1)
+ \bold v(C,C \cup D_0,\lambda \restriction D_0)\\
  &-\alpha \bold e(C',C'\cup D_1,\lambda\restriction D_1) \\
  &-\alpha \bold e(C'\cup D_1,C'\cup D_1\cup D_0,\lambda
\restriction D_0)\\
  &\ge \bold v(C',C'\cup D_1,\lambda\restriction D_1)
+ \bold v(C,C\cup D_0,\lambda \restriction D_0)\\
  &-\alpha \bold e(C',C'\cup D_1,\lambda\restriction D_1)
-\alpha \bold e(C,C\cup D_0,\lambda\restriction D_0)\\
  &= \bold w(C',C'\cup D_1,\lambda\restriction D_1) +
\bold w(C,C\cup D_0,\lambda \restriction D_0)\geq 0,
\endalign
$$
\mn
and the (strict) inequality holds by the irrationality of $\alpha$, i.e. by
\scite{1.9}(1). So actually $(C', B, \lambda)$ satisfies the requirements
on $C$, $\lambda_0$ thus giving contradiction to the minimality of $C$.
\bn
\ub{The if direction}:  

As the case $A = B$ is obvious, we can assume that the second half of
\scite{1.16} holds. 
So let $\lambda$ be as required in the second half of \scite{1.16}. 

Suppose $A <^* C \leq^* B$, and 
we shall prove that $\neg (A<^*_i C)$ thus finishing by Definition 
\scite{1.11}.
We shall show that $(A',\lambda') = (A,\lambda \restriction (C
\backslash A))$ satisfies $(*)_2[A,A,C,\lambda]$ from \scite{1.14}, 
thus $(ii)$ of \scite{1.14} fail hence $(i)$ of \scite{1.14} fails, 
i.e., $\neg(A <^*_i C)$ as required. Let
$D=\bigcup\{x/\lambda: x\in C\setminus A\}$,
so $D$ is a non empty $\lambda$-closed subset of $B\setminus A$. Hence
by the present assumption on $A$,
$B$, $\lambda$ we have $\bold w(A, A\cup D,\lambda\restriction D)>0$. Now

$$
\bold v(A,C,\lambda\restriction C) = |C/\lambda|= |D/\lambda|=
\bold v(A,D,\lambda\restriction D)
$$
\mn
and

$$
\bold e(A,C,\lambda \restriction C) \le \bold e(A,D,\lambda\restriction D)
$$
\mn
so $\bold w(A,C,\lambda\restriction C)\geq \bold w
(A,D,\lambda\restriction D)>0$ as requested. 
\hfill$\square_{\scite{1.16}}$\margincite{1.16}
\enddemo
\bigskip

\proclaim{\stag{1.17} Claim}:  1)  $\leq^*_i$ is transitive. \nl
2) $\leq^*_s$ is transitive. \nl
3)  For any $A\leq^* C$ for some $B$ we have $A\leq^*_i B \leq^*_s C$.
\nl
4)  If $A<^* B$ and $\neg(A\leq^*_s B)$ \ub{then} $A <^*_c B$ or there
is $C$ such that $A <^* C <^* B$, $\neg(A <^*_s C)$. \nl
5)  Smoothness holds (with $<^*_i$ instead of $<_i$ see [I,2.5,(4)]),
that is
\mr
\item "{$(a)$}"   if $A \le^* C \le^* M \in {\Cal K},A \le^* B
\le^* M$, $B \cap C=A$ \ub{then} $A <^*_c B \Rightarrow C <^*_c B \cup
C$ and $A \le^*_i C \Rightarrow B \le^*_i C \cup B$
\sn
\item "{$(b)$}"   if in addition $\nonforkin{C}{B}_{A}^{M}$ then 
$A <^*_c B \Leftrightarrow C \le^*_c B \cup C$ and $A \le^*_i C
\Leftrightarrow \le^*_i C \cup B$ and $A \le^*_s B \Leftrightarrow C
\le^*_s B \cup C$.
\ermn
6)  For $A<^* B$ we have $\neg(A\leq^*_sB)$ \ub{iff} $(\exists C)(A <^*_c
C\leq^* B)$. \nl
7)  If $A\leq^*B \leq^*C$ and $A\leq^*_sC$ \ub{then} $A \leq^*_s
B$. \nl
8) If $A_\ell \le^*_s B_\ell$ for $\ell =1,2,A_1 \le^* A_2,B_1 \le^*
B_2$ and $B_2 \backslash A_2 = B_1 \backslash A_1$ \ub{then}
$$
\xi(A_1,B_1) \ge \xi(A_2,B_2).
$$
9) In (8), equality holds iff $A_2,B_1$ are freely amalgamated over
$A_1$ inside $B_2$. \nl
10) If $A <^*_s B_\ell$ for $\ell =1,2$ and $B_1 <^*_i B_2$ \ub{then}
$\xi(A,B_1) > \xi(A,B_2)$. \nl
11) If $B_1 <^* B_2$ and for no $x \in B_1,y \in B_2 \backslash B_1$
is $\{x,y\}$ an edge of $B_2$ \ub{then} $B_1 <^*_s B_2$. \nl
12) If $A \le^* B \le^* C$ and $A \le^*_i C$ \ub{then} $B \le^*_i
C$. \nl
13) If $A <^*_{\text{pr}} B$ and $a \in B \backslash A$ \ub{then} $A
\cup \{a\} \le^*_i B$. \nl
14) If $A_1 <^*_{\text{pr}} B_1,A_1 \le^* A_2 \le^* B_2$ and $B_1
\le^* B_2$ \ub{then} $A_2 \le^*_s B_2$ or $A_2 <^*_{\text{pr}} B_2$.
\endproclaim
\bigskip

\demo{Proof}  1) So assume $A\leq^*_iB\leq^*_i C$ and we shall prove
$A \leq^*_i C$.  By \scite{1.14} it suffice to prove that clause (ii) there
holds with $A,C$ here standing for $A,B$ there. So assume
$A\leq^*A'<^*C$, $(A',C,\lambda) \in {\Cal T}$ and we shall prove that
$\bold w_{\lambda}(A',C) \le 0$, this suffice. Let $A'_1=A'\cap B$,
$A'_0=: B \cup A' \cup \cup\{x/\lambda:x \in B\}$, now 
as $A\leq^*_iB$ by \scite{1.14} + \scite{1.15}(2) we have
$\bold w_{\lambda}(A'_1,B) \leq 0$, and by \scite{1.15}(3) we have
$\bold w_\lambda(A',B \cup A') \le \bold w_\lambda(A'_1,B)$ and by
\scite{1.9}(4) we have $\bold w_\lambda(A',A'_0) \le \bold
w_\lambda(A',A' \cup B)$.  Those three inequalities together gives
$\bold w_\lambda(A',A'_0) \le 0$ and as $B \leq^*_i C$ by \scite{1.14} we have 
$\bold w_{\lambda}(A'_0,C)\le 0$.  By \scite{1.9}(2)(c) we have
$\bold w_{\lambda}(A',C) = \bold w_{\lambda}(A',A'_0) +
\bold w_{\lambda}(A'_0,C)$ and by the previous sentence the latter 
is $\le 0 + 0 = 0$, so $\bold w_{\lambda}(A,A')\leq 0$ as required.\nl
2) We use the condition from \scite{1.16}. So assume
$A_0 \leq^*_s A_1 \leq^*_s A_2$ and $\lambda_\ell$ witness 
$A_\ell \leq^*_s A_{\ell+1}$ (i.e. $(A_\ell,A_{\ell+1},\lambda_\ell)$
is as in \scite{1.16}). Let $\lambda$ be the equivalence relation on 
$A_2\setminus A_0$ such that for $x\in A_{\ell+1} \setminus A_\ell$ we
have $x/\lambda = x/\lambda_\ell$. 
Easily $(A_0,A_2,\lambda) \in {\Cal T}$.  Now, by \scite{1.9}(2)(c) and
\scite{1.16} the triple $(A_0, A_2, \lambda)$ satisfies
the second condition in \scite{1.16} so $A_0 \leq^*_s A_2$. \nl
3) Let $B$ be maximal such that $A \leq^*_i B\leq^*C$, such $B$ exists
as $C$ is finite \footnote{Actually the finiteness is not needed if for
possibly infinite $A,B$ we define $A \le^*_i B$ iff for every finite
$B' \le^* B$ there is a finite $B''$ such that 
$B' \le^* B'' \le^* B$, and $B''\cap A<^*_iB''$.} and
for $B=A$ we get $A \leq^*_i B \leq^*C$. Now if $B \leq^*_s C$ we are done,
otherwise by the definition of $\le^*_s$ in \scite{1.11}(4) there 
is $B'$ such that $B <^*_i B'\leq^* C$,
now by part ($1$) we have $A \leq^*_i B'\leq^* C$, contradicting the
maximality of $B$, so really $B \leq^*_s C$ and we are done.\nl
4) We assume $A <^* B$ and now if $A <^*_c B$ we are done hence we can assume
$\neg(A <^*_c B)$, clearly 
there is $\lambda$ such that $(A,B,
\lambda)\in {\Cal T}$ and $\bold w_\lambda(A,B) \ge 0$. 
So by the irrationality of $\alpha$ the inequality is strict 
and by \scite{1.13}
there is $C$ such that $A \leq^* C<^* B$, $C$ is 
$\lambda$-closed, $\bold w_\lambda(C,B) > 0$ and 
if $C'\subseteq B\setminus C$ is non-empty
$\lambda$-closed and $ \ne B \backslash C$ then 
$\bold w_\lambda(C, C\cup C') > 0 \and \bold w_\lambda(C\cup C',B)<0$.
So by \scite{1.16} + inspection, 
$C<^*_s B$, 
so by \scite{1.17}(2), $A\leq^*_s C \Rightarrow A^* \leq^*_s B$, but 
we know that $\neg(A <^*_s B)$ hence 
$\neg(A \leq^*_s C)$, so the second possibility in the conclusion
holds.\nl
5) \ub{Clause $(a)$}:  $A \le^*_c B \Rightarrow C <^*_c B \cup C$ and
$A \le^*_i C \Rightarrow B \le^*_i B \cup C$. \nl
[Why?  Note by our assumption $C <^* B \cup C$ and $B <^* B \cup C$.
The first desired conclusion is easier, 
so we prove the second hence assume $A \le^*_i C$.
If $B \le^* D <^* B \cup C$, and $(D,B \cup C,\lambda) \in {\Cal T}$ then
$A \le^* D\cap C <^* C$ so as $A \le^*_i C$, by the definition of $\leq_i$
we have $\bold w_{\lambda}(D\cap C,C) < 0$ hence (noting $C \setminus
D\cap C= B\cup C\setminus D)$ by Observation \scite{1.15}(3) we have
$\bold w_{\lambda}(D, B\cup C) \leq \bold w_{\lambda}(D \cap C,C) < 0$ hence
(by the Def. of $\leq^*_c$) $D \le^*_c B \cup C$. As this holds
for any such $D,\lambda$ by Def. \scite{1.11}(3) we have $B \le^*_i B\cup C$
as required.]
\bn
\ub{Clause $(b)$}:   If in addition $\nonforkin{C}{B}_{A}^{M}$ then $A <^*_c C
\Leftrightarrow B \le^*_c C \cup B$ and
$A\leq_i C \Leftrightarrow B\leq_i C\cup B$ and $A \le^*_s B
\Leftrightarrow B \le^*_s B \cup C$. \nl
[Why?  Immediate by \scite{1.15}(4), Definition \scite{1.11} and part
(a).]
\nl
6) The ``only if" direction can be prove by induction on $|B|$, 
using \scite{1.17}(4).  For the if direction assume that for 
some $C$, $A <^*_c C \le^*B$ and choose a minimal $C$ like that.  Now if 
$A\leq^*A^*<^*C$, and $\lambda_1$ is an equivalence relation on 
$C\setminus A^*$ then let $\lambda_0$ be an equivalence relation 
on $A^*\setminus A$ such that
$\bold w_{\lambda_0}(A,A^*) \ge 0$ (exists by the minimality of $C$) and let
$\lambda=\lambda_0\cup\lambda_1$ so $(A,C,\lambda) \in {\Cal T}$ and by
\scite{1.9}(2)(i) we have $\bold w_\lambda(A^*,C) = \bold w_\lambda(A,C)-
\bold w_{\lambda}(A,A^*)$; but as 
$A<^*_c C$ we have $\bold w_\lambda(A,C) < 0$, and by the choice of
$\lambda_0$ we have $\bold w_\lambda(A,A^*) \ge 0$ hence $\bold
w_{\lambda}(A^*,C)<0$ 
hence $\bold w_{\lambda_1}(A^*,C) = \bold w_{\lambda}
(A^*,C)< 0$. As $\lambda_1$ was any equivalence relation on
$C\setminus A^*$ by Def. \scite{1.11}(2) we have shown that 
$A^* <^*_c C$.  By the definition of $\leq^*_i$ (\scite{1.11}(3)),
as $A^*$ was arbitrary such that $A \leq^* A^* <^* C$ by Definition
\scite{1.11}(3) 
we get that $A<^*_i C$, hence by the definition of
$\le_s$ (\scite{1.11}(4)), we can deduce $\neg (A\leq^*_s B)$ as required.\nl
7) Immediate by Definition \scite{1.11}(4). \nl
8) It is enough to prove that
\mr
\item "{$\circledast$}"  if $\lambda \in \Xi(A_2,B_2)$ then $\lambda
\in \Xi(A_1,B_1)$ and $\bold w_\lambda(A_1,B_1) \ge \bold
w_\lambda(A_2,B_2)$.
\ermn
So assume $\lambda$ is an equivalence relation over $B_2 \backslash
A_2$ which is equal to $B_1 \backslash A_1$, now for every non-empty
$\lambda$-closed $C \subseteq B_1 \backslash A_1$ we have
\mr
\widestnumber\item{$(iii)$}  
\item "{$(i)$}"  $\bold v_\lambda(A_1,A_1 \cup C) = |C/\lambda| =
\bold v_\lambda(A_2,A_2 \cup C)$
\sn
\item "{$(ii)$}"  $\bold e_\lambda(A_1,A_1 \cup C) \le \bold
e_\lambda(A_2,A_2 \cup C)$ \nl
[as any edge in $e_\lambda(A_1,A_1 \cup C)$ belongs to
$e_\lambda(A_2,A_2 \cup C)$] hence
\sn
\item "{$(iii)$}"  $\bold w_\lambda(A_1,A_1 \cup C) \ge \bold
w_\lambda(A_2,A_2 \cup C)$.
\ermn
So by the definition of $\Xi(A_1,B_1)$ we have $\lambda \in
\Xi(A_2,B_2) \Rightarrow \lambda \in \Xi(A_1,B_1)$ and, moreover, the desired
inequality in $\circledast$ holds. \nl
9) If $\nonforkin{A_2}{B_1}_{A_1}^{B_2}$ in the proof of (8) we get
$e_\lambda(A_1,A_1 \cup C) = e_\lambda(A_2,A_2 \cup C)$ hence 
$\bold w_\lambda(A_1,A_1 \cup C) = \bold w_\lambda(A_2,A_2 \cup C)$,
in particular $\bold w_\lambda(A_1,B_1) = \bold w_\lambda(A_2,B_2)$.
Also now the proof of (8) gives $\lambda \in \Xi(A_1,B_1) \Rightarrow
\lambda \in \Xi(A_2,B_2)$ so trivially $\xi(A_1,B_1) = \xi(A_2,B_2)$.

If $\neg(\nonforkin{A_2}{B_1}_{A_1}^{B_2})$ then for every equivalence
relation $\lambda$ on $B_1 \backslash A_1 = B_2 \backslash A_2$ we
have
\mr
\item "{$(ii)^+$}"  $\bold e_\lambda(A_1,B_1) < \bold
e_\lambda(A_2,B_2)$ \nl
[Why?  As $e_\lambda(A_1,B_1)$ is a proper subset of
$e_\lambda(A_2,B_2)$ by our present assumption]
\nl
hence
\sn
\item "{$(iii)^+$}"  $\bold w_\lambda(A_1,B_1) > \bold
w_\lambda(A_2,B_2)$.
\ermn
As the number of such $\lambda$ is finite and as we have shown
$\Xi(A_2,B_2) \subseteq \Xi(A_1,B_1)$ we get $\xi(A_1,B_1) >
\xi(A_2,B_2)$.
\nl
10) This follows from $\circledast_1 + \circledast_2$ below and the
finiteness of $\Xi(A,B_2)$ recalling
Definition \scite{1.8}(3).
\mr
\item "{$\circledast_1$}"  $\lambda \in \Xi(A,B_2) \Rightarrow
\lambda \restriction (B_1 \backslash A) \in \Xi(A,B_1)$ \nl
[Why?  If $\lambda$ is an equivalence relation on $B_2 \backslash A$
and let $\lambda_1 = \lambda \restriction (B_1 \backslash A)$ so $\lambda_1$
is an equivalence relation on $B_1 \backslash A$ and for any non-empty
$\lambda_1$-closed $C_1 \subseteq B_1 \backslash A$, letting $C_2 =
\cup\{x/\lambda:x \in C_1\}$ we have
\nl
$\bold w_\lambda(A,A \cup C_1) \ge \bold w_\lambda(A,A \cup C_2)$ by 
\scite{1.9}(4) and the latter is positive because $\lambda \in
\Xi(A,B_2)$.]  \nl
And
\sn
\item "{$\circledast_2$}"  $\lambda \in \Xi (A,B_2) \Rightarrow \bold
w_\lambda(A,B_2) < \bold w_\lambda(A,B_1)$.
\ermn
Why?  Otherwise let $\lambda$ be from $\Xi(A,B_2)$ and let
$C_\lambda = \cup\{x/\lambda:x \in B_1 \backslash A\} \cup A$ so $B_1 \le^*
C_\lambda \le^* B_2$.
\bn
\ub{Case 1}:  $C_\lambda = B_2$.

So $\bold v_\lambda(A,B_1) = \bold v(A,B_1,
\lambda \restriction (B_1 \backslash A)) =
|(B_1 \backslash A)/\lambda| = |(B_2 \backslash A)/\lambda| = \bold
v_\lambda(A,B_2,\lambda)$.
As $B_1 <^*_i B_2$, by part (11) below, for some $x \in B_1,y \in B_2
\backslash B_1$ the pair $\{x,y\}$ is an edge so
$e(A,B_1,\lambda \restriction (B_1
\backslash A))$ is a proper subset of $e(A,B_2,\lambda)$ hence

$$
\bold e_\lambda(A,B_1) < \bold e_\lambda(A,B_2)
$$
hence

$$
\bold w_\lambda(A,B_1) > \bold w_\lambda(A,B_2)
$$
\mn
is as required.
\bn
\ub{Case 2}:  $C_\lambda \ne B_2$.

As in case 1, $\bold w_\lambda(A,B_1) \ge \bold
w_\lambda(A,C_\lambda)$.
Now $B_1 \le^* C_\lambda \le^* B_2$, (using the case assumption) and
$B_1 <^*_i B_2$ by an assumption so by part (12) below we have
$C_\lambda \le^*_i B_2$ hence $C_\lambda <^*_i B_2$ and 
by \scite{1.14} this implies
$\bold w_\lambda(C_\lambda,B_2) < 0$.   So $\bold w_\lambda(A,B_2) 
= \bold w_\lambda(A,C_\lambda) + \bold w_\lambda(C_\lambda,B_2) 
\le \bold w_\lambda(A,B_1) + \bold w_\lambda(C_\lambda,B_2) < 
\bold w_\lambda(A,B_1)$ as required.
\sn
11) Define $\lambda$, it is the equivalence relation with exactly
one class on $B_2 \backslash B_1$ so $(B_1,B_2,\lambda) \in {\Cal
T},\bold v_\lambda(B_1,B_2) = 1,\bold e_\lambda(B_1,B_2) = 0$ so
$\bold w_\lambda(B_1,B_2) \ge 0$ hence $\lambda \in \Xi(B_2,B_2)$
hence $B_1 <_s B_2$. \nl
12) By the Definition \scite{1.11}(3). \nl
13) Clearly $A \cup \{a\} \le^* B$ hence by part (3) for some $C$ we
have $A \cup \{a\} \le^*_i C \le^*_s B$.  If $C=B$ we are done,
otherwise $A <^*_s C$ by part (7) so we have $A <^*_s C <^*_s B$,
contradiction. \nl
14) Easy by Definition \scite{1.11}(6).
  \hfill$\square_{\scite{1.17}}$\margincite{1.17}
\enddemo
\newpage

\head {\S5 The probabilistic inequalities} \endhead  \resetall \sectno=5
 \spuriousreset
\bigskip

In this section we deal with probabilistic inequalities about the number
of extensions for the context ${\Cal M}^0_n$.

Note: the proof of almost simple niceness of ${\frak K}$ is in the next
section.                                                              
\demo{\stag{2.1} Context}   As in \S4, so 
$p_i=1/{i^\alpha}$, for $i>1$, $p_1=p_2$ (where $\alpha\in (0,1)_{\Bbb
R}$ irrational) and ${\Cal M}_n = {\Cal M}^0_n$ (i.e. only the graph).
\enddemo
\bn
We may in the definition demand $f$ to be just a one to one function.
\definition{\stag{2.2} Definition}   Let $\varepsilon > 0,k \in \Bbb N,
{\Cal M}_n \in {\Cal K}$ and $A<^*B$ be in ${\Cal K}_\infty$.  Assume
$f:A \hookrightarrow {\Cal M}_n$ is an embedding or just $f:A
\hookrightarrow [n]$ which means it is one to one.  Define

$$
\align
{\Cal G}^{\varepsilon,k}_{A,B}(f,{\Cal M}_n) := \bigl\{\bar g:&(1) \quad 
\bar g = \langle g_\ell:\ell<k\rangle, \\
   &(2) \quad f\subseteq g_\ell, g_\ell \text{ a one-to-one function
from } B \text{ into } |{\Cal M}_n| \\
  &(3) \quad g_\ell:B \hookrightarrow_f {\Cal M}_n,
\text{ for } \ell \leq k 
\text{ or just } g_\ell:B \hookrightarrow_A {\Cal M}_n \\
  &\qquad \,\,\text{ which means}: 
\{a,b\} \in \text{ Edge}(B) \backslash \text{ Edge}(A) \\
  &\qquad \,\, \Rightarrow \{g(b),g(b)\} \in \text{ Edge}({\Cal M}_n) \\
  &\qquad \,\, \text{( and } g \text{ is one to one extending } f) \\
  &(4) \quad \ell_1 \ne \ell_2 \Rightarrow \text{ Rang}(g_{\ell_1})
\cap \text{ Rang}(g_{\ell_2}) = \text{ Rang}(f) \\
  &(5) \quad [\ell < k \and  x \in B \backslash  A \and y \in A] 
\Rightarrow |g_\ell(x)-g_\ell(y)|\ge n^\varepsilon \bigr\}.
\endalign
$$
\mn
The size of this set has natural connection with the number of pairwise
disjoint extensions $g: B \hookrightarrow {\Cal M}_n$ of $f$, hence 
with the holding of $A<_s B$, see \scite{2.3} below. 
\enddefinition
\bn
\ub{\stag{2.3} Fact}:  For every 
$\varepsilon$ and $k$ and $A \leq^* B$ we have: 
\mr
\item "{$(*)$}"  for every $n,k$ and
$M \in {\Cal K}_n$ and one to one $f:A \hookrightarrow_A M_n$ we have: 
if ${\Cal G}^{\varepsilon,k}_{A, B}(f,M_n) = \emptyset$ \ub{then}

$$
\align
\max\{\ell:&\text{ there are }g_m:B \hookrightarrow_A M
\text{ for } m < \ell \text{ such that } f\subseteq g_m \text{ and} \\
  &[m_1< m_2 \Rightarrow \text{ Rang}(g_{m_1}) \cap \text{ Rang}
(g_{m_2})\subseteq \text{ Rang}(f)]\} \\
  & \le 2(|A|) n^\varepsilon + (k-1)
\endalign
$$
\endroster
\bigskip

\demo{Proof}   Assume that there are $g_m$ for $m <\ell^*$ where
$\ell^* > 2|A|n^\varepsilon + k-1$ as above. \nl
By renaming \wilog \, for some $\ell^{**} \le \ell^*$ we have 
$\text{Rang}(g_m) \backslash \text{ Rang}(f)$ if $m < \ell^{**}$ 
is with distance $\ge n^\varepsilon$ from Rang$(f)$ but if 
$\ell \in [\ell^{**},\ell^*]$ then
Rang$(g_\ell) \backslash \text{ Rang}(f)$ has distance 
$< n^\varepsilon$ to Rang$(f)$.  Recall that 
by one of our assumptions $\ell^{**} \le k-1$.  Now 
for each $x \in \text{ Rang}(f)$, there are $\le 2 n^\varepsilon$
numbers $m \in [\ell^{**},\ell^*)$ such that min$\{|x-g_m(y)|:
y \in B \setminus A\} \le n^\varepsilon$.
So by the demand on $\ell^{**}$ we have $\ell^* - \ell^{**} \le 
(|A|) \times (2n^\varepsilon) =  2|A| n^\varepsilon$ and as 
$\ell^{**} < k$ we are done.  \hfill$\square_{\scite{2.3}}$\margincite{2.3}
\enddemo
\bigskip

The following is central, it does not yet prove almost niceness but the
parallels (to \scite{2.4}) from \cite{ShSp:304}, \cite{BlSh:528} were
immediate, and here we see main additional difficulties we have which
is that we are looking for copies $B$
over $A$ but we have to take into account the distance, the closeness of
images of points in $B$ under embeddings into ${\Cal M}_n$. Now, in order to
prove \scite{2.4} we will have to look for different types of g's
which satisfies Condition (5) from the Def. of 
${\Cal G}^{\varepsilon,k}_{A,B}(f,{\Cal M}_n)$ and restricting ourselves to one
kind we will calculate the expected value of ``relevant part" of
${\Cal G}^{\varepsilon,1}_{A,B}(f,{\Cal M}_n)$ and we will show that it is
small enough.
\bigskip

\proclaim{\stag{2.4} Theorem}  Assume $A<^*B$ (so both in ${\Cal
K}_\infty$).  \ub{Then} a sufficient condition for 
\mr
\item "{$\bigotimes_1$}"   for every $\varepsilon >0$, for some 
$k \in \Bbb N$, for every random enough ${\Cal M}_n$ we have:
{\roster
\itemitem{ $(*)$ }  if $f:A \hookrightarrow [n]$ then 
${\Cal G}^{\varepsilon,k}_{A,B}(f,{\Cal M}_n)= \emptyset$ 
\endroster}
\ermn
is the following:
\mr
\item "{$\bigotimes_2$}"   $A <^*_a B$ 
(which by Definition \scite{1.11}(5) means $A <^* B \and 
\neg(A <_s B)$).
\endroster
\endproclaim
\bigskip

\remark{\stag{2.5} Remark}  From $\bigotimes_1$ we can conclude: 
for every $\varepsilon \in \Bbb R^+$ we have:
for every random enough ${\Cal M}_n$, for every 
$f:A \hookrightarrow_A {\Cal M}_n$, there cannot be 
$\geq n^\varepsilon$ extensions $g:B \hookrightarrow_A {\Cal M}_n$ of 
$f$ pairwise disjoint over $f$. \nl
For this, first choose $\varepsilon_1 < \varepsilon$.
Note that for any $k$ we have
${\Cal G}^{\varepsilon,k}_{A,B}(f,{\Cal M}_n)
\subseteq {\Cal G}^{\varepsilon_1,k}_{A,B}(f,{\Cal M}_n)$. Choose $k_1$
for $\varepsilon_1$ by \scite{2.3}. Then the number of pairwise disjoint
extensions $g:B \hookrightarrow_A {\Cal M}_n$ of 
$f$ is $\le 2(|A|) n^{\varepsilon_1} + (k_1-1)$. 
For sufficiently large $n$ this is $<n^\varepsilon$.
\endremark
\bigskip

\remark{\stag{2.6} Remark}   We think of $g$ extending $f$ such that
$g:B \hookrightarrow {\Cal M}_n$ that satisfies, for some constants 
$c_1$ and $c_2,c_2 > 2 c_1$:

$$
x \lambda y \Rightarrow |g(x)-g(y)|<c_1
$$
and

$$
[\{x,y\} \subseteq B \and  \{x,y\} \nsubseteq A \and  
\neg x \lambda y] \Rightarrow |g(x)-g(y)| \ge c_2.
$$
\mn
So the number of such $g \supset f$ is $\sim n^{|(B \backslash
A)/\lambda|} = n^{\bold v(A,B,\lambda)}$, the probability of each 
being an embedding, assuming 
$f$ is one, is $\sim n^{-\alpha \bold e(A,B,\lambda)}$, hence 
the expected value is $\sim n^{\bold w_\lambda(A,B)}$ ($\sim$ means 
``up to a constant''). So $A<^*_i B$ implies that usually there are few
such copies of $B$ over any copy of $A$, i.e. the expected value is
$<1$. In \cite{ShSp:304}, $\lambda$ is equality, here things are more
complicated.
\endremark
\bigskip

\demo{Proof of \scite{2.4}}  We assume $\otimes_2$; and we shall prove
$\otimes_1$.
\enddemo
\bn
\ub{Stage A}:

Without loss of generality $B$ is minimal, i.e., 
$A<^*B'<^*B \Rightarrow \neg(A<_a^*B')$ which by Definition
\scite{1.11}(5) means $A <^* B' <^* B \Rightarrow A <^*_s B'$.

This is because if there exists a $B'$ contradicting this then 
$\otimes_1$ for $A$ and $B'$ implies $\otimes_1$ for $A$ and 
$B$.  By \scite{1.17}(4) we can conclude $A <^*_c B$.  
Clearly if $\bigotimes_1$ holds for $\varepsilon_1$ it holds for every
$\varepsilon_2 \ge \varepsilon_1$ so \wilog \, $\varepsilon > 0$ is
small enough, in particular moreover, we choose positive reals $\bold
c,\varepsilon,\zeta$ such that:
\mr 
\item "{$\boxtimes_1$}"  $(a) \quad$ if $A \le^* A' \le^* B$ 
and $(A',B,\lambda) \in {\Cal T}$ \ub{then} 
$\bold w_\lambda(A',B) \ne 0 \Rightarrow$ \nl

\hskip45pt $\bold w_\lambda(A',B) \notin [-2 \varepsilon,2 \varepsilon]$, a
sufficient condition is $t - \alpha r \notin [-2 \varepsilon,+2
\varepsilon]$ \nl

\hskip45pt when $t \in \{0,1,\dotsc,|B \backslash A|\}$ and 
$r \in \{1 \dotsc,|\text{edge}(B)|\}$,
\sn
\item "{${{}}$}"  $(b) \quad 0 < \varepsilon < 1$ and $0 < \zeta <
\varepsilon/2$
\sn
\item "{${{}}$}"  $(c) \quad \bold c > |B \backslash A|$.
\endroster
\bn
\ub{Stage B}:

In order to prove $\otimes_1$ it suffices to prove that for 
some $\zeta >0$ we have:
\mr
\item "{$(*)_\zeta$}"   Exp$(|{\Cal G}^{\varepsilon,1}_{A,B}
(f,{\Cal M}_n)|) \le \frac{1}{n^\zeta}$ if $f:A \hookrightarrow [n]$.
\ermn
(remember that Exp stands for ``expected value'').
\sn
To see that choose $k$ such that $k \cdot \zeta >|A|$.  Now
assuming $f:A \hookrightarrow [n]$ we have (the conditions
(1)-(5) are from Definition \scite{2.2}:

$$
\align
\text{Exp}(&|{\Cal G}^{\varepsilon,k}_{A,B}(f,{\Cal M}_n)|) = \\
  &\Sigma \bigl\{\text{Prob} (\dsize \bigwedge_{\ell<k} g_\ell:B
\hookrightarrow_A {\Cal M}_n):
\bar g \text{ satisfies clauses (1), (2), (4) and (5)
(of Def. \scite{2.2})}\} = \\
  & \quad \text{(so by disjointness and independence)} \\
  &\Sigma \bigl\{\dsize \prod_{\ell<k} \text{Prob}(g_\ell:B 
\hookrightarrow_A {\Cal M}_n):\bar g \text{ satisfies clauses
(1),(2),(4) and (5)}\bigr\} \le \\
  & \quad \text{(so as the sum has just more terms, and all terms are non
negative)} \\
  &\Sigma \bigl\{\dsize \prod_{\ell<k} \text{Prob}(g_\ell:B
\hookrightarrow_A {\Cal M}_n):\bar g \text{ satisfies clauses (1),(2) 
and (5)}\bigr\} = \\
  &\quad \text{(as now we have no demand connecting } g_\ell,g_m 
\text{ for } \ell,m < k) \\
  &\Sigma\bigl\{\dsize \prod_{\ell<k} \text{Prob}(g_\ell:B
\hookrightarrow_A {\Cal M}_n):\dsize \bigwedge_{\ell<k} [g_\ell 
\text{ satisfies clauses (2) and (5)}]\bigr\} = \\
  &\quad \text{(so since now all } k \text{-element sequences of ``allowable"
extensions of } f \text{ are allowed)} \\
  &\dsize \prod_{\ell<k} \Sigma \bigl\{\text{Prob}(g_\ell:B
\hookrightarrow_A {\Cal M}_n):g_\ell \text{ satisfies clauses 
(2) and (5)}\bigr\} = \\
  &\dsize \prod_{\ell<k} \text{ Exp}(|{\Cal G}^{\varepsilon,1}_{A,B}
(f,{\Cal M}_n)|\big) \leq ({1\over{n^\zeta}})^k.
\endalign
$$
\mn
By a well known easy lemma, for each $f:A \rightarrow \{1,\ldots,n\}$,
it follows that

$$
\text{Prob}(|{\Cal G}^{\varepsilon,k}_{A,B}(f,{\Cal M}_n)| \ne
\emptyset) < \frac{1}{n^{k\cdot\zeta}}.
$$
\mn
Therefore

$$
\align
\text{Prob}(\text{there is } f:&A \hookrightarrow [n] \text{ such that }
|{\Cal G}^{\varepsilon,k}_{A,B}(f,{\Cal M}_n)|\ne \emptyset) < \\
  &\frac{1}{n^{k\cdot\zeta}}|\{f:f:A \rightarrow [n]\}| \le
\frac{n^{|A|}}{n^{k\cdot\zeta}} \longrightarrow 0.
\endalign
$$
\mn
(So if $A$ has many edges there is a waste and if $A$ has no edges then this
is the right choice of $k$).
\bn
\ub{Stage C}:

We now show, with some effort, that without loss of generality
Rang$(f)=[|A|](=\{1,2,\ldots,|A|\})$ - this is certainly helpful, though
it isn't really necessary but it clarifies the situation.

To see that let $f_0$ be a bijection between $A$ and $[|A|]$. Divide
the possible $g$'s to finitely many classes by their types (defined below)
(``finite'' means here: with a bound depending on $A$ and $B$ but not 
on $n$).
For each type we will prove that the expected value for $f$ under the
assumption $f:A \hra [n]$ is $\le$ Const$(A,B)\times$
[the expected value for $f_0$ as required under the assumption 
$f_0:A \hra [n]$].  We may continue to ignore ``$f_0$ is an embedding 
of $A$''.
\bn
\ub{Notation}:   For any one to one $g:B \hra [n]$ we define the 
linear order $<^*_g$ on $B$: $x<^*_g y$ iff
$g(x) < g(y)$.\nl
Let $A=\{a_1,\ldots,a_{m-1}\}$ where $0<i<j < m \Rightarrow 
a_i<^*_f a_j$ (that is, if $g \in {\Cal G}^{\varepsilon,1}_{A,B}
(f,{\Cal M}_n)$ then
$g(a_i)< g(a_j)$).  We add two ``pseudo-elements'' $a_0$ and
$a_m$ to $A$ and stipulate $f(a_0)=0$, $f(a_m)=n+1$ so if $g:B \hra
[n]$ is one to one and extends $f$ then
$a_0<^*_g b<^*_g a_m$ for all $b \in B$ and for each $\ell< m$, we let 
the set $\{b \in B \backslash A:a_\ell<^*_g b<^*_g a_{\ell+1}\}$ 
be listed by $b^g_{\ell,1} <^*_g b^g_{\ell,2} \ldots <^*_g 
b^g_{\ell,{m^g_\ell-1}}$ and we stipulate $b^g_{\ell,0} = a_\ell$, 
$b^g_{\ell,{m^g_\ell}} = a_{\ell+1}$.

For $f$, $n$ , $\varepsilon$, $A$, $B$ define

$$
\align
{\Cal G}^{1,\varepsilon}_{A,B}(f,[n]):=\{g:&g:B \hra [n],\text{ i.e., g is 
one-to-one and} \\
   &g \text{ satisfies clauses (2) and (5)}\}.
\endalign
$$
We define ${\frak p} = \text{ tp}^0(g) = \text{ tp}^0_{A,B}(g)$ to be 
the following information (where the number of possibilities has a 
bound which depends only on $|B|$, and remember $\lfloor \beta \rfloor$ 
is the largest integer $\le\beta$):
\mr
\item "{$(a)$}"   $\{(x,y):x,y \in B,g(x)< g(y)\}$, i.e., the linear
order $<^*_g$ so this fixes the $b^g_{\ell,j} -s$ and $m^g_\ell - s$
\sn
\item "{$(b)$}"   $\{\langle \ell,j,\lfloor\frac{g(b^g_{\ell,j+1})-
g(b^g_{\ell,j})}{f(a_{\ell+1}) - f(a_\ell)} \cdot (m^g_\ell + 8)^3
\rfloor \rangle:\ell<m \text{ and } j < m^g_\ell\}$
\ermn
(note \footnote{We do not care about the exact bound}: the number of
possibilities in clause (b) is at most $\dsize \prod_{\ell<m}
\dsize \prod_{j<m_\ell}(m^g_\ell+8)^3 \le (|A| +1)^{3|B \backslash A|}$).

Fixing the type ${\frak p}$, we know $m$, 
$\langle m_\ell:\ell<m\rangle$, of course, $f$ is given hence 
we know $\langle a_\ell:\ell=1,2,\ldots,m-1\rangle$ and $a_0$, $a_m$
too and $\langle b_{\ell,j}:\ell < m,j < m_\ell \rangle$, i.e., they
are the same for all the $g$'s of this type, so we shall omit the
superscript $g$; if confusion may arrive we can write
$m^{\frak p}_\ell,b^{\frak p}_{\ell,j}$.  In
addition for each $\ell<m$, when $m_\ell>1$, we can compute some
$x_\ell^*<x_\ell^{**}$ and $j_\ell<m_\ell$, depending only on $f$, $n$ and
${\frak p}$, such that if $g \in {\Cal G}^{\varepsilon,1}_{A,B}(f,[n])$ and 
tp$^0(g) = {\frak p}$ then 
\mr
\item "{$(*)$}"   $f(a_\ell)\leq g(b_{\ell,{j_\ell}})\le x_\ell^*<x_\ell^{**}
\le g(b_{\ell,{j_\ell+1}})\leq f(a_{\ell+1})$ and $x_\ell^{**}-
x_\ell^* \ge$ \nl

$\frac{1}{2m_\ell} \cdot (a_{\ell+1}-a_\ell -2)$ and $(x^*_\ell=
a_\ell \vee x^{**}_\ell = a_{\ell+1}) \Rightarrow a_{\ell+1} \le
a_\ell+2$.
\ermn
For any $a \in [n]$ let $\ell(a) = \ell_{\frak p}(a) \in \{0,\dotsc,m\}$ be
the unique $\ell$ such that for every $t \le m$ 
\mr
\item "{$(\alpha)$}"  $m > t > \ell(a) \Rightarrow a \le (x^*_t +
x^{**}_t)/2$
\sn
\item "{$(\beta)$}"  $t \le \ell(a) \Rightarrow 
(x^*_t + x^{**}_t)/2 < a$.
\ermn
We define a linear order $<_0$ on $[n]$:
\nl
$c <_0 d \text{ \ub{iff} } (c,d \in [n] \text{ and})$
\mr
\item "{$(i)$}"  $|c - a_{\ell(c)}| < |d - a_{\ell(d))}|$ \ub{or}
\sn
\item "{$(ii)$}"  $(|c - a_{\ell(c)}| = |d - a_{\ell(d)}|) \and c < d$.
\ermn
Let $\sigma$ be the unique permutation of $[n]$ such
that for $c,d \in [n]$ we have $(c <_0 d) \equiv (\sigma(c) < \sigma(d))$.
It is easy to check that $\sigma(a_1)=1,\sigma(a_2) =
2,\dotsc,\sigma(a_{m-1}) = m-1$.  We stipulate $\sigma(a_0) =
a_0,\sigma(a_m) = a_m$. 
\bn
\ub{\stag{2.7} Fact}:  1) There is some $\bold c = \bold c_{\frak p}>0$ 
such that:
\block
if $i,j\in\{a_1,\ldots,a_{m-1}\} \cup \dbcu_{\ell=0}^{m-1}
(a_\ell,x_\ell^*] \cup \dbcu_{\ell=0}^{m-1}[x_\ell^{**},a_{\ell+1})$
then $|\sigma(i) - \sigma(j)| \le \bold c \times|i-j|$.
\endblock
\nl
2) If $j \in [n]$ and $k \le n$ and
$\dsize \bigwedge_{\ell=0}^m |j-a_\ell|\geq k$ \nl
then $\dsize \bigwedge_{\ell=1}^{m-1} |\sigma(j)-\sigma(a_\ell)|\geq
k$ (actually $\ge 2m(k-1)+1$).
\bigskip

\demo{Proof}   1) By symmetry \wilog \, $i<j$, and 
by the triangle inequality \wilog \,:
\mr
\item "{$\circledast_0$}"   $\ell(i) = \ell(j) \and i+1=j$ \nl
\ub{or} $(i,j) = (x^*_\ell,x^{**}_\ell)$ for some $\ell$.
\ermn
First assume that the first case in $\circledast_0$ occurs, then
$\bold c \ge 2m$ suffice. \nl
[Why?  If $a_{\ell(i)} \le i$ then (as $\ell(i+1) = \ell(j) = \ell(i)$
we have) $a_{\ell(i)} \le i < i+1 = j \le 
x^*_{\ell(i)}$ so ($\ell(i) = \ell(j)$
and) $i <_0 j$ hence $|j-a_{\ell(j)}| = |i-a_{\ell(i)}|+1$, so

$$
\align
|\sigma(j)-\sigma(i)| = \sigma(j)-\sigma(i) &= |\{t:i \le_0 t <_0 j\}| \\
  &=|\{t:|t-a_{\ell(t)}| = |i-a_{\ell(i)}| \text{ and } i \le t \\
  &\qquad \qquad \text{ \ub{or} }  |t-a_{\ell(i)}|  
= |j-a_{\ell(j)}| = |i-a_{\ell(i)}| + 1 
\text{ and } t<j\}| \le 2m.
\endalign 
$$
\mn
If $a_{\ell(i)} > i$ then $j = i+1 \le a_{\ell(i)} = a_{\ell(j)}$ then the
proof is similar.] \nl
Second, assume that the second possibility in $\circledast_0$ occurs so
$(i,j) = (x^*_\ell,x^{**}_\ell)$, so $\ell(i) = \ell \and \ell(j) =
\ell +1$.  The first subcase is to assume that $i <_0 j \and (x^*_\ell
\ne a_\ell \and x^{**}_\ell \ne a_{\ell +1})$  hence $a_{\ell +1} -
x^{**}_\ell \ge x^*_\ell - a_\ell$ and $\sigma(i) < \sigma(j)$ so

$$
\align
|\sigma(j)-\sigma(i)| &= |\{t:t \in [n] \text{ and } i \le_0 t <_0 j\}| \\
  &\le |\{t \in [n]:|i-a_{\ell(i)}| \le |t-a_{\ell(t)}| \le
|j-a_{\ell(j)}|\}| \\
  &= |\{t \in [n]:|x^*_\ell -a_\ell| \le |t-a_{\ell(t)}| \le |x^{**}_\ell -
a_{\ell +1}|\}| \\
  &= \dsize \sum_{s \le m} |\{t \in [n]:|x^*-a_\ell| \le |t-a_s| \le
|x^{**}_\ell- a_{\ell +1}| \text{ and } \ell(t)=\ell\}| \\
  &\le 2m \cdot |\{r:|x^*_\ell -a_\ell| \le r \le |x^{**}_\ell -
a_{\ell +1}|\}| \\
  &= 2m((a_{\ell +1}- x^{**}_\ell +1) - (x^*_\ell - a_\ell)) \\
  &\le 2m(a_{\ell +1} - x^{**}_\ell) \le 2m(a_{\ell +1} - a_\ell -1) \le
(2m)(2m_\ell)(x^{**}_\ell - x^*_\ell) +2m \\
  &= 4m m_\ell|j-i|+2m
\endalign
$$
\mn
so if $\bold c \ge 4m m_\ell +2m$ for $\ell < m$ this is O.K., too. 
The second subcase is $j <_0 i \and (x^*_\ell \ne a_\ell \and
x^{**}_\ell \ne a_{\ell +1})$.  The proof is similar and the 
cases $x^*_\ell = a_\ell,x^{**}_\ell = a_{\ell +1}$ are easier.
\nl
2) If $j \in \{a_1,\dotsc,a_{m-1}\}$ this is trivial as $k$ is
necessarily zero; so assume not hence $\sigma(a_\ell) < \sigma(j)$ for
$\ell =1,\dotsc,m-1$.  First consider $a_{\ell(j)} < j$ and $\ell \le m$ then
$|\sigma(j) - \sigma(a_\ell)| = \sigma(j) - \sigma(a_\ell) = |\{t \in
[n]:a_\ell \le_0 t <_0 j\}|$.
\nl
Now the set $\{t:a_\ell \le_0 t <_0 j\}$ includes
$a_\ell,a_{\ell(t)}+1,\dotsc,j-1$ hence has $\ge j-a_{\ell(t)}$ members so
$|\sigma(j)-\sigma(a_\ell)| \ge j-a_{\ell(t)} = |j-a_{\ell(t)}| \ge k$.  
Second consider $a_{\ell(j)} > j$, the 
proof is similar.  Together we are done. 
${{}}$ \hfill$\square_{\scite{2.7}}$\margincite{2.7}
\enddemo
\bigskip

\demo{Continuation of the proof of \scite{2.4}}   
Let ${\Cal G}^{\varepsilon,1,{\frak p}}_{A,B}(f,{\Cal M}_n) = \{g \in
{\Cal G}^{\varepsilon,1}_{A,B}(f,{\Cal M}_n):\text{tp}^0 (g) = {\frak p}\}$.

Now $\sigma$ induces a 1-to-1 function from 

$$
{\Cal G}^{\varepsilon,1,{\frak p}}_{A,B}(f,[n]) :=
\{g \in {\Cal G}^{\varepsilon,1}_{A,B}(f,[n]):\text{tp}^0(g) = {\frak p}\}
$$
\mn
into ${\Cal G}^{\varepsilon,1}_{A,B}(f_0,[n])$, 
that is $g\mapsto g\circ\sigma^{-1}$.
[Why? Because $\sigma$ is a permutation of $[n]$ and $f\circ \sigma^{-1}=f_0$
and \scite{2.7}(2).]

We note
\mr
\item "{$(*)$}"  $\text{Prob}(g \in {\Cal G}^{\varepsilon,1,{\frak p}}_{A,B}
(f,{\Cal M}_n)) \le \bold c^1_{\frak p} \cdot \text{Prob}(g \circ 
\sigma^{-1}\in {\Cal G}^{\varepsilon,1}_{A,B}(f_0,{\Cal M}_n))$
\ermn
for some constant $\bold c^1_{\frak p}$. \nl
(By \scite{2.7}(1); of course we may choose to consider also
the probability of non edges but this is easier as there is a constant upper
bound $(1-p_1) < 1$.)
So by the definitions of ${\Cal G}^{\varepsilon,1}_{A,B}(f_0,{\Cal M}_n),
{\Cal G}^{\varepsilon,1,{\frak p}}_{A,B}(f,[n])$,$f_0$, and $(*)$ we have

$$
\align
\text{Exp}(|{\Cal G}^{\varepsilon,1,{\frak p}}_{A,B}(f,{\Cal M}_n)|) &= \Sigma
\{\text{Prob}(g\in {\Cal G}^{\varepsilon,1}_{A,B}(f,{\Cal M}_n)):g \in 
{\Cal G}^{\varepsilon,1,{\frak p}}_{A,B}(f,[n])\} \\
  &\le \bold c^1_{\frak p} \Sigma\{\text{Prob}(g \circ \sigma^{-1} \in 
{\Cal G}^{\varepsilon,1}_{A,B}(f_0,{\Cal M}_n)):g \in
{\Cal G}^{\varepsilon,1,{\frak p}}_{A,B}(f,[n])\} \\
  &\le \bold c^1_{\frak p} \Sigma\{\text{Prob}(g'\in 
{\Cal G}^{\varepsilon,1}_{A,B}(f_0,{\Cal M}_n)):
g'\in {\Cal G}^{\varepsilon,1,{\frak p}}_{A,B}(f_0,[n])\} \\
  &\le \bold c^1_{\frak p} 
\text{ Exp}(|{\Cal G}^{\varepsilon,1}_{A,B}(f_0,{\Cal M}_n)|).
\endalign
$$
\mn
Also

$$
{\Cal G}^{\varepsilon,1}_{A,B}(f,[n]) =
\bigcup\{{\Cal G}^{\varepsilon,1,{\frak p}}_{A,B}(f,[n]):
{\frak p} \text{ is a tp}^0 \text{-type}\}
$$
\mn
and the number of tp$^0$--types has a bound not depending on $n$, say
$\bold c^2$.  Hence letting $\bold c^1 = \Sigma\{\bold c^1_{\frak p}:{\frak p}$
is a tp$^0$-type$\}$ it is constant and we have

$$
\align
\text{Exp}(|{\Cal G}^{\varepsilon,1}_{A,B}(f,{\Cal M}_n)|) &= 
\underset{\frak p} {}\to \Sigma \,
\text{ Exp}(|{\Cal G}^{\varepsilon,1,{\frak p}}_{A,B}(f,{\Cal M}_n)|) \\
  &\le \underset{\frak p} {}\to \Sigma \bold c^1_{\frak p} \text{ Exp}
(|{\Cal G}^{\varepsilon,1}_{A,B}(f_0,{\Cal M}_n)|) =
 \bold c^1 \text{ Exp}(|{\Cal G}^{\varepsilon,1}_{A,B}
(f_0,{\Cal M}_n)|).
\endalign
$$
\mn
So the expected value of ${\Cal G}^{\varepsilon,1}_{A,B}(f,{\Cal M}_n)$ 
is $\le$ a constant times the 
expected value of ${\Cal G}^{\varepsilon,1}_{A,B}(f_0,{\Cal M}_n)$.
Therefore it is enough to prove that the expected
value of the second value is small (i.e., $\leq n^{-\xi}$ for some
$\xi \in \Bbb R^{>0}$).\nl
Note: we are using heavily that we are dealing with expected values that are
additive, so dependencies between various $g$'s are irrelevant.
\enddemo
\bn
\ub{Stage D}:

By the previous stage it is enough for a fix $f_0$ to bound
Exp$({\Cal G}^{\varepsilon,1}_{A,B}(f_0,{\Cal M}_n))$ where $f_0$ is a
one-to-one mapping from $A$ onto $\{1,\dotsc,|A|\}$.

In this step we prove that the expected value of 
${\Cal G}^{\varepsilon,1}_{(A,B)}(f_0,{\Cal M}_n)$ has an appropriate
upper bound.

First we will define a new type of $g$, tp$^1(g)$ for $g:B \hra_A {\Cal M}_n$
(assuming $g\supseteq f_0$).  Now we let tp$^1(g)$ includes the following
information: $(\bar b^g,\sigma[g])$ is defined in $(\alpha), (\beta)$ below
\mr
\item "{$(\alpha)$}"   let 
$\bar b^g = \langle b^g_1,\dotsc,b^g_m \rangle$ list $B \backslash
A$ such that $g(b^g_1) < g(b^g_2) < \ldots < g(b^g_m)$.
\ermn
Stipulate $b_0\in A$ being $f_0(b_0) = |A|$, so if $A \ne \emptyset$
then $f_0(b_0)$ is the maximal element
in Rang$(f_0)$, and stipulating $b_{m+1}\in B$
with $g(b_{m+1})=n+1$.
\mr
\item "{$(\beta)$}"   essentially the order between $r_\ell =
r_\ell[g] =: g(b^g_{\ell+1}) - g(b^g_\ell)$ for $\ell\leq m$, i.e., 
the truth value of each statement $r_{\ell_1} < r_{\ell_2}$, 
but pedantically the value of $\sigma[g]$ defined below.
\ermn
Let $\sigma = \sigma[g] = \sigma(\text{tp}^1(g))$ be 
a permutation on $\{0,1,\ldots,m-1\}$ (without $m$, as below
$r_m=n+1-\sum\limits_{\ell<m} r_\ell-(m-1)$ so no need to 
say what is $r_m$) such that 
$r_{\sigma(0)} \le r_{\sigma(1)} \le \ldots \le 
r_{\sigma(m-1)}$ and under those restraints, is the first such
permutation by lexicographic order (by fixing $\sigma$ we restrict
the order of the distances among the vertexes). Note that if we use the first
version of $(\beta)$ when $r_\ell = r_i$ we lose a degree
of freedom, this of course does not really matter.
Let $\bar r[g] = \langle r_\ell[g]:\ell<m\rangle$.
Let a tp$^1$-type mean tp$^1(g)$ for some $g$ and for a tp$^1$-type
${\frak p},(\bar b^{\frak p},\sigma[{\frak p}])$ is well defined and
we let $<_{\frak p}$ be the linear order on $B \backslash A$ in which $1
\le \ell < k \le m \Rightarrow b^{\frak p}_\ell < b^{\frak p}_k$ and let

$$
\align
\Gamma_*[{\frak p}] = \Gamma^n_*[{\frak p}] := \{\bar r:&\bar r=\langle
r_0,\ldots,r_{m-1}\rangle,r_{\sigma(0)}\le
r_{\sigma(1)}\le \ldots \le r_{\sigma(m-1)}, \\
   &\text{each } r_\ell \text{ is an integer } >0 \text{ and} \\
  &r_0 \ge n^\varepsilon \text{ and } |A|+
\dsize \sum_{\ell<m} r_\ell\le n\}
\endalign
$$

$$
\align
\Gamma[{\frak p}] = \Gamma^n[{\frak p}] = \{\bar r \in
   \Gamma_*[{\frak p}]:&\text{ for every } \ell < m-1 \text{ we have} \\
   & r_{\sigma(\ell)} = r_{\sigma(\ell+1)}
\Rightarrow \sigma(\ell) < \sigma(\ell+1)\}
\endalign
$$

$$
\Gamma = \Gamma^n = \dbcu_{\frak p} \Gamma^n[{\frak p}]
$$

$$
{\Cal G}^{1,\varepsilon,{\frak p}}_{A,B}(f_0,{\Cal M}_n) =
\{g \in {\Cal G}^{1,\varepsilon}_{A,B}(f_0,{\Cal M}_n):\text{tp}^1(g)
= {\frak p}\}
$$

$$
{\Cal G}^{1,\varepsilon,{\frak p}}_{A,B}(f_0,[n]) =
\{g:\text{ for some } M \in {\Cal K}_n \text{ we have } g \in 
{\Cal G}^{1,\varepsilon,{\frak p}}_{A,B}(f_0,M)\}.
$$
\mn
The $\Gamma[{\frak p}]$'s are necessarily pairwise disjoint and
$\Gamma^n[{\frak p}] \subseteq \Gamma^n_*[{\frak p}] \subseteq \Gamma^n$.  
Note that $r_{\sigma(\ell)} = r_{\sigma(\ell +1)}
\Rightarrow \sigma(\ell) < \sigma(\ell +1)$ is to be consistent with
the way we chose $\sigma[g]$.

Note, for each $\bar b^* = \langle b^*_1,\dotsc,\bar b^*_m)$ listing
$B \backslash A$ and $\bar r \in \Gamma^n$ there is one and
only one $g = g_{\bar b^*,\bar r} = g^n_{\bar b^*,\bar r} \in
{\Cal G}^{\varepsilon,1}_{A,B}(f_0,[n])$ such that 
$\bar b^g = \bar b^* \and \bar r(g) = \bar r$ and
${\Cal G}^{\varepsilon,1}_{A,B}(f_0,[n]) = \{g^n_{\bar b^*,\bar r}:
\bar r \in \Gamma^n$ and $\bar b^* = \langle b^*_1,\dotsc,b^*_m
\rangle$ list $B \backslash A\}$. \nl
[Why? Let $g^n_{\bar b^*,\bar r}(b^*_\ell) = |A| + \Sigma\{r_m: m<\ell\}$.  
So $g$ is one to one and, moreover, $g$ preserves the order $<_{\frak p}$ (as
each $r_m$ is an integer $>0$) also $a \in A \and b \in B
\backslash A \Rightarrow g(b) - g(a) \ge n^\varepsilon$ as the worst
case is $a = |A| \and b = b^*_1$ so $g(b) - g(a) = (|A| + r_0) -
|A| = r_0 \ge n^\varepsilon$ by one of the demands in the
definition of $\Gamma$ is $``r_0 \ge n^\varepsilon"$.]
\bn
\ub{Stage E}:

By stage D we fix a tp$^1$-type ${\frak p}$ and will bound
Exp$({\Cal G}^{1,\varepsilon,{\frak p}}_{A,B}
(f_0,{\Cal M}_n))$.  As $\bar b^{\frak p}$
is constant we rename it as $\langle b_1,\dotsc,b_m \rangle$ and let
$\sigma = \sigma[{\frak p}]$. \nl
In order to compute the probability that $g$ is a graph embedding, we
introduce the following notation. Let

$$
\align
E_B = \{(c,d):&c \in B,d \in B,\{c,d\} \nsubseteq A \text{ and either }
c \in A \and  d \in B \setminus A \text{ or} \\
   &\bigvee\{c=b_{\ell_1} \wedge d=b_{\ell_2}:1 \le \ell_1<\ell_2 \le m\}\},
\endalign
$$

$$
E^0_B = \{(c,d)\in E_B:\{c,d\} \in \text{ edge}(B) \text{ and for no
}\ell,(c,d) = (b_\ell,b_{\ell +1})\},
$$

$$
E^1_B = |\{(c,d) \in E_B:(c,d) \notin \text{ edge}(B)\}
$$
\mn
and

$$
E^2_B = E_B \setminus E^0_B \backslash E^1_B.
$$
\mn
Note considering whether $g \in {\Cal G}^{1,\varepsilon,
{\frak p}}_{A,B}(f_0,[n])$ belong to ${\Cal
G}^{1,\varepsilon}_{A,B}(f_0,{\Cal M}_n)$, we have to
consider all unordered pairs $\{c,d\} \subseteq B,\nsubseteq A$, but
$E_B$ contains exactly one ordered pair among $(c,d),(d,c)$.
For $(c,d)\in E_B$ let $w(c,d)$ be defined as follows: if $(c,d)=(b_{\ell_1},
b_{\ell_2})$ and $1\leq\ell_1<\ell_2\leq m$ then $w(c,d)=\{\ell_1,\ell_1+1,
\ldots,\ell_2-1\}$ and if $c\in A$, $d=b_{\ell_2}$ then $w(c,d)=\{0,1,2,
\ldots,\ell_2-1\}$. So $w(c,d)$ is a non empty subset of $\{0,\ldots,m-1\}$,
in fact an interval of it. Thus $\sum\limits_{\ell\in
w(c,d)} r_\ell$  essentially is the distance between $g(c)$ 
and $g(d)$ (only if $c\in A$ we should add
$g(b_0)-g(c)$ but it does not matter as this number is $\le |A|$, 
it contributes just a little to the first
inequality on $\beta^*$ below).  Let $k(c,d)$ be the 
$k \in w(c,d)$ such that $(\sigma[g])(k)$ is maximal, that is $r_k$
maximal and if there are several such candidates, the last such $k$.
So $k(c,d)$ is computable from ${\frak p}= \text{ tp}^1(g)$.  

As we are fixing the 
(new) type ${\frak p}$, hence $\sigma = \sigma[{\frak p}]$, we have:
\nl
(where we assume $f_0:A \hra {\Cal M}_n$)

$$
\align
\beta^*_{\frak p} &=: \text{ Exp}(|\{g\in 
{\Cal G}^{\varepsilon,1}_{A,B}(f_0,{\Cal M}_n):
\text{tp}^1(g) = 
{\frak p} \text{ hence } \bar r[g]\in\Gamma^n[{\frak p}]\}|) \\
  & = \dsize \sum_{\bar r \in\Gamma[{\frak p}]} \text{ Exp}
(|\{g_{\bar r} \in {\Cal G}^{\varepsilon,1}_{A,B}(f_0,{\Cal M}_n):
\text{tp}^1(g) = {\frak p} \text{ and } \bar r[g]=\bar r\}|,  \\
   & = \dsize \sum_{\bar r \in\Gamma[{\frak p}]}
\text{Prob}(g_{\bar r}\in
{\Cal G}^{\varepsilon,1}_{A,B}(f_0,{\Cal M}_n))
\endalign
$$
\mn
and this implies, letting $\beta'_{\frak p} = 1/2^{\alpha|E^2_\beta|}$
(because $p_1$ is $p_2$ and not $\frac{1}{1^\alpha}$ so essentially
can be replaced by 1):

$$
\beta^*_{\frak p} \le \beta'_{\frak p} \cdot
\dsize \sum_{\bar r \in\Gamma[{\frak p}]}
\dsize \prod_{(c,d)\in E^0_B}\bigr(\frac{1}{(\sum_{\ell\in
w(c,d)} r_\ell)^\alpha}\bigl)
\le \beta'_{\frak p} \dsize \sum_{\bar r\in\Gamma[{\frak p}]}
\big(\dsize \prod_{(c,d)\in E^0_B}
\frac{1}{(r_{k(c,d)})^\alpha}\big).
$$
\mn
This means, $\beta^*_{\frak p} \le \beta'_{\frak p} \cdot
\dsize \sum_{\bar r \in\Gamma[{\frak p}]} 
\dsize \prod_{\ell\in [0,m)}\frac{1}{r^{\alpha\cdot j(\ell)}_\ell}$ 
where we let
\mr
\item "{$\circledast_1$}"  $j(\ell)=:|\{(c,d) \in E^0_B:k(c,d) =\ell\}|$
so neccessarily $j(m)=0$; actually we should write $j(\ell,{\frak p}))$.
\ermn
For $\bar r \in\Gamma$ let $\bar s(\bar r):=\langle
s_\ell(r_\ell):\ell<m\rangle$ where $s_\ell(r_\ell)
= \lfloor\log_2(r_\ell) \rfloor$ and if $\bar g \in {\Cal
G}^{1,\varepsilon,{\frak p}}_{A,B}(f_0,[n])$ let $\bar s(g) = \bar
s(\bar r)$ when $g = g_{\bar b^{\frak p},\bar r}$.  \nl
Define

$$
\align
S =:\{\bar s:&\bar s = \langle s_\ell:\ell<m\rangle,0 \le s_{\sigma(0)}\le
s_{\sigma(1)} \le \ldots \le \lfloor\log_2(n)\rfloor  \\
  &\text{ each } s_\ell \text{ is of course an integer} \\   
  &\text{and } s_{0} \ge \lfloor\log_2(n^\varepsilon) \rfloor\}.
\endalign
$$
\mn
Why we move from $\bar r$ to $\bar s$?  As there are not too many
possible sequences $\bar s$, i.e., $\le (\log_2 n)^m$, but they are
enough in giving the bounds. \nl
So by the above we have:

$$
\beta^*_{\frak p} \le \beta'_{\frak p} \cdot \dsize \sum_{\bar s\in S} \, 
\dsize \sum_{\bar r \in\Gamma[{\frak p}],\bar s(\bar r)=
\bar s}(\dsize \prod_{\ell\in[0,m)}(\frac{1}{r_\ell^{\alpha\cdot
j(\ell)}})) = \beta'_{\frak p} \cdot 
\dsize \sum_{\bar s\in S}\beta_{{\frak p},\bar s}
$$
\mn
where $\beta_{{\frak p},\bar s} := 
\dsize \sum_{\bar r\in\Gamma[{\frak p}],\bar s
(\bar r) = \bar s}(\dsize \prod_{\ell\in[0,m)}
(\frac{1}{r_\ell^{\alpha\cdot j(\ell)}}\big))$.  Now

$$
\lfloor\log_2(r_\ell)\rfloor = 
s_\ell(r_\ell)\leq\log_2(\gamma_\ell)<s_\ell(r_\ell)+1
$$
\mn
and hence

$$
\frac{1}{r_\ell^{\alpha\cdot j(\ell)}}=\frac{1}{2^{(\log_2(r_\ell))
\cdot\alpha\cdot j(\ell)}}\leq
\frac{1}{2^{s_\ell(r_\ell)\cdot\alpha
\cdot j(\ell)}}
$$
\mn
and

$$
[\bar r \in\Gamma[{\frak p}] \and \bar{s}(\bar r)=\bar{s}]
\Rightarrow \dsize \prod_{\ell\in [0,m)}\frac{1}{(r_\ell)^{\alpha\cdot
j(\ell)}} \le \dsize \prod_{\ell\in [0,m)}\frac{1}{2^{s_\ell(r_\ell)
\cdot\alpha\cdot j(\ell)}}.
$$
\mn
Note that it is a common bound to all 
$\bar r \in\{\bar r\in\Gamma[{\frak p}]:
\bar{s}(\bar r) = \bar{s}\}$.  \nl
Hence

$$
\beta_{{\frak p},\bar s} \le 
|\big\{\bar r \in\Gamma[{\frak p}]:\bar s(\bar r)=\bar s\big\}|
\cdot\prod_{\ell\in [0,m)}\big(\frac{1}{2^{s_\ell(r_\ell)
\cdot \alpha\cdot j(\ell)}}\big).
$$
Now

$$
\align
|\{\bar r \in\Gamma[{\frak p}]:\bar{s}(\bar r)=\bar{s}\}| &\le
|\{\bar r: \bar r =\langle r_\ell:\ell<m \rangle
\text{ and } \lfloor\log_2(r_\ell)\rfloor = s_\ell\}| \\
  &= \dsize \prod_{\ell<m}|\{r:\lfloor\log_2(r) \rfloor =
s_\ell\}| \\
  &= \dsize \prod_{\ell<m}|\{r:2^{s_\ell} \le r < 2^{s_\ell+1}\}|\} \\
  &= \dsize \prod_{\ell<m}(2^{s_\ell+1}-2^{s_\ell})
\endalign
$$
\mn
and hence

$$
\align
\beta_{{\frak p},\bar s} &\le \dsize \prod_{\ell\in
[0,m)}(2^{s_\ell+1}-2^{s_\ell})
\cdot \dsize \prod_{\ell\in[0,m)}(\frac{1}{2^{s_\ell\cdot\alpha\cdot
j(\ell)}}) \\
  &= \dsize \prod_{\ell\in[0,m)}(2^{s_\ell}) \cdot
\dsize \prod_{\ell\in [0,m)}(\frac{1}{2^{s_\ell\cdot\alpha\cdot 
j(\ell)}}) \\
  &= \dsize \prod_{\ell\in[0,m)}(2^{s_\ell(1-\alpha \cdot j(\ell))}) \\
  &= 2^{\underset {\ell \in[0,m)} {}\to \Sigma s_\ell(1-\alpha\cdot j(\ell))}.
\endalign
$$
\mn
Choose $\bar s^*\in S$ which is best in the following sense: the 
last expression is maximal and if
there is more than one candidate, demand further that
$|\text{Rang}(\bar s^*) \backslash \{0,\lfloor
\log_2(n^\varepsilon)\rfloor,\lfloor \log_2 n \rfloor\}|$ is
minimal (among the candidates). Note that for each
$\bar s\in S$, hence also for this $\bar s^*$, every $s^*_\ell$ is restricted
to the interval $[0,\lfoot\log_2(n)\rfoot]$ when $\ell>0$ and to the interval
$[\lfoot\log_2(n^\vep)\rfoot,\lfoot\log_2(n)\rfoot]$ 
when $\ell=0$ and, of course, $s_{\sigma(\ell)} \le 
s_{\sigma(\ell+1)}$. We claim that
\mr
\item "{$\circledast_2$}"  every $s^*_\ell$ belongs to the set 
$\{0,\lfoot\log_2(n^\varepsilon)\rfoot,\lfoot\log_2(n)\rfoot\}$. 
\ermn
[Why? Let $E = E_{\frak p}$ be an equivalence relation on $[0,m)$ such that

$$
{\ell}E{k} \quad \text{  if and only if } \quad s^*_\ell=s^*_k.
$$
\mn
For $\langle \sigma(\ell):\ell<m\rangle$ this equivalence relation 
$E$ is convex (i.e. if $\sigma(\ell_1)< \sigma(\ell_2) <
\sigma(\ell_3)$ and $\sigma(\ell_1) E \sigma(\ell_3)$ then 
$\sigma(\ell_1) E \sigma(\ell_2)$). \nl
[Why?  Because if $\sigma(\ell_1)< \sigma(\ell_2)< \sigma(\ell_3)$ then 
$[\text{log}_2(r_{\sigma(\ell_1)})] \leq 
[\text{log}_2(r_{\sigma(\ell_1)})] \leq
[\text{log}_2(r_{\sigma(\ell_3)})]$ but $\sigma(\ell_1)E \sigma(\ell_3)$
which means that $\lfoot \text{log}_2(r_{\sigma(\ell_3)})\rfoot=
\lfoot \text{log}_2(r_{\sigma(\ell_1)})\rfoot$ so 
we obtain $\lfoot \text{log}_2(r_{\sigma(\ell_1)})\rfoot =
\lfoot \text{log}_2(r_{\sigma(\ell_2)})\rfoot$ by using
inequality above which means $\sigma(\ell_1)$ E $\sigma(\ell_2)$.] \nl

As was said above we are dealing with $tp^1(g) = {\frak p}$ and so we
fix the order on the distances
$r_{\ell}$'s by $\sigma$ and we want that under this order $E$ 
``will act nicely" which means $E$ will satisfies convexity on 
$\langle \sigma(\ell):\ell < m \rangle$.  So each equivalence class 
can be identify with one number that belongs to the 
interval $[0,\text{log}_2(n^{\varepsilon})]$ or the interval $[\text{log}
(n^{\varepsilon})],[\text{log}_2(n)]$ so we actually have natural 
order on the set of equivalence classes. \nl

Suppose that some $s^*_\ell$ is not one of the 
$\{0,\lfoot\text{log}_2(n^{\vep})\rfoot,\lfoot\text{log}_2(n)\rfoot\}$. 
\nl
Now take the last $E$-equivalence class 
with a value that is different from our three desired values and move
this block up or down to get a better expression. 
That is, let $\ell$ be maximal such that $s^*_\ell \notin \{0,\lfoot
\text{log}_2(n^\varepsilon)\rfoot,\lfoot\text{log}_2(n)\rfoot\}$, and
for $i \in \{1,-1\}$ define $\bar s^i$ by $s'_k$ is $s_k$ if $\neg(k E
\ell)$ and $s'_k = s_k +i$ if $k E \ell$.  Clearly $\bar s^1,\bar
s^{-1} \in S$, (and we do not care that it does not correspond to some
$g$ with tp$^1(g) = {\frak p}$)
but the expression above for $\bar s^1$ or for $\bar
s^{-1}$ is strictly smaller than the one for $\bar s^*$ except when $1
- \alpha \cdot \Sigma\{j(k):k \in \ell/E\} = 0$ but 
then $\alpha = 1/\Sigma\{j(k):k \in \ell/E\}$ so 
$\alpha$ is rational, contradiction; so $\circledast_2$ holds.]

So the number we are estimating is bounded from above by (more see
\scite{2.8} below)
\mr
\item "{$\circledast_3$}"  $\beta^*_{\frak p} \le 
\dsize \sum_{\bar s \in S} \beta_{{\frak p},\bar s} \le \dsize
\sum_{\bar s \in S} 2^{\Sigma\{s_\ell(1 - \alpha \cdot j(\ell))\}} \le
|S| \times 2^{\Sigma\{s^*_\ell(1 - \alpha \cdot j(\ell))\}} \le
(\log_2(n)+1)^m \times 2^{\Sigma\{s^*_\ell(1 - \alpha \cdot j(\ell))\}}$.
\ermn
The third inequality follows from the choice of $\bar s^*$ and the
second inequality follows from the 
fact that $\beta_{{\frak p},\bar s} \leq 
2^{\Sigma_{\ell \in [0,m)}s_\ell(1-\alpha(j(\ell))}$. 
For the fourth inequality we have to note that $|S| \leq (\log_2(n))^m$.
\bn
\ub{Stage F}:   Now we show that for the 
$\bar{s}^*$ (chosen above) $\beta'_{{\frak p},\bar{s}^*} \times
2^{\Sigma\{s^*_\ell(1 - \alpha \cdot j(\ell))\}} \le
\frac{1}{n^\zeta}$ for the appropriate $\zeta>0$ (this will suffice by
$\circledast_3$ above). \nl
Define two equivalence relations on 
$B \backslash  A = \{b_1,\dotsc,b_{m}\}$:
\mr
\item   $b_i\lambda_1 b_j$ \ub{if and only if} 
we have $i=j$ or $i<j \and \dsize \bigwedge_{i\le\ell\le j-1} 
s^*_\ell \in \{0,\lfoot\log_2(n^\vep)\rfoot\}$ or
$j < i \and \dsize \bigwedge_{j\le\ell\le i-1}
s^*_\ell \in \{0,\lfoot\log_2(n^\varepsilon)\rfoot\}$; because of
$\circledast_2$ this is an equivalence relation which by its
definition is convex
for the natural ordering $<_{\frak p}$ on $\{b_1\ldots b_{m}\}$ \nl
(that is each $i/E$ is a convex subset of $\{1,\dotsc,m\}$ (recall
that we have agreed that $\bar b^{\frak p} = \langle b_1,\dotsc,b_m \rangle$))
\sn
\item  $b_i \lambda_0 b_j$ \ub{if and only if}  we have $i=j$ or 
$i < j \and \dsize \bigwedge_{i\le\ell\le j-1}s^*_\ell =0$ or 
$j<i \and \dsize \bigwedge_{j\le\ell\le i-1} s^*_\ell=0$; 
this is an equivalence relation, convex for the natural ordering 
on $\{b^{\frak p}_1 \ldots b^{\frak p}_{m}\}$ and it refines $\lambda_1$
\sn
\item  Let $C_1<C_2<\ldots < C_t$ be the $\lambda_1$-equivalence classes,
ordered naturally by $<_{\frak p}$.
\sn
\item   Let $\lambda_2 = \lambda_1 \restriction (\dbcu_{\ell =2}^t C_\ell)$.
\endroster
\bn
\ub{Case 1}:   $t>1,\bold w_{\lambda_2}(A\cup C_1,B)<0$
and $s^*_0 = \lfoot\log_2(n^\varepsilon)\rfoot$. \nl
So for $\bar r \in \Gamma[{\frak p}]$ or just $\bar r \in
\Gamma'_*[{\frak p}]$ and $g=g_{\bar b^{\frak p},\bar r}$ 
with $s(g) = \bar s^*$, the set
$\{g(b_\ell):\ell \in C_1\}$ is ``near 0, i.e., near $f_0(A)$" 
(that is of distance $\sim
n^\varepsilon$, as $s^*_0 = \lfoot\log_2(n^\varepsilon)\rfoot$) but
$\dbcu^t_{\ell=2} C_\ell$ is not empty (as $t>1$) and the nodes in
$\{g(b):b \in \dbcu^t_{\ell=2} C_\ell\}$ are quite ``far" from zero.
\mn
Note:  $s^*_\ell = \lfoot\log_2(n)\rfoot$ \ub{iff} $s^*_\ell \notin
\{0,\lfoot\log_2(n^\varepsilon)\rfoot\}$ \ub{iff} for some
$r\in \{1,\ldots, t-1\}$ we have:

$$
\{\ell,\ell+1\} = \{\text{the maximal }i \in C_r, \text{ the minimal }
j\in C_{r+1}\}.
$$
\mn
Hence $\{\ell < m:s^*_\ell = \lfoot\log_2 n \rfoot\}$ 
has exactly $t-1$ members which is $\bold v_{\lambda_2}
(A \cup C_1, B)$; also recalling the definition of $j(\ell)$ in
$\circledast_1$ we get $\sum \{j(\ell):\ell < m$ and 
$s^*_\ell = \lfoot\log_2 n \rfoot\}$
is exactly $\bold e_{\lambda_1}(A \cup C_1,B)$ hence together
$\Sigma\{1-\alpha j(\ell):
\ell <m$ and $s^*_{\ell} = \lfoot\log_2 n \rfoot\}$ is 
equal to $\bold w_{\lambda_2}(A \cup C_1,B)$.

Recalling from stage A that $\bold c,\zeta \in \Bbb R^+$ depend on 
$A$, $B$ only and $\bold c \ge \Sigma\{1-\alpha \cdot j(\ell):
\ell$ satisfies $s^*_\ell \in 
\{0,\lfoot\log_2(n^{\varepsilon})\rfoot\}$ because $\bold c \ge |B
\backslash A|$ and we have

$$
\align
\beta'_{{\frak p},\bar s^*} &\le 2^{\sum_\ell s^*_\ell\cdot(1-\alpha \cdot
j(\ell))} \\
  &= 2^{\Sigma\{s^*_\ell\cdot(1-\alpha\cdot j(\ell)):\ell\text{ satisfies }
s^*_\ell\in\{0,\lfoot{\log}_2(n^\varepsilon)\rfoot\}\}} \\
  &\times 2^{\Sigma\{s^*_\ell\cdot(1-\alpha\cdot j(\ell)):
\ell \text{ satisfies } s^*_\ell = \lfoot{\log}_2(n)\rfoot\}} \\
  &\le n^{\bold c \varepsilon} \times n^{(t-1) - 
\alpha \bold e_{\lambda_2}(A\cup C_1,B)} \\
  &= n^{\bold c\cdot \varepsilon} \times 
n^{\bold w_{\lambda_2}(A\cup C_1,B)}\le\frac{1}{n^{ 2 \zeta}},
\endalign
$$
\mn
and we are done (by the assumption of the case 
$\bold w_{\lambda_2}(A\cup C_1,B) <0$, 
hence as by $\boxtimes_1$ from Stage A,  
$\varepsilon>0$ is small enough such that
$- \bold w_{\lambda_2}(A\cup C_1,B)/2 > |\bold c| \cdot \varepsilon,
\zeta \le (- \bold w_{\lambda_2}(A\cup C_1, B)/2$ all is O.K.).
\bn
\ub{Case 2}:   $s_0^* = \lfoot\log_2(n)\rfoot$ (here we use
$\lambda_1$ for the relevant $g$'s, all members of $\{g(b):b \in B
\backslash A\}$ are far from zero).
\sn
Now $\bold w_{\lambda_1}(A,B)<0$ (as $A <_c^* B$ by stage A; see
Definition \scite{1.11}(2)) and the proof is similar to that of case 1,
but here we will have $\bold v_{\lambda_1}(A,B)=t$ by the definition
of $\lambda_1$.
\bn
\ub{Case 3}:   $t>1$ and $\bold w_{\lambda_2}(A \cup C_1,B) \ge 0,
s_0^* = \lfoot\log_2(n^\varepsilon)\rfoot$.
\mn
We will try to get a contradiction to the minimality of the set $B$ from
stage A, and let us try $B'=A\cup C_1$.

Let $\lambda$ be an equivalence relation on $C_1$, by \scite{1.16} it
suffice to prove that $\bold w_\lambda(A,A\cup C_1)<0$.
\sn

Let $\lambda^*$ be the equivalence relation on $B \backslash A$ defined by:
$C_2,\ldots,C_t$ are equivalence classes of $\lambda^*$, no $x\in C_1$ and
$y \in \dbcu^t_{\ell=2} C_\ell$ are equivalent, and $\lambda^*
\restriction C_1=\lambda$, clearly $(A,B,\lambda^*) \in {\Cal T}$.
Now (first inequality as $A <^*_c B$ which holds by Stage A) we have

$$
\align
0> \bold w_{\lambda^*}(A,B) &= |(B \backslash A)/\lambda^*|
-\alpha \cdot \bold e_{\lambda^*}(A,B) \\
  &= |(B \backslash C_1 \backslash A)/\lambda^*| +
|C_1/\lambda^*| - \alpha \cdot \bold e_{\lambda^*}(A,B) \\
  &= |(B \backslash (C_1 \cup A))/\lambda_2| \\
  &+ |C_1/\lambda| - \alpha \cdot (\bold e_\lambda(A,A \cup C_1) +
\bold e_{\lambda_2}(A\cup C_1,B)) \\
  &= [|(B \backslash (C_1\cup A))/\lambda_2| - 
\alpha \cdot \bold e_{\lambda_2}(A\cup C_1,B)] \\
  &+ [|C_1/\lambda| - \alpha\cdot \bold e_\lambda(A,A\cup C_1)] \\
  &= \bold w_{\lambda_2}(A \cup C_1,B) + \bold w_\lambda(A,A\cup C_1).
\endalign
$$
\mn
Therefore $0 > \bold w_{\lambda_2}(A\cup C_1,B) +
\bold w_\lambda(A,A\cup C_1)$. By the assumption of the present 
case 3 we have $0\leq \bold w_{\lambda_2}(A\cup C_1,B)$ hence 
$0 > \bold w_\lambda(A,A\cup C_1)$. As $\lambda$ was any equivalence 
relation on $C_1$ we get $A<_c^*A\cup C_1$ so by \scite{1.17}(6), 
$\neg(A <^*_s A\cup C_1)$, a contradiction as desired (for this case).
\bn
\ub{Case 4}:  $s_0^* = \lfoot\log_2(n^\varepsilon)\rfoot$ and $t=1$.

This is like case 2, using $\lambda_0$ and replacing $n$ by $n^\varepsilon$.
That is, let $C'_1 < \ldots < C'_{t'}$ be the $\lambda_0$-equivalence
classes ordered naturally.  Clearly $s^*_\ell =
\lfoot\log_2(n^\varepsilon) \rfoot$ iff $s^*_\ell \ne 0$ iff for $r
\in \{0,\dotsc,t'-1\}$ we have $\{\ell,\ell +1\} = \{\text{the maximal
} i \in C'_r$, the minimal $j \in C'_{r+1}\}$ where we stipulate $C_0
= \{|A|\}$.  Hence $\{\ell:\ell < m$ and $s^*_\ell
= \lfoot\log_2(n^\varepsilon)\rfoot\}$ has exactly $t'$ members which
is $\bold v_{\lambda_0}(A,B)$ and $\Sigma\{j(\ell):s^*_\ell =
\lfoot\log_2(n^\varepsilon) \rfoot\}$ is exactly $\bold
e_{\lambda_0}(A,B)$, see the definition of $j(\ell)$ in $\circledast_1$.

So if $0 < \zeta \le \varepsilon \times |-\bold w_{\lambda_0}(A,B)|$ then

$$
\beta'_{{\frak p},\bar s^*} = (n^\varepsilon)^{t'-\bold e_{\lambda_0}(A,B)} =
n^{\varepsilon(\bold w_{\lambda_2}(A,B))} < 1/n^{\varepsilon-(2
\varepsilon)} < 1/n^{2 \zeta}.
$$
\mn
By Stage A this is enough.

The four cases cover all possibilities (remembering $s_0^*\geq
\lfoot\log_2(n^\varepsilon)\rfoot$).
So we have finished getting a bound on 
$\beta_{{\frak p},\bar{s}^*}$ of the form $1/n^{\zeta + \varepsilon}$ 
hence $\le 1/n^\zeta$ (see end of stage E) and 
as said in stages B, C this suffices. \hfill$\square_{\scite{2.4}}$\margincite{2.4}
\bigskip

\remark{\stag{2.8} Remark}  1) Actually the situation is even 
better as $\dsize \sum_{\bar s\in S} \beta'_{{\frak p},\bar s} \leq \bold c 
\beta'_{{\frak p},\bar s^*}$ by the formula for the sum of a geometric series 
(induction on $m$). \nl
2) We can get logarithmic bound instead $n^\varepsilon$, a simple way
is to divide the problem according to $\min\{\lfoot\log_2(|g(x)-g(y)\rfoot):x
\in A,y \in B\setminus A\}$. \nl
3) Similarly we can get better bounds in \scite{2.9}, \ub{but} all this is
not needed to our main purpose.
\endremark
\bn
By \scite{2.4} we have sufficient conditions for (given $A \le^* B$)
``every $f:A \hra {\Cal M}_n$ has few pairwise disjoint extension
to $g:B \hra {\Cal M}_n$".  Now we try to get a dual, a sufficient
condition for: (given $A \le^* B$) for every random enough 
${\Cal M}_n$, every $f:A \hra {\Cal M}_n$ has ``many" pairwise disjoint
extensions to $g:A \hra {\Cal M}_n$.
\bigskip

\proclaim{\stag{2.9} Lemma}   Assume 
\mr
\item "{$(A)$}"   $(A,B,\lambda) \in {\Cal T}$,
\sn
\item "{$(B)$}"   $(\forall B')[A <^* B' \le^* B \and B' \text{ is } \lambda
\text{-closed} \Rightarrow \bold w_{\lambda}(A,B') > 0]$ \nl
(recall that ``$B'$ is $\lambda$-closed'' means 
$x \lambda y \and x \in B' \Rightarrow y \in B'$).
\ermn
\ub{Then} there is $\zeta \in \Bbb R^+$, in fact we can let

$$
\zeta =: \min\{\bold w_{\lambda \restriction B'}(A,B'):A \subseteq B'
\subseteq B \text{ and } B' \text{ is } \lambda \text{-closed}\},
$$
\mn
such that:
\mr
\item "{$\otimes$}"    for every small enough $\varepsilon > 0$, for
random enough 
${\Cal M}_n$, for every $f:A \hra {\Cal M}_n$ and $k$ with 
$0<k<k+n^{1-\varepsilon} < n$, there are $\ge n^{(1-\varepsilon)
\cdot\zeta}$ pairwise disjoint extensions $g$ of $f$ satisfying 
{\roster
\itemitem{ $(i)$ }   $g:B \hra {\Cal M}_n$,
\sn
\itemitem{ $(ii)$ }   $g(B \backslash A) \subseteq
[k,k+n^{1-\varepsilon})$.
\endroster}
\endroster
\endproclaim
\bigskip

\demo{Proof}   Let the 
$\lambda$-equivalence classes be $B_1\ldots B_{\ell(*)}$, and
$B_\ell=\{b_{\ell,1}\ldots b_{\ell,m_\ell}\}$. \nl
Without loss of generality Rang$(f) \cap [k,k + 
\frac{1}{|A|+1} n^{1-\varepsilon}) = \emptyset$ as for this we can
replace $k$ by some $k' \in [k,k+n^{1-\varepsilon}(1-\frac{1}{|A|+1}))$.
\nl
For $\ell=1,\ldots, \ell(*)$ let $k_\ell=k+\frac{n^{1-\varepsilon}}{2\ell(*)
+1}\cdot\frac{2\ell-1}{|A|+1}$ and $I_\ell=[k_\ell,k_\ell +
\frac{n^{1-\varepsilon}}{2\ell(*)+1} \cdot \frac{1}{|A|+1})$ and 
$J_\ell=\{i\in I_\ell: m_\ell \text{ divides }(i-k_\ell)$ and $i +
m_\ell \in I_\ell\}$.
(The point is that for some positive constant $\bold c,|J_\ell| \ge \bold c
n^{1 - \varepsilon}$ and $j_1 \in I_{\ell_1} \and j_2 \in I_{\ell_2} 
\and \ell_1 \ne \ell_2 \Rightarrow |j_1 -j_2| \ge 
\bold c n^{1-\varepsilon} +|B|$ and $j_1 \in \text{ Rang}(f) \and j_2 \in
\dbcu_\ell I_\ell \rightarrow |j_1-j_2| \ge \bold c n^{1- \varepsilon}$ and
$j_1,j_2 \in J_\ell$ and $i_1,i_2 \in \{1,\dotsc,m_\ell\} \and j_1+i_1 =
j_2+i_2 \Rightarrow (j_1,i_1) = (j_2,i_2))$, e.g.,
for $\bold c = 1/(2^{\ell(*)+1}) \times (|B| +1)$, (the $|B|$ an overkill).

Now if $\bar j=\langle j_1 \ldots j_{\ell(*)}\rangle \in
\dsize \prod_\ell J_\ell$ and we let $g_{\bar j}$ be the function with 
domain $B$ satisfying $g_{\bar j}(b_{\ell,m}) = j_\ell + m$ and 
$g_{\bar j}\restriction A=f$, then from the random ${\Cal M}_n$ we can
compute the set $Y = Y[{\Cal M}_n] =: \{\bar j:g_{\bar j}$ is an
embedding of $B$ into ${\Cal M}_n\}$.  
As said above, each $g_{\bar j}$ is a one-to-one mapping of the
desired kind and $g_{\bar j}(b_{\ell,m}) \in I_\ell$.

Note that 
\mr
\item "{$\circledast$}"  $(\text{Rang}(g_{\bar j_1}) 
\backslash \text{ Rang}(f))$ 
is disjoint from $(\text{Rang}(g_{\bar j_2}) \backslash \text{ Rang}(f))$
\ub{if and only if} \nl
$\text{Rang}(\bar j_1)$ is disjoint from $\text{Rang}(\bar j_2)$.
\ermn
Reflecting, this is similar to \cite{ShSp:304}, \cite{BlSh:528} and is
a particular case
of \cite{Sh:550}; informally the transformation of \cite{ShSp:304},
\cite{BlSh:528} corresponding to non $\lambda$-equivalent pairs
consists of
\mr
\widestnumber\item{$(iii)$}
\item "{$(i)$}"   the probability of an edge is $\in [\bold c_1
\frac{1}{(n^{1-\varepsilon})^\alpha},\bold c_2 
\frac{1}{(n^{1-\varepsilon})^\alpha}]$ here for some positive 
constants $\bold c_1< \bold c_2$ (and not as there $\frac{1}{n^\alpha}$),
as $|J_\ell|$ is $\sim n^{1-\varepsilon}$,
\sn
\item "{$(ii)$}"   the number of nodes is now $n^{1-\vep}$ (and not $n$),
factoring by $\lambda$,
\sn
\item "{$(iii)$}"   here we have several candidates for edges for 
each pair $j_1$, $j_2$, as we are thinking of equivalence classes 
as nodes, 
\sn
\item "{$(iv)$}"   we have random unary predicates coming from 
edges inside a $\lambda$--equivalence class and edges from a point in
$A$ to a point in $B \backslash A$.
\ermn
However, (i), (iii) and (iv) are insignificant changes and (ii) is o.k. if 
$\varepsilon$ is small enough.

Formally, for every graph $M \in {\Cal K}_n$ we define a model 
$N[M]$ (this depends on $f$, $A$, $B$, $\langle B_\ell:
\ell=1,..., \ell(*)\rangle$ and the $k$).
Its set of elements is

$$
\dbcu_{\ell=1}^{\ell(*)} R^{N[M]}_\ell,
$$
\mn
where
\mr
\item "{$(a)$}"   $R^{N[M]}_\ell=R^n_\ell$ is a unary 
predicate, $R^{N[M]}_\ell$ is
$$
\{i:1 \le i \le \frac{n^{1-\varepsilon}}{2\ell(*)+1} \cdot
\frac{1}{|A|+1} \cdot \frac{1}{m_\ell}\}\times\{\ell\}
$$
(so is constant), a member $i$ of $J_\ell$ is represented by
$(i-k_\ell)/m_\ell$. 
We may say that the set $R^n_\ell$ 
consists of representatives of the set $J_\ell$
\sn
\item "{$(b)$}"   if 
$e=\{b_1,b_2 \}$, $b_1=b_{\ell_1,j_1}$, $b_2=b_{\ell_2,j_2}$
with $\ell_1<\ell_2$ from $[\ell(*)]$ (so 
$j_1\in [m_{\ell_1}]$, $j_2\in [m_{\ell_2}])$ \ub{then}
$Q^{N[M]}_e$ is the following binary predicate:
$$
\align
Q^{N[M]}_e=\big\{\{(i_1,\ell_1),(i_2,\ell_2)\}:\{&k_{\ell_1}+m_{\ell_1}\times
i_1+j_1, k_{\ell_2}+m_{\ell_2}\times i_2+j_2\} \\
  &\text{is an edge of }M \big\}
\endalign
$$
(i.e. it is symmetric and trivially irreflexive) (i.e., (b) 
is talking 
about long edges; edges between two disjoint equivalence classes and we
represent it as a binary predicate),
\sn
\item "{$(c)$}"   if $e=\{a,b\}$, $a\in A$, $\ell_1\in [\ell(*)],
j \in [m_{\ell_1}],b=b_{\ell_1,j}$ \ub{then} 
$Q_e$ is a unary predicate of $N[M]$:
$$
Q^{N[M]}_e=\big\{(i,\ell_1):\{f(a),k_{\ell_1}+m_{\ell_1}\times
i+j\} \text{ is an edge of }M\big\},
$$
[i.e., (c) is talking about edges from $f(A)$ to $g(B \setminus A)$ and this
is represented by unary predicate.] 
\sn
\item "{$(d)$}"   if $e=\{b_{\ell,j_1},b_{\ell,j_2}\}$, $\ell\in [\ell(*)]$,
$j_1<j_2$ in $[m_\ell]$, \ub{then} $Q_e$ is a unary predicate
$$
\align
Q^{N[M]}_e=\big\{(i,\ell):&\{k_\ell+m_\ell\times i+j_1,k_\ell+m_\ell\times i
+ j_2\} \\
  &\text{ is an edge of }M\big\}.
\endalign
$$\nl
[i.e., (d) is talking about the ``short edges" and so it is
represented as a unary predicate.]
\ermn
Now let ${\Cal K}^*_n = \{N[M]:M\in {\Cal K}_n\}$ and 
let ${\Cal N}_n$ vary on ${\Cal K}^*_n$ with 

$$
\text{Prob}_{\mu^*_n}({\Cal N}_n=N) = \text{Prob}_{\mu_n}(N[{\Cal M}_n]=N).
$$
\mn
Now reflecting, the relations $Q_e$ are drawn independently, moreover all the
instances are drawn independently and the probabilities are essentially
constant.  More formally, for $e=\{b_{\ell_1,j_1},b_{\ell_2,j_2}\}
\subseteq B$ as in clause (b),
$\ell_1<\ell_2$ in $[\ell(*)]$ and $x\in R_{\ell_1}$, $y\in R_{\ell_2}$ the
event

$$
{\Cal E}(x,y,e) =: [{\Cal N}_n \models Q_e(x,y)]
$$
\mn
has probability $p^e_{(x,y)}$ which satisfies

$$
\align
p^e_{x,y} &= \text{Prob}_{\mu_n}(\{k_{\ell_1} + 
m_{\ell_1}\times x+j_1,k_{\ell_2} +
m_{\ell_2}\times y+j_2\} \text{ an edge of } {\Cal M}_n) \\
  &= \big(|(k_{\ell_1}+m_{\ell_1}\times x+j_1)-(k_{\ell_2}+ m_{\ell_2}\times
y+j_2)|\big)^{-\alpha}.
\endalign
$$
\mn
Now

$$
\align
\frac{n^{1-\varepsilon}}{(2\ell(*)+1)\times |B|} &\le
|k_{\ell_1}-m_{\ell_1}\times x+j_1-(k_{\ell_2}+m_{\ell_2}\times
y+j_2)| \\
  &\le n^{1-\varepsilon}.
\endalign
$$
\mn
[Why?  We first evaluate the lower bound.   The worst case is 
$\ell_1$,$\ell_2= \ell_1+1$. 
The difference is at least min$\{|g(w)-g(y):w \in B_{\ell_2}$ and $y
\in B_{\ell_1}\} \ge \min\{|w-y|:w \in I_{\ell_1},y \in I_{\ell_2}\}
\ge \frac{1}{|2 \ell(*)+1)\times |B|} \ge n^\varepsilon$ (see
beginning).
\nl
The upper bound: any two members of $g(z),z \in B \backslash A$ are in
$[k,k+n^\varepsilon)$ hence for some constants $\bold c_1 < \bold c_2$ 
(i.e., numbers depending on $A,B$ but not on $n$) 

$$
\bold c_1\cdot (n^{1-\varepsilon})^{-\alpha} \le 
p^e_{x,y}\leq \bold c_2\cdot(n^{1-\varepsilon})^{-\alpha}.
$$
\mn
Similarly, if $e=\{a,b_{\ell,j}\}$, $a\in A$, $\ell\in [\ell(*)]$, $j\in
[m_\ell]$ and $x\in R^n_\ell$ the event

$$
{\Cal E}(x,e) =: {\Cal N}_n\models Q_e(x)
$$
\mn
has probability $p^e_x$ which for some $\varepsilon(e)\in
[0,\varepsilon]_{\Bbb R}$ and constants $\bold c_3$, $\bold c_4$ 
($>0$), satisfies 

$$
\bold c_3(n^{1 - \varepsilon})n^{-\alpha} 
\le p^e_x \le \bold c_4\cdot (n^{1-\varepsilon(e)})^{- \alpha}.
$$
\mn
Lastly, for $e=\{b_{\ell,j_1},b_{\ell,j_2}\}$, $\ell\in [\ell(*)]$, $j_1<j_2$
in $[m_\ell]$ and $x\in R^n_\ell$ the event

$$
{\Cal E}(x,e) =: {\Cal N}_n\models Q_e(x)
$$
\mn
has probability $p^e_x$ which is constant $\in (0,1)_{\Bbb R}$: it is

$$
p_{|j_1-j_2|} =\cases
|j_1-j_2|^{-\alpha} &\text{ if } \quad |j_1-j_2|>1 \\
 \frac{1}{2^\alpha}        &\text{ if } \quad |j_1-j_2|=1.
\endcases
$$
\mn
Now these events ${\Cal E}(x,y,e),{\Cal E}(x,e)$ are independent.
This is a particular case of \cite{Sh:550}, but if we like to deduce 
it from \cite{BlSh:528}, first draw the unary predicates. 
There we know with large probability good lower and upper 
bounds to the number of elements in

$$
\align
R^{n,*}_\ell = \{x\in R^n_\ell:&\text{if } e=\{b_{\ell,j_1},b_{\ell,j_2}\},\
j_1<j_2 \text{ in }[m_\ell] \\
 &\text{or } e=\{a,b_{\ell,j}\},\ j\in [m_\ell] \\
 &\text{then } e \in \text{ edge}(B) \Rightarrow {\Cal N}_n \models Q_e(x)\}.
\endalign
$$
\mn
Now the rest of the drawing and what we need is like \cite{BlSh:528} only
having $\ell(*)$ sorts of elements, which are not exactly of the same size: we
can throw some nodes to equalize.  \hfill$\square_{\scite{2.9}}$\margincite{2.9}

Now, \scite{2.4}, \scite{2.9} are 
enough for proving $<^*_i=<_i$, $<^*_s=<_s$, weakly nice 
and similar things. But we need more.
\enddemo
\bigskip

\proclaim{\stag{2.10} Lemma}  Assume
\mr
\item "{$(A)$}"   $(A,B,\lambda) \in {\Cal T}$,
\sn
\item "{$(B)$}"   $\xi = \bold w_\lambda(A,B)>0$,
\sn
\item "{$(C)$}"   if $A<^*C<^*B$, and $C$ is $\lambda$-closed 
then $\bold w_\lambda(C,B)< 0$ 
(hence neccessarily $\xi\in (0,1)_{\Bbb R}$ and $C\neq\emptyset
\Rightarrow  \bold w_\lambda(A,C)>0$ and even $\bold w_\lambda(A,C) > \xi$).
\ermn
\ub{Then} for every $\varepsilon\in \Bbb R^+$, for 
every random enough ${\Cal M}_n$, for every $f:A \hra {\Cal M}_n$ we
have
\mr
\item "{$(a)$}"   the number of $g:B \hra {\Cal M}_n$ 
extending $f$ is at least $n^{\xi-\varepsilon}$
\sn
\item "{$(b)$}"   also the maximal number of pairwise 
disjoint extension $g:B\hra {\Cal M}_n$ of $f$ is at least this
number.
\ermn
\endproclaim
\bigskip

\remark{\stag{2.11} Remark}   1) We can 
get reasonably much better bound (see \cite{ShSp:304}, \cite{BlSh:528}
and \cite{Sh:550}) but this suffices. \nl
2) In the most interesting cases of \scite{2.10} we have $A
<^*_{\text{pr}} B$ but it applies to more cases.
\endremark
\bigskip

\demo{Proof}   Repeat the proof of \scite{2.9} 
noting that for ($\bar{j},g_{\bar j}$) as there, for 
random enough ${\Cal M}_n$:
\mr
\item "{$(*)$}"   for some $k^*=k^*(A,B)$, for every $x$ (recalling $J
= \dsize \prod_\ell J_\ell$)
$$
|\{\bar{j}\in J:x \in \text{ Rang}(g_j)\setminus \text{ Rang}(f)\}|<k^*.
$$
\ermn
Why $(*)$?  By \scite{2.4} and \scite{1.17}(13).  
\hfill$\square_{\scite{2.10}}$\margincite{2.10}
\enddemo
\bigskip

\proclaim{\stag{2.12} Claim}  Assume
\mr
\item "{$(A)$}"   $(A,B,\lambda) \in {\Cal T}$ and
\sn
\item "{$(B)$}"   if $C\subseteq B\setminus A$ is 
non empty and $\lambda$--closed then $\bold w_\lambda (A,A\cup C)>0$.
\ermn
\ub{Then} for some $\varepsilon_0 \in \Bbb R^+$ for every $\varepsilon\in
(0,\varepsilon_0)_{\Bbb R}$, every 
random enough ${\Cal M}_n$ for every $f:A \hra {\Cal M}_n$ we have
\mr
\item "{$(a)$}"   the number of $g: B \hra {\Cal M}_n$ 
extending $f$ is at least $n^{\bold w_\lambda (A,B)-\varepsilon}$;
moreover
\sn
\item "{$(b)$}"   for every $X\subseteq [n]$, $|X|\leq n^{\varepsilon_0-
\varepsilon}$, the number of $g:B \hra {\Cal M}_n$ 
extending $f$ with ${\text{\rm Rang\/}}(g)\cap X \subseteq 
{\text{\rm Rang\/}}(f)$ is at 
least $n^{\bold w_\lambda(A, B)-\varepsilon}$. 
\endroster
\endproclaim
\bigskip

\remark{\stag{2.13} Remark}   1) By \scite{1.16} the statement

\block
``for some $\lambda$ the hypothesis of \scite{2.12} holds"
\endblock
is equivalent to ``$A \le^*_s B$''. \nl
2) The affinity of this claim to being nice (see [I,\S2]) should be
clear. \nl
3) If $|X| \ge n^{1-\varepsilon}$ we can demand Rang$(g) \subseteq X$
but no need arises.
\endremark
\bigskip

\demo{Proof}   Choose a sequence 
$\bar B=\langle B_\ell:\ell \le k\rangle$ such that:
\mr
\widestnumber\item{$(iii)$}  
\item "{$(i)$}"     $A=B_0<^* B_1<^* \ldots<^* B_k = B$,
\sn
\item "{$(ii)$}"    each $B_\ell$ is $\lambda$--closed,
\sn
\item "{$(iii)$}"   $\xi_\ell = \bold w(B_\ell,B_{\ell+1},
\lambda\restriction(B_{\ell+1} \setminus B_\ell))>0$,
\sn
\item "{$(iv)$}"    if $C\subseteq B_{\ell+1}\setminus B_\ell$ is non empty
$\lambda$--closed then $\bold w(B_\ell,B_\ell \cup C,
\lambda\restriction C)>0$,
\sn
\item "{$(v)$}"     if $B_\ell<^* C<^* B_{\ell+1}$ and 
$C$ is $\lambda$--closed \ub{then} $\bold w_\lambda(C,B_{\ell+1}) < 0$.
\ermn
[Why $\bar B$ exists?   As $\{\bar B':\bar B' 
\text{ satisfies } (i)-(iv)\}$ is
not empty (as $\langle A,B\rangle$ belongs to it) and every member has length
$<|B\setminus A|+1$, so there is $\bar B = \langle B_\ell:\ell \le k
\rangle$ in the family of maximal length.  If it fails, clause (v), say
for $m$ and we have a $\lambda$-closed set $C$ 
satisfying $B_m <^* C <^* B_{m+1}$ such that
$\bold w_\lambda(C,B_{m+1}) \ge 0$ hence $\bold w_\lambda(C,B_{m+1}) >
0$, and \wilog \, $C$ is maximal under those conditions.  Let $\bar B' =
\langle B'_\ell:\ell \le k+1 \rangle,B'_\ell$ is $B_\ell$ if $\ell \le
k$ and is $B_{\ell -1}$ if $\ell > k$ and is $C$ if $\ell=k$.  It
is easy to check that $\bar B'$ satisfies clauses (i)-(iv), too. 
So we got a contradiction to the maximal length.]

Now, by \scite{1.9}(2)(c):

$$
\xi =: \bold w_\lambda(A,B) = \dsize \sum_{\ell<k}
\bold w(B_\ell,B_{\ell+1},\lambda\restriction (B_{\ell+1}\setminus B_\ell)),
$$
\mn
and let

$$
\varepsilon_0 = \min\{\bold w_\lambda(B_\ell,B_{\ell+1}): \ell<k\}
$$
\mn
(it is $>0$). For each $(B_\ell,B_{\ell+1})$ we can apply \scite{2.10} (for
disjointness to $X$ use clause (b) of \scite{2.10}).
\enddemo
\bigskip

\demo{Detailed Proof of $(b)$}   Let $\varepsilon\in 
(0,\varepsilon_0)$, let ${\Cal M}_n$ be random enough and let 
$f:A \hra {\Cal M}_n$. Let $X\subseteq [n]$, $|X|\leq
n^{\varepsilon_0-\varepsilon}$. For $\ell\in\{0,\ldots,k\}$ we let

$$
\align
{\Cal F}_\ell = \big\{g:&g \text{ is an embedding of }B_\ell \text{ into }
{\Cal M}_n \text{ extending }f \\
   &\text{and such that Rang}(g)\cap X \subseteq \text{ Rang}(f)\big\}.
\endalign
$$
\mn
Clearly $|{\Cal F}_0|=1$. 
Now let $\varepsilon_1 = \varepsilon/k$ and we 
prove by induction on $\ell$ that $|{\Cal F}_\ell|\geq
n^{\underset {i<\ell} {}\to \Sigma \bold w_\lambda(B_i,B_{i+1})
-\ell \varepsilon(1)}$. \nl
For $\ell=0$ this is clear; for $\ell+1$ it is enough to 
prove that for each $f^* \in {\Cal F}_\ell$ the set

$$
{\Cal F}_{\ell+1,f^*} = \{g\in {\Cal F}_{\ell+1}:
g\restriction B_\ell = f^*\}
$$
\mn
has $\geq n^{\bold w_\lambda(B_\ell,B_{\ell+1}) -
\varepsilon(1)}$ members (noting that ${\Cal F}_{\ell+1}$ is the disjoint 
union of $\{{\Cal F}_{\ell+1,f^*}: f^*\in {\Cal F}_\ell\}$).
Now clearly

$$
\align
{\Cal F}_{\ell+1,f^*} = \big\{g:&g \text{ is an embedding of }
B_{\ell+1} \text{ into } {\Cal M}_n \text{ extending } f^* \\
  &\text{and such that Rang}(g)\cap X \subseteq \text{ Rang}(f^*)\big\}
\endalign
$$
(as Rang$(f^*)\cap X \subseteq \text{ Rang}(f)$ as $f^*\in 
{\Cal F}_\ell$). Now the
inequality follows by clause (b) of 
\scite{2.10}. \hfill$\square_{\scite{2.12}}$\margincite{2.12}
\enddemo
\bigskip

\proclaim{\stag{2.14} Claim}   Assume $A<_{pr}^* B$ and 
$k=|B\setminus A|$, and we let $\xi = \xi(A,B)$ that is

$$
\align
\xi = \max\{\bold w_\lambda(A,B):&(A,B,\lambda) \in {\Cal T}
\text{ and} \\
 &\text{for every }\lambda \text{-closed non empty } C \subseteq B
\setminus A \\
  &\text{we have } \bold w(A,A \cup C,\lambda\restriction C)>0\}.
\endalign
$$
\ub{Then} for every 
$\varepsilon \in \Bbb R^+$ for every random enough ${\Cal M}_n$ for
every $f:A \hra {\Cal M}_n$ we have
\mr
\item "{$(*)$}"   the number of $g:B \hra {\Cal M}_n$ extending 
$f$ is at most $n^{\xi+\varepsilon}$.
\endroster
\endproclaim
\bigskip

\demo{Proof}   Like the proof of \scite{2.4} and \scite{2.9}.

In detail, for this it 
is enough to prove that for each $f:A \rightarrow [n]$ the
probability of failure is $< 1/n^k$ for any $k$.  By \scite{2.4} (or
more fully, \scite{3.2} below) it is enough to prove that: for every
$\varepsilon \in \Bbb R^+$ for every random enough ${\Cal M}_n$, for
every $f:A \hookrightarrow {\Cal M}_n$
\mr
\item "{$(*)'$}"  there are no $n^{\xi + \varepsilon}$ extensions $g:B
\hookrightarrow {\Cal M}_n$ of $f$ with Rang$(g) \backslash \text{
Rang}(f)$ pairwise disjoint.
\ermn
We can consider ${\Cal G}_{A,B}(f,{\Cal M}_n) = \{g:g$ a one to one function
from $B$ to $[n]$ extending $f$ and mapping every edge of $B$ not
$\subseteq A$ to an edge of ${\Cal M}_n\}$.
As in stage C in the proof of \scite{2.4}, 
for each $f:A\rightarrow [n]$, we can divide the
extensions $g: [B] \rightarrow [n]$ of $f$ to types (by tp$^0$), and the
set of types is computable from $A$, $B$, so it is enough to restrict
ourselves to one type, and as there, it is enough to deal with the
case Rang$(f)=\{1,\ldots,|A|\}$ (of course, $h^*$ changes too).  This
just simplifies the computations.

Again as in the proof of \scite{2.4}, stage D, we can fix $B\setminus
A = \{b_1,\ldots,b_m\}$, stipulate $b_0= f^{-1}(|A|)$, and having
fixed $\bar \nu$ we fix
tp$^1(g)$ (i.e., the permutation $\sigma= \sigma(g)$ which essentially is
the order of $\langle g(b_{\ell+1})-g(b_\ell):\ell \rangle$) and we can fix
$s_\ell = \lfoot\log_2(g(b_{\ell+1})-g(b_\ell))\rfoot$ 
(again partition to some power of $\log_2(n)$ cases) as in \scite{2.4},
this disappears and we take more cases.  
No harm in increasing the probability.  So the situation is similar
enough to \cite{BlSh:528}.

Let for $\ell_1 < \ell_2$, $k(\ell_1,\ell_2)$ be the $k$ such that
$\ell_1 \le k < \ell_2$ and tp$^1(g)$ says that $\sigma(k)$ is maximal
equivalently $s_k$ is and $k$ is maximal among those satisfying this
condition if we have some candidates.

So it is enough to prove 
\mr
\item "{$(a)$}"   assume: $s_\ell \leq \log_2(n)$, for 
$\ell<m$, $X_\ell$ is a set of $2^{s_{\ell-1}}$  elements 
for $\ell=1, \ldots, m$,$X_0=\{1, \ldots, |A|\}$,
pairwise disjoint for simplicity, we draw a graph on
$X = \dbcu_{\ell<m} X_\ell$, flipping a coin for each edge
independently, the probability of ``$\{x,y\}$ is an edge'' is:
$$
\text{zero} \quad \text{ \ub{if} } \dsize \bigvee_i \{x, y\}\subseteq
X_\ell
$$

$$
\frac{1}{2^{s_{k(\ell_1, \ell_2)}\alpha}} \text{ \ub{if} } x \in
X_{\ell_1}, y\in X_{\ell_2}, \ell_1< \ell_2.
$$
We have to show: with probability $\geq 1-\frac{1}{n^{|B|+1}}$ the
number of embeddings $g \in {\Cal G}_{A,B}(f,M)$ satisfying $\ell \in
\{1,\dotsc,m\} \Rightarrow g(b_\ell)\in X_\ell$, is $< n^{\xi+\varepsilon}$.
\ermn
Of course, we can discard the cases $\dsize \bigvee_\ell 2^{s_\ell} <
n^{\varepsilon/2}$ (as then we can apply \scite{2.4} to 
$(A\cup \{b_\ell\}, B)$). The rest should be clear.
\enddemo
\bn
\ub{Alternative proof}:  Let 
$m = \lfoot\log(n)\rfoot+1$, and $\eta_i \in {}^{\{0,\dotsc,m-1\}}2$
be the binary representation of $i$ for $i\in [n]$ that is $i = \dsize
\sum_{\ell < m} \eta(\ell)2^\ell$.  
For $\eta \ne \nu \in {}^{\{0,\dotsc,m-1\}}2$ we let $\ell_1(\eta,\nu)
= \text{ Min}\{\ell < m:\eta(\ell)\ne \nu(\ell)\}$ and
$\ell_2(\eta,\nu) = \text{ Max}\{\ell:\ell_1(\eta,\nu) < \ell \le m$
and if $\ell_1(\eta;\nu) + 2\le m$ then $(\forall t)(\ell_1(\eta,\nu)
+ 1 < t < \ell \Rightarrow \eta(t) = 1 - \eta(\ell_1(\eta,\nu)+1) \and
\nu(t) = \nu(\ell_1(\eta,\nu)+1)\}$.
For each one to one $g: B\rightarrow [n]$ let

$$
\align
w_g =: \{\ell_1(\eta_{g(b_1)},\eta_{g(b_2)}):b_1 \ne b_2 \in B\} \cup
\{\ell_2(\eta_{g(b_1)},\eta_{g(b_2)}):&b_2 \ne b_2 \in B \\
  &\text{and } \ell_1(g(b_1),g(b_2)) + 1 < m \text{ and} \\
  &\ell_2(\eta_{g(b_1)},\eta_{g(b_2)}) < m\}
\endalign
$$
\mn
$\bar \nu^g = \langle \eta_{g(b)} \restriction w_g:b \in B \rangle$;
this is a little more information than we need.
\nl
Easily $|w_g| \le 2|B|$.  As there are $< (\log n)^{2 \cdot |B|}$ 
cases, and $(\log n)^k << n^\varepsilon$ for $n$ large enough, it
suffices to prove $(*)$ when we restrict ourselves to one $\bar \nu$ as
$\bar \nu^g$.  Fix $\varepsilon > 0$, and let $n$ be large enough.

Lastly, let ${\Cal G}_{\bar \nu}(f,{\Cal M}_n) = \{g:g$ is a
one-to-one function from $B$ into ${\Cal M}_n$, it extends $f$ and 
$\bar \nu_g = \bar \nu$  and $g$ maps edges of $B$
not $\subseteq A$ to edges of ${\Cal M}_n\}$ and let

$$
\align
{\Cal G}_{\bar \nu}(f,[n]) = \{g:&\text{ for some } 
M \in {\Cal K}_n \text{ into
which } f \text{ is an embedding of} \\
  &A, \text{ we have } f \subseteq g:B \hookrightarrow M \text{ and }
\bar \nu_g = \bar \nu\}
\endalign
$$
\mn
Now we try to bound $\beta_{\bar \nu} =: \text{ Exp}(|{\Cal G}^{\xi +
\varepsilon}_{\bar \nu}(f,{\Cal M}_n)|)$ where

$$
\align
{\Cal G}^{\xi + \varepsilon}_{\bar \nu}(f,{\Cal M}_n) = 
\{\bar g:&\bar g = \langle g_i:i < n^{\xi + \varepsilon} \rangle, \\
  &\text{ satisfies } g_i \in {\Cal G}_{\bar \nu}(f,{\Cal M}_n) \text{
and} \\
  &i < j \Rightarrow \text{ Rang}(g_i) \cap \text{ Rang}(g_1) = \text{
Rang}(f)\}. 
\endalign
$$

$$
{\Cal G}^{\xi + \varepsilon}_{\bar \nu}(f,[n]) \text{ is defined
similarly}.
$$
\mn
Let $\bold p_{\bar \nu} = \text{ max}\{\text{Prob}(g$ embed $B$ into
${\Cal M}_n$ over $A):g \in {\Cal G}_{\bar \nu}(f,[n])\}$. \nl
Clearly

$$
\align
\beta_{\bar \nu} &\le \Sigma\{\text{Prob}(\dsize \bigwedge_i g_i 
\in {\Cal G}_{A,B}(f,{\Cal M}_n)):\bar g \in {\Cal G}^{\xi +
\varepsilon}_{\bar \nu}(f,[n])\} \\
  &= \Sigma\{\dsize \prod_i (\text{Prob}(g_i \in {\Cal
  G}_{A,B}(f,{\Cal M}_n)):\bar g \in {\Cal G}^{\xi + 
\varepsilon}_{\bar \nu}(f,[n])\} \\
  &= |{\Cal G}^{\xi + \varepsilon}_{\bar \nu}(f,[n]| \cdot (\bold p_{\bar
\nu})^{n^{\xi + \varepsilon}} \\
  &\le |{\Cal G}_{\bar \nu}(f,[n])^{n^{\xi + \varepsilon}}|
(\bold p_{\bar \nu})^{n^{\xi + \varepsilon}} 
\endalign
$$
\mn
Recall that for simplicity we 
assume \wilog \, Rang$(f) = \{1,\dotsc,|A|\}$. \nl
Let $B \backslash A = \{b_1,\dotsc,b_k\}$ and let $b_0$ be such that $f(b_0) =
|A|$ and if $A = \emptyset$ we stipulate such $b_0$.  
Clearly there are $\langle \ell^1_{i,j}:i < j \le k
\rangle,\langle \ell^2_{i,j}:i < j \le k \rangle$ such that

$$
g \in {\Cal G}_{\bar \nu}(f,[n]) \Rightarrow \ell^1_{i,j} =
\ell_1(\eta_{g(b_i)},\eta_{g(b_j)}) \and \ell^2_{i,j} =
\ell_2(\eta_{g(b_1)},\eta_{g(b_2)}).
$$
\mn
Now note that for every $i \in \{1,\dotsc,k\}$ and $g_* \in {\Cal
G}_{\bar \nu}(f,[n])$

$$
|\{g(b_\ell):g_* \restriction \{b_0,\dotsc,b_{\ell-1}\} \subseteq g
\in {\Cal G}_\nu(f,[n])\}|
$$
\mn
is $\le 2^{m-\ell^{2_{i-1,i}}}$.
\mn  
Hence

$$
|{\Cal G}_\nu(f,[n])| \le \dsize \prod^k_{\ell =1} 2^{m-\ell^2_{i-1,i}} =
2^{\Sigma\{(m-\ell^2_{i,i +1}):i=0,\dotsc,k-1\}}.
$$
\mn
What about $\bold p_{\bar \nu}$?  For each $c \in \text{ Edge}(B)
\backslash \text{ Edge}(A)$ let $e = \{a(e),b(e)\}$ such that $b(e) =
b_{i_2(e)}$ and $(a(e) \in A \and i_0(e) =0) \vee (a(e) = b_{i_0(e)} \and i_0(e) <
i_2(e))$.  Let $i_1(e)$ be the $i \in \{i_0(e),i_0(e) +
1,\dotsc,i_2(e)-1\}$ such that $\ell^2_{i,i+1}$ is minimal, and among
those, $i$ is minimal.

For $i \in \{0,\dotsc,k-1\}$ let $j(i) = \{e:e \in \text{ Edge}(B)
\backslash \text{ Edge}(A)$ and $i_1(e) = i\}$.

The rest should be clear being similar to \scite{2.4}.
\newpage

\head {\S6 The conclusion} \endhead  \resetall \sectno=6
 \spuriousreset
\bigskip

{\bf Comment}:  In this section it is shown that $<^*_i$ and
$<^*_s$ (introduced in \S4) agree with the $<_i$ and $<_s$ of [I,\S1]
by using the  probabilistic information  from \S5. Then it is
proven that main  context ${\Cal M}_n$ is simply nice (hence simply almost
nice) and it satisfies  the $0-1$ law.
\bigskip

\demo{\stag{3.1} Context}  As in \S4 and \S5, so 
$p_i=1/{i^\alpha}$, for $i>1$, $p_1=p_2$ (where $\alpha
\in(0,1)_{\Bbb R}$ irrational) and ${\Cal M}_n = {\Cal M}^0_n$ 
(only the graph) and $\leq_i$, $\leq_s$, cl are as defined \S1. (So
${\Cal K}_\infty = {\Cal K}$ by \scite{1.4}).
\enddemo
\bigskip

Note that actually the section has two parts of distinct flavours: in
\scite{3.2} - \scite{3.6} we use the probabilistic information from \S5 to
show that the definition of $<_x$ from [I,\S1] and of $<^*_x$ from \S4
give the same relation. But to actually prove almost niceness, 
we need more work on the relations $\leq^*_x$ defined in \S4;
this is done in \scite{3.9}, \scite{3.10}, \scite{3.11}. Lastly we
put everything together.

The argument in \scite{3.2} - \scite{3.6} parallels that in
\cite{BlSh:528} which is more hidden in 
\cite{ShSp:304}. The most delicate step is to establish clauses
(A)($\delta$) and ($\varepsilon$) of Definition [I,2.13,(1)]
(almost simply nice). For this, we consider 
$f:A \hra {\Cal M}_n$ and try to extend $f$ to 
$g:B \hra {\Cal M}_n$ where $A\leq_s B$ such that 
Rang$(g)$ and $c \ell^k(f(A),{\Cal M}_n)$ are ``freely amalgamated''
over Rang$(f)$.  The key facts have been established in Section 5. If
$\zeta = \bold w(A,B,\lambda)$ we have shown (Claim \scite{2.12}) that 
for every
$\varepsilon>0$ for every random enough ${\Cal M}_n$, there are 
$\ge n^{\zeta-\varepsilon}$ embeddings of $B$ into ${\Cal M}_n$
extending $f$. 
But we also show (using \scite{2.14}) that for each obstruction to 
free amalgamation there is a $\zeta'<\zeta$ such that for every 
$\varepsilon_1>0$ such that the number of embeddings satisfying 
this obstruction is $<n^{\zeta'+\varepsilon_1}$, 
where $\zeta' = \bold w(A,B',\lambda)$ (for some $B'$ exemplifying
the obstruction) with $\zeta'+\alpha\leq\zeta$. So if 
$\alpha>\varepsilon+\varepsilon_1$ we overcome the 
obstruction. The details of this computation for various 
kinds of obstructions are carried out in proving Claim \scite{3.6}.
\bigskip

\proclaim{\stag{3.2} Claim}  Assume $A<^* B$. 
\ub{Then} the following are equivalent:
\mr
\item "{$(i)$}"     $A<^*_i B$ (i.e. from Definition \scite{1.11}(3)),
\sn
\item "{$(ii)$}"    it is not true that: for some $\varepsilon$, for 
every random enough ${\Cal M}_n$ for every $f:A \hra {\Cal M}_n$, 
the number of $g: B\hra {\Cal M}_n$ extending $f$ is $\geq n^\varepsilon$,
\sn
\item "{$(iii)$}"   for every $\varepsilon \in \Bbb R^+$ for 
every random enough ${\Cal M}_n$ for every $f:A \hra {\Cal M}_n$ the 
number of $g:B \hra {\Cal M}_n$ extending $f$ is $<n^\varepsilon$
(this is the definition of $A <_i B$ in [I,\S1]).
\endroster
\endproclaim
\bigskip

\demo{Proof}   We shall use the well known finite 
$\Delta$-system lemma: if $f_i:B\rightarrow[n]$ is one to one for 
$i<k$ then for some $w \subseteq \{0,\ldots,k-1\}$, $|w| \ge
\frac{1}{|B \backslash A|^2} k^{1/2^{|B|}}$, and $A'\subseteq B$ and
$f^*$ we have: 
$\dsize \bigwedge_{i\in w} f_i \restriction A'=f^*$ and 
$\langle \text{ Rang}(f_i\restriction (B\setminus A'):i\in w \rangle$ 
are pairwise disjoint (so if the $f_i$'s are pairwise distinct then
$B\setminus A'\neq\emptyset$).

We use freely Fact \scite{1.2}.  First, clearly $(iii) \Rightarrow
(ii)$. \nl
Second, if $\neg(i)$, i.e., $\neg(A<^*_i B)$ then by \scite{1.14} 
(equivalence of first and last possibilities + \scite{1.13}(1)) there are $A'$,
$\lambda$ as there, that is such that:
\block
$A\leq^* A'<^* B$ and $(A',B,\lambda)\in {\Cal T}$ and
if $C\subseteq B\setminus A'$ is non empty $\lambda$-closed
then $\bold w(A',A'\cup C,\lambda\restriction C)>0$ (see \scite{1.14}).
\endblock
So $(A',B,\lambda)$ satisfies the assumptions of \scite{2.9} which gives
$\neg(ii)$, i.e., we have proved 
$\neg (i) \Rightarrow \neg(ii)$ hence $(ii) \Rightarrow (i)$.

Lastly to prove $(i) \Rightarrow (iii)$ assume $\neg(iii)$ 
(and we shall prove $\neg(i)$).  So for some $\varepsilon \in \Bbb R^+$:
\mr
\item "{$(*)_1$}"   $0 < \lim \underset {n \rightarrow \infty} {}\to
\sup \text{Prob}$ (for some $f:A \hra {\Cal M}_n$, the number of 
$g:B \hra {\Cal M}_n$ extending $f$ is $\geq n^\varepsilon)$.
\ermn
But the paragraph of this proof showing $\neg(i) \Rightarrow \neg(ii)$ 
shows also that from $(*)_1$ we can deduce that for some $\zeta \in \Bbb R^+$
\mr
\item "{$(*)_2$}"   $0 < \underset {n\rightarrow\infty} {}\to {\lim \sup} 
\text{ Prob}$ (for some $A',A \le^* A'<^* B$, and $f':A' \hra {\Cal
M}_n$ there are $\geq n^\zeta$ functions $g:B \hra {\Cal M}_n$ which 
are pairwise disjoint extensions of $f')$.
\ermn
So for some $A'$, $A\leq^*A'<^*B$ and
\mr
\item "{$(*)_3$}"   $0 < \lim \underset {n\rightarrow \infty} {}\to
\sup \text{  Prob}$ (for some $f': A' \hra {\Cal M}_n$ there are 
$\geq n^\zeta$ functions $g: B \hra {\Cal M}_n$ which are pairwise 
disjoint extensions of $f^*)$.
\ermn
By \scite{2.4} (and \scite{2.3}, \scite{2.2}) we have $\neg(A'<^*_a B)$,
which by Definition \scite{1.11}(5) means that
$A' <^*_s B$ which (by Definition \scite{1.11}(4)) 
implies $\neg(A<^*_i B)$ so $\neg$(i) holds as required.
\hfill$\square_{\scite{3.2}}$\margincite{3.2}
\enddemo
\bigskip

\proclaim{\stag{3.3} Claim}  For $A<^*B \in {\Cal K}_\infty$, the 
following conditions are equivalent:
\mr
\widestnumber\item{$(iii)$}
\item "{$(i)$}"     $A<^*_s B$,
\sn
\item "{$(ii)$}"    it is not true that: ``for every 
$\varepsilon \in \Bbb R^+$, for every random enough ${\Cal M}_n$ for every 
$f:A \hra {\Cal M}_n$, there are no $n^\varepsilon$ pairwise 
disjoint extensions $g:B \hra {\Cal M}_n$ of $f$'',
\sn
\item "{$(iii)$}"   for some $\varepsilon \in \Bbb R^+$ for 
every random enough ${\Cal M}_n$ for every $f:A \hra {\Cal M}_n$ 
there are $\geq n^\varepsilon$ pairwise disjoint extensions
$g:B \hra {\Cal M}_n$ of $f$.
\endroster
\endproclaim
\bigskip

\demo{Proof}  Reflection shows that $(iii) \Rightarrow (ii)$.

If $\neg(i)$, i.e., $\neg(A<^*_s B)$ then by Definition
\scite{1.11}(4) for some $B',A<^*_i B' \le B$, hence by \scite{3.2} 
easily $\neg(ii)$, so $(ii) \Rightarrow (i)$.

Lastly, it suffices to prove $(i) \Rightarrow (iii)$.  Now by (i) and 
\scite{1.16} for some $\lambda$ the assumptions 
of \scite{2.9} holds hence the conclusion 
which give clause (iii).  \hfill$\square_{\scite{3.3}}$\margincite{3.3}
\enddemo
\bigskip

\demo{\stag{3.4} Conclusion}  1) $<^*_s=<_s$ and $<^*_i=<_i$, and 
${\frak K}$ is weakly nice where $<_s,<_i$ are defined in
[I,1.4(4),1.4(5)] hence $<^*_{pr} = \le_{pr}$. \nl
2) $({\Cal K},c \ell)$ is as required in [I,\S2], and the $\leq_i$,
$\leq_s$ defined in [I,\S2] are the same as those defined in [I,\S1]
for our context, of course when for $A \le B \in {\Cal K}_\infty$ we
let $c \ell(A,B)$ e minimal $A'$ such that $A \le A' \le^*_s B$. \nl
3) Also ${\frak K}$ (that is $({\Cal K},c \ell)$) is transitive local 
transparent and smooth, \nl
(see [I,2.2(3),2.3(2),2.5(5),2.5(4)]).
\enddemo
\bigskip

\demo{Proof}   1) $<^*_s=<_s$ and $<^*_i=<_i$ by \scite{3.2}, 
\scite{3.3} and see Definition in [I,\S1]. \nl
Lastly, ${\frak K}$ being weakly nice follows from \scite{3.3}, see
Definition in [I,\S1]. \nl
2) By [I,2.6]. \nl
3) By [I,2.6] the transitive local and transparent follows
(see clauses $(\delta), (\varepsilon), (\zeta)$ there).  As for
smoothness, use \scite{1.17}(5).
\hfill$\square_{\scite{3.4}}$\margincite{3.4}
\enddemo
\bn
Note that we are in a ``nice'' case, in particular
no successor function.  Toward proving it we 
characterize ``simply good''.
\proclaim{\stag{3.6} Claim}  If $A \leq_s^* B$ and $k,t \in \Bbb N$
satisfies $k + |B| \le t$, 
\ub{then} for every random enough ${\Cal M}_n$, for every 
$f:A \hra {\Cal M}_n$, we can find $g:B \hra {\Cal M}_n$ extending $f$
such that:
\mr
\widestnumber\item{$(iii)$}
\item "{$(i)$}"     ${\text{\rm Rang\/}}(g) \cap c \ell^t
({\text{\rm Rang\/}}(f),{\Cal M}_n) = \,{\text{\rm Rang\/}}(f)$,
\sn
\item "{$(ii)$}"    $\nonforkin{{\text{\rm Rang\/}}(g)}
{c \ell^t({\text{\rm Rang\/}}(f),{\Cal M}_n)}
_{{\text{\rm Rang\/}}(f)}^{{\Cal M}_n}$,
\sn
\item "{$(iii)$}"   $c \ell^k({\text{\rm Rang\/}}(g),{\Cal M}_n)
\subseteq \, {\text{\rm Rang\/}}(g) \cup c \ell^k 
({\text{\rm Rang\/}}(f),{\Cal M}_n)$.
\endroster
\endproclaim
\bigskip

\remark{\stag{3.7} Remark}  1) The reader may compare this 
with \scite{4.4}(2) below where we have successor function. \nl
2) Note that we can in clauses (i), (ii) of \scite{3.6} replace $t$ by
$k$ - this just demands less.  We shall use this freely.  Have we put
$t$ in the second appearance of $k$ in clause (iii) of \scite{3.6} the
loss would not be great: just in [I], we should systematically use
[I2.12(2)] instead of [I2.12(1)].
\endremark
\bigskip

\demo{Proof}   We prove this by induction on 
$|B\setminus A|$, but by the character of the desired 
conclusion, if $A<_s^* B<_s^* C$, to prove it for the pair $(A,C)$ 
it suffices to prove it for the pairs $(A,B)$ and $(B,C)$.
Also, if $B=A$ the statement is trivial (because we can take $f=g$). 
So, without loss of generality, $A <_{pr}^* B$ (see Definition
\scite{1.11}(6)). \nl
Let $\lambda$ be such that $(A,B,\lambda)\in {\Cal T}$ and for every
$\lambda$-closed $C\subseteq B\setminus A$ we have 
$\bold w_\lambda(A,A \cup C)>0$ and 

$$
\align
\xi =: \bold w_\lambda(A,B) = \max\{&\bold w_{\lambda_1}(A,B):
(A,B,\lambda_1) \in {\Cal T} \text{ satisfies}: \\
  &\text{for every } \lambda_1 \text{--closed non-empty } 
C \subseteq B \setminus A \\
  &\text{we have } \bold w_{\lambda_1}(A,A\cup C)> 0\}.
\endalign
$$
\mn
Choose $\varepsilon\in \Bbb R^+$ small enough.  The
requirements on $\varepsilon$, $k(*)$ will be clear by the end of the 
argument.

Let ${\Cal M}_n$ be random enough, and $f:A \hra {\Cal M}_n$. 
Now by \scite{3.2} and the definition of $c \ell^t$ we have $(*)$
and by \scite{2.10} we have $(*)_1$ where 
\mr
\item "{$(*)$}"     $|c \ell^t(f(A),{\Cal M}_n)| \le 
n^{\varepsilon/k(*)}$,
\sn
\item "{$(*)_1$}"   $|{\Cal G}|\geq n^{\xi-\varepsilon/2}$,
\ermn
where

$$
{\Cal G} = \{g:g \text{ extends } f \text{ to an embedding of } B 
\text{ into } {\Cal M}_n\}.
$$
\mn
Recall that
\mr
\item "{$\circledast_1$}"  if $A' \subseteq M \in {\Cal K}$ and
$a \in c \ell^k(A',M)$ \ub{then} for
some $C$ we have $C \subseteq c \ell^k(A',M),|C| \le k,a \in C$
and $c \ell^k(C \cap A',C) = C$ 
\ermn
(by the Definition of $c \ell^k$, see [I,\S1]).
\nl
We intend to find $g \in {\Cal G}$ satisfying the 
requirements in the claim. Now $g$
being an embedding of $B$ into ${\Cal M}_n$ extending $f$ 
follows from $g \in {\Cal G}$.  So it
is enough to prove that $<n^{\xi-\varepsilon}$ members $g$ of ${\Cal G}$ 
fail clause (i) and similarly for clauses (ii) and (iii).

More specifically, let ${\Cal G}^1 = \{g \in {\Cal G}:g(B) \cap c
\ell^t(f(A),{\Cal M}_n) \ne A\},{\Cal G}^2 = \{g \in {\Cal G}:g
\notin {\Cal G}^1$ but clause (ii) fails for $g\}$ and 
${\Cal G}^3 = \{g \in {\Cal G}$: clause (iii) fails for $g$
but $g \notin {\Cal G}^1 \cup {\Cal G}^2\}$.
So clearly it is
enough to prove ${\Cal G} \nsubseteq {\Cal G}^1 \cup {\Cal G}^2 \cup
{\Cal G}^3$ (because (i) fails for 
$g \Rightarrow g \in {\Cal G}^1$, (ii) fails
for $g \Rightarrow g \in {\Cal G}^2 \vee g \in {\Cal G}^1$ and 
(iii) fails for $g \Rightarrow g \in {\Cal G}^3 \vee g \in {\Cal G}^2 \vee 
g \in {\Cal G}^1$.
\enddemo
\bn
\ub{On the number of $g \in {\Cal G}^1$}:  For $a \in B \backslash A$
and $x \in c \ell^t(f(A),{\Cal M}_n)$ let ${\Cal G}^2_{a,x} =
\{g \in {\Cal G}^2:g(a) = x\}$, so ${\Cal G}^2 = \cup\{{\Cal G}^2_{a,x}:a \in
B \backslash A$ and 
$x \in c \ell^t(f(A),{\Cal M}_n)\}$, and by \scite{1.17}(3)
clearly $A \cup \{a\} \le_i B$ (as $A <_{pr} B$).  The rest is as in
the proof for ${\Cal G}^2$ below, only easier.
\bn
\ub{On the number of $g \in {\Cal G}^2$}:

If $g \in {\Cal G}^2$ then for some 

$$
x_g \in c \ell^t(\text{Rang}(f),{\Cal M}_n) \setminus 
\text{ Rang}(f) \text{ and } y\in B \setminus A
$$
\mn
we have: $\{x_g,g(y)\}$ is an edge of ${\Cal M}_n$. 
Note $x_g \notin g(B)$ as $g \in {\Cal G}_2$.

We now form a 
new structure $B^2=B\cup\{x^*\}$, ($x^*\notin B$), such that $g\cup\{
\langle x^*, x_g\rangle\}:B^2 \hra {\Cal M}_n$ and let 
$A^2 = B^2\restriction(A\cup \{x^*\})$. 
Now up to isomorphism over $B$ there is a finite number (i.e., with
a bound not depending on $n$) of such $B^2$, say $\langle B^2_j:
j < j^*_2 \rangle$.

For $x \in c \ell^{k,m}(\text{Rang}(f),{\Cal M}_n)$ and $j<j^*$ let

$$
\align
{\Cal G}^2_{j,x} =:\{g:&g \text{ is an embedding of } B^2_j
\text{ into } {\Cal M}_n \text{ extending } f \\
   &\text{ and satisfying } g(x^*)=x\}
\endalign
$$

$$
{\Cal G}^2_j =: \dbcu_{x \in c \ell^{k,m}(f(A),{\Cal M}_n)} \, {\Cal
G}^2_{j,x}. 
$$
\mn
So:
\mr
\item "{$(*)_2$}"   if $g \in {\Cal G}^2$ then

$$
g \in \dbcu_{j} \{\{g'\restriction B:g' \in {\Cal G}^2_{j, x}\}:j<j^*_2
\text{ and } x \in c \ell^{k,m}(f(A),{\Cal M}_n)\}.
$$
\ermn
Now, if $\neg(A^2_j <_s B^2_j)$ then as $A <^*_{pr} B$ 
easily $A^2_j<_i B^2_j$ so by \scite{3.2} using $(*)$ (with
$\varepsilon/2 - \varepsilon/k(*)$ here standing for $\varepsilon$ in
(iii) there) we have
\mr
\item "{$(*)_3$}"   if $\neg(A^2_j<_s B^2_j)$ then 
$|{\Cal G}^2_j| \leq n^{\varepsilon/2}$.
\ermn
[Why?  As $A^2_j<_i B^2_j$ on the one hand for each $x \in c
\ell^t(\text{Rang}(f),{\Cal M}_n)$
by \scite{3.2} the 
number of  $g:B^2_j
\hra {\Cal M}_n$ extending $f \cup \{\langle x^*,x \rangle\}:
A^2_j \hra {\Cal M}_n$ is $<n^{\vep /k(*)}$ and on the other hand the
number of candidates for $x$ is 
$\le |c \ell^t({\text{\rm Rang\/}}(f),{\Cal M}_n)| \le 
n^{\frac{\vep}{k(*)}}$. 
So $|{\Cal G}^2_j| \le n^{\varepsilon/k(*)} \times
n^{\varepsilon /k(*)} \le n^{2 \varepsilon /k(*)} \le
n^{\frac{\vep}{2}}$.]

If $A^2_j<_s B^2_j$, then by \scite{1.17}(14) 
still $A^2_j<_{pr} B^2_j$ and letting

$$
\align
\xi^2_j =: \max\{\bold w_\lambda(A^2_j,B^2_j):&(A^2_j,B^2_j,\lambda)
\in {\Cal T} \text{ and for every} \\
  &\lambda \text{-closed non-empty } C \subseteq B^2_j\setminus A^2_j \\
  &\text{we have } \bold w(A^2_j,A^2_j\cup C,\lambda\restriction C)>0\}
\endalign
$$
\mn
clearly $\xi^2_j<\xi-2\varepsilon$ (as we retain the ``old" edges, 
and by at least one
we actually enlarge the number of edges but we keep the number of
``vertex", i.e., equivalence classes see \scite{1.17}(9)).

So, by \scite{2.14},
\mr
\item "{$(*)_4$}"   if $A^2_j <^*_{pr} B^2_j$ then 
$|{\Cal G}^2_j|\leq n^{\xi-2\varepsilon}$.
\ermn
As $\xi -2 \varepsilon > \varepsilon$ by $(*)_3 + (*)_4$, multiplying
by $j^*$, as $n$ is large enough
\mr
\item "{$(*)_5$}"   the number of $g \in {\Cal G}  \backslash {\Cal G}^1$ 
failing clause (ii) of \scite{3.6} is $\leq n^{\xi-\varepsilon}$.
\endroster
\bn
\ub{On the number of $g \in {\Cal G}^3$}:

Next assume $g \in {\Cal G}^3$.  So there are 
$A^+$, $B^+$, $C$, $g^+$ such that
\mr
\item "{$\otimes$}"   $A \le_i A^+ \le_s B^+$, $B \le B^+$, 
$B \cap A^+=A$, $|B^+| \le |B| + k,|A^+| \le |A| +k,C \nsubseteq B
\cup A^+,B^+ \setminus  B
\subseteq C\subseteq B^+$,$c \ell^k(C \cap B,B^+) \supseteq C$ 
and $g \subseteq g^+$, $g^+:
B^+ \hra {\Cal M}_n$, $g^+(A^+)\subseteq c \ell^t(f(A),{\Cal M}_n)$.
\ermn
[Why?  As $g \in {\Cal G}^3$ there is 
$y_g \in c \ell^k(g(B),{\Cal M}_n)$ such
that $y_g \notin g(B)$ and moreover $y_g \notin c \ell^k(f(A),{\Cal
M}_n)$.  By the first statement (and $\circledast_1$ above) there 
is $C^* \subseteq c \ell^k
(g(B),{\Cal M}_n)$ with $\le k$ elements such
that $y_g \in C^*$ and $C^* \cap g(B) \le_i C^*$.
Let $B^* = g(B) \cup C^* \le {\Cal M}_n$.  Let 
$B^+,g^+$ be such that $B \le B^+ \in {\Cal K},g \subseteq g^+,g^+$ an
isomorphism from $B^+$ onto $B^*$, and let
$C = g^{-1}(C^*)$.  Lastly, choose $A^+$ such that $A' \le_i A^+
\le_s B^+$, clearly it exists by \scite{1.17}(2).  
Now $|A^+| \le |B| + |C| \le t$ by the assumptions on $A,B,k,t$ 
hence $g^+(A^+) \subseteq c
\ell^t(f(A),{\Cal M}_n)$ but as $g \in {\Cal G}^3$ we
have $g \notin {\Cal G}^1$ hence
$A = g(B) \cap c \ell^t(f(A),{\Cal M}_n)$ so 
we have $A^+ \cap B = A$.  
Also $C \nsubseteq B \cup A^+$, otherwise, as $g \notin {\Cal G}^2,g
\notin {\Cal G}^1$ we have $\nonforkin {B}{A^+}_{A}^{B^+}$ hence
$\nonfork{C \cap B}{C \cap A^+}_{C \cap A}$ but as $C \cap B <^*_i C$
so by smoothness (e.g. \scite{1.17}(5))
we get $C \cap A <^*_i C \cap A^+$ hence $C \cap A^+ \subseteq c
\ell^k(A,B^+)$ hence $C^* \backslash g(B) = g^+(C \backslash B)
\subseteq g^+(C \cap A^+) \subseteq c \ell^k(f(A),{\Cal M}_n)$
hence $y_g \in c \ell^k(f(A),{\Cal M}_n)$, contradiction.
So $\otimes$ holds.]

Now if $\lambda \in \Xi(A^+,B^+)$ (see \scite{1.8}(2)), as $A
\subseteq A^+ \and B \backslash A \subseteq B^+ \backslash A^+$ easily
$\lambda\restriction (B \setminus A) \in \Xi(A,B)$, 
and $\bold w(A^+,B^+,\lambda) < \bold w(A^+,B \cup A^+) \le 
\bold w(A,B,\lambda 
\restriction (B\setminus A))$, see \scite{1.17}(8) noting that $A^+
\cup B <^*_i B^+$.  \nl
Let $\{(A^+_j,B^+_j):j<j^*_3\rangle$ list the possible
$(A^+,B^+,\lambda)$ up to isomorphism over $B$ as described above. Let

$$
{\Cal G}^3_{j,h} =: \{g \in {\Cal G}:g \text{ embeds } B^+_j \text{ into }
{\Cal M}_n \text{ extending }f \text{ and moreover } h\}
$$
\mn
for any $h \in {\Cal H}^3_j =: \{h:h:A^+_j \hra {\Cal M}_n$ extending
$f\}$, so $h$ necessarily satisfies 
$h(A^+_j) \subseteq c \ell^k(f(A),{\Cal M}_n) \subseteq c
\ell^t(f(A),{\Cal M}_n)$.
\nl
Now easily (for random enough ${\Cal M}_n$)
\mr
\item "{$(*)_6$}"   if $g \in {\Cal G}^3$ then
$$
g \in \dbcu_{j<j^*_3} \, \dbcu_{h \in {\Cal H}^3_j} 
\{g'\restriction B:g' \in {\Cal G}^3_{j,h}\},
$$
\sn
\item "{$(*)_7$}"  $|{\Cal G}^3_{j,h}|< n^{\xi-2\varepsilon}$ for each
$h \in {\Cal H}^3_j$
\sn
\item "{$(*)_8$}"  $|{\Cal H}^3_j| < |c \ell^k(f(A),{\Cal M}_n)|^k \le
|c \ell^t(f(A),{\Cal M}_n)|^k < n^{\varepsilon,k}$.
\ermn
Together
\mr
\item "{$(*)_9$}"   the number of $g \in {\Cal G}^3$ is
$<n^{\xi-\varepsilon}$.
\hfill$\square_{\scite{3.6}}$\margincite{3.6}
\endroster
\bigskip

\demo{\stag{3.8} Conclusion}  If $A<^*_s B$ and $B_0 \subseteq B$ 
and $k \in \Bbb N$ \ub{then} the tuple $(B,A,B_0,k)$ is simply good 
(see Definition [I,2.12(1)]).
\enddemo
\bigskip

\demo{Proof}   Read \scite{3.6} and Definition [I,2.12(1)].
\hfill$\square_{\scite{3.8}}$\margincite{3.8}
\enddemo
\bn
\centerline {$* \qquad * \qquad *$}
\bn
Toward simple niceness the ``only'' thing left is the universal part, i.e.,
Definition [I,2.13(1)(A)]. \nl
The following claims \scite{3.9}, \scite{3.10} do not use \S5 and have
nothing to do with probability; 
they are the crucial step for proving the satisfaction of
Definition [I,2.13(1)(A)] in our case; claim 
\scite{3.9} is a sufficient condition for goodness (by \scite{3.8}). 
Our preceding the actual proof (of \scite{3.11}) by the two claims
(\scite{3.9},\scite{3.10}) and separating them
is for clarity, though it has a bad effect on the bound; 
see also - \scite{3.9} using ``$c \ell^{k,m}(\bar a b,M)$" instead 
$c \ell^k(\bar a b,M)$ when $k'<k$ may improve the bound.
\proclaim{\stag{3.9} Claim}   For every $k$ and $\ell$ (from $\Bbb N$) 
there are natural numbers
$t=t(k, \ell)$ and $k^*(k,\ell)\geq t$, $k$ such that for 
any $k^*\geq k^*(k, \ell)$ we have:
\mr
\item "{$(*)$}"   if $m^\otimes\in \Bbb N$ and $M \in {\Cal K}$, $\bar a
\in {}^{\ell\geq}M$, $b\in M$ \ub{then}
{\roster
\itemitem{ $\otimes$ }  the set
$$
\align
{\Cal R} =: \{(c,d):&d \in c \ell^k(\bar a b,M) \setminus 
c \ell^{k^*,m^\otimes+k}(\bar a,M) \text{ and} \\
  &c\in c \ell^{k^*,m^\otimes}(\bar a,M) \text{ and } \{c,d\}
\text{ is an edge of }M\}
\endalign
$$
has less than $t$ members.
\endroster}
\endroster
\endproclaim
\bigskip

\demo{Proof}  If $k=0$ this is trivial so \wilog \, $k>0$.  
Choose $\varepsilon \in \Bbb R^+$ small enough such that
\mr
\item "{$(*)_1$}"  $C_0 <^* C_1 \and (C_0,C_1,\lambda) \in {\Cal T} \and 
|C_1|\leq k \Rightarrow \bold w_\lambda(C_0,C_1)\notin 
[-\varepsilon,\varepsilon]$
\ermn
(in fact we can restrict ourselves to the case $C_0<^*_i C_1$)).

Choose $\bold c\in \Bbb R^+$ large enough such that
\mr
\item "{$(*)_2$}"  $(C_0,C_1,\lambda) \in 
{\Cal T},|(C_1 \setminus C_0)/\lambda|\leq
k \Rightarrow \bold w_\lambda(C_0,C_1)\leq \bold c$.
\ermn
(so actually $\bold c=k$ is enough). 
Choose $t_1>0$ such that $t_1 > \bold c /\varepsilon$ and $t_1 > 2$. \nl
Choose $t_2\geq 2^{2^{t_1+k+\ell}}$ (overkill, we mainly need to apply
twice the $\Delta$-system lemma; but note that in the proof of
\scite{3.10} below we will use Ramsey Theorem). 
Lastly, choose $t> k^2 t_2$ and let $k^* \in \Bbb N$ be large enough 
which actually means that $k^*>k \and k^*\ge (k+1) \times t_2$ so
$k^*(k,\ell) =: (k+1) \times t_2$ is O.K.

Suppose we have $m^\otimes$, $M$, $\bar a$, $b$ as in $(*)$ but such that
the set ${\Cal R}$ has at least $t$ members. Let $(c_i,d_i)\in {\Cal R}$ for
$i<t$ be pairwise distinct
\footnote{Note: we do not require the
$d_i$'s to be distinct; though if $w = \{i:d_i=d^*\}$ has $\geq k'>
\frac{1}{\alpha}$ elements then $d^*\in c \ell^{k'}
(c \ell^{k^*,m^\otimes+k}(\bar a, M))$}. 

As $d_i\in c \ell^{k}(\bar a b,{\Cal M}_n)$, we
can choose for each $i<t$ a set $C_i\leq M$ such that:
\mr
\widestnumber\item{$(iii)$}
\item "{$(i)$}"     $C_i\leq M$,
\sn
\item "{$(ii)$}"  $|C_i|\leq k$,
\sn
\item "{$(iii)$}"  $d_i\in C_i$,
\sn
\item "{$(iv)$}"    $C_i\restriction (C_i\cap(\bar ab))<_i C_i$.
\ermn
For each $i<t$, as $C_i\cap c \ell^{k^*,m^\otimes+k}(\bar{a},M)$ is 
a proper subset of $C_i$ (this is witnessed by $d_i$, 
i.e., as $d_i \in C_i \setminus
(C_i\cap c \ell^{k^*,m^\otimes+k}(\bar{a},M_n)))$ clearly this set has 
$<k$ elements and hence for some $k[i]<k$ we have
\mr
\item "{$(v)$}"   $C_i\cap c \ell^{k^*,m^\otimes+k[i]+1}(\bar a,M)
\subseteq c \ell^{k^*,m^\otimes+k[i]}(\bar a,M)$.
\ermn
So without loss of generality
\mr
\item "{$(vi)$}"    $i<t/k^2 \Rightarrow k[i]=k[0] \and |C_i|=|C_0| =
k' \le k$,
\ermn
remember $t_2 < t/k^2$; also
\mr
\item "{$(vii)$}"  $b \in C_i$. \nl
[Why?  If not then by clause (iv) we have $(C_i \cap \bar a) <_i C_i$,
hence $d_i \in C_i \subseteq c \ell^k(\bar a,M) \subseteq c
\ell^{k^*,m^\otimes +k}(\bar a,M)$, contradiction.]
\ermn
As $k^*\geq t_2\times (k+1)$ (by the assumption on $k^*$), 
clearly $|\dbcu_{i<t_2} C_i| \cup |\{c_i:i < t_2\}| \le 
\dsize \sum_{i<t_2} |C_i| + t_2 \le
\dsize \sum_{i < t_2} k + t_2 \le t_2 \times (k+1) \le k^*$ and we define

$$
D = \dbcu_{i<t_2} C_i \cup \{c_i:i < t_2\}
$$

$$
D'=D \cap c \ell^{k^*,m^\otimes+k[0]}(\bar a,M).
$$
\mn
So by the previous sentence we have $|D'| \le |D| \le k^*$. \nl
Now
\mr
\item "{$\otimes_0$}"  $D' <_s D$. \nl
[Why?  As otherwise ``there is $D''$ such that $D'<_i D'' \le_s D$ so as
$|D''|\leq |D| \leq k^*$ clearly $D'' \subseteq 
c \ell^{k^*,m^\otimes+k[0]+1}(\bar a b,M)$; contradiction.]  

\ermn
So we can choose $\lambda\in \Xi(D',D)$ (see Definition \scite{1.8}(2)).
Let $C_i=\{d_{i,s}:s<k'\}$, with $d_{i,0}=d_i$ and recalling (vii) also
$b \neq d_i \Rightarrow b=d_{i,1}$ and with no repetitions.

Clearly $d_{i,0} = d_i \notin D'$.  
By the finite $\Delta$-system lemma for some $S_0,S_1,S_2\subseteq
\{0,\ldots,k'-1\}$ and $u\subseteq\{0,\ldots,t_2-1\}$ with $\geq 
t_1$ elements we have:
\mr
\item "{$\oplus_1(a)$}"   $\lambda' =: \{(s_1,s_2):d_{i,s_1}
\lambda d_{i,s_2}\}$ is the same for all $i\in u$ and $S_0 =
\{0,\dotsc,k'\}  \backslash \text{ Dom}(\lambda')$ so $d_{i,s} \in D'
\Rightarrow i \in S_0$
\sn
\item "{$(b)$}"   for each $j < \ell g(\bar{a})+1$, and $s<k'$, 
the truth value of $d_{i,s}=(\bar ab)_j$ is the same for all 
$i \in u$ and $s \in S_0 
= \{0,\dotsc,k'-1\} \backslash \text{ Dom}(\lambda')$ so
$d_{i,s} \in D' \Rightarrow i \in S_0$ 
\sn
\item "{$(c)$}"   $d_{i_1,s_1} = d_{i_2,s_2} \Rightarrow 
s_1=s_2$ for $i_1,i_2 \in u$ 
\sn
\item "{$(d)$}"   $d_{i_1,s} = d_{i_2,s} \Leftrightarrow s\in S_1$ for $i_1
\neq i_2 \in u$
\sn
\item "{$(e)$}"   $d_{i_1,s_1} \lambda d_{i_2,s_2} \Rightarrow d_{i_1,s_1}
\lambda d_{i_1,s_2} \and d_{i_1,s_2} \lambda d_{i_2,s_2}$ for $i_1 \ne
i_2 \in u$
\sn
\item "{$(f)$}"   $d_{i_1,s} \lambda d_{i_2,s} \Leftrightarrow s\in S_2$ for
$i_1 \neq i_2 \in u$: so $i \in u \and s \in S_2 \Rightarrow d_{i,s}
\notin D'$
\sn
\item "{$(g)$}"   the statement
$b = d_{i,0}$ has the same truth value for all $i\in u$
\sn
\item "{$[(h)$}"   in \S7 we also demand $(d_{i_1,s_1} Sd_{i_1,s_2}) =
(d_{i_2,s_1} Sd_{i_1,s_2})$ where $S$ is the successor relation].
\ermn
Now we necessarily have:
\mr
\item "{$\oplus_2$}"   $0 \notin S_2$ (i.e., 
$\lambda \restriction \{d_i:i\in u\}$ is equality).
\ermn
[Why?  Otherwise, letting $X = d_i/\lambda$ for any $i \in u$, 
the triple $(D',D' \cup X,\lambda \restriction X) \in
{\Cal T}$ has weight

$$
\align
\bold w(D',D'\cup X,\lambda\restriction X) &= \bold v(D',D' \cup X,
\lambda\restriction X) \\
  &- \alpha \bold e(D',D'\cup X,\lambda \restriction X) \\
  &= 1 - \alpha \times|\{e:e \text{ an edge of } M 
\text{ with one end in} \\
  & \quad D' \text{ and the other in } X\}|.
\endalign
$$
\mn
Now as $c_i \in c \ell^{k^*,m^\otimes}(\bar{a},M)$ clearly $c_i \in D'$
and the pairs
$\{c_i,d_i\} \in \text{ edge}(M)$ are distinct for different $i$
clearly the number above is \nl
$\le 1-\alpha\times|\{(c_i, d_i):i\in u\}| = 
1-\alpha \times|u|=1 - \alpha \times t_1 < 0$;
contradiction to $\lambda\in\Xi(D',D)$.]

Let $D_0= \bar a \cup \bigcup \{d_{i,s}/\lambda:s \in S_2 \text{ and } i \in
u\}$, clearly $D_0$ is $\lambda$-closed subset of $D$ though not
necessarily $\subseteq \text{ Dom}(\lambda) = D \backslash D'$ because
of $\bar a$.
\mr
\item "{$\oplus_3$}"   $b = d_{i,1}$ and $1\in S_1 \backslash S_0$ and
$0 \notin S_0 \cup S_1 \cup S_2$ and $S_1 \backslash S_0 \subseteq
S_2$ (hence $b \in D_0$). \nl
[Why?  The first two clauses hold as $b\in C_i$, $b\in\{d_{i,0},d_{i,1}\}$
and $\oplus_2$ and (g) of $\oplus_1$.  The last clause holds by
$\oplus_1$(d),(f) and the ``hence $b \in D_0$" by the definition of
$D_0,S_1 \backslash \text{ Dom}(\lambda') \subseteq S_2$ and the first clause.
Also $0 \notin S_0 \cup S_1 \cup S_2$ should be clear.] 
\sn
\item "{$\oplus_4$}"   For each $i\in u$ we have 
$\bold w_\lambda(C_i\cap D_0, C_i)<0$.
\nl
[Why? As $C_i \cap (\bar a b) \subseteq C_i\cap D_0$ by clauses (b) + (f) of
$\oplus_1$ and by monotonicity of $<_i$ we have 
$C_i\restriction (C_i\cap \bar{a}b) <_i C_i \Rightarrow
C_i\cap D_0\leq_i C_i$ but $d_{i,0} = d_i \in C_i \backslash C_i \cap D_0$.]
\ermn
Hence
\mr
\item "{$\oplus_5$}"   $\bold w_\lambda(C_i\cap D_0,C_i) \le -
\varepsilon$ for $i \in u$. 
\nl
[Why?  See the choice of $\varepsilon$.]
\ermn
Let
$$
D_1=: D'\cup \bigcup \{d_{i,s}/\lambda:i \in u,s < k'\} = 
D' \cup D_0 \cup \{d_{i,s}/\lambda:i \in u,s < k' \and s \notin S_2\}
$$
so clearly $D_1$ is $\lambda$-closed subset of $D$ including $D'$ but
$D_1\neq D'$ as $i \in u \Rightarrow d_i \in D_1$ by $\oplus_2$.
Also clearly
\mr
\item "{$\oplus_6$}"  $D'\subseteq D' \cup D_0\subseteq D_1\subseteq D$
and $D_0$, $D_1$ are $\lambda$-closed.
\ermn
So, as we know $\lambda\in \Xi(D',D)$, we have
\mr
\item "{$\oplus_7$}"   $\bold w_\lambda(D', D_1) > 0$.
\ermn
Now:

$$
\align
\bold w_\lambda(D', D_1) &= \bold w_\lambda(D',\bigcup\{x/\lambda: 
x \in \dbcu_{i\in u} C_i\setminus D'\}\cup D') \\
  &= \bold w_\lambda(D',D'\cup D_0) + \bold w_\lambda(D'\cup D_0,
D'\cup D_0 \\
  &\hskip 85pt \cup\bigcup\{d_{i,s}/\lambda:i \in u,s<k',s \notin S_0
  \cup S_2\}) \\
  &\quad \text{ [now by \scite{6.7A} below with } B_i = \{d_{i,s}:s <
k',s \notin S_2 \cup S_0\} \text{ and} \\
  &\quad \quad B^+_i = \cup\{d_{i,s}/\lambda:s < k',s \notin S_2\}] \\  
  &\le \bold w_\lambda(D',D'\cup D_0) + \dsize \sum_{i\in u}
\bold w_\lambda(D'\cup D_0,D'\cup D_0\cup\{d_{i,s}:s<k',s\notin S_2
\cup S_0\}) \\
  &\quad \text{ [now as } C_i = \{d_{i,s}:s < k'\} \text{ and } d_{i,s} \in
D' \cup D_0 \text{ if } s < k',s \in S_0 \cup S_2,i \in u] \\
  &\le \bold w_\lambda(D',D' \cup D_0) + \dsize \sum_{i \in u} \bold
w_\lambda(D' \cup D_0,D' \cup D_0 \cup C_i) \\
  &\quad \text{ [now as } \bold w_\lambda(A_1, B_1) \le \bold
w_\lambda(A,B) \\
  &\quad \text{ when } A \le A_1 \le B_1,A \le B \le B_1,B_1
\setminus A_1=B\setminus A) \text{ by \scite{1.15}(3)}] \\
  &\le \bold w_\lambda(D'\cap D_0,D_0)+
\dsize \sum_{i\in u} \bold w_\lambda(C_i\cap D_0, C_i) \\
  &\quad \text{ [now by the choice of } \bold c,D_0  
\text{ i.e., } (*)_2 \text{ and the choice of }\varepsilon,u + \oplus_5
\text{ respectively]} \\
  &\le \bold c + |u|\times (-\varepsilon) = \bold c-t_1\varepsilon<0,
\endalign
$$
\mn
contradicting the choice of $\lambda$, i.e., 
$\oplus_7$.  \hfill$\square_{\scite{3.9}}$\margincite{3.9}
\enddemo
\bigskip

\demo{\stag{6.7A} Observation}  Assume
\mr
\item "{$(a)$}"  $A \le^* A \cup B_i \le^* A \cup B^+_i \le^* B$ for $i \in u$
\sn
\item "{$(b)$}"   $B \backslash A$ is the disjoint union of $\langle
B^+_i:i \in u \rangle$
\sn
\item "{$(c)$}"  $\lambda$ is an equivalence relation on $B \backslash
A$
\sn
\item "{$(d)$}"  each $B^+_i$ is $\lambda$-closed
\sn
\item "{$(e)$}"  $B^+_i = \cup\{x/\lambda:x \in B_i \backslash A\}$
\ermn
\ub{Then} $\bold w_\lambda(A,B) \ge \Sigma\{\bold w_\lambda(A,B_i):i \in u\}$.
\enddemo
\bigskip

\demo{Proof}  By clause (b) + (d)

$$
\bold v_\lambda(A,B) = \Sigma\{\bold v_\lambda(A,A \cup B^+_i):i \in
u\} = \Sigma\{\bold v_\lambda(A,A \cup B_i):i \in u\}
$$
\mn
and by clause (b) the set $e_\lambda(A,B)$ contains the disjoint union of
$\langle e_\lambda(A,B_i):i \in u \rangle$.  Together the result
follows.
\hfill$\square_{\scite{6.7A}}$\margincite{6.7A}
\enddemo
\bigskip

\proclaim{\stag{3.10} Claim}  For every $k$, $m$ and $\ell$ from 
$\Bbb N$ for some $m^*=m^*(k, \ell, m)$, for any 
$k^*\geq k^*(k,\ell)$ (the function $k^*(k,l)$ is the one from 
claim \scite{3.9}) satisfying $k^* \geq k \times m^*$  we have
\mr
\item "{$(*)$}"   if $M \in {\Cal K},\bar a \in{}^{\ell \geq}M$ and
$b \in M \setminus c \ell^{k^*,m^*}(\bar a,M)$
\ub{then} for some $m^\otimes\leq m^*-m$ we have
$$
c \ell^{k}(\bar a b,M) \cap c \ell^{k^*,m^\otimes+m}(\bar a,M)
\subseteq c \ell^{k^*,m^\otimes}(\bar a,M).
$$
\endroster
\endproclaim
\bigskip

\demo{Proof}   For $k=0$ this is trivial so \wilog \, $k>0$.
Let $t= t(k,\ell)$ be as in the previous claim \scite{3.9}.
Choose $m^*$ such that, e.g., 
$\lfoot m^*/(km) \rfoot \rightarrow (t+5)^2_{2^{k!+\ell}}$ in the
usual notation in Ramsey theory.   We could get more reasonable
bounds but no need as for now.  Remember that $k^*(k, \ell)$ is from
\scite{3.9} and $k^*$ is any natural number $\geq k^*(k, \ell)$ such that
$k^*\geq km^*$.

If the conclusion fails, then the set 

$$
Z =: \{j\leq m^* - k:c \ell^k(\bar a b,M) \cap c \ell^{k^*,j+1}
(\bar a,M) \nsubseteq c \ell^{k^*,j}(\bar a,M)\}
$$
\mn
satisfies:

$$
j\leq m^*-m-k \Rightarrow Z\cap [j,j+m) \ne \emptyset.
$$
\mn
Hence $|Z|\geq (m^*-m-k)/m$. 
\nl
For $j\in Z$ there are $C_j\leq M$ and $d_j$ such that

\block
$|C_j|\leq k$ and $(C_j\cap(\bar ab)) <^*_i C_j$, and
$d_j\in C_j \cap c \ell^{k^*,j+1}(\bar a,M)\setminus 
c \ell^{k^*,j}(\bar a,M)$.
\endblock
\sn
Now we use the same argument as in the proof of \scite{3.9}.
\nl
As $d_j \in C_j \cap c \ell^{k^*,j+1}(\bar a,M) \setminus
c \ell^{k^*,j}(\bar a,M)$ we will get that $C_j \cap c 
\ell^{k^*,j}(\bar a,M)$ is
a proper subset of $C_j \cap c \ell^{k^*,j+1}(\bar a,M)$(witness by $d_j$) 
so $|C_j \cap c \ell^{k^*,j}(\bar a,M)|< 
|C_j \cap c \ell^{k^*,j+1}(\bar a,M)| \leq k$
so $|C_j \cap c \ell^{k^*,j}(\bar a,M)|<k$.  Hence for 
some $k_j \in \{1,\dotsc,k\}$ we have 
$C_j \cap c \ell^{k^*,m^*-k_j+1}(\bar a,M) \subseteq c \ell^{k^*,m^*-k_j}
(\bar a,M)$ hence for some $k' \in \{1,\dotsc,k\}$ we have $|Z'| \geq 
(m^*-m-k)/(mk)$ where $Z' = \{j \in Z:k_j = k'\}$.

Let $C_j=\{d_{j, s}: s<s_j\leq k\}$ with $d_{j,0}=d_j$ and no
repetitions. 
We can find $s^* \le k$ and $S_1,S_0 \subseteq \{0,\dotsc,s^*-1\}$ and
$u \subseteq Z'$ satisfying $|u|=t+5$  such that
(because of the partition relation):
\mr
\item "{$(a)$}"   $i \in u \Rightarrow s_j=s^*$,
\sn
\item "{$(b)$}"   for each $j < \ell g(\bar{a})+1$ and $s < s^*$ the 
truth value of $d_{i,s} = (\bar a b)_j$ is the same for all 
$i\in u$, 
\sn
\item "{$(c)$}"   if $i\neq j$ are from $u$ then $|i-j|>k+1$,
i.e., the $C_i$'s for $i \in u$ are quite far from each other
\sn
\item "{$(d)$}"   the truth value of 
``$\{d_{i,s_1}, d_{i,s_2}\}$ is an edge'' is the same for all 
$i\in u$, 
\sn
\item "{$(e)$}"   for all $i_0<i_1$ from $u$:
$$
[d_{i_0,s} \in c \ell^{k^*,i_1}(\bar{a},M)] \Leftrightarrow s\in S_0,
$$
\sn
\item "{$(f)$}"   for all $i_0<i_1$ from $u$:
$$
d_{i_1,s} \in c \ell^{k^*,i_0}(\bar{a},M) \Leftrightarrow s\in S_1,
$$
\sn
\item "{$(g)$}"   for each $s<s^*$, the sequence $\langle d_{i, s}:
i\in u\rangle$ is constant or with no repetition,
\sn
\item "{$(h)$}"   if $d_{i_1,s_1}= d_{i_2, s_2}$ 
then $d_{i_1, s_1}= d_{i_1,s_2}= d_{i_2,s_2}$, moreover, $s_1 = s_2$
(recalling that $\langle d_{i,s}:s < s_j \rangle$ is with no repetitions).
\ermn
Now let $i(*)$ be, e.g., the third element of the set $u$ and

$$
B_1 =: C_{i(*)}\cap c \ell^{k^*,\min(u)}(\bar a,M), \text{ and }
 B_2 =: C_{i(*)}\cap c \ell^{k^*,\max(u)}(\bar a,M).
$$
\mn
So
\mr
\item "{$\circledast_1$}"   $B_1 <^* B_2 \le^* C_{i(*)}$ 
(note: that $B_1\neq B_2$ because $d_{i(*)}\in B_2 \setminus B_1$) and
\sn
\item "{$\circledast_2$}"   $(\bar a b)\cap B_2\subseteq B_1$
by clause (b) and
\sn
\item "{$\circledast_3$}"   there is no edge in 
$(C_{i(*)}\setminus B_2)\times (B_2\setminus B_1)$.
\ermn
Why?  Toward contradiction assume that this fails.
Let the edge be $\{d_{i(*),s_1},d_{i(*), s_2}\}$ with

$$
d_{i(*),s_1}\in C_{i(*)} \setminus B_2 \text{ and }
d_{i(*),s_2}\in B_2\setminus B_1;
$$
\mn
hence
\mr
\item "{$(*)_1$}"   $d_{i(*),s_1}\in C_{i(*)} \setminus 
c \ell^{k^*,\max(u)}(\bar a,M)$ and \nl
$d_{i(*),s_2}\in c \ell^{k^*,\max(u)}(\bar a, M)\setminus
c \ell^{k^*,\min(u)}(\bar a,M)$ \nl
and $\{d_{i(*),s_1},d_{i(*),s_2}\}$ is an edge. 
\ermn
Hence by clause (d)
\mr
\item "{$(*)_2$}"  $\{d_{i,s_1},d_{i,s_2}\}$ is an edge for every $i \in u$
\ermn
and by clauses $(e),(f)$ we have
\mr
\item "{$(*)_3$}"   if $i_0<i_1<i_2$ are in $u$ then $d_{i_1,s_2} \notin
c \ell^{k^*,i_0}(\bar a,M)$ and $d_{i_1,s_2} \in c \ell^{k^*,i_2}(\bar
a,M)$ and $d_{i_1,s_1} \notin c \ell^{k^*,i_2}(\bar a,M)$ 
\ermn
and so necessarily
\mr
\item "{$(*)_4$}"   $\langle d_{i,s_2}:i \in u \rangle$ is with no 
repetitions. \nl
[Why?  By clause (g) and $(*)_3$.]
\ermn  
So the set of edges 
$\{\{d_{i,s_1},d_{i,s_2}\}:i\in u \text{ but } |u\cap i|\geq 2 
\text{ and } |u\setminus i|\geq 2\}$ 
contradicts \scite{3.9} using $m^\otimes =\max(u)-k$ there
(and our choice of parameters and $C_i \subseteq c \ell^k(\bar a
b,M)$).  So $\circledast_3$ holds.
\sn
As $C_{i(*)}\restriction (\bar ab)<_i C_{i(*)}$ and
 $B_2\cap(\bar ab) \subseteq B_1$ (by $\circledast_2$), 
clearly $C_{i(*)}\cap \bar a b \subseteq
C_{i(*)}\setminus (B_2\setminus B_1)\subset C_{i(*)}$, 
the strict $\subset$ as $d_{i(*)} \in C_{i(*)} \cap (c
\ell^{k^*,i(*)+1}(\bar a b,M) \backslash c \ell^{k^*,i(*)}(\bar a
b,M)) \subseteq B_2\setminus B_1$.
But, as stated above, $\nonforkin{C_{i(*)}\setminus (B_2\setminus
B_1)}{B_2}_{B_1}^{}$, hence by the previous sentence (and smoothness, see
\scite{1.17}(5)) we get $B_1<^*_i B_2$; also 
$|B_2|\leq |C_{i(*)}|\leq k\leq k^*$.
By their definitions, $B_1\subseteq c \ell^{k^*,\min(u)}(\bar{a},M)$, but
$B_1\leq^*_i B_2$, $|B_2|\leq k \le k^*$ and hence 
$B_2 \subseteq c \ell^{k^*,2^{\text{nd}}\text{member of }u}(\bar a,M)$. 
Contradiction to the choice of $d_{i(*)}$. 
\hfill$\square_{\scite{3.10}}$\margincite{3.10}
\enddemo
\bigskip

\proclaim{\stag{3.11} Lemma}  For every $k$, $m$ and $\ell$ (from
$\Bbb N$), for some $m^*$, $k^*$ and $t^*$ we have:
\mr
\item "{$(*)$}"   if $M \in {\Cal K},\bar a \in {}^{\ell\geq }M$ and 
$b\in M\setminus c \ell^{k^*,m^*}(\bar a,M)$ \ub{then} for 
some $m^\otimes\leq m^*-m$ and $B$ we have
{\roster
\itemitem{ $(i)$ }     $|B|\leq t^*$,
\sn
\itemitem{ $(ii)$ }   $\bar a \subseteq B \subseteq c \ell^k(B,M)\subseteq
c \ell^{k^*,m^\otimes} (\bar a,M)$,
\sn
\itemitem{ $(iii)$ }   $c \ell^{k^*,m^\otimes+m}(\bar a,M),
(c \ell^k(\bar a b,M) \setminus c \ell^{k^*,m^\otimes +m}
(\bar a,M))\cup B$ are free over $B$ inside $M$,
\sn
\itemitem{ $(iv)$ }    $B \leq_s^* B^* =: M \restriction
((c \ell^{k^*}(\bar{a}b,M) \setminus c \ell^{k^*,m^\otimes+m}(B,M))
\cup B)$.
\endroster}
\endroster
\endproclaim
\bigskip

\remark{\stag{6.9areM} Remark}   Clearly this will finish the proof of
simply nice. 
\endremark
\bn
\ub{\stag{6.9B} Comments}:  Let us describe the proof below. \nl
1)  In the proof we apply the last two claims.  By them we
arrive to the following situation: inside $ c \ell^k(\bar a b,M)$ we have
$B \le B^*,|B| \le t^*$ and there is no ``small" $D$ such that $B
<^*_i D \le B^*$ and we have to show that $B <^*_s B^*$, a kind of
compactness lemma. \nl
2) Note that for each $d \in c \ell^k(\bar a b,M)$ there is $C_d
\subseteq c \ell^k(\bar a b,M)$ witnessing it, i.e., $C_d \cap (\bar a
b) \le_i C_d,d \in C_d,|C_d| \le k$.  To prove the statement above we
choose an increasing sequence $\langle D_i:i \le i(*) \rangle$ of
subsets of $B^*,D_0 = B \cup \{b\},|D_i|$ has an a priori bound,
$D_{i+1}$ ``large" enough.  So by our assumption toward contradiction
$B <^*_s D_{i(*)}$ hence there is
$\lambda \in \Xi(B,D_{i(*)})$, \wilog, $B^* = B \cup \bigcup\{C_d:d
\in D_{i(*)}\}$.  For each $i < i(*)$ we try to ``lift" $\lambda
\restriction (D_i \backslash B)$ to $\lambda^+ \in \Xi(B,B^*)$, a
failure will show that we could have put elements satisfying some
conditions in $D_{i+1}$ so
we had done so.  As this occurs for every $i < i(*)$, by weight
computations we get contradiction.
\bigskip

\demo{Proof}  Without loss of generality $k>0$.  Let $t = t(k,\ell),
k^*(k,\ell)$ be as required in \scite{3.9} (for our given $k$, $\ell$).

Choose $m(1)=t\times (m+1) +k +2$ and let $t^*=t+\ell +k$.

Choose $m^*$ as in \scite{3.10} for $k$ (given in \scite{3.11}), $m(1)$
(chosen above) and $\ell$ (given in \scite{3.11}), i.e., $m^* =
m^*(k,m(1),\ell)$ . Let $\varepsilon^* \in \Bbb R^{>0}$ be such that

$$
(A',B',\lambda)\in {\Cal T} \and  |B'| \le k \and A'\neq B' \Rightarrow
\bold w_\lambda(A',B') \notin (-\varepsilon^*,\varepsilon^*).
$$
\mn
Let $i(*)> \frac{1}{\varepsilon^*}$. Define inductively $k^*_i$ for
$i\leq i(*)$ as follows

$$
k^*_0 = \max\{k^*(k,\ell),m \times k,m^* \times t^* +1\}
$$

$$
k^*_{i+1}= 2^{2^{k^*_i}}
$$
\mn
and lastly let

$$
k^*= k\times k^*_{i(*)}.
$$
\mn
We shall prove that $m^*$, $k^*$, $t^*$ are as required in
\scite{3.11}.  So let 
$M$, $\bar a$, $b$ be as in the assumption of $(*)$
of \scite{3.11}.  So $M\in {\Cal K}$, $\bar a \in {}^{\ell\geq}M$ and 
$b \in M \setminus c \ell^{k^*,m^*}(\bar a,M)$, but this means that the
assumption of $(*)$ in \scite{3.10} holds for $k$, $m(1)$, $\ell$,
so we can apply it (i.e., as $m^* = m^*(k,m(1),\ell),k^* \ge
k^*(k,\ell)$
where $k^*(k,\ell)$ is from \scite{3.9} and $k^* \ge k \times m^*$ as
$k^* \ge k^*_{i(*)} > k^*_0 \ge m^* \times k$).
So for some $r\leq m^*-m(1)$ we have
\mr
\item "{$\oplus_1$}"   $c \ell^k(\bar a b,M) \cap c \ell^{k^*,r+m(1)}(\bar a,
M)\subseteq c \ell^{k^*,r}(\bar a,M)$.
\ermn
Let us define

$$
\align
{\Cal R} = \{(c,d):&d \in c \ell^k(\bar a b,M) \setminus c \ell^{k^*,r+m(1)}
(\bar a,M) \text{ and} \\
  &c \in c \ell^{k^*,r+m(1)-k}(\bar a,M) \text{ and} \\
  &\{c,d\} \text{ is an edge of }M\}.
\endalign
$$
\mn
How many members does ${\Cal R}$ have?  By \scite{3.9} (with 
$r + m(1) - k$ here standing for $m^\otimes$ there as $k^* \ge
k^*(k,\ell)$) at most $t$ members.  But by $\oplus_1$
above

$$
\align
{\Cal R} = \{(c,d):&d \in c \ell^k(\bar a b,M) \setminus c \ell^{k^*,r}
(\bar a,M) \text{ and} \\
  &c \in c \ell^{k^*,r+m(1)-k}(\bar a,M) \text{ and} \\
  &\{c,d\} \text{ is an edge of }M\}.
\endalign
$$
\mn
But $t \times (m+1) + 1 < m(1)-k$ by the choice of $m(1)$ (and, of course,
$c \ell^{k^*,i}(\bar a,M)$ increase with $i$) hence for some
$m^\otimes \in \{r+1, \ldots, r+m(1)-k-m\}$ we have
\mr
\item "{$\oplus_2$}"   $(c,d)\in {\Cal R} \Rightarrow c \notin 
c \ell^{k^*,m^\otimes+m}(\bar a,M) \setminus c \ell^{k^*,m^\otimes-1}
(\bar a,M)$.
\ermn
So
\mr
\item "{$\oplus_3$}"   $r \le m^\otimes -1 < m^\otimes + m \le r+m(1)-k$.
\ermn
Let

$$
B =: \{c \in c \ell^{k^*,m^\otimes-1}(\bar a,M):\text{ for some } d 
\text{ we have }(c,d)\in {\Cal R}\} \cup \bar a.
$$
\mn
So by the above $B=\{c \in c \ell^{k^*,m^\otimes+m}(\bar a,M):(\exists
d)((c,d)\in {\Cal R})\} \cup \bar a$.

Let us check the demands (i) - (iv) of $(*)$ of \scite{3.11}, remember that
we are defining $B^* = (c \ell^k(\bar a b,M) \setminus c \ell^{k^*,m^\otimes+m}
(\bar a,M))\cup B$, that is the submodel of $M$ with this set of elements.
\enddemo
\bn
\ub{Clause (i)}:  $|B|\leq t^*$.

As said above, $|{\Cal R}| \le t$, hence clearly $|B|\leq t + \ell g
(\bar a)\leq t+\ell \le t^*$.
\bn
\ub{Clause (ii)}:  $\bar a \subseteq B \subseteq c \ell^k(B,M)
\subseteq c \ell^{k^*,m^\otimes}(\bar a,M)$.

As by its definition
$B \subseteq c \ell^{k^*, m^\otimes -1}(\bar a, M)$, and $k\leq k^*$
clearly $c \ell^k(B,M) \subseteq c \ell^{k^*, m^\otimes}(\bar a, M)$ and
$B \subseteq c \ell^k(B,M)$ always and $\bar a\subseteq B$ by its
definition.
\bn
\ub{Clause (iii)}:

Clearly 

$$
\align
B &= c \ell^{k^*,m^\otimes+m}(\bar a,M) \cap ((c \ell^k(\bar a b,M) 
\setminus c \ell^{k^*,m^\otimes+m}(\bar a,M)) \cup B) \\
  &= (c \ell^{k^*,m^\otimes+m} (\bar a,M)) \cap B^*.
\endalign
$$
\mn
Now the ``no edges'' holds by the definitions of $B$ and ${\Cal R}$.
\bn
\ub{Clause (iv)}:  $B \le^*_s B^*$.

Clearly $B \subseteq B^*$ by the definition of $B^*$ before the proof
of clause (i).
Toward contradiction assume $\neg(B\leq_s^* B^*)$ then \scite{1.17}(2)
for some $D$, $B <_i D \leq B^*$, and choose such $D$ with minimal
number of elements. Note that
as $B\subseteq c \ell^{k^*, m^\otimes- 1}(\bar a, M)$ and $B^*\cap
c \ell^{k^*, m^\otimes+m}(\bar a, M)=B$, necessarily 
$|D|>k^*$ (and $B <^* D \le B^*$). 
For every $d \in D\setminus B$, as $d \in B^*$ clearly $d \in c
\ell^k(\bar a b,M)$ hence  there is a 
set $C_d\leq M$, $|C_d|\leq k$ such that
$C_d\restriction (\bar ab) \leq_i C_d$, $d \in C_d$; note that
$C_d \subseteq c \ell^k(\bar a b,M)$ by the definition of $c \ell^k$,
hence by the choice of $B^*$ and $m^\otimes$ and $\oplus_1$ we have
$C_d\subseteq B^* \cup c \ell^{k^*,m^\otimes-1}(\bar a,M)$.
Let $C^\prime_d = C_d\cap (B\cup\{b\})$, $C^{\prime\prime}_d 
= C_d\cap B^*$. 
Clearly $C_d \cap (\bar a b) \le C'_d \le C^{\prime\prime}_d \leq C_d$ hence
$C'_d\leq_i C_d$.  Now by clause (iii),
$\nonforkin{C^{\prime\prime}_d}{C^\prime_d\cup(C_d\setminus
C^{\prime\prime}_d)}_{C^\prime_d}^{M}$ hence (by smoothness) we have
$C^\prime_{d_{i}} \le_i C^{\prime\prime}_{d_{i}}$ . Of course, 
$|C^{\prime\prime}_d|\leq |C_d|\leq k$.  For $d \in B$ let 
$C_d = C'_d = C''_d = \{d\}$.

We now choose a set $D_i$ by induction on $i\leq i(*)$, such
that (letting $C^{**}_i = \dbcu_{d\in D_i} C''_d$):
\mr
\item "{$(a)$}"   $D_0 = B \cup\{b\}$
\sn
\item "{$(b)$}"   $j<i \Rightarrow D_j\subseteq D_i \subseteq D$
\sn
\item "{$(c)$}"  $|D_i| \leq k^*_i$
\sn
\item "{$(d)$}"   if $\lambda$ is  an equivalence relation on
$C^{**}_i \setminus B$ and for some $d\in D\setminus D_i$ one of
the clauses below holds \ub{then} there is such $d\in D_{i+1}$ where
{\roster
\itemitem{ $\otimes^1_{\lambda,d}$ }  for some  $x\in C''_d
\setminus C^{\ast\ast}_i$, there are no $y\in 
C^{\prime\prime}_d\cap C^{\ast\ast}_i$, $j^*\in \Bbb N$ and $\langle
y_j: j\leq j^*\rangle$ such that $y_j \in C''_d,y_{j^*} = x,y_0 =y,
\{y_j, y_{j+1}\}$ an edge of $M$, (actually an empty
case, i.e., never occurs see $(*)_{14}$ below)
\sn
\itemitem{ $\otimes^2_{\lambda, d}$ }  there are $x\in
C^{\prime\prime}_d\setminus C^{\ast\ast}_i$, $y\in 
(C^{**}_i\setminus C^{\prime\prime}_d)\cup B$ and
$y' \in C''_d \cap C^{**}_i$ such that $\{x,y\}$ is an edge of $M$ and
$y$ is connected by a path $\langle y_0,\dotsc,y_1 \rangle$ inside
$C''_d$ to $x$ so $x=y_j,y=y_0$ and $[y_i \in C^{**}_i \equiv i=0]$ and
$\neg(y' \lambda y)$
\sn
\itemitem{ $\otimes^3_{\lambda,d}$ }  there is an edge 
$\{x_1,x_2\}$ of $M$ such that we have:
\sn
\itemitem{ ${{}}$ }  $\quad (A) \quad \{x_1,x_2\}
\subseteq C^{\prime\prime}_d$ 
\sn
\itemitem{ ${{}}$ } $\quad (B) \quad \{x_1,x_2\}$ is disjoint to $C^{**}_i$
\sn
\itemitem{ ${{}}$ }  $\quad (C) \quad$ for 
$s\in \{1,2\}$ there is a path $\langle
y_{s,0},\ldots,y_{s,j_s}\rangle$ in \nl

$\qquad \qquad \quad C^{\prime\prime}_d,y_{s, j_s}=x_s,
[y_{s, j}\in C^{\ast\ast}_i\equiv j=0]$ and 
$\neg (y_{1, 0}\lambda y_{2,0})$
\endroster}
\item "{$(e)$}"   if $\lambda$ is an equivalence relation on
$C^{\ast\ast}_i\setminus B$ to which clause (d) does not apply but
there are $d_1$, $d_2\in D$ satisfying one of the following then we can find 
such $d_1$, $d_2\in D_{i+1}$ 
{\roster
\itemitem{ $\otimes^4_{\lambda,d_1,d_2}$ }  for some 
$x_1\in C_{d_1}'' \setminus C^{**}_i$, $x_2\in C_{d_2}^{''} 
\setminus C^{**}_i$ and $y_1\in C^{\prime\prime}_{d_1}\cap
C^{\ast\ast}_i,y_2 \in C''_{d_2} \cap C^{**}_i$ we have: 
for $s = 1,2$ there is a path $\langle
y_{s, 0},\ldots, y_{s,j_s} \rangle$ in $C^{\prime\prime}_{d_s},y_{s,j_s}=x$, 
$y_{s, 0}=y_s$, $[y_{s, j}\in C^{\ast\ast}_i \Leftrightarrow j=0]$
and: $x_1 = x_2 \and \neg(y_1 \lambda y_2)$
\sn
\itemitem{ $\otimes^5_{\lambda,d_1,d_2}$ }  for some $x_1,x_2,y_1,y_2$
as in $\otimes^4_{\lambda, d_1, d_2}$
we have: $\neg(y_1\lambda y_2)$ and $\{x_1, x_2\}$ an edge.
\endroster}
\ermn
So $|D_{i(*)}|\leq k^*/k$ (by the choice of $k^*$, $i(*)$ and clause
(c)), hence $C^{**}_{i(*)} =: \dbcu_{d\in D_{i(*)}} C''_d$ 
has $\leq k^*$ members, $\bar a b \subseteq B \cup \{b\} \subseteq D_0
\subseteq C^{**}_{i(*)} \subseteq c \ell^k(\bar a b,M)$ 
and $C^{**}_{i(*)} \cap
c \ell^{k^*, m^\otimes+m}(\bar a, M)= B\subseteq 
c \ell^{k^*,m^\otimes-1}(\bar a, M)$ hence necessarily 
$B\leq_s C^{**}_{i(*)}$ hence there is $\lambda \in
\Xi(B,C^{**}_{i(*)})$. 
Let $\lambda_i=\lambda\restriction (C^{**}_i \setminus B)$. \nl
Now
\mr
\item "{$\boxdot$}"  $(B,C^{**}_i,\lambda_i) \in \Xi(B,C^{**}_i)$. \nl
[Why?  Easy.]
\endroster
\bn
\ub{Case 1}:  For some $i$, (d) and (e) are vacuous for $\lambda_i$.

Let $\lambda^*_i$ be the set of pairs $(x,y)$ from 
$C^{**}\setminus B$ where $C^{**}= \bigcup\limits_{d\in D}
C^{\prime\prime}_d$ which satisfies $(\alpha)$ or $(\beta)$ where
\mr
\item "{$(\alpha)$}"   $x,y \in C^{**}_{i}\setminus B$ and $x
\lambda_i \beta$
\sn
\item "{$(\beta)$}"   for some $d \in D$ we have 
$x \in C^{**} \setminus C^{**}_i,x \in 
C^{\prime\prime}_d$, $y\in C^{**}_i \cap C^{\prime\prime}_d$ and there
is a sequence $\langle y_j: j\leq j^*\rangle$, $j^*\geq 1$ such that
$y_j \in C''_d,y_0=y,\{y_j,y_{j+1}\}$ an 
edge of $M$ and $[j>0 \Rightarrow y_j\notin C^{**}_i]$.
\ermn
This in general is not an equivalence relation. \nl
Let $C^\otimes = \{x:$ for some $(x_1,x_2) \in \lambda^*_i$ we have $x
\in \{x_1,x_2\}\}$

$$
\align
\lambda^+_i = \{(x_1,x_2):&\text{for some } y_1,y_2 \in D_i \text{ we
have} \\
  &y_1 \lambda y_2,(x_1,y_1) \in \lambda^*_i,(x_2,y_2) \in
\lambda^*_i\}.
\endalign
$$
\mn
Now
\mr
\item "{$(*)_1$}"  $\lambda^+_i$ is a set of pairs from $C^\otimes$
with $\lambda^+_i \restriction D_i = \lambda_i$
\sn
\item "{$(*)_2$}"  $x \in C^\otimes \Rightarrow (x,x) \in \lambda^+_i$
\nl
[why?  read $(\alpha)$ or $(\beta)$ and the choice of $\lambda^+_i$]
\sn
\item "{$(*)_3$}"   for every $x \in C^\otimes$ for some $y \in
C^{**}_i$ we have $x \lambda^*_i y$ \nl
[why?  read the choice of $\lambda^+_i,\lambda^*_i$.]
\sn
\item "{$(*)_4$}"  $\lambda^+_i$ is a symmetric relation on $C^\otimes$
\nl
[why? read the definition of $\lambda^+_i$ recalling $\lambda$ is
symmetric.]
\sn
\item "{$(*)_5$}"  $\lambda^+_i$ is transitive \nl
[why?  looking at the choice of $\lambda^*_i$ this is reduced to the
case excluded in $(*)_6$ below]
\sn
\item "{$(*)_6$}"  if $(x,y_1),(x,y_2) \in \lambda^*_i,\{y_1,y_2\}
\subseteq D_i,x \notin D_i$, then $y_1 \lambda y_2$. \nl
[Why?  Because clause (e) in the choice of $D_{i+1}$ is vacuous. 
More fully, otherwise possibility $\otimes^4_{\lambda,d_1,d_2}$ holds for
$\lambda_i$.]
\sn
\item "{$(*)_7$}"   for every $x \in C^{**}\setminus C^{**}_i$, 
clause $(\beta)$ apply to $x \in C^\otimes$, i.e., $C^\otimes_{**} =
C^{**}$ \nl
[why?   as $x \in C^{**}$ there is $d \in D$ such that $x \in C''_d$,
hence by $\otimes^1_{\lambda,d}$ of clause (d) of the choice of 
$D_{i+1}$ holds for $x$ hence is not vacuous contradicting the
assumption on $i$ in the present case]
\sn
\item "{$(*)_8$}"  $\lambda^+_i$ is an equivalence relation on $C^{**}
\backslash B$ \nl
[why?  its domain is $C^{**} \backslash B$ by $(*)_7$, it is an
equivalence relation on its domain by $(*)_1 + (*)_2 + (*)_4 + (*)_5$.]  
\nl
Also
\sn
\item "{$(*)_9$}"   $\lambda^+_i \restriction C^{**}_i = \lambda_i$
\nl
[why? by the choice of $\lambda^+_i$ that is by $(*)_1$]
\sn
\item "{$(*)_{10}$}"   every $\lambda^+_i$-equivalence class is 
represented in $C^{**}_i$ \nl
[why?  by the choice of $\lambda^+_i$ and $\lambda^*_i$]
\sn
\item "{$(*)_{11}$}"   if $x_1$, $x_2\in C^{**}\setminus B$ and $\neg(x_1
\lambda^+_i x_2)$ but $\{x_1, x_2\}$ is an edge \ub{then} 
$\{x_1,x_2\} \subseteq C^{**}_i$. \nl
\ermn
[Why $(*)_{11}$ holds?  
Assume $\{x_1,x_2\}$ is a counterexample, so $\{x_1,x_2\}
\nsubseteq C^{**}_i$, so without loss of generality 
$x_1 \notin C^{**}_i$.  Now for $\ell = 1,2$ if
$x_\ell \notin C^{**}_i$ then we can choose $d_\ell \in D_i$ and $y_\ell \in
C''_{d_\ell} \cap C^{**}_i$ such that $d$ witnesses that
$(x_\ell,y_\ell) \in \lambda^*_i$ that is, as in clause $(\beta)$
there is a path $\langle y_{\ell,0},\dotsc,y_{\ell,j_\ell} \rangle$
such that $y_{\ell,0} = y_\ell,y_{\ell,0} = x_\ell$ and $(j>0
\Rightarrow y_{\ell,j} \notin C^{**}_i)$.
\nl
We separate to cases:
\mr
\item "{$(A)$}"   $x_1,x_2 \notin C^{**}_i,d_1 = d_2$.  \nl
This case can't happen as $\otimes^3_{\lambda,d_1}$ of clause (d) is vacuous
\sn
\item "{$(B)$}"   $x_1,x_2 \notin C^{**}_i,d_1 \ne d_2$. \nl
In this case by the vacuousness of $\otimes^5_{\lambda_i,d_1,d_2}$ of 
clause (e) we get contradiction 
\sn
\item "{$(C)$}"    $x_1 \in C''_d$ and $x_2 \in C^{**}_i$ \nl
by the vacuousness of $\otimes^2_{\lambda_i,d_1}$ of clause (d).
\ermn
Together we have proved $(*)_{11}$.]

As $\lambda_i \in \Xi(B, C^{**}_i)$ by
$(*)_8 + (*)_9 + (*)_{10} + (*)_{11}$ and $\boxdot$, easily 
$\lambda^+_i\in \Xi(B, C^{**})$,
hence (see \scite{1.16}) $B <^*_s C^{**}$, so as $B\subseteq
D\subseteq C^{**}$ we have $B <^*_s D$, the desired contradiction.
\bn
\ub{Case 2}:   For every $i<i(*)$, at least one of the clauses (d), (e)
is non vacuous for $\lambda_i$.

Let $\bold w_i = \bold w_{\lambda_i}(B,C^{**}_i)$. 
For each $i$ let $\langle d_{i,j}:j< j_i\rangle$ list 
$D_{i+1}\setminus D_i$, such that: if clause (d) applies to 
$\lambda_i$ then $d_{i, 0}$ form a witness and if clause (e) 
applies to $\lambda_i$ then $d_{i, 0}$, $d_{i, 1}$ form
a witness. For $j\leq j_i$ let $C^{**}_{i, j}= C^{**}_i \cup
\dbcu_{s<j} C^{\prime\prime}_{d_{i, s}}$, so $C^{**}_{i, 0} =
C^{**}_i$, $C^{**}_{i, j_i}= C^{**}_{i+1}$.
Let $\bold w_{i,j} = \bold w_{\lambda_i}(B,C^{**}_{i,j})$.

So it suffice to prove:
\mr
\item "{$(A)$}"   $\bold w_{i,j}\geq \bold w_{i, j+1}$
\sn
\item "{$(B)$}"   $\bold w_{i,0} - \varepsilon^*\geq \bold w_{i,1}$ 
or $\bold w_{i,1}-\varepsilon^* \geq \bold w_{i,2}$.
\ermn
Let $i< i(*)$, $j< j_i$.

Clearly $C^{**}_{i, j+1}\setminus C^{**}_{i, j} \subseteq
C^{\prime\prime}_{d_{i, j}}\subseteq C^{**}_{i, j+1}$, let

$$
A_{i,j} = \{x \in C^{\prime\prime}_{d_{i,j}}:x \in B \text{ or } x/ \lambda 
\text{ is not disjoint to } C^{**}_{i, j}\}.
$$
\mn
Clearly $A_{i,j} \backslash B$ is 
$(\lambda\restriction C''_{d_{i,j}})$-closed 
hence $A_{i, j} \leq^* C^{\prime\prime}_{d_{i, j}}$,
$C^{\prime\prime}_{d_{i, j}}\setminus A_{i, j}$ is disjoint to
$C^{**}_{i,j}$ and
$C^\prime_{d_{i,j}}= C_{d_{i,j}} \cap (B \cup \{b\})
\subseteq C^{**}_{i, j}$, and $C^\prime_{d_{i, j}} 
\subseteq C^{\prime\prime}_{d_{i, j}}$ hence 
$C^\prime_{d_{i, j}}\subseteq A_{i,j},A_{i, j} \le^* C''_{d_{i, j}}$, 
but $C'_{d_{i,j}} \leq_i
C^{\prime\prime}_{d_{i,j}}$ so $A_{i, j}\leq^*_i C''_{d_{i,j}}$ 
(the $\leq^*$ in this sentence serves \S7 where we say, 
``repeat the proof of \scite{3.14}").

Clearly
\mr
\item "{$(*)_{12}$}"   $\bold w_{i,j+1} = \bold w_{i,j} + 
\bold w_\lambda(A_{i, j},C''_{d_{i, j}})- \alpha \bold e^1_{i,j} - 
\alpha \bold e^2_{i,j}$ where

$$
\align
\bold e^1_{i,j} = |\{\{x, y\}:&\{x,y\} \text{ an edge of } M, 
\{x,y\} \subseteq A_{i,j}, \\
  &\neg(x\lambda y) \text{ but } \{x,y\} \nsubseteq C^{**}_{i, j}\}|
\endalign
$$

$$
\align
\bold e^2_{i,j} = |\{\{x,y\}:&\{x,y\} \text{ an edge of } M,x \in
C''_{d_{i,j}} \setminus C^{**}_{i,j}, \\
  & y \in C^{**}_{i,j} \setminus C''_{d_{i,j}}
\text{ but } \neg(x\lambda y)\}|.
\endalign
$$
\ermn
Note
\mr
\item "{$(*)_{13}$}"   $\bold w_\lambda(A_{i,j},
C''_{d_{i,j}})$ can be zero if $A_{i, j}= C''_{d_{i,j}}$ 
and is $\leq -\varepsilon^*$ otherwise.\nl
[Why?   As $A_{i,j}\le^*_i C^{\prime\prime}_{d_{i, j}}$.]
\sn
\item "{$(*)_{14}$}"   in clause (d), $\otimes^1_{\lambda,d}$ never
occurs \nl
[Why?  If $x \in C''_d$ is as there, let $Y = \{y \in C''_d:y,x$ are
connected in $M \restriction C''_d\}$.  So $x \in Y \subseteq C''_d,Y
\cap C^{**}_i = \emptyset$ hence $C'_d = C''_d \cap (B \cup \{b\}) =
C''_d \cap C^{**}_0 \subseteq C^{**}_i$.  Hence $(C''_d \backslash Y)
<^*_i C''_d$, but the equivalence relation $\{(y',y''):y',y'' \in Y\}$
exemplify that this fail.]
\endroster
\bigskip

\demo{Proof of (A)}  

Easy by $(*)_{12}$, because $\bold w_\lambda(A_{i,j},
C^{\prime\prime}_{d_{i,j}})\leq 0$ holds by $(*)_{13}$, 
$-\alpha \bold e^1_{i, j}\leq 0$, and $-\alpha \bold e^2_{i,j} \le 0$ as 
$\bold e^1_{i,j}$, $\bold e^2_{i,j}$ are natural numbers.
\enddemo
\bigskip

\demo{Proof of (B)}

It suffice to prove that 
$\bold w_{i,0} \neq \bold w_{i,1}$ or $\bold w_{i, 1}\neq
\bold w_{i,2}$ (as inequality implies the right order (by clause (A)) and
the difference is $\geq \varepsilon^*$ by definition of $\varepsilon^*$ (if
$\bold w_\lambda(A_{i,1},C^{\prime\prime}_{d_{i, j}})\neq 0$) and $\geq
\alpha$ (if $\bold e^1_{i, j}\neq 0$ or $\bold e^2_{i, j}\neq 0$)). But if
$\bold w_{i, 0}= \bold w_{i, 1}$ recalling $(*)_{14}$ easily clause 
(d) does not apply to
$\lambda_i$, and if $\bold w_{i, 0} = \bold w_{i,1} 
= \bold w_{i, 2}$ also clause (e) does not apply.

So (A),(B) holds so we are done proving case 2 hence the claim.
\hfill$\square_{\scite{3.11}}$\margincite{3.11}
\enddemo
\bigskip

\remark{\stag{6.9A} Remark}  
\mr
\item "{$(a)$}"   We could use smaller $k^*$ by building a 
tree $\langle (D_t,D^+_t,C_t,\lambda_t): t\in T\rangle$, $T$ a 
finite tree with a
root $\Lambda$, $D_\Lambda =\emptyset$, $D^+_\Lambda = B\cup\{b\}$,
for each $t$ we have $\lambda_t$ an equivalence relation on $C_t
\backslash B$ and $C_t = \cup\{C''_d:d \in D_t\} \cup B,
s \in \text{ suc}_T(t)\Rightarrow D^+_t = D_s$ and $D^+_t \backslash
D_t$ is $\{d\}$ or $\{d_1,d_2\}$ which witness clause (d)
or clause (e) for $(D_t,\lambda_t)$ when $t \ne \Lambda$ and
$$
\align
\{(D_s,\lambda_s):&s\in \text{ suc}_T(t)\} = \{(D^+_t,\lambda):
\lambda \restriction D_t = \lambda_t,\lambda \\
  &\text{ an equivalence relation on } D^+_t\setminus B\}.
\endalign
$$
\item "{$(b)$}"   We can make 
the argument separated that is prove as a separate claim that is for any $k$
and $\ell$ there is $k^*$ such that: if $A$, $B\subseteq
M \in {\Cal K}$, $|B|$, $|A|\leq \ell$, $B\subseteq B^*$, 
$c \ell^k(A, M)\setminus c \ell^k(B, M)\subseteq B^*\setminus
B \subseteq c \ell^k(A, M)$ and $(\forall C)(B\subseteq C\subseteq B^*
\wedge |C|\leq k^* \Rightarrow B<_s C)$ \ub{then} $B <_s B^*$.\nl
This is a kind of compactness.
\endroster
\endremark
\bigskip

\demo{\stag{3.13} Conclusion}  Requirements (A) of [I,{2.13(1)] 
and even (B) + (C) of [I,2.13(3)] hold.
\enddemo
\bigskip

\demo{Proof}   Requirement (B) of [I,2.13(3)] holds by
\scite{3.8}.   Requirement (A) of [I,2.13(2)] holds by
\scite{3.11} (and the previous sentence).  \hfill$\square_{\scite{3.13}}$\margincite{3.13}
\enddemo
\bigskip

\demo{\stag{3.14} Conclusion}  
\mr
\item "{$(a)$}"  ${\frak K}$ is smooth and 
transitive and local and transparent
\sn
\item "{$(b)$}"   ${\frak K}$ is simply nice (hence simply almost
nice)
\sn
\item "{$(c)$}"   ${\frak K}$ satisfies the 0-1 law.
\endroster
\enddemo
\bigskip

\demo{Proof}   1) By \scite{3.4}. \nl
2) By \scite{3.13} we know that ${\frak K}$ is simply nice. \nl
3) By \scite{1.2} we know that for each $k$, for every random enough
${\Cal M}_n,c \ell^k(\emptyset,{\Cal M}_n)$ is empty. Hence by 
[I,2.19(1)] we get the desired conclusion.
 \hfill$\square_{\scite{3.14}}$\margincite{3.14}
\enddemo
\newpage

\head {\S7 Random graphs with the successor function} \endhead  \resetall \sectno=7
 \spuriousreset
\bigskip

Recall that ${\Cal M}^1_n$ is defined like ${\Cal M}^0_n$ (called
${\Cal M}_n$ earlier) but we expand it by interpreting the two-place
predicate $S$ by the relation $S_n = \{(\ell,\ell +1):\ell \in
\{1,\dotsc,n-1\}$.  Trivally, ${\Cal M}_n^1$ fails the 0-1 law; (e.g.,
``the first two elements (by $S$) are connected by an edge" has 
probability $p_1$.  So now
${\Cal K} = \{M:M = (|M|,R,S),(|M|,R)$ a graph, $S$ an anti-symmetric
irreflexive two-place relation satisying $x_1 S x_2 \wedge y_1 S y_2
\Rightarrow (x_1 = y_2 \equiv x_2 = y_2)$.  Still we prove the convergence
law.  We may remedy the failure by replacing the relation $S_n$ by
$S'_n =: S_n \cup \{(n,1)\}$, call this random model ${\Cal M}^{0,5}_n$.
\bigskip

\proclaim{\stag{4.1} Theorem}  Let ${\frak K}$ be the 0-1 context as follows:
\block
$\alpha\in (0,1)_{\Bbb R}$ be irrational, \nl
$p_\ell=\frac{1}{\ell^\alpha}$ for $\ell>1$, \nl
$p_1=\frac{1}{2^\alpha}$ (we can omit this) and ${\Cal M}_n =
{\Cal M}^1_n$
\endblock
\ub{Then} ${\frak K}$ satisfies the convergence law.
\endproclaim
\bigskip

\demo{Proof}   This is proved later.
\enddemo
\bigskip

\remark{\stag{4.2} Remark}  If the probability of ``$\{i, j\}$ is an edge'' is 
$\frac{1}{n^{\alpha}}+\frac{1}{|i-j+1|^\alpha}$, then
the same conclusion holds, see [I,\S3].
\endremark
\bigskip

\proclaim{\stag{4.3} Theorem}  Let 
$\alpha\in (0,1)_{\Bbb R}$ be irrational, $p_i=1/i^\alpha$ for $i>1$,
$p_1=\frac{1}{2^\alpha}$. Let ${\Cal M}_n = {\Cal M}^{0.5}_n$ 
be $([n],R,S'_n)$ where $R$ is a graph on $[n]$ chosen randomly: \nl
$\{i,j\}\in R$ has probability $p_k$ where $k\in\{1,\ldots,n-1\}$ is minimal
such that $n$ divides $|i-j-k|$, the choices are independent for 
distinct edges and, lastly, the predicate $S$ is interpreted by the relation

$$
S'_n =: \{(x,y):x\in [n], y\in [n] \text{ and } x+1=y \,{\text{\rm
mod\/}} \, n\}.
$$
\mn
So a 0-1 context ${\frak K}$ is defined.\nl
\ub{Then}
\mr
\item "{$(a)$}"   ${\frak K}$ is smooth, transitive, local and
transparent,
\sn
\item "{$(b)$}"   ${\frak K}$ is simply almost nice 
\sn
\item "{$(c)$}"   for every $k$ for every random enough
${\Cal M}_n$ we have: $c \ell^k(\emptyset,{\Cal M}_n) =
\emptyset$ hence $\langle {\Cal M}_n \restriction c
\ell^k(\emptyset,{\Cal M}_n):n<\omega\rangle$ satisfies the 0-1 law
\sn
\item "{$(d)$}"   Consequently ${\frak K}$ satisfies the 0-1 law.
\endroster
\endproclaim
\bigskip

\demo{Proof}   We repeat the proof for ${\Cal M}_n = {\Cal M}^0_n$ in
\S4, \S5, \S6 with some changes. Let ${\Cal K}$ be the
class of finite models $(X,R,S)$ such that: $R$ is a symmetric irreflexive
two--place relation on $X$, $S$ is an irreflexive antisymmetric two--place
relation satisfying

$$
(\forall x_0,y_0)(\forall x_1,y_1)(S(x_0,y_0) \wedge S(x_1,y_1)
\Rightarrow (x_0=x_1 \equiv y_0=y_1)).
$$
\mn
We use below unary predicates $P_f$, $P_\ell$ for the proof of
\scite{4.1} later so at present 
$P_f(x)\equiv P_\ell(x)\equiv \text{ false}$.

Let ${\Cal K}'$ be the set of $(X,R,S) \in {\Cal K}$ 
with no $S$--cycle.  ${\Cal K}_n$ is the set of possible values 
of ${\Cal M}^{0.5}_n$ so every member $M$ of ${\Cal K}_n$ is an $S^M$-cycle
but ${\Cal K}_\infty = {\Cal K}'$, see below.

Easily ${\Cal K}_n,{\Cal K}' \subseteq {\Cal K}$ and
 by the proof of \scite{1.2}:
\mr
\item "{$(*)$}"   for each $A \in {\Cal K}'$ 
for some $\bold c \in \Bbb R^+$ for every random enough ${\Cal M}_n$
there are $\geq \bold c \times n$ pairwise disjoint embeddings of $A$
into ${\Cal M}_n$.
\ermn
[Why?  Let $A=\{a_\ell: \ell< k\}$, and without loss of generality
$a_\ell S a_m \rightarrow \ell+1=m$ and let

$$
u = \{\ell< k:\ell=0 \text{ or } \ell>0 \and \neg a_{\ell-1} S a_\ell\},
$$
\mn
and for $r < \lfoot n/2k \rfoot$ let ${\Cal E}^n_r$ be the event
$a_\ell \mapsto 2rk+\ell + |\{m:m \le \ell \and  m \in u\}| 
\text{ is an embedding of } A \text{ into }{\Cal M}_n$.
\mn
Note that for $S$ this is always an embedding and those embeddings
have pairwise disjoint ranges.]
\nl
Trivially
\mr
\item "{$(*)$}"  if $A \in {\Cal K} \backslash {\Cal K}'$ 
and $M \in {\Cal K}_n$ \ub{then} $A$ cannot be embedded into $M$.
\ermn
Hence ${\Cal K}' = {\Cal K}_\infty \subseteq {\Cal K}$, and 
hence (the parallel of) \scite{1.3}, \scite{1.4} hold, but in 
\scite{1.4} we replace simple by simple$^*$ which is defined in
[I,2.21]. \nl
And so replace [I,2.17] by [I,2.23]. 
If $A\subseteq B$ we let

$$
scl^0(A,B) = A
$$

$$
\align
scl^1(A,B) = \{x\in B:&x \in A \text{ or }(\exists y\in A)(S^B(x,y)
\vee S^B(y,x)) \\
  &\text{ or (for the proof of \scite{4.1} below) } P_f(x)\vee P_\ell(x)\}
\endalign
$$

$$
scl^{k+1}(A,B) = scl^1(scl^k(A,B))
$$

$$
rcl^0(A,B) = A
$$

$$
r c \ell^1(A,B) = \{x \in B:x \in A \text{ or } (\exists y \in
A)(S^B(x,y) \vee S^B(y,x)\}
$$

$$
r c \ell^{k+1}(A,B) = r c \ell^1(r c \ell^k(A,B)).
$$
\mn
So for our present context $({\Cal M}^{0,5}_n)$ we have $sc
\ell^k(A,B) = rc \ell^k(A,B)$ but not so for the proof of \scite{4.1}.
Lastly, let

$$
scl(A,B)= \dbcu_{k \in \Bbb N} scl^k(A,B)
$$

$$
rcl(A,B)= \dbcu_{k \in \Bbb N} rcl^k(A,B).
$$
\mn
Let $A\leq B$ means $A \subseteq B \in {\Cal K}'$ (submodel) so $M \in
{\Cal K}_n \Rightarrow \neg(M \le M)$.

Let $A\leq^* B$ means $A\leq B \and scl(A,B)=A$.  Clearly it is a 
partial order on ${\Cal K}_\infty$.

Let $E_A$ be the finest equivalence relation on $A$ such that
$S^A(x,y) \Rightarrow xEy$.

Lastly, let $A\leq^{**} B$ means: $A\leq B$ and $E_A=E_B\restriction
A$ and $x\in scl (A,B)\setminus A  \and y \in scl(A,B)
\Rightarrow \{x, y\}$ not an edge. Clearly it is a partial order on
${\Cal K}'$.

Now we define ${\Cal T}$ (instead of definition \scite{4.6}(1)):

$$
\align
{\Cal T} = \{(A,B,\lambda):&A<^* B \in {\Cal K}_\infty,\lambda
\text{ an equivalence relation on } B \setminus A \\
  &\text{ which } E_B \text{ refines, which means that } 
x \in  B \setminus A \and  y\in B \\
  &\and (S(x,y) \vee S(y,x)) \Rightarrow  y\in (x/\lambda)\}.
\endalign
$$
\mn
Note: generally in this version cases of $S$ counts as edges (even more
so).  \ub{But} in the definition of $\bold e(A,B,\lambda)$ below they do not
count as any $x S^B y \Rightarrow \{x, y\}\subseteq A\ \vee\ (x\lambda y)$
and $x\lambda y \rightarrow \{x, y\}\subseteq B\setminus A$.

We define $\lambda$-closed as in definition \scite{4.6}(2) (but $A\leq^*
B$ has already been defined).
So $(A, B, \lambda)\in {\Cal T} \Rightarrow scl(A,B)=A$ and $A\leq^* B
\Rightarrow A \leq^{**}B$. Now for applying \S2 we are interested in
$\leq$, but in applying \S4 - \S6 only in $<^*$, so these
sections are written trying to have this in mind.

Now $\bold v,\bold e,\bold w$ are defined for 
$A\leq^* B$ as in \scite{1.7}, \scite{1.8} and $A<^*_x B$ are 
defined as in \scite{1.11}, and the parallel to \scite{1.9} and 
\scite{1.14} - \scite{1.17} hold.

In the proof of \scite{2.4} we should be careful to preserve $\pm S$ in the
relevant cases, in particular in the definition of 
${\Cal G}^{\varepsilon,k}_{A,B}(f,{\Cal M}_n)$ in condition (1) 
(all in Definition \scite{2.2}) and of ${\Cal G}^{\varepsilon,1}_{A,
B}(f,[n])$ in stage C of the proof of \scite{2.4},
choose $m^*$ large enough and in the types $(\text{tp}^0,\text{tp}^1$
mentioned there) of $g$ 
fix $g(b)+m^*{\Bbb Z}$ for each $b\in B$. So there is
$m^{\otimes}\in\{0,\ldots,m^*-1\}$ such that for no $b$ is 
$g(b)=m^\otimes \text{ mod } m^*$, and we move blocks 

$$
\{m^*\cdot i+m^\otimes+1,m^*\cdot i+m^\otimes+2,\ldots,m^*\cdot(i+1)+
m^\otimes\}
$$
\mn
together. And concerning $f_0$ we just ask

$$
\text{rang}(f_0) \subseteq \{i:1\leq i\leq m^*|A|
\text{ or } n- m^*|A|\leq i\leq n\}
$$
\mn
so there are some possibilities for $f_0$ (but there is a bound on the
number not depending on $n$).
The other cases are even less serious.

In \scite{2.9}, on each $\lambda$-equivalence class we should preserve
$S$ and we change $g_{\bar j}$ as we have for \scite{4.1} (i.e., in
proving $(*)$ above).
Note that in \scite{2.9} - \scite{2.12} we have $(A, B, \lambda)\in
{\Cal T} \Rightarrow A\leq^* B$ and $\lambda$-closed implies $S$-closed in
relevant places.  Now \scite{3.4} 
has to be rephrased by \scite{4.4}(1) below. Also we have to 
change somewhat \scite{3.6} because \scite{3.6} says we
have the nice case, whereas here we only have the almost nice case, so
we replace it by \scite{4.4}(2) below.
\enddemo
\bigskip

\proclaim{\stag{4.4} Claim}  1) For $A$, $B\in {\Cal K}_\infty$ we have
\mr
\item "{$(a)$}"  $A <_s B \Leftrightarrow A <^*_s B$,
\sn
\item "{$(b)$}"   $A <_i B \Leftrightarrow A \leq B \and scl(A,B) 
\leq^*_i B$.
\ermn
2) If $A<^*_s B$ and $k \in \Bbb N$ \ub{then} for every 
random enough ${\Cal M}_n$ and $f:A \hra {\Cal M}_n$ we can find 
$g:B \hra {\Cal M}_n$ extending $f$ and such that
\mr
\widestnumber\item{$(iii)$}
\item "{$(i)$}"   ${\text{\rm Rang\/}}(g) \cap c \ell^{k}
({\text{\rm Rang\/}}(f),{\Cal M}_n) = \,{\text{\rm Rang\/}}(f)$,
\sn
\item "{$(ii)$}"   $rc \ell^{k+1}(f(A),{\Cal M}_n),rc \ell^{k+1}(g(B
\backslash A),{\Cal M}_n)$ are disjoint
\sn
\item "{$(iii)$}"  letting $B^+ = f(A) \cup rc \ell^k(g(B \backslash
A),{\Cal M}_n)$ we have $g(B) \leq^{**} B^+$
\sn
\item "{$(iv)$}"   $\nonforkin{B^+}{c \ell^k (f(A),
{\Cal M}_n)}_{f(A)}^{{\Cal M}_n}$
\sn
\item "{$(v)$}"  if $|C|\leq k$, $C \subseteq {\Cal M}_n$ and $C\restriction
(C\cap g(B))<_i C$ then $C \setminus c \ell^k(f(A),{\Cal M}_n)
\subseteq \text{ rcl}(C\cap g(B) \backslash A,C)$.
\endroster
\endproclaim
\bigskip

\demo{Proof}   Like \scite{3.6}. \nl
1a) By the parallel of \scite{3.2} and the definitions of $<_s,<^*_s$.
\nl
1b)  By the parallel of \scite{3.3} of the definitions of $<_i,<^*_i$.
\nl
2) (i), (iv) and (v) follows from \scite{3.6}.  (iii) follows from
definition of $B^+$ and \scite{3.6}. \nl
(Could use $k,k'$ as in [I,2.12,2.13]!) 
\hfill$\square_{\scite{4.4}}$\margincite{4.4}
\enddemo
\bn
Similarly we replace \scite{3.8} by \scite{4.5} below as first
approximation (see \scite{4.7} later).
\proclaim{\stag{4.5} Claim}  Assume
\mr
\item "{$(a)$}"   $A <_s B$
\sn
\item "{$(b)$}"   $B_0 \subseteq B$, and (for \scite{4.1}) for 
simplicity if ${\Cal M}_n \models \exists x P_f(x) \vee 
(\exists y) P_\ell(y)$ then $(\exists x,y \in B_0)(P_f(x) \wedge 
P_\ell(y))$
\sn
\item "{$(c)$}"   if $A'\subseteq A,|A'| \le k^2$ and 
$x\in c \ell^k(B_0\cup A',B),x \notin A$
and the $S$-component of $x$ in $B$ is $\{x_0, \ldots, x_j\}$ with
$S^B(x_\ell, x_{\ell+1})$ and $x=x_{i(*)}$ \ub{then} 
$(i(*)>k)$ and $(i(*)< j-k)$.
\ermn
\ub{Then} $(B, A, B_0, k)$ is simply good (see Definition [I,2.12]).
\endproclaim
\bigskip

\remark{\stag{7.4newA} Remark}  What is the reason for assumption (c)?

Assume $f:A \hra {\Cal M}_n$ and $A <^*_s B$ and $B_0 \subseteq B$ and
we would like to find $g:B \hra {\Cal M}_n$ extending $A$ such that $c
\ell^k(g(B_0),{\Cal M}_n) \subseteq g(B) \cup c \ell^k(f(A),{\Cal
M}_n)$.  A problematic point which does not arise earlier is: suppose $b \in
B_0 \backslash A$ and $sc \ell^k(\{g(b)\},{\Cal M}_n)$ include an
element $c \in {\Cal M}_n \backslash g(B)$.  Also if $b \in g(B \cap
c \ell^{k(1)}(g(\bar a b,{\Cal M}_n)$ and $c \in r c
\ell^{k(2)}(\{b_1\},{\Cal M}_n)$ and $k(1) + k(2) \le k$ we have a
similar problem.  The aim of clause (c) is to exclude those cases.

The need of something like $A'$ arises in the possibility of $C$; can
be avoid if we use smoothness.
\endremark
\bigskip

\demo{Proof}   Let ${\Cal M}_n$ be random enough and $f:A \hra {\Cal M}_n$. By
\scite{4.4}(2) applied to large enough $k^*$, we can 
find an extension $g: B\hra {\Cal M}_n$ of $f$, and
$B^+$ as there.
Now from the demands in Definition [I,2.12(1)]: clause (i) there holds
by \scite{4.4}(2)(ii) + (iii).  \nl
Clause (ii) there holds by \scite{4.4}(2)(iv) and
obvious monotonicity property of $\nonforkin{}{}_{}^{}$.
So our problem is to show that clause (iii) of Definition 
[I,2.12(1)] holds, i.e., 
\mr
\item "{$(*)$}"   $c \ell^k(g(B_0),{\Cal M}_n) \subseteq g(B) \cup 
c \ell^k(f(A),{\Cal M}_n)$.
\ermn
Toward contradiction suppose $d \in c \ell^k(g(B_0),{\Cal M}_n) \setminus
(g(B) \cup c \ell^k(f(A),{\Cal M}_n))$, so for some $C\subseteq {\Cal M}_n$
we have $|C|\leq k$, $d \in C$ and let $C_1=: C\restriction (g(B_0)\cap
C)$ we have $C_1 <_i C$.

So also $(g(B)\cap C)\leq_i C$, hence by \scite{4.4}(2)(v) we have 
$C \subseteq \text{ rcl}(g(B) \cap C \backslash f(A,C) 
\cup c \ell^k(f(A),{\Cal M}_n)$ and $x = x_j$ which means that
$C \subseteq B^+\cup c \ell^k(f(A),{\Cal M}_n)$.  Hence by clause (ii)
of \scite{4.4}(2) clearly let $C_2 \subseteq |C_2|\leq k^2$ be such that

$$
C_2 \subseteq f(A)
$$

$$
C \cap c \ell^k(f(A),{\Cal M}_n)\subseteq c \ell^k(C_2,{\Cal M}_n)
$$

$$
g(A) \cap C \subseteq C_2.
$$
\mn
Let $B'_0 = f(B_0) \cap C$.
\nl
Let $C_0 = (B'_0) \cup C_2,C^+ = C \cup C_2$ so as $B'_0 \le_i C$ also
$C_0 \le_i C^+$.  By clauses (iv) + (v) of \scite{4.4}(2) we get
\mr
\item "{$\circledast$}"  $C_3 =: C^+ \cap (B \cup c \ell^k(f(A),{\Cal
M}_n)) \le^{**} C^+$. 
\ermn
Clearly $C_0 \le C_3 \le C^+$ and as by $\circledast$ and claim
\scite{4.7} below we get $C_0 \le_i C_3$, also $|C_3| \le |C^+| \le
|C| + |C_2| = k+k^2$.  As $C_3 \le^{**} C^+$ and $C^+ \backslash C_3
\subseteq B^+ \backslash B$, clearly
\mr
\item "{$\circledast$}"  if $c \in C^+ \backslash C_3$ and $d \in C_3$
then $\{c,d\}$ is not an edge of ${\Cal M}_n$.
\ermn
Now this implies that $sc \ell(C_3,C^+) \le_s C^+$ (using the choice
of $B^+$) but as $C_0 \le_i C^+,C_0 \le C_3 \le C^+$ clearly $C_3
\le_i C^+$.  So necessarily $sc \ell(C_3,C^+) = C^+$.

But $d \in C^+ \backslash C_3$ we have $C_3 \ne C^+$ hence for some
$d' \in C_3$ and $d'' \in C^+$ are $S$-nbs.  But this contradicts
clause (c).  \hfill$\square_{\scite{4.5}}$\margincite{4.5}
\enddemo
\bigskip


\demo{\stag{4.6} Fact}  If $A \le^* B\leq^{**} C$ and
$A\leq^*_i C$ \ub{then} $A\leq^*_i B$.
\enddemo
\bigskip

\demo{Proof}  Remember that, e.g., $E_C$ is the closure of $S^C$ to 
an equivalence relation, and similarly $E_B$.
If $\neg(A\leq^*_i B)$ then we can find $A'$ such that $A\leq^*
A' <^*_s B$, hence there is $\lambda \in \Xi(A', B)$ (see
\scite{1.8}(2)).  Now clearly scl$(A',B)=A'$, and let

$$
\align
A^+ = A'\cup \{x\in C:&(x/E_C) \cap A' \ne \emptyset \text{
equivalently} \\
  &(x/E_C) \cap A' \ne \emptyset \and  (x/E_C)\cap B\subseteq A'\},
\endalign
$$
\mn
so $A^+\cap B=A'\cap B=A'$.  We
define $\lambda'$, an equivalence relation on $C\setminus A^+$, by:

$$
x \lambda' y \Leftrightarrow (\exists x'\in x/E_C)(\exists y'\in y /
E_C)(x'\lambda y').
$$
\mn
As $E_C \restriction B = E_B$ (by the definition of $<^{**}$),
clearly $\lambda'\restriction (B\setminus A^+)=\lambda$,
and every equivalence class of $\lambda'$ has member in $B \setminus
A'$ (as $B \leq^{**} C$), so $\bold v(A^+,C,\lambda') = 
\bold v(A',B,\lambda)$.   Also every edge of $C$ not included in $A^+$
is an edge of $B$ not included in $A'$ hence 
$\bold e(A^+,C,\lambda') = \bold e(A',B,\lambda)$, together 
$\bold w(A^+, C, \lambda')= \bold w(A', B, \lambda)$.  Moreover, if 
$D'\subseteq C\setminus A^+$ is $\lambda'$-closed then
$D = D' \cap B$ is $\lambda$-closed and
$\bold w_\lambda (A^+, A^+ \cup D') = \bold w_\lambda (A',A'\cup D)>0$. 
So $(A^+,C,\lambda') \in \Xi(A^+,C)$ so 
$A^+ <_s C$ which contradicts $A \le^*_i C$.  \hfill$\square_{\scite{4.6}}$\margincite{4.6}
\enddemo
\bn
\centerline {$\ast\qquad\qquad\ast\qquad\qquad\ast$}
\bn
The rest (generalizing \scite{3.9}, \scite{3.10}, \scite{3.11}, 
\scite{3.13}) is similar to the original except the following.
In \scite{3.9} - \scite{3.11} instead ``$\{x, y\}$ an edge" we should say
``$\{x, y\}$ an edge or $xSy$ or $ySx$".
In \scite{3.9}, \scite{3.10}, \scite{3.11} $M\in {\Cal K}'$ 
(rather than $M \in {\Cal K}$).
In the proof of \scite{3.9}, in $\oplus_1$ we add
\mr
\item "{$(h)$}"   the truth value of $d_{i,s_i},Sd_{i,s_2}$ do not 
depend on $i$;
\sn
\item "{$(i)$}"   if $d_{i_1,s_1},Sd_{i_2,s_2}$ then $d_{i_1,s_1},
Sd_{i_1,s_2}$ (so necessarily $i_1=i_2 \vee \{s_1,s_2\}\subseteq
S_1$).
\ermn
Also in the proof of $\oplus_4$ (inside the proof of \scite{3.9}), $d_i$
is not in $\{d_{i, s}: s\in S_2\}$ which is closed under $S^{C_i}$,
hence (see clause (i)), $C_i\cap D_0 <^* C_i$ but, as there, $C_i\cap
D_0 \leq_i C_i$. Now check the inequality.

In the proof of \scite{3.10} note that after (h), really $B_1 <^* B_2
\leq^* C_{i(*)}$ by clauses (f)+(g) there, as
$b'S^M b'' \and b' \in c \ell^{k^*,\ell},\ell(\bar a,M) \Rightarrow b'' \in c 
\ell^{k^*,\ell+1}(\bar a,M)$.
\mn
In claim \scite{3.11}, easily $B <^* B^*$ (as in the addition to
\scite{3.10}); also in the proof, ``$\{x, y\}$ is an edge'' means $xRy
\vee xSy \vee ySx$; in $(*)_2$ there 
$\bold w_\lambda(A_{i, j},C''_{d_{i, j}})=0$ holds iff $A_{i, j} \leq^{**}
C^{\prime \prime}_{d_{i, j}}$, and be more careful in $(*)_4$ and the
proof of clause (B). However there is a gap: \scite{3.11} does not
give clause (c) of \scite{4.5}. 
For this we can use the ``simply$^*$ almost nice", i.e.,
use [I,2.20 - 2.24].  So the parallel to \scite{3.8} is:
\proclaim{\stag{4.7} Claim}   Assume
\mr
\item "{$(a)$}"   $A <_s B$, $\bar c$ list $A$,
\sn
\item "{$(b)$}"  $B_0\subseteq B_1= B$
\sn
\item "{$(c)$}"   if $x\in B_1$ and $\{x_0, \ldots, x_j\}$ is its
$S^B$-component, $\dsize \bigwedge_{i<j} S^B(x_i,x_{i+1})$,
$x=x_{i(*)}$ then $P_f(x_0)\vee (i(*)>k)$ and $P_\ell(x_j)\vee (i(*)<
j-k)$
\sn
\item "{$(d)$}"   $A \leq A'\leq N$, $B\leq N$, 
$\nonforkin{B}{A'}_{A}^{N}$, and $c \ell^k(B_0, N)\subseteq B_1\cup A'$.
\ermn
\ub{Then} for some $\psi(\bar x)$ (of size depending on $k,
\ell g(\bar x) = \ell g(\bar c)$ only), the sequence $(B,\bar c,\psi(\bar
x),\langle B_0, B_1\rangle, k, k)$ is simply$^*$ good.
\endproclaim
\bigskip

\demo{Proof}  Let $\psi(\bar x)$ say 
exactly which quantifier free types over $A$ of 
$\leq k$ elements are realized in $A'$.

Let ${\Cal M}_n$ be random enough and $f:A \hra {\Cal M}_n$ be such that

$$
{\Cal M}_n \restriction c \ell^k(f(A),{\Cal M}_n)\models \psi(f(\bar c)).
$$
\mn
By \scite{4.4}(2) we find an extension $g:B \hra {\Cal M}_n$ of $f$ and $B^+$
as there. \nl
We continue as in \scite{4.5}'s proof till we conclude $C\subseteq
B^+ \subseteq B^+ \cup c \ell^k(f(A),{\Cal M}_n)$ (and including it),
but we do not define $C_2$. Instead we note that (by the choice of 
$\psi$) there is a embedding $h$, from some $C_{-2}\leq A'$ onto 
$C_2 =: C\cap c \ell^k(A),{\Cal M}_n)$ such that 
$h \restriction (C_{-2}\cap A) = f\restriction (C_{-2}\cap A)$. 
As $\nonforkin{B}{A'}_{A}^{N}$ and $\nonforkin{g(B)}{c \ell^k(f(A),
{\Cal M}_n)}_{f(A)}^{{\Cal M}_n}$ clearly $g'=: g\cup h$ embedd 
$B\cup C_{-2}$ (i.e., $N\restriction (B\cup C_{-2}))$ onto
$g(B)\cup C_2 \subseteq B^+ \cup (C\cap c \ell^k(f(A),{\Cal M}_n))$).

Clearly $g(B)\cup C_2 \leq^{**} B^+\cup C_2$ (as
$\nonforkin{B^+}{f(A)\cup C_2}_{f(A)}^{{\Cal M}_n}$ which holds by
\scite{4.4}(2)(iii)) so (by $\otimes$ from the proof of \scite{4.4})

$$
C\cap (g(B)\cup C_2\leq^{**} C.
$$
\mn
As also $C\cap g(B_0)\leq_i C$, $C\cap g(B_0)\subseteq C\cap (g(B)\cup
C_2)\subseteq C$, by \scite{4.6} we know that $C\cap g(B_0)\leq_i C\cap
(g(B)\cup C_2)$. By assumption (c) of \scite{4.7} (instead (iii) of
\scite{4.4}) we finish as in the proof of \scite{4.4}.

Now \scite{3.11} fit well with \scite{4.7}. So we have finished proving
Theorem \scite{4.3}.  \hfill$\square_{\scite{4.3}}$\margincite{4.3}
\enddemo
\bigskip

\demo{Proof of Theorem \scite{4.1}}  Like the 
proof of \scite{4.3}, but in the vocabulary we have also two unary 
predicates $P_f$, $P_\ell$, and we replace ${\Cal K},{\Cal K}'$ by

$$
\align
{\Cal K} = \{(X,R,S,P_f,P_\ell):&(X,R,S) \in {\Cal K}'
\text{ from the proof of \scite{4.3}}, \\
  &|P_f|\leq 1,|P_\ell|\leq 1,\text{ and } S(x,y) 
\Rightarrow \neg P_f(x) \wedge \neg P_\ell(y) \text{ and} \\
  &(P_f(x) \wedge P_\ell(y)\wedge(x,y) \text{ are } S
  \text{-connected)}) \Rightarrow  X \text{ is } S \text{-connected}\}.
\endalign
$$

$$
{\Cal K}' = \{(X,R,S,P_f,P_\ell) \in {\Cal K}:\text{ no } S 
\text{-cycle in }(X,S)\}.
$$
\mn
We can prove, as in \scite{4.3}, that
\mr
\item "{$\otimes$}"   ${\frak K}$ is simply$^*$ almost nice.
\ermn
The difference in the proof with \scite{4.3} is that here for positive
$k \in \Bbb N$, for random enough ${\Cal M}_n,c \ell^{k,m}
(\emptyset,{\Cal M}_n)$ is not empty, so we get only convergence in 
\scite{1.3} (in fact, e.g., 

$$
(\exists x,y,z)(S(x,y)\wedge P_f(x)\wedge R(x,z))
$$
\mn
has probability $p_2$).  Still applying [I,2.19] we need
\mr
\item "{$\otimes_1$}"   for every f.o. $\varphi$ the sequence
$$
\langle \text{prob}({\Cal M}_n \restriction c \ell^k(\emptyset,{\Cal
M}_n) \models \varphi):n \in \Bbb Nn\rangle
$$
converges.
\ermn
For this it suffices to prove
\mr
\item "{$\otimes_2(a)$}"  for every $\varepsilon \in \Bbb R^{>0}$ and
$k \in \Bbb Nn$, for some $m\in \Bbb Nn$, for every $n$ large enough 
$1-\varepsilon\leq \text{prob}({\Cal E}^n_k)$, where ${\Cal E}^n_k$ 
is the event
$$
e \ell^k(\emptyset,{\Cal M}_n) \subseteq \{1,\ldots,m\} \cup
\{n-m+1,\ldots,n-1,n\},
$$
\sn
\item "{$(b)$}"   for every first order $\varphi$, for some $k$
assuming ${\Cal E}^n_k$ occurs and for random enough ${\Cal M}_n$ 
(in particular, $n>2k+1$),
the satisfaction of ${\Cal M}_n \restriction c \ell^k(\emptyset,{\Cal
M}_n) \models \varphi$ depends only on the isomorphism type of 
${\Cal M}_n \restriction(\{1,\ldots,m\}\cup\{n-m+1,\ldots,n\})$,
\sn
\item "{$(c)$}"   for all $n>2k$, for every $N\in {\Cal K}$ with 
$2k$ elements, the probability $\text{prob}(({\Cal M}_n \restriction
(\{1,\ldots,m\}\cup\{n-m+1,\ldots,n\}) \cong N)$ does not depend on 
$n$ or at least, as a function of $n$, it converge.
\ermn
Now in $\otimes_2$, clauses (b) and (c) are immediate. For proving
(a), we show by induction on $\ell$ that
\mr
\item "{$\otimes_3$}"   for every $\varepsilon$ for some $m$, for every $n$
large enough and $w \subseteq [n]$ with $k-\ell$ elements
$$
1-\varepsilon\leq \text{ prob} \cases
\text{if } &{\Cal M}_n \restriction w<_i A \leq {\Cal M}_n,|A|\leq k \\
\text{then } &(\forall x\in A\setminus w)(\exists y\in w)[|x-y|\leq m]
\endcases
$$
\ermn
(Note: this is for a fixed $w$; if we say ``for every $w$'', this is a
different matter).

If you have read the proof of \scite{2.4} this should be clear.
\hfill$\square_{\scite{4.1}}$\margincite{4.1}
\enddemo
\newpage

     \shlhetal 

\nocite{ignore-this-bibtex-warning} 
\newpage
    
REFERENCES.  
\bibliographystyle{lit-plain}
\bibliography{lista,listb,listx,listf,liste}

\def\germ{\frak} \def\scr{\cal} \ifx\documentclass\undefinedcs
  \def\bf{\fam\bffam\tenbf}\def\rm{\fam0\tenrm}\fi 
  \def\defaultdefine#1#2{\expandafter\ifx\csname#1\endcsname\relax
  \expandafter\def\csname#1\endcsname{#2}\fi} \defaultdefine{Bbb}{\bf}
  \defaultdefine{frak}{\bf} \defaultdefine{mathfrak}{\frak}
  \defaultdefine{mathbb}{\bf} \defaultdefine{mathcal}{\cal}
  \defaultdefine{beth}{BETH}\defaultdefine{cal}{\bf} \def\bbfI{{\Bbb I}}
  \def\mbox{\hbox} \def\text{\hbox} \def\om{\omega} \def\Cal#1{{\bf #1}}
  \def\pcf{pcf} \defaultdefine{cf}{cf} \defaultdefine{reals}{{\Bbb R}}
  \defaultdefine{real}{{\Bbb R}} \def\restriction{{|}} \def\club{CLUB}
  \def\w{\omega} \def\exist{\exists} \def\se{{\germ se}} \def\bb{{\bf b}}
  \def\equivalence{\equiv} \let\lt< \let\gt> \def\implies{\Rightarrow}
\begin{thebibliography}{ShSp 304}
\makeatletter \renewcommand{\@biblabel}[1]{[#1]} \makeatother
\def\eprintfootnotetext{References of the form {\tt math.XX/$\cdots$}
 refer to the {\tt xxx.lanl.gov} archive}
\ifx\documentstyle\undefinedcontrolsequence
   \def\anyfootnote{\footnote{*}}
   \else\def\anyfootnote{\footnote}\fi
\def\eprintfn{\ifEprint\anyfootnote{\eprintfootnotetext}\fi\Eprintfalse }
\newif\ifEprint  \Eprinttrue

\bibitem[BlSh 528]{BlSh:528}John~T. Baldwin and Saharon Shelah.
\newblock {Randomness and Semigenericity}.
\newblock {\em {Transactions of the American Mathematical Society}}, {\bf
  349}:1359--1376, 1997.
\newblock math.LO/9607226.

\bibitem[Sh 550]{Sh:550}Saharon Shelah.
\newblock {0--1 laws}.
\newblock {\em {Preprint}}.
\newblock math.LO/9804154.

\bibitem[Sh 467]{Sh:467}Saharon Shelah.
\newblock {Zero-one laws for graphs with edge probabilities decaying with
  distance. Part I}.
\newblock {\em {Fundamenta Mathematicae}}, {\bf 175}:195--239, 2002.
\newblock math.LO/9606226.

\bibitem[ShSp 304]{ShSp:304}Saharon Shelah and Joel Spencer.
\newblock {Zero-one laws for sparse random graphs}.
\newblock {\em {Journal of the American Mathematical Society}}, {\bf
  1}:97--115, 1988.

\end{thebibliography}

\enddocument